\documentclass[a4paper,12pt]{amsart}

\usepackage{amsthm}
\usepackage{amsmath}
\usepackage{amsfonts}
\usepackage{amssymb}
\usepackage{times}
\usepackage{graphicx}
\usepackage{rotating}

\theoremstyle{definition}
\newtheorem{defn}{\indent\bf Definition}
\newtheorem{rem}[defn]{\indent\bf Remark}
\newtheorem{ex}[defn]{\indent\bf Example}
\newtheorem{obs}[defn]{\indent\bf Observation}

\theoremstyle{plain}
\newtheorem{lemma}[defn]{\indent\bf Lemma}
\newtheorem{prop}[defn]{\indent\bf Proposition}
\newtheorem{thm}[defn]{\indent\bf Theorem}
\newtheorem{cor}[defn]{\indent\bf Corollary}

\begin{document}

\title[Type II$_1$ von Neumann algebra representations of Hecke 
operators]{Type II$_1$ von Neumann  algebra representations of Hecke 
operators on Maass forms and the Ramanujan-Petersson conjectures}
\author[Florin R\u adulescu]{Florin R\u adulescu${}^*$
  \\ \\ 
Dipartimento di Matematica\\ Universita degli Studi di Roma ``Tor Vergata''}
\dedicatory{Dedicated to Professor Dan Virgil Voiculescu on the occasion of his 60'th anniversary}

\maketitle 

\thispagestyle{empty}

\def\tilde{\widetilde}
\def\a{\alpha}
\def\T{\theta}
\def\PSL{\mathop{\rm PSL}\nolimits}
\def\SL{\mathop{\rm SL}\nolimits}
\def\PGL{\mathop{\rm PGL}}
\def\Per{\mathop{\rm Per}}
\def\GL{\mathop{\rm GL}}
\def\Out{\mathop{\rm Out}}
\def\Int{\mathop{\rm Int}}
\def\Aut{\mathop{\rm Aut}}
\def\ind{\mathop{\rm ind}}
\def\card{\mathop{\rm card}}
\def\d{{\rm d}}
\def\Z{\mathbb Z}
\def\R{\mathbb R}
\def\cR{{\mathcal R}}
\def\tPsi{{\tilde{\Psi}}}
\def\Q{\mathbb Q}
\def\N{\mathbb N}
\def\C{\mathbb C}
\def\bH{\mathbb H}
\def\Y{{\mathcal Y}}
\def\L{{\mathcal L}}
\def\G{{\mathcal G}}
\def\U{{\mathcal U}}
\def\QC{{\mathcal Q}}
\def\H{{\mathcal H}}
\def\I{{\mathcal I}}
\def\A{{\mathcal A}}
\def\S{{\mathcal S}}
\def\O{{\mathcal O}}
\def\V{{\mathcal V}}
\def\D{{\mathcal D}}
\def\B{{\mathcal B}}
\def\K{{\mathcal K}}
\def\cC{{\mathcal C}}
\def\cR{{\mathcal R}}
\def\cX{{\mathcal X}}
\def\ptimes{\mathop{\boxtimes}\limits}
\def\potimes{\mathop{\otimes}\limits}

\begin{abstract}
Classical Hecke operators on Maass 
forms are unitarily equivalent, up to a commuting, operatorial, phase, to completely positive
maps on II$_1$ factors, associated to a pair of isomorphic
subfactors, and an intertwining unitary. This  representation is obtained through a quantized representation of the Hecke operators. The Hecke operators  act on the Berezin's quantization, deformation algebra of the fundamental domain of $\PSL(2,\Z)$ in the upper halfplane. The Hecke operators are inheriting from the ambient, non-commutative algebra on which they act, a rich structure of matrix inequalities. Using this construction we obtain  that,  for every prime $p$, the essential spectrum of the classical Hecke operator  $T_p$ is contained in the interval $[-2\sqrt p, 2\sqrt p]$, predicted by the Ramanujan Petersson conjectures. In particular, given an open interval containing $[-2\sqrt p, 2\sqrt p]$,  there are at most a finite number of possible exceptional eigenvalues lying outside this interval.
The main tool for obtaining this  representation of the Hecke operators (unitarily
equivalent to the classical representation, up to a commuting phase) is a Schurr  type, positive  "square root" of the state on $\PGL(2,\Q)$,
measuring the displacement of fundamental domain of $\PSL(2,\Z)$ in $\mathbb H$, by translations in $\PGL(2,\Q)$. The "square root" is obtained from the matrix coefficients
of the discrete series representations of $\PSL(2,\R)$ restricted to $\PGL(2, \Q)$. The methods in this paper may also be applied to any finite index, modular  subgroup $\Gamma_0(p^n)$, $n\geq 1$, of  $\PSL(2,\Z)$. In this case the essential norm of the Hecke operator is equal to the norm of the corresponding convolution operator on the cosets Hilbert space $\ell^2((\Gamma_0(p^n))\backslash \PGL(2,\Z[1/p])$.
\end{abstract}

\renewcommand{\thefootnote}{}
\footnotetext{${}^*$ Member of the Institute of  Mathematics ``S. Stoilow" of the Romanian Academy}
\footnotetext{${}^*$  
Supported in part by PRIN-MIUR, and EU network Noncommutative Geometry MRTN-CT-2006-0031962,  PN-II-ID-PCE-2012-4-0201 and a personal grant from WBS Holding, Bucharest}

\section*{}

\section*{Introduction}

In this paper we obtain an operator algebra representation for the classical 
Hecke operators. We prove that the classical operators admit a "quantized" representation, 
 to which they are  unitarily equivalent, up to commuting phase.
The "quantized" Hecke operators  act on the noncommuative von Neumann algebra associated to the
$\PSL(2,\Z)$-equivariant,  Berezin's quantization deformation of the upper halfplane ([Ra1]). Using matrix positivity properties, inherent to operator algebra structures, we  deduce various properties for  the Hecke operators on Maass forms,
e.g. we  compute the essential spectrum.

Some of the results in this paper are valid in a more general setting. We start with countable a discrete group $G$ with an almost normal subgroup $\Gamma$, such
that the set $\S$ of finite index subgroups of the form $\Gamma_\sigma = \Gamma \cap \sigma \Gamma \sigma^{-1}$, $\sigma$ in $G$,
generates a downward directed, modular lattice, with respect to inclusion.

The Hecke algebra $\H_0=\H_0(\Gamma,G)$  of double cosets of $\Gamma$ in $G$ has a canonical representation,
called left regular representation, acting by left convolution 
on $\ell^2(\Gamma/G)$ (see [BC]). Our basic object will be von Neumann algebra $\H$ , the closure of $\H_0$,
in the weak operator  topology  on the bounded linear operators $B(\ell^2(\Gamma/G)$, that are acting on $\ell^2(\Gamma/G)$. We will refer to the von Neumann algebra $\H$ as to the reduced von Neumann Hecke algebra (as customary  in operator algebra). When taking the norm closure of $\H_0$ in $B(\ell^2(\Gamma/G)$ we obtain the reduced $C^\ast$ reduced Hecke algebra. We denote this $C^\ast$-algebra by $\H_{\rm {red}}$. 

Our main assumption is that there exists a (projective) unitary representation $\pi$
of $G$ on $\ell^2(\Gamma)$, extending the left regular representation  (projective, when a group 2-cocycle is present) of $\Gamma$ on 
$\ell^2(\Gamma)$. 
This assumption  implies in particular that $[\Gamma:\Gamma_{\sigma}]=[\Gamma:\Gamma_{\sigma^{-1}}]$
for all $\sigma$ in $G$ (see the paper [Ra7] for a more general setting). 

This assumption is equivalent to the existence of an   isometric embedding  of the Hilbert  spaces having as orthonormal basis the left (respectively right) cosets in $G$, of  the subgroup $\Gamma$, into the Hilbert space associated to the type II$_1$ von Neumann algebra $\L(G)$ associated to the discrete group $G$. We require that this embedding transforms cosets concatenation into algebra multiplication and and we require  that the $\ast$-operation on $\L(G)$ moves the image (through the embedding) of a  left coset $\Gamma \sigma$ into the corresponding  image of the  right coset $\sigma^{-1}\Gamma$, for all $\sigma$ in $G$.

The above mentioned embedding  is constructed out of the data given by the matrix coefficients of the representation $\pi$. 
 
   More precisely, we  let
  $\cC(G,\Gamma)$ be the linear space spanned by all sets of the form $[\sigma_1\Gamma\sigma_2]$,
$\sigma_1, \sigma_2 \in G$, subject to the obvious relation that 
$$\sum_i [\sigma_1^i\Gamma\sigma_2^i]=\sum_j [\theta_1^j \Gamma \theta_2^j],$$
 whenever
$\sigma_\varepsilon^i, \theta_\varepsilon^j$, $\varepsilon = 1,2$ are elements in $G$, 
such that the sets $(\sigma_1^i\Gamma\sigma_2^i)_i$ and respectively $(\theta_1^j\Gamma \theta_2^j)_j$
are disjoint,  and 
$$\bigcup_i  \sigma_1^i\Gamma\sigma_2^i = \bigcup _j \theta_1^j \Gamma \theta_2^j.$$
of equal union. The adjoint map $\ast$ on $\cC(G,\Gamma)$ is defined by mapping $[\sigma_1\Gamma\sigma_2]$  into
$[\sigma_2^{-1}\Gamma\sigma_1^{-1}]$. In particular   the $*$ operation maps $[\sigma_1\Gamma]$ into $[\Gamma{\sigma_1^{-1}}]$.

Let $\C(\Gamma/G)$, (respectively $\C(G/\Gamma)$) be the vector space having as basis the left (respectively right) cosets of $\Gamma$ in $G$. There exists a canonical pairing  $\C(G/\Gamma)\times \C(\Gamma/G)\to \cC(G,\Gamma)$,  mapping $[\sigma_1\Gamma] \times[\Gamma \sigma_2]$ into $[\sigma_1\Gamma\sigma_2]$
$\sigma_1, \sigma_2 \in G$ (this is what we call coset concatenation).
This map obviously factors   to  $\C(\Gamma/G)\otimes_{\H_0} \C(G/\Gamma)\to \cC(G,\Gamma)$,  and hence gives another way to define  the multiplication on $\H_0=\H_0(\Gamma,G)$.

A representation of $\cC(G,\Gamma)$ into a II$_1$ factor $M$ with trace $\tau$ is an isometric embedding of the Hilbert spaces $\ell^2(\Gamma\setminus G)$,
$\ell^2(G/\Gamma)$  into the standard Hilbert space $L^2(M,\tau)$ associated to $M$ and $\tau$ via the GNS representation. This embedding  should be compatible with the $*$ operation, and  should transform the concatenation 
$[\sigma_1\Gamma]\times[\Gamma{\sigma_2}]=[\sigma_1\Gamma\sigma_2]$ into the algebra product in $M$.

 Let  $\L(G,\varepsilon)$, $\L(\Gamma,\varepsilon)$  be the finite von Neumann algebras, with cocycle $\varepsilon$ associated to the discrete groups $G, \Gamma$ (see e.g.  [Su] for definitions). Here $\varepsilon$ is the two cocycle on $G$ associated with the projective representation $\pi$ of $G$ considered above.
 
From the matrix coefficients of the representation $\pi$, we construct the representation  $t$ of $\cC(G,\Gamma)$ into the von Neumann II$_1$ factor $\L(G,\varepsilon)=
\overline{\C(G,\varepsilon)}^w$ (by $\overline{\ \cdot\  } ^w$ we designate the closure in the weak operator topology), associated with the group $G$ and with the 2-group cocycle $\varepsilon$.
In this representation, the cosets $[\Gamma\sigma]$ are mapped into a family $t^{\Gamma\sigma}\in \ell^2(\Gamma\sigma)\cap \L(G,\varepsilon)$, $\sigma\in G$. 

The formula for $t^{\Gamma\sigma}$, $\sigma\in G$ depends on
 the matrix algebra coefficients of the representation $\pi$ with respect to
the unit vector $I$ in $\ell^2(\Gamma)$ corresponding to the identity element of $\Gamma$.
More precisely,
$$
t^{\Gamma\sigma} = \sum _{\theta \in \Gamma\sigma} \overline{\langle\pi(\theta)I,I\rangle}\theta,
$$
where $I \in \ell^2(\Gamma) $ is the vector corresponding to the neutral element of $\Gamma$.

Even more general, if $A$ is a subset of $G$ we define 
$$t^A= \sum _{\theta \in A} \overline{\langle\pi(\theta)I,I\rangle}\theta.
$$

The main property of the elements $t^{\Gamma\sigma}\in \L(\Gamma,\varepsilon) $, $\sigma \in G$, is that, 
with respect to the adjoint and  multiplication operation on $\L(\Gamma,\varepsilon)$, we have 
$$(t^{\Gamma\sigma_1})^{\ast}t^{\Gamma\sigma_2}=t^{\sigma_1\Gamma}t^{\Gamma\sigma_2}= 
t^{\sigma_1\Gamma\sigma_2},\quad \sigma_1, \sigma_2 \in G.
$$
This defines the representation $t$ of $\cC(G,\Gamma)$ into $\L(G,\varepsilon)$.

 Our construction also shows that we may chose
$(t^{\Gamma\sigma})_{\sigma\in G}$ so they constitute a basis of $\L(G,\varepsilon)$ as a module over
$\L(\Gamma,\varepsilon)$ (a Pimsner-Popa basis ([PP]) for $\L(\Gamma,\varepsilon) \subseteq \L(G,\varepsilon)$).
Moreover, $t^{\Gamma\sigma}$ is supported in $\ell^2(\Gamma\sigma)$.

By using this representation we construct a $*$-algebra representation of the Hecke  algebra $\H_0$, mapping a double coset $[\Gamma\sigma\Gamma]$ into a 
 into a completely positive map $\Psi_{[\Gamma\sigma\Gamma]}$ on the von Neumann algebra associated with $G$. We will refer to the completely positive maps $\Psi_{[\Gamma\sigma\Gamma]}$ by calling them "quantized" or abstract Hecke operators. This is because,  when $\Gamma$ is $\PSL_2(\Z)$, these completely positive maps are proven to be unitarily equivalent, via the quantization representation, to the classical Hecke operators for $G = \PGL_2(\Z[\frac1{p}])$, (modulo a positive phase operator, commuting to the Laplacian). These "quantized' Hecke operators  are constructed, as described bellow,  by 
using the representation $t$ of $\cC(G,\Gamma)$.
  
 Let  $\sigma  \in G$, and let  $[\Gamma \sigma \Gamma]$  be the corresponding double coset. According to the previous definition for general subsets $A$ of $G$, the operator $t^{\Gamma\sigma\Gamma}$ is simply $\sum\limits_{[\Gamma{\sigma s]}\subseteq[\Gamma\sigma\Gamma]} t^{\Gamma{\sigma s}}$ where $s\in \Gamma$ runs over a system of representatives for cosets of $\Gamma_{\sigma}\subseteq\Gamma$. Let $E^{\L(G,\varepsilon)}_{\L(\Gamma,\varepsilon)}$ be
the canonical conditional expectation from $\L(G,\varepsilon)\to\L(\Gamma,\varepsilon)$
(the conditional  is the linear, positive  map on $\L(G,\varepsilon$, killing all $g$ with $g$ not in $\Gamma$, extended then by linearity and continuity to $\L(G,\varepsilon)$ .

The abstract Hecke operators are constructed as follows.
The abstract Hecke operator $\Psi_{[\Gamma\sigma\Gamma]}$ is the completely positive, unital operator on $\L(\Gamma,\varepsilon)$ (extendable to $\ell^2(\Gamma)$), defined by the formula
$$
\Psi_{[\Gamma\sigma\Gamma]}(x) = [\Gamma:\Gamma_{\sigma}] E^{\L(G,\varepsilon)}_{\L(\Gamma,\varepsilon)} (t^{\Gamma\sigma\Gamma}x (t^{\Gamma\sigma\Gamma})^{\ast}), \quad x \in \L(\Gamma,\varepsilon).  \leqno(1)
$$

In general, throughout the paper we will use the convention that $\Psi_{[\Gamma\sigma\Gamma]}$ is the non-normalized operator corresponding to the coset $[\Gamma\sigma\Gamma]$, while $\tilde{\Psi}_{[\Gamma\sigma\Gamma]}$ will stand for the normalized Hecke operator: $$\tilde{\Psi}_{[\Gamma\sigma\Gamma]}=\frac {1}{[\Gamma : \Gamma_{\sigma}]}\Psi_{[\Gamma\sigma\Gamma]},$$ so that $\tilde{\Psi}_{[\Gamma\sigma\Gamma]}(1)=1$.

The abstract Hecke operators are canonically determined by   the representation $t$ of $\cC(G, \Gamma)$ that
we described above (or equivalently since $t$ is computed  from the coefficients of $\pi$ by the representation $\pi$ of $G$). 

We will prove in Appendix 4 and in Example \ref{Heckeas2} that this new representation of the Hecke algebra corresponds to a new method of constructing Hecke algebra representations. One starts with a representation of  the groupoid $(G\times G^{\rm op}) \rtimes K$ on a Hilbert space $V$  ($K$ is the profinite completion of $\Gamma$). By restricting to 
$\Gamma \times\Gamma$ invariant vectors in $V$, one obtains a new representation of the Hecke algebra associated to $\Gamma \subseteq G$. In Example \ref{Heckeas2}, we prove that the above construction is a particular realization of this new model for the Hecke operators.

  In particular, the family $\Psi_{[\Gamma\sigma\Gamma]}, \sigma\in G,$ forms a hypergroup (see e.g. [Ve]) of completely positive maps (that is the product of any two elements in the family is a linear combination of elements in the family, with positive coefficients). Then formula (1) is a Stinespring dilation of the hypergroup $\Psi_{[\Gamma\sigma\Gamma]}, \sigma\in G$. 
  
  Indeed, recall that in quantum dynamics ([Bel], [Bh], [Par], [Ar]), for a semigroup of unital completely positive maps $\Phi_n$, $n\geq 0$, $n\in \Z$, on a II$_1$ von Neumann algebra $A$, on finds a larger II$_1$ von Neumann algebra $B$, a semigroup of endomorphisms $\rho^n$, $n\in\N$ of $B$, such that $\rho^n (B)$ is increasing with $n$ and such that if $E=E^B_A$ is the conditional expectation from $B$ onto $A$, then
  $$ \Phi_n(x)=E(\rho^n (x)),\quad x \in B.$$
  When $\rho$ is inner, that is, if there exists a unitary $u$ in $B$ such that  $\rho(x)= uxu^{\ast}$, this is analogous to formula (1).

A generalized form of the Ramanujan Petersson conjectures can be formulated as follows:
\vskip6pt
{\bf Generalized Ramanujan Petersson Conjecture.}
{\it Let  $G$ be a countable discrete group with an almost normal subgroup $\Gamma$. Let $t$ be  a representation of $\cC(G,\Gamma)$ with the properties outlined above (equivalently $t$ is defined by using the matrix coefficients of unitary representation $\pi$ of $G$ extending the left regular representation). 

The statement of the conjecture is that  
the $\ast$-algebra homeomorphism
$$[\Gamma \sigma \Gamma]\rightarrow  \Psi_{[\Gamma\sigma\Gamma]}, \quad\sigma\in G,$$
 from the Hecke algebra $\H_0=\H_0(\Gamma,G)$  into the bounded linear maps on $\L(\Gamma,\varepsilon)$ (extended by continuity to the to bounded linear maps on $\ell^2(\Gamma)$), has a continuous  extension (with respect to the weak operator topology on $\H$) from 
 $\H$ into~$B(\ell^2(\Gamma))$.}
\vskip6pt

We will prove that for $\Gamma= \PSL_2(\Z)$ this corresponds to the classical Ramanujan Petersson conjectures for Hecke operators on Maass wave forms.
Our main result is the following
\vskip6pt

{\bf Theorem.} 
{\it Let $p$ be a prime number. Let $G=\PGL_2(\Z[\frac1{p}])\supseteq\Gamma= \PSL_2(\Z)$. Let $\pi$ be the representation $\pi_{13}|_G$, where $\pi_{13}$ is  the 13-th projective unitary representation  in the discrete series of unitary representations of $\PSL(2,\R)$. For a double coset $[\Gamma\sigma\Gamma]$ let $\Psi_{[\Gamma\sigma\Gamma]}$ be the completely positive map constructed as above.
Let $\Pi_Q$ be the projection from $B(\ell^2(\Gamma))$  onto the Calkin algebra
(see e.g. [Do] for the definition of Calkin algebra)
 $$Q(\ell^2(\Gamma)= B(\ell^2(\Gamma))/K(\ell^2(\Gamma)).$$ 
 Then 
the $\ast$-algebra homeomorphism
$$[\Gamma \sigma \Gamma]\rightarrow  \Pi_Q(\Psi_{[\Gamma\sigma\Gamma]}),\quad \sigma\in G,$$
 from the Hecke algebra $\H_0=\H_0(\Gamma,G)$  into $Q(\ell^2(\Gamma)$ has a continuous  extension (with respect to the weak operator topology on $\H$) from 
 $\H$ into  $Q(\ell^2(\Gamma))$.
 
Moreover, the operators $\Psi_{[\Gamma\sigma\Gamma]}$, $\sigma\in G$, are unitarily equivalent, up to a
commuting phase to the classical Hecke operators on Maass wave forms,  corresponding to the cosets $[\Gamma\sigma\Gamma]$.
This implies that the essential spectrum of the classical operators $T_n=T_{[\Gamma\sigma_{p^n}\Gamma]}$
coincides with the spectrum in the representation of the Hecke algebra on $\ell^2(\Gamma\setminus G)$.
The spectrum in this  last representation  coincides with the spectrum  predicted by the Ramanujan-Petersson Conjectures. 

This result holds true for finite index modular subgroups of $\PSL(2,\Bbb Z)$, the essential norm of the corresponding Hecke operators is then equal to the norm of the
corresponding convolutor in the reduced C$^\ast$ Hecke algebra $\mathcal H_{\rm red}$.
} 
\vskip6pt

We explain bellow the reformulation of this result in classical terms.
Recall that the classical Hecke operators are acting on $L^2(F,\frac{\d\overline{z}\d z}{({\rm Im}\, z)^2})$,
where $F$ is a fundamental domain for the action of $\PSL_2(\Z)$ in the upper half plane $\mathbb H$. Let $n$ be an natural number. The
classical Hecke operator, corresponding to the sum of double cosets in matrices of determinant
$n$ is given by the formula
$$
T(n) f(z) = \sum_{{\rm ad}=n\atop b=0,1,\ldots,d-1} f\left(\frac{az+b}{d}\right)
$$
and the normalized version
$$
\tilde{\tilde{T}}(n)=\frac{1}{\sqrt{n}} T_n.
$$

The Ramanujan-Petersson conjecture states that if $c(p)$ are the eigenvalues
for a common eigenvector $\xi \neq 0$, for all the $\tilde{\tilde{T}}_p$'s, then $c(p)\in [-2,2]$ for 
all primes $p$ (see [Hej]). This corresponds, when working with the non-normalized Hecke operator $T_p$, to the fact that the eigenvalues should be in the interval $[-2\sqrt p, 2\sqrt p]$ (it is well known see e.g. [Hej] that it is sufficient to verify the conjecture for $n$ a prime number).

It is well known (going back to Hecke's and Peterssons's work (see e.g. [Krieg])) that the Hecke operators give a $\ast$-algebra representation for the Hecke algebra associated to $G=\PGL_2(\Q)\supseteq\Gamma= \PSL_2(\Z)$. As  formulated above,  the conjecture is equivalent to the   continuity, with respect to the weak operator topology on the Hecke algebra, of the linear application  mapping a double coset (which is labelled by $n$ -the determinant),
in the Hecke algebra $\H_0( \PSL_2(\Z),\PGL_2(\Q))$,  into the Hecke operator $T(n).$

 The conjecture thus makes sense  in the more general setting of a group $G$, an almost normal subgroup $\Gamma$ and $\pi$ a projective unitary representation of $G$ on $l^2(\Gamma)$ extending the left regular representation (with cocycle) of $\Gamma$. 
The Hecke operators are replaced by the operators in formula (1), and the Ramanujan Petersson conjectured estimates are equivalent to conjecturing 
 the continuity of the linear application which maps $[\Gamma\sigma\Gamma] $ into the completely positive map $\Psi_{[\Gamma\sigma\Gamma]} $. This is equivalent, by what we explained in the preceding paragraph to the classical case for
 $\PGL_2(\Z[\frac1{p}])\supseteq \PSL_2(\Z)$.
 
 We prove therefore that this continuity holds, when replacing the Hecke operators with their image in the Calkin algebra, and thus prove that
 the essential spectrum of the Hecke operator sits in the predicted interval ($[-2\sqrt p, 2\sqrt p]$). Therefore our main result implies the following:
 \vskip6pt
 
 {\bf Corollary.} {\it
 For every prime $p$ the essential spectrum of the classical Hecke operator  $T_p$ is contained in the interval $[-2\sqrt p, 2\sqrt p]$, predicted by the Ramanujan-Petersson conjectures. In particular, given an open interval containing $[-2\sqrt p, 2\sqrt p]$,  there are at most a finite number of possible exceptional eigenvalues lying outside this interval.}
\vskip6pt

Our result shows that the representation of the Hecke algebra into  completely positive maps have a canonical extension to $\cC(G, \Gamma)$.
Hence their knowledge is  relevant for the determination of the eigenvalues.

The fact that the classical Hecke operators are unitarily equivalent to the abstract
Hecke operators in formula (1) is outlined bellow.

First we describe a more abstract setting. Let $G$  be a discrete countable group and $\Gamma$ an almost normal subgroup, with the modular family of subgroups described above.
Assume that $H$  is a  Hilbert space acted unitarily by $G$, with a rich family of $\Gamma$ fixed vectors.
We denote by $H^{\Gamma_\sigma}$  the Hilbert space  of vectors in $H$ fixed by the subgroup $\Gamma_\sigma$, $\sigma \in G$. Then 
the Hecke operator $v \to T_\sigma(v) = \sum s_i \sigma v$ (where $\Gamma = \bigcup s_i \Gamma_\sigma$ is the decomposition into right cosets of the group $\Gamma$)
is obtained by composing the maps in the following diagram
$$
\begin{array}{ccc}
H^{\Gamma_{\sigma^{-1}}} & \mathop{\longrightarrow}\limits^{\sigma} &  H^{\Gamma_\sigma}\\
\phantom{aaaa}_{\rm inc}\nwarrow &  &\hspace{-0.3cm} 
\swarrow {}_P \\[4pt]
&  H^{\Gamma}&
\end{array}
$$
where $P$ is the orthogonal projection from $H^{\Gamma_\sigma}$ onto $H^\Gamma$. Thus 
$$T_{\sigma}v = [\Gamma:\Gamma_{\sigma}]P(\sigma v),\quad \sigma\in G,\ v \in H^\Gamma.$$

The commutant algebras 
$\{\Gamma\}', \{\Gamma_\sigma\}', \{\Gamma_{\sigma^{-1}}\}'$ in $B(\ell^2(\Gamma))$, 
are II$_1$ factors, so there is a canonical conditional expectation 
$E= E_{\{\Gamma\}'}^{\{\Gamma_\sigma\}'}$
from $ \{\Gamma_\sigma\}'$ onto $ \{\Gamma\}'$, which plays the role of the projection $P$.

In particular, if we let $\Gamma$ act on $\ell^2(\Gamma)$ (eventually with a cocycle $\varepsilon$)
and $\pi$ a unitary representation of $G$ on $\ell^2(\Gamma)$  with cocycle  $\varepsilon$, extending the left regular representation to $G$,  then
the following diagram (with $E=E_{\{\Gamma\}'}^{\{\Gamma_\sigma\}'}$, the canonical conditional expectation from $\{\Gamma_\sigma\}'$ onto $\{\Gamma\}'$)
$$
\begin{array}{ccc}
\{\Gamma_{\sigma^{-1}}\}' & \mathop{\longrightarrow}\limits^{{\rm Ad}\,\pi(\sigma)} &  \{\Gamma_\sigma\}'\\
\phantom{aaaa}_{\rm inc}\nwarrow &  &\hspace{-0.3cm} \swarrow {}_E \\[4pt]
&  \{\Gamma\}'&
\end{array}
$$
for $\sigma$ in $G$, yields a Hecke operator, $\Psi_{\sigma}=\Psi_{[\Gamma\sigma\Gamma]}$, defined by the formula:
$$
\Psi_{\sigma}{(X)} = [\Gamma:\Gamma_{\sigma}] E_{\{\pi(\Gamma)\}'}^{\{\pi(\Gamma_\sigma)\}'} (\pi(\sigma)X\pi(\sigma)^*)
=\sum_{i=1}^n\pi(s_i\sigma)(X)\pi(s_i\sigma)^*, \leqno(2)
$$
where $\Gamma = \bigcup s_i \Gamma_\sigma$ is the coset decomposition of $\Gamma$ with respect to $\Gamma_\sigma$.

The equivalence of the two representations of the Hecke operators is based on the following theorem of V.F.R.\ Jones
(see e.g.\ [GHJ]).
Let $M$ be the factor generated by the image of $\PSL_2(\Z)$
through the discrete series representation $\pi_{13}$ of $\PSL_2(\R)$.
Then as proven in ([GHJ]), $M$ is unitarily equivalent to the factor
$\L(\PSL_2(\Z),\varepsilon)$ associated to the left regular representation of the
discrete group $\PSL_2(\Z)$.
Thus in the case of $G = \PGL_2(\Z[\frac1{p}])$, $\Gamma = \PSL_2(\Z)$, the left regular representation
of $\Gamma$ on $\ell^2(\Gamma)$, with cocycle $\varepsilon$, is equivalent by [GHJ], with the restriction to 
$\Gamma$ of the 13-th element $\pi_{13}$ in the discrete series representation of $\PSL_2(\R)$.

The Hilbert space $H_{13}$ of $\pi_{13}$ is the space $H^2(\bH, {\rm d}\nu_{13})$, where ${\rm d}\nu_{13} = 
({\rm Im}\, z)^{13-2} {\rm d}z {\rm d}{\overline z}$, and $\pi$ acts by left translations via 
M\"obius transforms, corrected by the factor $J(g,z)^{13} = (cz+d)^{-13}$, $z\in \bH$,
$g = \left( \begin{array}{cc} a & b \\ c & d \end{array} \right)$ in $\PSL_2(\R)$.

The operators in $B(H_{13})$ (the bounded linear operators on $H_{13}$), by using Berezin's quantization method ([Be]), are represented 
by reproducing kernels  $k(\overline{z},\eta)$, $z,\eta \in \bH$, which are analytic
functions on $\eta$ and antianalytic functions of $z$, subject to certain growth condition ([Ra1]). 
Then  $\{\pi_{13}(\Gamma)\}'$ consists of kernels $k$ such that $k(\overline{\gamma z}, \gamma \eta) = k(\overline{z},\eta)$, for all  $\gamma \in \Gamma$, 
$z,\eta \in \bH$. The action of $ \Psi_{\sigma}$, $\sigma \in G$, on the operator $X$ with kernel $k_X$ gives un operator with kernel given by the kernel
$$
z,\eta \to  \sum k_X(\overline{s_i \sigma z}, s_i \sigma \eta).
$$

The completely positive maps $ \Psi_{\sigma}$, $\sigma \in G$   may be looked at as a quantization of the classical Hecke operators, as they are acting on the algebra of a  deformation quantization of their classical domain. If we restrict to the diagonal we get the classical Hecke operators.  By the theory of the  Berezin transform ([Be])
(which is in fact the same as the Selberg transform) we know that the comparison between the  kernel itself and its restriction to the diagonal
is given by an invertible phase, e.g. a positive transformation - the Berezin transform.

This allows to prove that the operators in (1) and (2) are equivalent
(up to a commuting, operatorial phase).
Hence the analysis of the spectrum of the classical Hecke operators is reduced to the analysis
of the operators in formula (1).

To analyze the essential spectrum of the operators in formula (1) we compute the values of the positive states on the image of the Hecke algebra, in  the Calkin algebra (we will refer to such states as to essential states, since they determine the essential spectrum).
The states are then generically  averages, over points in $\Gamma$, distributed in cosets of modular subgroups. Thus, when passing to the Calkin algebra, equalities of the type $g_1 \gamma_0 g_2 = \gamma$, are replaced by equalities on  average, with respect to the measure induced on the profinite completion of $\Gamma$, by the supports of the finite sets of points in $\Gamma$.

Identifying the corresponding states amounts, at least when reducing to the case of limits of  finite sets of averages that converge to the identity, (i.e. the averaging points sit inside a  family of normal subgroups shrinking to the identity) to the study of the  the space of
conjugation orbits in the group, viewed as infinite measure spaces, with the counting measure on sets of orbits. Fortunately for $\PSL(2, \Z)$ this can be done exactly.

We let $G\times G^{\text {op}}$ act as a groupoid (by left and right action) on $\Gamma$ and thus on $\ell^2(\Gamma)$ as partial isometries. Let $K=\PSL_2(\mathbf Z_p)$, be the profinite completion of $\Gamma$ with respect to the modular family $\S$, with $\mathbf Z_p$ the $p$-adic integers. Denote by  $\mu_p$ the Haar measure on $K$. The algebra $C(K)$ of continuos functions on the profinite completion of $\Gamma$ is contained in $\ell^{\infty }(\Gamma)$ and thus acts on $\ell^2(\Gamma)$.\  Hence we can construct the (groupoid) reduced and maximal $C^*$-crossed product algebra
$$
\A = C^*_{\rm red} ((G \times G^{\rm op})\rtimes C(K)), \quad \A_{\text{max}}=C^*((G \times G^{\rm op})\rtimes C(K)).
$$
To construct the reduced crossed product algebra we use the canonical trace $\tau_p$ on the algebraic crossed product
$(G \times G^{\rm op})\rtimes C(K)$ induced by the $G \times G^{\rm op}$ invariant measure
$\mu_p$ on $K$.

We have a covariant representation of the crossed product
$$C^* ((G \times G^{\rm op})\rtimes C(K))$$ which comes from the embedding of $C(K)$ into
$\B(\ell^2(\Gamma))$ described above, and by representing elements in $(G \times G^{\rm op})$ as left or right convolutors. By $\B$ we denote the $C^*$ algebra that is the image of this representation.

Our main tool is a local version for the group $G=\PGL(2, \Z[1/p])$ of the Akeman- Ostrand
result ([AO], [Oz]). Indeed we prove that the image of $\B$ in the Calkin algebra (the quotient modulo the compact operators) is the reduced $C^*$- algebra product. The result is:

\vskip6pt

{\bf Theorem.} {\it Let $p$ be a prime number and let $G$ be the group $\PGL(2,\Z[\frac1p])$,
 $\Gamma=
\PSL(2,\Z)$. 
Let $\A_0=\B/K(\ell^2(\Gamma))$ be the projection in the Calkin algebra of the algebra $\B$ considered above (generated by left and right convolutors and by $C(K)$ acting on $\B\ell^2(\Gamma)$). Then $\A_0$ is isomorphic to the $C^*$-algebra $\A$, the reduced groupoid crossed product of $G \times G^{\rm op}$ acting on $K$, with respect to the invariant Haar measure on K. This remains valid if instead $C^*(G)$ we use the skewed $C^*$-algebra by  the canonical 2-group cocycle on $PSL(2,\Z)$.
}

\vskip6pt
  
Using this, we  prove  that the map $[\Gamma\sigma\Gamma]\to\Pi_{Q(\ell^2(\Gamma))}([\Psi_{\Gamma\sigma\Gamma}])$
is preserving the essential states, and hence is continuous with respect to the reduced Hecke algebra topology
on $\H$. Hence it follows that the Ramanujan-Petterson estimate holds true for the essential spectrum  in the case $G=\PGL_2(\Z[\frac1{p}])$.

Our methods also allows to derive matrix inequalities on eigenvalues for Hecke operators. This inequalities are encoded in the fact  the linear map on the reduced $C^\ast$- Hecke algebra multiplying a double coset by the corresponding (normalized) eigenvalue is  a completely positive map on the Hecke algebra $\H$.

Assume the completely positive maps ${\Psi}_\a$ in formula (1), where $\a$ runs over the
space of double cosets of $G$ have a joint eigenvector $\xi\neq 0$, and
denote by $\tilde{c}(\a)$ the corresponding eigenvalue.

The above description allows one to prove the following

\vskip6pt
{\bf Theorem.}
{\it The map on the Hecke algebra
that maps a coset $\a=[\Gamma \sigma\Gamma]$ into $$\Phi_{\tilde{c}}(\a) = \tilde{c}(\a)[\Gamma
\sigma \Gamma]$$
extends to a completely positive map on the reduced von Neumann algebra of the Hecke algebra.

In particular, this proves that the sequence $(\tilde{c}(\a))_{\a \in \Gamma \setminus G/\Gamma}$
is a completely positive multiplier for the Hecke $C^*$-algebra of $\Gamma$ in $G$.
}
\vskip6pt

This information encodes positive definiteness for various matrices whose coefficients are
linear combinations of the Hecke operators eigenvalues $\tilde{c}(\a)$'s.

In fact, the representation we obtained for the Hecke operators, through the completely positive maps  ${\Psi}_\a$ encodes a stronger positivity
result, based on the complete positivity of the bilinear form of $\H$
$$
(a,b) \to \tau_{\L(G)} (\xi^*a\xi b^*).
$$

This happens because the type II$_1$ representations encodes an action of 
$\H \otimes \H$. The Hecke operators on Maass form only the ``diagonal" part of
this action.

Another consequence of our representation for the Hecke operators  is the following; let $\A (G,\Gamma)$ the free ${}^* -\C$-algebra generated by all the cosets   $[\Gamma\sigma]$, $\sigma \in  G$, and their adjoints
($[\Gamma\sigma]^\ast=[\sigma^{-1} \Gamma]$, subject to $$
\sum[\sigma_1^i\Gamma][\Gamma\sigma_2^i]=\sum[\theta_1^j\Gamma][\Gamma \theta_2^j]
$$
if $\sigma_s^i, \theta_r^j$ are elements of $G$, and the disjoint union $\sigma_1^i \Gamma \sigma_2^i$
is equal to the disjoint union of $\theta_1^j \Gamma \theta_2^j$. Note that the above relation corresponds exactly to the fact that the Hecke algebra
of double cosets is a a subalgebra of $\A (G,\Gamma)$, by the trivial embedding of a double coset into the formal sum of its left or right cosets (using representatives).
Then we have (see Appendix 2).

\vskip6pt
{\bf Theorem.}
{\it The ${}^*- \C$-algebra $\A (G,\Gamma)$  admits at least one unital $C^*$ algebra representation.}

 \vskip6pt

Note that the Hecke  algebra operator represention in formula (1) admits an extension to the algebra $\A (G,\Gamma)$ ([Ra5]), and the content of the Ramanujan Petersson conjecture can be viewed a s a conjecture on the representations of $\A (G,\Gamma)$.

\vskip6pt

The author is indebted to Professors F. Boca, A. Figa-Talamanca, A. Gorodnik, R. Grigorchuk  D. Hejhal,  N. Monod,  H. Moscovici,  R. Nest.,   Lizhen Ji, P. Sarnak,  G. Skandalis,  Tim Steger  and L. Zsido and to the anonymous  referee for a first version of this paper for
several discussions regarding topics related to the
subject of this paper. The author is particularly indebted to Professor N. Ozawa  for several comments on this paper and for providing him his personal notes for a seminary at the University of Tokyo on the content of this paper (see [Ra3]). The author is specially thanking to Professor  S. Neshveyev for very pertinent questions on arguments in the proofs. The author is specially indebted to Professor A. Gorodnik for inviting him to the University of Bristol, and for pointing out that calculating states of equidistributed points would not be sufficient to determine the states corresponding to singular measures with respect Haar measure. 
The author is indebted to Professor Ovidiu P\u as\u arescu for pointing him out
the relation between Loeb measures and essential states.
 The author also thanks his formers colleagues at the University of Iowa, for the warm supporting environment, during first attempts toward this work, several years ago.

\section{
Hecke operators and Hilbert spaces}

In this chapter we present known facts about Hecke operators, from the point of view of  
Hecke operators as orthogonal projections composed with translation operators. This point of view
is particularly relevant when dealing with finite von Neumann algebras, since in that case the projections are
conditional expectations between von Neumann algebras. This representation of the Hecke operators as conditional expectations
unveils an operators system structure on the Hilbert space of cosets, which in turn determines the structure
of the Hecke algebra. 

Let $G$ be a discrete group and $\Gamma$ an almost normal subgroup.

We assume that the modular set  $\S$ generated by all finite index subgroups $\Gamma_\sigma$ of the
form
$\Gamma_\sigma = \sigma \Gamma{\sigma^{-1}}\cap\Gamma$, $\sigma\in\Gamma$ has the modular property, that
is for any $\sigma_1,\sigma_2$ in $G$ there exists $\sigma_3$ in $G$, such that $\Gamma_{\sigma_1} \cap
\Gamma_{\sigma_2}\supseteq \Gamma_{\sigma_3}$. Later we will also need the assumption that the indices
$[\Gamma:\Gamma_\sigma]$ and $[\Gamma: \Gamma_{\sigma^{-1}}]$ are equal.

We introduce following type of unitary representations of the group $\Gamma$.

\begin{defn}
An adelic Hilbert space representation of the group $G$, consists of 
the following data. Let $\V$ be a topological vector space, acted by $G$, and let $H\subseteq \V$ be a dense Hilbert space unitarily  
acted by $G$ (this is not the Hilbert space of the adelic Hilbert space representation).

For $\Gamma_\sigma\in\S$, we denote by $\V^{\Gamma_\sigma}$ the set of vectors in $\V$ fixed by
$\Gamma_\sigma$. We assume that we are given a family of Hilbert space $H^{\Gamma_\sigma}$ for
$\Gamma_\sigma$ in $\S$ with the following properties:

1) For all $\Gamma_1\subseteq \Gamma_0$, for $\Gamma_1, \Gamma_0$ in $\S$ then
$$
H^{\Gamma_0} = H^{\Gamma_1} \cap \V^{\Gamma_0}.
$$

2) The Hilbert space norm on $H^{\Gamma_\sigma}$, for all $\Gamma_\sigma$ in $\S$ has the property  that if
$\Gamma_{\sigma_1}\subseteq \Gamma_{\sigma_0}$ then the inclusion 
$H^{\Gamma_{\sigma_0}}\subseteq H^{\Gamma_{\sigma_1}}$ is isometric.

3) Note that if $v\in \V^{\Gamma_\sigma}$ then $\sigma_1 v$ is invariant by the group 
$\sigma_1 \Gamma_{\sigma} \sigma_1^{-1}$ and thus by 
$\sigma_1 \Gamma_\sigma \sigma_1^{-1} \cap \Gamma = \Gamma_{\sigma \sigma_1}\cap \Gamma_{\sigma_1}$,
which by modularity contains some subgroup $\Gamma_{\sigma_2}\in \S$. Thus
$\sigma_1 (\V^{\Gamma_{\sigma}})$ is contained in $\V^{\Gamma_{\sigma_2}}$ and consequently $\sigma(H^{\Gamma_\sigma})$ is
contained in $H^{\Gamma_{\sigma_2}}$. 

In particular, the group $G$ acts on the reunion of all the spaces $H^ {\Gamma_{\sigma}}, \sigma \in G$ and
$$
\sigma H^{\Gamma}\sigma^{-1}= H^{\Gamma_\sigma}.
$$

Thus $G$ acts on $H^{\rm ad}=\bigcup_{\Gamma_\sigma \in \rho} H^{\Gamma_\sigma}$ and the inductive
limit of Hilbert spaces (since all the inclusions are isometric) carries a natural inductive limit
Hilbert space pre-norm. Let $\overline{H}^{\rm ad}$ be Hilbert space completion of $H^{\rm ad}$.

We assume that $G$ acts unitarily on $\overline{H}^{\rm ad}$. We will refer to   the Hilbert space $\overline{H}^{\rm ad}$ as to   the adelic Hilbert space.

The following axiom will not be used, although it holds true in all examples. It relates the Hilbert space
$H$ with the Hilbert spaces $H^{\Gamma_\sigma}$, $\Gamma_\sigma\in \S$.

4) We assume that there exist $\langle\; ,\,\rangle$ a pairing  between a dense subspace of $H$ and the Hilbert
space $H^\Gamma$ such that  for all $\Gamma_\sigma\in \S$, and $v,w\in  H^{\Gamma_\sigma}$, such
 there exists a vector $\xi$ in $\V$, such that 
$v=\sum\limits_{\gamma\in\Gamma_\sigma}\gamma\xi$, for the topology on $\V$,
then $$\langle v,w\rangle_{H^{\Gamma_\sigma}}=\frac{1}{[\Gamma:\Gamma_\sigma]}\langle\xi,w\rangle.$$
\end{defn}

In the following we describe the orthogonal projection from $H^{\Gamma_\sigma}$
onto $H^\Gamma$. This will then be used to define an abstract Hecke operator.

\begin{defn}
Fix $\Gamma_{\sigma_0}\supseteq\Gamma_{\sigma_1}$ two subgroups in $\S$ and denote by $P_{H^{\Gamma_{\sigma_0}}}$
the orthogonal projection from $\overline{H}^{\rm ad}$
onto $H^{\Gamma_{\sigma_0}}$ and by $P_{H^{\Gamma_{\sigma_0}}}^{H_{\Gamma_{\sigma_1}}}$
the restriction of $P_{H^{\Gamma_{\sigma_0}}}$ to $H^{\Gamma_{\sigma_1}}$
(which is the same as the orthogonal projection from $H^{\Gamma_{\sigma_1}}$ onto $H^{\Gamma_{\sigma_0}}$).

When $\Gamma_{\sigma_0}=\Gamma$, we denote, the above projection, simply by $P^{H_{\Gamma_{\sigma_1}}}$.
\end{defn}

The projection $P_{H^{\Gamma_\sigma}}$
has the following property

\begin{lemma}
For all $v$ in $H^{\rm ad}$, $a$ in $\Gamma_\sigma$,
$P_{H^{\Gamma_\sigma}}$ has the property 
$P_{H^{\Gamma_\sigma}}(av)=P_{H^{\Gamma_\sigma}}(v)$.
To give a suggestive description of this property we will write 
$P_{H^{\Gamma_\sigma}}([\Gamma_\sigma]v)=P_{H^{\Gamma_\sigma}}(v)$.
\end{lemma}

\begin{proof}
Indeed for all $w\in H^{\Gamma_\sigma}$ we have
$$\langle P_{H^{\Gamma_\sigma}}(av),w\rangle_{H^{\Gamma_\sigma}}=
\langle av,w\rangle_{H^{\Gamma_\sigma}}=
\langle v,a^{-1}w\rangle_{H^{\Gamma_\sigma}}=
\langle v,w\rangle_{H^{\Gamma_\sigma}}.\qedhere
$$
\end{proof}

The following proposition is almost contained in Sarnak [Sa 1].

\begin{prop} Let $\Gamma_\sigma$ in $\S$  and let $(s_i)_{i=1}^n$ (where $n$ is 
the index $[\Gamma:\Gamma_\sigma]$) be a system of right coset representatives for $\Gamma_\sigma$
in $\Gamma$ (that is $\Gamma = \bigcup\limits_{i=1}^n s_i \Gamma_\sigma$).
Define $Q_\sigma : \V \to \V$ by the formula
$Q_\sigma v = \frac1{n}\Big( \sum\limits_{i=1}^n s_i v\Big)$, $v\in V$.

Then $Q_\sigma|_{H^{\Gamma_\sigma}}$ is the orthogonal projection from 
$H^{\Gamma_\sigma}$ onto $H^{\Gamma}$. 
\end{prop}

\begin{proof}
First, we note that indeed $Q_\sigma$ is a projection from $H^{\Gamma_\sigma}$ onto $H^{\Gamma}$.
Indeed, for all $\gamma \in \Gamma$, and for every $i$ in $\{1,2,\ldots,n\}$ there exists $\theta_i(\gamma)$
an element in $\Gamma_\sigma$ and $\pi_\gamma$ a permutation of $\{1,2,\ldots,n\}$ such that
$$
\gamma s_i = s_{\pi_\gamma(i)} \theta_\gamma(i).
$$

Hence for all $v$ in $\V^{\Gamma_\sigma}$
(by the argument in [Sa1]), for $v$ in $\V^{\Gamma_\sigma}$
$$
\gamma(Q_{\sigma} v)=\frac{1}{n}\sum_{i=1}^{n}\gamma s_1 v=
\frac{1}{n}\sum_{i=1}^{n} s_{\pi_{\gamma}(i)}\theta_i(\gamma)v=
\frac{1}{n}\sum_{i=1}^{n} s_1 v=Q_{\sigma}(v).
$$

(This holds true  since $\theta v = v$ for all $\theta$ in $\Gamma_\sigma$.) Since
$Q_\sigma$ is obviously the identity when restricted to $\V^{\Gamma}$  it follows that $Q_\sigma$ is a projection onto $H^{\Gamma}$.

The complete the proof we have to show that $Q_\sigma$ is indeed an orthogonal projection, i.e., that
the adjoint of $Q_\sigma$ is equal to $Q_\sigma$.

For $v,w$ in $H^{\Gamma_\sigma}$ we have 
\begin{gather*}
\langle Q^\sigma v,w\rangle_{H^{\Gamma_\sigma}}=
\frac{1}{n}\sum_{i=1}^{n} \langle s_i v,w\rangle_{H^{\Gamma_\sigma}}=\\
=\frac{1}{n}\sum_{i=1}^{n} \langle v, s_i^{-1}w\rangle_{H^{\Gamma_\sigma}}
=\frac{1}{n}\sum_{i=1}^{n} \langle v, P_{H^{\Gamma_\sigma}} (s_i^{-1}w)\rangle.
\end{gather*}
Hence for $w$ in $H^{\Gamma_\sigma}$
$$
Q_\sigma w = \frac{1}{n}\sum_{i=1}^{n} P^{H^{\Gamma_\sigma}} (s_i^{-1}w)
$$
and by using the notation in the previous lemma we have 
$$
(Q_\sigma)^* (w) = \frac{1}{n}\sum_{i=1}^{n} P^{H^{\Gamma_\sigma}} ([\Gamma_\sigma]s_i^{-1}w).
$$
But $\Gamma = \bigcup s_i \Gamma_\sigma$ and hence $\Gamma = \bigcup \Gamma_\sigma (s_i)^{-1}$ 
and hence we can arrange by taking appropriate representatives for the right cosets of $\Gamma_\sigma$
that 
$$
(Q_\sigma)^* (w) = \frac{1}{n}\sum P^{H^{\Gamma_\sigma}} (s_iw)=
\frac{1}{n} P^{H^{\Gamma_\sigma}} \bigg( \sum_{i=1}^{n} s_iw\bigg).
$$
Since $\sum\limits_{i=1}^{n} s_iw$ is already in $ H^\Gamma$ this is further equal
to 
$$
\sum_{i=1}^{n} s_iw = \theta_\sigma(w).
$$
Thus $Q_\sigma$ is a selfadjoint projection. We note as a consequence of the previous proof that
$P_{H^{\Gamma}}(s \sigma v)=P(\sigma v)$ for all $v$ in $H^{\Gamma}$, $s$ in $\Gamma$, $\sigma$ in $G$. 
Indeed in this case $\sigma v$ is in $H^{\Gamma_\sigma}$ and hence
$$
P_{H^\Gamma}(\sigma v) = P_{H^\Gamma}([\Gamma]\sigma v).\qedhere
$$
\end{proof}

As a corollary, we have the following equivalent description of the Hecke operator.

\begin{prop}
Fix $\sigma$ in $G$. Let $T_{[\Gamma\sigma\Gamma]}=T_\sigma : H^\Gamma \to H^\Gamma$ be the abstract Hecke operator, defined
by the formula 
$$
T_\sigma v = \sum_{i=1}^{n} s_i \sigma v,\quad v \in H^\Gamma,
$$
where $s_i$ is a system of representatives for right cosets for $\Gamma_\sigma$ in $\Gamma$
(that is $\Gamma = \bigcup s_i \Gamma_\sigma$) 

Let $P_{H^\Gamma}^{H^{\Gamma_\sigma}}$ be the orthogonal projection from 
$H^{\Gamma_\sigma}$ onto $H^{\Gamma}$ and note that $\sigma o$ belongs to $H^{\Gamma_\sigma}$.
Then $$T_{\sigma}v = [\Gamma: \Gamma_\sigma]P_{H^\Gamma}^{H^{\Gamma_\sigma}}(\sigma v) = [\Gamma: \Gamma_\sigma]
P_{H^\Gamma}^{H^{\Gamma_\sigma}} ([\Gamma\sigma\Gamma]v),$$
(where the last term of the equality is rather a notation to suggest that it doesn't depend on which
element in the coset we choose: that is $P_{H^\Gamma}^{H^{\Gamma_\sigma}}(\sigma v)=
P_{H^\Gamma}^{H^{\Gamma_\sigma}} (\sigma_1 v)$ for all $\sigma_1$ in $\Gamma\sigma\Gamma$;
also $\Gamma_{\sigma_1} = \Gamma_\sigma$ if $\sigma_1 = \gamma_1 \sigma \gamma_2$).
\end{prop}

\begin{proof} 
This is a direct consequence of the last proposition and of the remark afterwords.
\end{proof}

\begin{cor}
The composition of the arrows in the following diagram gives the Hecke operator. Let $\sigma$ in
$G$. The diagram is
$$
\begin{array}{ccc}
H^{\Gamma_{\sigma^{-1}}} & \mathop{\longrightarrow}\limits^{\sigma} &  H^{\Gamma_\sigma}\\
\phantom{aaaa}_{\rm inc}\nwarrow & \circlearrowright \phantom{a}&\hspace{-0.3cm} \swarrow {}_{P_{H_{\Gamma}}^{H^{\Gamma_\sigma}}} \\[4pt]
&  H^{\Gamma}&
\end{array}.
$$
To get the non-normalized Hecke operator we have to multiply $P_{H_{\Gamma}}^{H^{\Gamma_\sigma}}$ by\break 
 $[\Gamma: \Gamma_\sigma].
$
\end{cor}

Bellow, we present some basic examples of this construction.
The first example corresponds to the induced $C^*$-Hecke algebra ([BC]) which also assigns
a canonical norm on the Hecke algebra (the reduced $C^*$-algebra norm).

\begin{ex}
Let $\V$ consist of the function on the discrete group $G$, and let $G$ act on $\V$ by left translation. We let
$H = \ell^2(G)$ and define $H^\Gamma$ as $\ell^2(\Gamma/G) \subseteq \V^\Gamma$
(since cosets of $\Gamma$ are $\Gamma$-invariant functions).
\end{ex}

We define the $\ell^2$ norm of cosets of $[\Gamma]$ to be equal to 1, and then for smaller cosets, we renormalize that
scalar product on $\ell^2(\Gamma_\sigma\setminus G)$ by the factor 
$\frac{1}{[\Gamma:\Gamma_\sigma]}$. Hence the canonical map $\ell^2(\Gamma/ G)\hookrightarrow
 \ell^2(\Gamma_\sigma/ G)$ becomes an isometry.

In this setting $s_i(\sigma\Gamma)$ is the set $s_i\sigma\Gamma$ which decomposes as a union smaller cosets. Hence 
 for the Hecke operator we have the formula
$( T_\sigma) [\sigma_1\Gamma] = \sum [s_i \sigma\sigma_1 \Gamma]$.

This means that in this representation the Hecke operator $ T_\sigma$ coincides with the
multiplication by $[\Gamma\sigma\Gamma]$ in the Hecke algebra ([Krieg]).

Thus the $\C$-algebra generated by the Hecke operators coincides with the Hecke algebra $\H_0$ of double
cosets.
Recall ([BC]) that if $[\Gamma:\Gamma_\sigma]= [\Gamma:\Gamma_{\sigma^{-1}}]$ for all $\Gamma_\sigma$ in $\S$,
then the vector state $\langle \cdot  [\Gamma],[\Gamma]\rangle$ is a trace on $\H_0$ and the reduced 
$C^*$-Hecke algebra $\H_{\rm{red}}$ is the closure of $\H_0$ in the topology induced by the GNS construction corresponding to this state.
(Thus $\H \subseteq B(\ell^2(\Gamma/ G))$ is the weak operator topology 
closure of the $*$-algebra $\H_0$.)

Recall ([Krieg]) that in the case 
$G = \PGL_2(\Z[\frac1{p}])$, $\Gamma=\PSL_2(\Z)$ and
$\sigma_{p^n}= \left(\begin{array}{cc} p^n & 0 \\ 0 & 1 \end{array} \right)$, $n$ a positive integer, 
then if $\chi_n=[\Gamma_{\sigma_{p^n}} \Gamma]$, the cosets $\chi_n$ generate the Hecke algebra and
are selfadjoint. The relations for the elements $\chi_n$ are as follows
$$
\chi_1 \chi_n =
\left\{\begin{array}{ll}
\chi_2+(p+1){\rm Id} & \hbox{if } n=1, \\
\chi_{n+1}+p\chi_{n-1} & \hbox{if } n\geq 2.
\end{array}\right.
$$
and the value of the state $\langle \cdot [\Gamma], [\Gamma]\rangle$ on $\chi_n$ is 0 unless $n=0$, when
the value is 1.

By comparing with [Py], we see that these are exactly the relations verified by the elements of the
radial algebra of a free group with $N = \frac{p+1}{2}$ generators.

We can define polynomials $t_n(\lambda)$ by the recurrence relations above
$$
t_1{(\lambda)} t_n(\lambda) = \left\{\begin{array}{ll}
t_2(\lambda)+2N & \hbox{if } n=1, \\
t_{n+1}(\lambda)+(2N-1)t_{n-1}(\lambda) & \hbox{if } n\geq 2.
\end{array}\right.
$$

Let $\varphi_\lambda$ be the character of the $*$-algebra $\H_0$ define by requiring $\varphi_\lambda(\chi_1)=\lambda$
(and thus $\varphi_\lambda(\chi_n) = t_n(\lambda)$). 
It turns out ([Py]) that $\varphi_\lambda$ is positive for $\lambda$ in $[-2N,2N]=
[-(p+1),(p+1)]$. Moreover  if $\lambda$ is in the interval $[-2\omega,2\omega]$, where $\omega=\sqrt{p}$, then $\varphi_\lambda$
is a state on the reduced $C^*$-algebra (it is actually a positive definite function on $F_N$ and it is affiliated
with the left regular representation). Thus the spectrum of $\chi_1$ in the reduced $C^*$-algebra
is equal to $[-2\omega,2\omega] = [-2\sqrt{p},2\sqrt{p}]$ and thus $\|\chi_1\|=2\sqrt{p}$.

In particular, the norm of $[\Gamma\sigma_p\Gamma]$ in the reduced $C^*$-Hecke algebra is equal to ~$2\sqrt{p}$.

It is thus natural, in view of this example to formulate a generalized Ramanujan-Petterson conjecture 
as follows.

\begin{defn}
{\bf Generalized  Ramanujan-Petersson conjecture for an adelic representation} of a discrete group $G$, containing
an almost normal subgroup $\Gamma$, such that the subgroups $\Gamma_\sigma = \sigma \Gamma {\sigma^{-1}}
\cap\Gamma$ generate a modular family and $[\Gamma:\Gamma_\sigma]= [\Gamma:\Gamma_\sigma^{-1}]$
(and thus $[\Gamma\sigma\Gamma]=[\Gamma\sigma^{-1}\Gamma]$) for all $\sigma$ in $G$.
For all $\sigma$ in G, let $T_{[\Gamma\sigma\Gamma]}=T_\sigma$  be the corresponding Hecke operator
acting on~$H^\Gamma$.

 The claim of the conjecture is that $\|T_\sigma\| = \|[\Gamma\sigma\Gamma]\|$, where the norm of $[\Gamma\sigma\Gamma]$
is calculated in the reduced $C^*$-Hecke algebra of double cosets of $\Gamma$ in~$G$.

Equivalently, for any adelic representation of $G$ on $H^{\rm ad}$ (as in the sense of Definition 1)
the $\Gamma$-equivalent states of $G$ from this representation are weak limits of $\Gamma$-invariant
states of $G$ derived from the left regular representation of the Hecke algebra.
\end{defn}

\begin{proof} (of the equivalence of the two statements). Indeed a $\Gamma$-equivariant state of $G$
is of the form $\varphi(g)=\langle g v,v\rangle$, where $v$ is $H^\Gamma$. On the other hand,
the Hecke algebra is the center of the algebra generalized
by $G$.

Indeed, if $v,w$ are two vectors in $H^\Gamma$ such that $\langle T_\sigma v,w\rangle=0$ for all 
$\sigma$ in $G$ then $\langle [\Gamma\sigma\Gamma]v,w\rangle=0$ for all $\sigma$ and thus
 $\langle gv,w\rangle=0$ for all $g$ in $G$. 
\end{proof}

\begin{rem}
In the case of $G=\PSL_2(\Z[\frac{1}{p}])$, $\Gamma = \PSL_2(\Z)$, the positives states on $\H_0$, 
are $\varphi_\lambda$, $\lambda \in [-(p+1),(p+1)]$.

In general, a positive state on $\H_0$ is not necessary a positive state on $G$ (see [Ha])
but in the case of $\PSL_2(\Z[\frac{1}{p}])$ all such states are positive definite on $G$, and
hence cannot he excluded a priori ([Lu]).

We now describe a second example, related to operators algebra. The
essential data here is a projective
unitary representation $\pi$ (with cocycle $\varepsilon$) which extends to $G$ the left regular representation
with cocycle $\varepsilon$ of $\Gamma$, on the Hilbert space $\ell^2(\Gamma)$. We assume that $\pi$ acts  on the same Hilbert space as the left regular representation.
\end{rem}

\begin{ex} Let $G,\Gamma$ as above; $\pi$ a (projective) unitary representation of $G$ on $H=\ell^2(\Gamma)$
extending the left regular representation.
Then, let $\V = B(H)$, let $G$ act on $\V$ by ${\rm Ad}(\pi(g))$. Note that even if $\pi$ may be  a projective representation, ${\rm Ad}(\pi)$ is an actual representation. Then $\V^\Gamma=\{\pi(\Gamma)\}'\cong
\cR(\Gamma)$ (the commutant). Let $H=L^2(\L(G),\tau)\cong \ell^2(G)$. Hence $H^\Gamma= \ell^2(\L(\Gamma))\cong \ell^2(\Gamma)$
and naturally $H^{\Gamma_\sigma} = L^2(\L(\Gamma_\sigma),\tau)'$, where $\L(\Gamma_\sigma)'$ is
endowed with the normalized trace $\tau$. 
Here if the representation $\pi$ is effectively projective, then we consider the skewed version of $\L(G)$.

Then clearly, $T_\sigma = \Psi_\sigma$ is a map from $\ell^2(\Gamma)$ into $\ell^2(\Gamma)$
induced by the map on $(\L(\Gamma))'$ given by the formula:
$$
\Psi_\sigma(X)= [\Gamma:\Gamma_\sigma]E(\pi(\sigma)X\pi(\sigma)^{-1})=\sum\pi(s\sigma)X\pi(\sigma s)^{-1})
$$
for $x$ in $\L(\Gamma)$ and where $s$ runs over a system of representatives of  left cosets of $\Gamma_\sigma$ in $\Gamma$. Note that $\Phi_\sigma$  is a completely positive  map.
\end{ex}

The classical setting also fits into this pattern:

\begin{ex}
Classical setting of Hecke operators acting on Maass forms. Let $G=\PGL_2(\Z[\frac{1}{p}])$, $\Gamma = \PSL_2(\Z)$. The group $G$ acts naturally on the upper halfplane $\bH$ by Moebius transforms.  The topological vector space  $\V$ is, in this example, the space of measurable  functions on $\bH$, and
$G$ acts on a function $f$ by mapping it into $gf(z)=f(g^{-1}z), z \in \bH$ and
$H^\Gamma = L^2(F_{\Gamma},\nu_0)$, $H^{\Gamma_\sigma} = L^2(F_{\Gamma_\sigma}), \frac {1}{[\Gamma:\Gamma_0]}\nu_0)$, (where
$F_{\Gamma_\sigma}$ is a fundamental domain for the action of the discrete group
$\Gamma_\sigma$ on the upper half plane $\bH$,$\sigma \in G$).
 Here, $T_\sigma f(z) = \sum f(s_i \sigma_z),\ z\in \bH $, with $s_i$ a system of representatives of left cosets of $\Gamma_\sigma$ in $\Gamma$.
Let  $\sigma_{p^n} = \left(\begin{array}{cc} p^n & 0 \\ 0 & 1 \end{array} \right),\ n\in\N$. Then the Hecke operator, $T_{\sigma_p}(f)(z)$,  has the form
 $\sum\limits_{d=0}^{p-1} f\left(\frac{z+d}{p}\right)+f(pz),\ z \in \bH$.
\end{ex}

In the  next chapter we explain  why Example 11 is equivalent  to Example 10 in the case
$G=\PGL_2(\Z[\frac{1}{p}])$, $\Gamma = \PSL_2(\Z)$.

Of course, the Hecke operators acting on automorphic forms are another example of this
setting.

\section{Abstract Hecke operators on II$_1$ factors}

In this section we introduce the abstract Hecke operators, associated
with a pair of isomorphic subfactors, of equal indices, of a given factor $M$.

In the case $M=\L(\PSL_2(\Z), \varepsilon)$ we prove that with a suitable choice 
of the unitary implementing the isomorphism, one recovers the classical Hecke operators
acting on Maass forms. This isomorphism is based on the Berezin's quantization
of the upper half plane introduced in [Ra1], [Ra2].

First, we introduce the definition of an abstract Hecke operator.

\begin{defn} Let $M$ be a type II$_1$ factor and let $P_0$, $P_1$ be two
subfactors of finite equal indices. 

Let $\theta : P_0 \to P_1$ be a von Neumann algebras isomorphism. Let $U$
be a unitary in $\U(L^2(M))$, that implements $\theta$, that is
$U p U^* = \theta(p)$ for all $p$ in $P_0$.
Since $P_0$, $P_1$ have equal indices there always exists such a unitary, which
is unique up to left multiplication by a unitary in $P'_1$.
Then $UP'_0 U^{\ast} = P'_1$ and hence we can define $\tilde{\Psi}_U$ as the composition
of the following diagram:
$$
\begin{array}{ccc}
P'_0 & \stackrel{{\rm Ad}\, U}{\longrightarrow} &  P'_1\\[3pt]
\phantom{aaaa}_{\rm inc}\nwarrow  & & \hskip-4pt\swarrow   {}_{E}\\[3pt]
 & M'  &
 \end{array}
$$
that is $\tilde{\Psi}_U(x)= E_{M'}^{P'_1}(Ux U^*)$, $x\in M'$,
where $E=E_{M'}^{P'_1}$ is the canonical conditional expectation. We will use in the sequel the notation ${\Psi}_U=[M:P_0]\tilde{\Psi}_U$, for the non-normalized version.
\end{defn}

\begin{rem}
If $\theta$ can be extended to an automorphism $\tilde{\theta}$ of $M$,
then we can choose $U$ such that $UxU^* = \tilde\theta(x)$,
for $x \in M$ and hence in this case it follows that $UM'U^* = M'$ and
hence $\tilde{\Psi}_U(x)$ is simply $UxU^*$, $x \in M'$, that
is $\tilde{\Psi}_U$ is an automorphism of $M'$.
\end{rem}

To get a more exact description of $\tilde{\Psi}_U$ in the case of group von Neumann algebras, we need
a more precise formula for the conditional expectation $E$ in the case of
$\Gamma_1 \subseteq \Gamma$ a subgroup of a discrete group of finite index.

\begin{lemma}
Let $\Gamma$ be a discrete group and let $\Gamma_1$ be a discrete subgroup of finite index.

Let $\Gamma_1$ act on $\ell^2(\Gamma)$, and let $\L(\Gamma_1)'$ be the commutant
of $\L(\Gamma_1)$ in $B(\ell^2(\Gamma))$.
Then the conditional expectation $E_{\L(\Gamma)'}^{\L(\Gamma_1)'}$ from 
$\L(\Gamma_1)'$ onto $\L(\Gamma)'$ is defined by following formula:
choose $(s_i)^n_{i=1}$ be a system of representatives for right cosets for
$\Gamma_1$ in $\Gamma$ (that is $\Gamma = \bigcup_{i=1}^n s_i\Gamma_1$ disjointly).

Denote by $L_{s_i}$ the operator of left convolution with $s_i$ acting on $\ell^2(\Gamma)$.
Then 
$$
E_{\L(\Gamma)'}^{\L(\Gamma_1)'}(x) = \frac{1}{n}\sum_{i=1}^n L_{s_i} x L_{s_i}^*, \quad
x \in \L(\Gamma_1)'.
$$
\end{lemma}

This formula is reminiscent of the average formula in the definition of a double coset
action on Maass forms.

\begin{proof}
The lemma is certainly well known for specialists in von Neumann algebras
although we could not find a citation. For the sake of completeness we include 
the proof.

The proof is identical to the argument used for proving that Hecke
operators are mapping $\PSL_2(\Z)$-invariant functions into $\PSL_2(\Z)$-invariant functions.

For every $\gamma$ in $\Gamma$ there exists a permutation $\pi_\gamma$ of $\{1,2,\ldots,n\}$
such that
$$
\gamma {s_i} = s_{\pi_\gamma(i)} \theta_i(\gamma), \quad i=1,2,\ldots,n.
$$
Here $\pi_\gamma(i)$ is uniquely determined by the requirement that the element
$\theta_i(\gamma)=s^{-1}_{\pi_\gamma(i)}\gamma {s_i}$ belongs to $\Gamma_1$.

We denote, for $x$ in $\L(\Gamma_1)'$, by $E(x)$ the expression
$$
E(x) = \frac{1}{n} \sum_{i=1}^n L_{s_i} x L_{s_i}^*.
$$
We have that for all $\gamma$ in $\Gamma$
$$
L_\gamma E(x)L_\gamma = \frac{1}{n} \sum_{i=1}^n L_{\gamma s_i} x L_{\gamma s_i}^*
= \frac{1}{n} \sum_{i=1}^n L_{s_{\pi_\gamma(i)}}L_{\theta_i(\gamma)} x 
L_{\theta_i(\gamma)}^* L_{s_{\pi_\gamma(i)}}.
$$
Since $x$ belongs to $\L(\Gamma_1)'$, and $\theta_i(\gamma)$ belongs to $\Gamma_1$, it follows
that $L_{\theta_i(\gamma)} x L_{\theta_i(\gamma)}^*=x$ and hence that
$$
L_\gamma E(x)L_\gamma^* = \frac{1}{n} \sum_{i=1}^n L_{s_{\pi_\gamma(i)}} x L_{s_{\pi_\gamma(i)}}^*
=E(x).
$$
Hence $E(x)$ belongs to $\L(\Gamma_1)'$ for all $x$ in $\L(\Gamma_1)'$. Moreover, it is
obvious that $E$ is positive and $E(x)=x$ for $x$ in $\L(\Gamma_1)'$.
Hence $E$ is the conditional expectation 
$E_{\L(\Gamma)'}^{\L(\Gamma_1)'}$. The fact that $E$ is selfadjoint was proved in the previous chapter, in Proposition 4. This completes the proof.
\end{proof}

Using this lemma we can conclude the unitary equivalence of the abstract
Hecke operators (in the case of $\Gamma=\PSL_2(\Z)$), for a specific choice of the
unitary $U$ coming from a representation of $\PSL_2(\R)$, with the classical Hecke operators
on Maass forms. This has been observed in [Ra2], and we recall the argument for the
comfort of the reader.

The analytic discrete series $\pi_n$, $n\geq 2$ of representations of $\PSL_2(\R)$
is realized by considering the Hilbert space $H_n = H^2(\mathbb{H}, \d  \mu_n)$
of analytic square summable functions on the upper half plane $\mathbb{H}
= \{z\in\C \mid {\rm Im}\, z> 0\}$ with respect to the measure 
$\d \mu_n = ({\rm Im}\, z)^{n-2} \d \overline{z} \d z$. For $g=\begin{pmatrix}
a& b \\ c & d \end{pmatrix}$ in $\PSL_2(\R)$, with the standard action on $\mathbb{H}$,
and automorphy factor $j(g,z)=(cz + d)$, $z \in \mathbb{H}$, the formula for the representation is 
$$
\pi(g) f(z) = f(g^{-1}z) j(g,z)^{-n}, \quad f\in H_n, \, z\in\mathbb{H}.
$$

For odd $n$ this corresponds to a projective, unitary representation of 
$\PSL_2(\R)$ (the author thanks to the anonymous referee of a first submitted version of the paper, who reminded to the author this detail). We denote the 2-cocycle  corresponding to the projective representation by $\varepsilon$, and note that it only takes the values $\pm 1$.

As a particular case of the results in [GHJ], the space $H_{13}$ is unitarily equivalent
to $\ell^2(\PSL_2(\Z))$ by a unitary isomorphism that transforms $\pi_{13}(\gamma)$ for
$\gamma$ in $\Gamma = \PSL_2(\Z)$ into the unitary operator of left convolution (with cocycle $\varepsilon$) with
$\gamma$ on $\ell^2(\Gamma)$.

Another way to rephrase this is to say that the Hilbert space $H_{13}$ contains
a  cyclic vector $\xi$ for $\pi_{13}(\Gamma)$ such that $\langle \pi_{13}(\gamma) \xi, \xi \rangle$ is $0$, for $\gamma \in \Gamma$ with the exception
of the case $\gamma = {\rm e}$.

In [Ra1] we proved that the commutant $\A_{13} = \{\pi_{13}(\Gamma)\}'\subseteq B(H_{13})$
(which is thus isomorphic to $\L(\PSL_2(\Z), \varepsilon)$, the $\varepsilon$ skewed, II$_1$ factor, associated to the discrete group $\PSL_2(\Z)$), can be described as the space of bivariant
kernels $k:\bH \times \bH \to \C$ (subject to a growth condition depending on the Hilbert space) that are analytic in the first variable and
anti-analytic in the second variable, and that are $\Gamma$-bivariant, that
is $k(\gamma z, \gamma \eta) = k(z,\eta)$ for all $\gamma$ in $\Gamma$, $z,\eta \in \bH$. 
The growth condition an the kernel $k$ is requiring that $k$ generates a bounded operator $X_k$ on $H_{13}$, via
the reproducing kernel formula
$$
(X_kf)(z) = \int_{\bH} k(z,\eta) f(\eta) \d \mu_{13}(\eta)
$$
for $z\in \bH$, $f$ in $H_{13}$. It is obvious that $X_k$ commutes with 
$\{\pi_{13}(\Gamma)\}$, and thus belongs to $\A_{13}$, because of the
$\Gamma$-invariance of the kernel.

The uniform norm of $X$ is difficult to compute, but the trace in $A_{13}$ of $X_k$ ($X_k$ is
an element in the type II$_1$ factor $\A_{13}$) is given by the formula
$$
\tau_{\A_{13}}(X_k) = \frac{1}{\mu_0(F)} \int_F k(z,z)\d \mu_0(z).
$$
Hence the $L^2$-norm of $X_k$, that is $\tau_{\A_{13}}(X_k^*X_k)^{1/2}$, is
given by the formula
$$
\tau_{\A_{13}}(X_k^*X_k)^{1/2} = \frac{1}{\mu(F)}\int_\bH \int_F |k(z,\eta)|^2|\d(z,\eta)|^{13}
\d \mu_0(z)\d \mu_0(\eta).
$$
Here $d(z,\eta) = \frac{|\overline{z} - \eta|^2}{{\rm Im}\, z\ {\rm Im}\,\eta}$ for $z,\eta \in \bH$
is the cosine of the hyperbolic distance from $z$ to $\eta$.

In [Ra1] it was proven that $L^2(\A_{13},\tau)$ is isomorphic to the Hilbert space of functions
on $F$, with scalar product formula
$$
\ll f,g \gg _{13} = \langle f,B_{13}(\Delta)\rangle_{L^2(F)},
$$
(the Selberg-Berezin transform [Be]), where $B_{13}(\Delta)$ is a positive, injective, selfadjoint operator, a well determined function of  the $G$-invariant
Laplacian $\Delta$, which therefore commutes  with all the Hecke operators. In fact $B_{13}(\Delta)$ is determined by the point pair invariant function $d(z,\eta)^{13}$ on $\bH\times\bH$.

The unitary map $\Phi_{13}$ from $L^2(\A_{13},\tau)$ into the space of functions on $F$ with Hilbert scalar product induced by $B_{13}(\Delta)$, is simply the
restriction of $k$ to the diagonal. If $g$ is an element in $\PSL_2(\R)$ and $X_k$ is an
element in $\A_{13}$ represented by the kernel $k$, then $\pi_{13}(\gamma) X_k
\pi_{13}^{-1}(\gamma)$ is represented by the kernel $\alpha_g(k)$ defined by the formula
$$
\alpha_g(k)(z,\eta)=k(g^{-1}z, g^{-1}\eta), \quad z,\eta \in \bH. \leqno(1)
$$
With these identifications we can prove the following proposition

Note that we are using here representations of  $\PGL_2(\Q_+)$, the quotient of $GL_2(\Q)_+$ by
its center.

\begin{prop}\label{identification}
Let $\Gamma = \PSL_2(\Z)$. Let $\Gamma\sigma\Gamma$ in $\PGL_2(\Q)$
be a double coset of $\Gamma$ in $\PGL_2(\Q)$,
where $\sigma\in\PGL_2(\Q)$. Then the classical Hecke operator associated 
to $\sigma$, is defined, by using a system of representatives $(s_i)_{i=1}^n$
for right cosets of $\Gamma_\sigma = \Gamma \cap \sigma \Gamma \sigma^{-1}$
in $\Gamma = \PSL_2(\Z)$, by the following formula: For $f$ a $\Gamma$-invariant
function on $\bH$, 
$$
(\tilde{T}_\sigma f)(z) = \frac{1}{n} \sum_{i=1}^n f((s_i\sigma)^{-1} z), 
\quad z\in \bH.
$$

Let $\tilde{\Psi}_\sigma(x)=E^{\{\pi_{13}(\Gamma_\sigma\}'}_{\{\pi_{13}(\Gamma\}'}
(\pi_{13} (\sigma) x \pi_{13} (\sigma)^*)$
be the abstract Hecke operator associated, to $\L(\Gamma_{\sigma^{-1}},\varepsilon)$, $\L(\Gamma_{\sigma},\varepsilon)$,
and the isomorphism $\theta_\sigma(x) = \sigma x \sigma^{-1}$, $x \in \L(\Gamma_{\sigma^{-1}},\varepsilon)$,
and unitary $U_\sigma = \pi_{13}(\sigma)$.

Then $\tilde{\Psi}_\sigma$ is unitarily equivalent to $\tilde{T_\sigma}$, up to a scalar phase, 
$B_{13}(\Delta)$, on $L^2(F,\d \nu_0)$. Since $B_{13}(\Delta)$ commutes with all Hecke operators
on $L^2(F, \nu_0)$, $\tilde{\Psi}_\sigma$ and $\tilde{T}_\sigma$ have the same eigenvalues, and the
eigenvectors are the same in the correspondence given by the restriction of the Berezin's bivariant kernels representing elements in the algebra to the diagonal.
\end{prop}

\begin{proof}
For the sake of completeness we verify that $\tilde{T}_\sigma$ maps $\Gamma$-invariant functions
into $\Gamma$-invariant functions.

Since $s_i$ was a system of representatives for right cosets of $\Gamma_\sigma$ in $\Gamma$, 
that is $\Gamma=\bigcup_{i=1}^n s_i \Gamma_\sigma$ as a disjoint union,
it follows that for every $\gamma$ in $\Gamma$, there exists a permutation
$\pi_\gamma$ of $\{1,2,\ldots,n\}$ such that
$$
\gamma {s_i} = s_{\pi_{\gamma}(i)}\theta_i(\gamma),
$$
with
$$
\theta_i(\gamma) = s_{\pi_\gamma(i)}^{-1} \gamma {s_i}
$$
belonging to $\Gamma_\sigma$.

Hence for all $i=1,2,\ldots,n$ 
$$
\gamma {s_i}\sigma = s_{\pi_{\gamma}(i)}\theta_i(\gamma) \sigma = 
s_{\pi_{\gamma}(i)}\sigma (\sigma^{-1} \theta_i (\gamma)\sigma).
$$
Note that $\theta_i(\gamma)$ belongs to $\Gamma_\sigma = \Gamma \cap \sigma
\Gamma \sigma^{-1}$ and hence that $\sigma^{-1} \theta_i(\gamma) \sigma$
belongs to $\Gamma_{\sigma^{-1}}=\Gamma \cap \sigma^{-1}
\Gamma \sigma \subseteq \Gamma$.

As a consequence, if $f$ is a $\Gamma$-invariant function on $\bH$, then
for $z\in\bH$, we have 
$$
(\tilde{T}_\sigma f)(\gamma^{-1}z) =
\frac{1}{n} \sum_{i=1}^n f((s\sigma)^{-1}\gamma^{-1}z) =
\frac{1}{n} \sum_{i=1}^n f((\gamma s_i \sigma)^{-1}z) =
$$
$$
=
\frac{1}{n} \sum_{i=1}^n f((s_{\pi_{\gamma}(i)} \theta_i (\gamma) \sigma)^{-1} z)
= \frac{1}{n} \sum_{i=1}^n f((s_{\pi_{\gamma}(i)}\sigma \cdot
(\sigma^{-1} \theta_i (\gamma) \sigma))^{-1}z)=
$$
$$
= \frac{1}{n} \sum_{i=1}^n f((\sigma^{-1} \theta_i (\gamma) \sigma)^{-1}
 (s_{\pi_{\gamma}(i)}\sigma)^{-1}z),
$$
but $f$ is $\Gamma$-invariant,  $\sigma^{-1} \theta_i (\gamma) \sigma$
belongs to $\Gamma$ and hence this is equal to 
$$
\frac{1}{n} \sum_{i=1}^n f((s_{\pi_{\gamma}(i)}\sigma)^{-1}z) = \tilde{T}_\sigma f(z).
$$
Hence $\tilde{T}_\sigma f$ is a $\Gamma$-invariant function on $\bH$.

The abstract Hecke operator associated to the unitary $U_\sigma = \pi_{13}(\sigma)$
is defined for $x$ in $\{\pi_{13}(\Gamma)\}'$, by the formula
$$
\tilde{\Psi_\sigma}(x) = 
E^{\{\pi_{13}(\Gamma_\sigma\}'}_{\{\pi_{13}(\Gamma\}'}
(U_\sigma x U_\sigma^*)
= 
\sum_{i=1}^n \frac{1}{n} \pi_{13}(s_i) U_\sigma x U_\sigma^* \pi_{13}(s_i)^*,
$$
where $s_i$ are a system of right representatives
for $\Gamma_\sigma$ in $\Gamma$ (that is $\Gamma = \bigcup s_i\Gamma_\sigma)$.
Because $\pi_{13}(s_i)U_\sigma = \pi_{13}(s_i \sigma)$, if $x$ is represented by a kernel
$k$, then by formula (1), we get that $\tilde{\Psi}_{\sigma}(x)$ is represented by the
kernel
$$\frac{1}{n}\sum_{i=1}^n k(s_i\sigma)^{-1} z, (s_i\sigma)^{-1}\eta),$$ 
$z,\eta\in \bH$.
If we identify $L^2 (\A_{13},\tau)$ with the Hilbert space
$L^2(F,\d \mu_0)$ with scalar product $$\ll f,g \gg = \langle f,B_{13}(\Delta)g\rangle_{L^2(F)},
$$
then, in this identification, $\tilde{\Psi}_\sigma$ will thus map a function $f$ in $L^2(F,{\rm d}\mu_0)$ into 
the function $$\tilde{\Psi}_\sigma (f)(z) = \frac{1}{n}\sum_{i=1}^n f((\sigma s_i)^{-1}z).$$
But this is exactly the Hecke operator $\tilde{T}_\sigma$, at least as a linear map.
The structure of eigenvector, eigenvalues and the selfadjointness is unchanged by the
new scalar product, since $B_{13}(\Delta)$ has zero kernel, and as a function of the invariant Laplacian,
commutes with all Hecke operators. \end{proof}

\section{Explicit description of the abstract Hecke operator\\ in the subgroup case}

In this section we assume that $\Gamma$ is a discrete subgroup and 
let $\Gamma_0$, $\Gamma_1$ be two isomorphic  subgroups 
of equal, finite index. Let $\theta$ be an isomorphism  between 
$\Gamma_0$, $\Gamma_1$ and let $U$ be a unitary in $B(\ell^2(\Gamma))$ that implements
$\theta$ (we can always find such a unitary since the subgroups have equal
index).
For $\gamma$ in $\Gamma$ we denote by $L_\gamma$, $R_\gamma$ the operators of left
and respectively right convolution on $\ell^2(\Gamma)$ by $\Gamma$.
All the statements in the chapters that follows are also valid in the presence of a  group 2-cocycle $\epsilon$ on the group $G$, which restricts to the group $\Gamma$ (see Appendix 1). We will assume that all the partial automorphisms of $\Gamma$, Ad\ $ \sigma $, $\sigma \in G$ are $\varepsilon$ preserving. This certainly happens in the case 
 $G=\PGL_2(\Z[\frac{1}{p}])$, $\Gamma = \PSL_2(\Z)$, since $\varepsilon$ is in this case  canonical. 
 In the appendix 1 we will provide an alternative approach for the case with cocycle.
 
For $m$ in $\ell^2(\Gamma)$, we denote by $L_m$, $R_m$
the (eventually unbounded) operator of left (respectively right) convolution on 
$\ell^2(\Gamma)$ with $m$.

By $\L(\Gamma)$ and respectively $\cR(\Gamma)$, we denote the algebra of
left (respectively right) bounded convolutors on $\ell^2(\Gamma)$. $\L(\Gamma)$ is
then the type II$_1$ factor associated with $\Gamma$. When a two cocycle 
$\varepsilon$ on $\Gamma$ is given we will use instead the notation 
$\L(\Gamma,\varepsilon )$ and respectively $\cR(\Gamma, \varepsilon)$

Recall that the anti-linear involution operator $J: \ell^2(\Gamma)\to \ell^2(\Gamma)$,
defined by $J x = x^*$, $x \in \ell^2(\Gamma)$ has the property that $J \L(\Gamma)J = \cR(\Gamma)$
and $J L_m J = R_m^*$.

We have that $R_a R_b = R_{ba}$, $a,b \in \ell^2(\Gamma)$
and $\Phi(L_x) = R_{x^*} = J L_x J$ is a $*$ isomorphism from $\L(\Gamma)$ onto
$\cR(\Gamma)$. 

Moreover, for the von Neumann algebra of a group, the conjugation map ${}^{\overline{\;\cdot \;}}$ which
maps $\sum\limits_{\gamma \in \Gamma} a_\gamma \gamma$ into 
$\sum\limits_{\gamma \in \Gamma} \overline{a}_\gamma \gamma$ ($a_\gamma \in \C, \gamma \in \Gamma$), is an antilinear isomorphism 
of von Neumann algebras (from $\L(\Gamma)$ onto $\L(\Gamma)$).

Now, if $U$ is a unitary implementing $\theta$ that is, $UL_\gamma U^*=L_{\theta (\gamma)}$  for
$\gamma$ in $\Gamma_0$, we obtain an expression for 
$$
\tilde{\Psi}_U(R_x) = E^{R(\Gamma_1)}_{R(\Gamma)} (UR_x U^*).
$$

We will transfer, via the canonical antiisomorphism $\Phi(L_x) = R_{x^*}$, this map to a
completely positive map  on $\L(\Gamma)$.
The ingredients for the explicit expression of $\tilde{\Psi}_U$ are the unit vectors
$t_i = U s_i$, where $(s_i) \in \Gamma \subseteq \ell^2(\Gamma)$
$i=1,2,\ldots$, $n=[\Gamma: \Gamma_0]$, is a system
of representatives for left cosets of $\Gamma_0$ in $\Gamma$.
Since $\text {Ad\ }U$ maps $\L(\Gamma_0)$ into $\L(\Gamma_1)$ and since $\{s_i\}_{i=1}^n$ are
a Pimsner--Popa basis [PP] for $\L(\Gamma_0) \subseteq \L(\Gamma)$ it follows that
$(t_i)_{i=1}^n$ are Pimsner--Popa basis for $\L(\Gamma_1) \subseteq \L(\Gamma)$.
More precisely, this is equivalent to the fact that  $\ell^2(\Gamma)$ is the orthogonal sum of the subspaces
$\ell^2(\Gamma_1)t_i$ and to the property that $\langle \gamma_1 t_i, \gamma_2 t_i \rangle_{\ell^2(\Gamma)}$
is equal to zero unless $\gamma_1 = \gamma_2$.

The properties of $t_i$ relative to $\L(\Gamma_1)$
can be also expressed by saying that $\tau(\gamma t_i t_j^*)$ is zero unless
$i=j$ and $\gamma$ is the identity. Equivalently, $E_{\L(\Gamma_1)}(t_i t_j^*)$
is zero unless $i=j$ and in this case $E_{\L(\Gamma_1)}(t_i t_i^*)=1$.

To prove the result we need first a lemma, which gives a tool for calculating 
conditional expectations from elements in $\L(\Gamma_1)'$ onto $\cR(\Gamma) = \L(\Gamma)'$.

\begin{lemma} Let $\Gamma$ be a discrete group and let $\Gamma_1$ be a subgroup
of finite index. Let $a,b$ two vectors in $\ell^2(\Gamma_1)$, that are left $\Gamma_1$
orthonormal, that is $E_{\L(\Gamma_1)}(aa^*) = E_{\L(\Gamma_1)}(bb^*)=1$.
Fix an element $m$ in $\L(\Gamma_1)$ and consider the operator $V_{ab}^m$ acting on 
$\ell^2(\Gamma)$, with initial space $\overline{\ell^2(\Gamma_1)a}$ and range contained
in $\overline{\ell^2(\Gamma_1)b}$ given by the formula
$$
V_{ab}^m (\gamma_1a) = \gamma_1 m b.
$$
Then $V_{ab}^m$ belongs to $\L(\Gamma_1)'$ and
$
E^{\L(\Gamma_1)'}_{\L(\Gamma)'}  (V_{ab}^m) =([\Gamma:\Gamma_1])^{-1} R_{a^*mb}
$
(here the product $a^*mb$ is computed in $\L(\Gamma)$). 
\end{lemma}

\begin{proof}
Let $V_a$ (respectively $V_b$) be the partial isometries with initial
space $\ell^2(\Gamma_1)$ and range $\overline{\ell^2(\Gamma_1)a}$
and $\overline{\ell^2(\Gamma_1)b}$ respectively.

Note that $V_a$, $V_b$ are partial isometries because $a,b$ are left
orthonormal with respect to $\L(\Gamma_1)$. Indeed, the relation
$E_{\L(\Gamma_1)}(aa^*)=1$ implies that for $\gamma \in \Gamma_1$,
$\tau_{\L(\Gamma_1)}(\gamma aa^*)$ is zero unless $\gamma$ is the identity
and hence $\langle \gamma_1a,\gamma_2a\rangle_{\ell^2(\Gamma)} = 
\tau(\gamma_2^{-1} \gamma_1 aa^*)$
is zero unless $\gamma_1=\gamma_2$. Similarly for $V_b$.

If $e$ is the projection from $\ell^2(\Gamma)$ onto $\ell^2(\Gamma_1)$
then $e\in \L(\Gamma_1)'$ and 
$$V_a=R_ae \quad \hbox{and} \quad V_b = R_be.
$$
Clearly, being an isometry $V_a^*$ is the partial isometry that maps $\gamma_1 a$ 
into $\gamma_1$ for $\gamma_1$ in~$\Gamma_1$.
Consequently,
$$
V_{ab}^m = V_b R_m V_a^* = R_b e R_m e R_a^*
$$
But if we use the map on $B(\ell^2(\Gamma))$ mapping $x$  into $Jx^*J$ then
$\cR(\Gamma)$ is mapped into $\L(\Gamma)$, $\L(\Gamma_1)'$ is
mapped into $J\L(\Gamma_1)'J$ and $JeJ=e$.
The inclusion $\cR(\Gamma) \subseteq \L(\Gamma_1)'$ is mapped into the
first step of the Jones basic construction for $\L(\Gamma_1)\subseteq \L(\Gamma)$.
Hence $e$ commutes with $R_m$ and $E^{\L(\Gamma_1)'}_ {\L(\Gamma)'} (e) = \frac{1}{[\Gamma:\Gamma_1]}$.

Thus $V_{ab}^m = R_b e R_m R_a^*$ and since $R_b, R_m, R_a^*$  all belong to
$\L(\Gamma)'$, it follows that
$$
E^{\L(\Gamma_1)'}_{\L(\Gamma)'} (V_{ab}^m) = \frac{1}{[\Gamma:\Gamma_1]} R_b R_m R_a^*
$$
which is further equal to 
$$
\frac{1}{[\Gamma:\Gamma_1]} R_{a^*mb}.\qedhere
$$
\end{proof}

As an exemplification we note the following corollary, which is certainly known to specialists.
We include its proof for completeness.

\vskip6pt

{\bf Corollary.} {\it Let $t$ in $\ell^2(\Gamma)$ be left orthonormal with respect to $\Gamma_1$
(that is\break $E_{\L(\Gamma_1)} (tt^*) =1$). Let $(s_i)_{i=1}^n$ be a system of right representatives
for $\Gamma_1$ in $\Gamma$, that is $\Gamma$  is the disjoint union of $s_i\Gamma_1$.

Denote by $P_{[s_i\Gamma t]}$ the projection onto the space
$\overline{{\rm Sp}\,s_i \Gamma t}$. Then 
$$
\sum P_{[s_i\Gamma t]} = R_{t^*t}.
$$
If we use the map $J \cdot J$ we get in $\L(\Gamma)$ that in $\L(\Gamma)$,
if $r_\alpha$ is a system of representatives for left cosets of $\Gamma_1$ 
in $\Gamma$ (that is $\Gamma=\bigcup \Gamma_1 r_\alpha$) then
$$
\sum_{\a=1}^n P_{[t^* \Gamma_1 r_\a]} = t^*t.
$$
}
\vskip6pt

\begin{proof}
The projection $P_{[\Gamma_1 t]}$ clearly belongs to $\L(\Gamma_1)'$
since it is invariant to left multiplication by $\Gamma_1$. In the
terminology of the previous lemma we have~that
$$
p = P_{[\Gamma_1 t]} = V_{tt}^1
$$
and hence
$$
E^{\L(\Gamma_1)'}_{\L(\Gamma)'}(p) = \frac{1}{[\Gamma:\Gamma_1]} R_{t^*t}.
$$
Now by Lemma 3, since $s_i$ is a system of right representatives for $\Gamma_1$ in
$\Gamma$ it follows that 
$$
E^{\L(\Gamma_1)'}_{\L(\Gamma)'} (p) = \frac{1}{[\Gamma:\Gamma_1]} \sum_{i=1}^n L_{s_i}pL_{s_i}^*.
$$
But $L_{s_i}P_{[\Gamma_1 t]}L_{s_i}^* = P_{[s_i\Gamma_1 t]}$.
Hence
$$
 \sum_{i=1}^n P_{[s_i\Gamma_1 t]} = R_{t^*t}.
 $$
 If we apply the conjugation map $J \cdot J$, the space $s_i\Gamma t$ gets mapped 
 into $J(s_i\Gamma t) = t^* \Gamma s_i^{-1}$ and
 $J(R_{t^*t})J = L_{t^*t}$. But $(s_i^{-1})_{i=1}^n$ is
 a system of left representatives for $\Gamma_1$ in $\Gamma$ and
 the result follows. \end{proof}
 
 We can now prove the main result of this section, which gives a concrete
 expression for the completely positive map
 $\tilde{\Psi}_U(R_x)=E^{\L(\Gamma_1)'}_{\L(\Gamma)'} (U R_x U^*).$
We will also describe this map as an operator from $\L(\Gamma)$ into $\L(\Gamma)$.

\begin{thm} Let $\Gamma$ be a discrete subgroup and let 
$\Gamma_0, \Gamma_1$ be two isomorphic subgroups of equal finite
index. Let $\theta$ be an isomorphism from $\Gamma_0$ onto $\Gamma_1$ and assume that
$U$ is a unitary in $B(\ell^2(\Gamma))$ that implements $\theta$,
that is $U L_{\gamma_0} = L_{\theta(\gamma_0)} U$, for $\gamma_0$ in~$\Gamma_0$.

Let $\tilde{\Psi}_U : \L(\Gamma)' \to \L(\Gamma)'$ be the corresponding completely positive 
map, defined by the formula 
$$\Psi_U(x)=
[\Gamma:\Gamma_0]\tilde{\Psi}_U(x)= [\Gamma:\Gamma_0]E^{\L(\Gamma_1)'}_{ \L(\Gamma)'} (U x U^*), \quad x\in\L(\Gamma_1)'.
$$

Let $(s_i)_{i=1}^n$, with $n=[\Gamma :\Gamma_0]=[\Gamma : \Gamma_1]$
be a system of representatives for left cosets for $\Gamma_0$ in $\Gamma$, that
is $\Gamma=\bigcup \Gamma_0 s_i$. Let $t_i = U(s_i)$, $i=1,2,\ldots,n$, which
as we observed before have the property that
$E_{\L(\Gamma_1)}(t_it_j^*)=\delta_{ij}$. Then
$$
[\Gamma:\Gamma_0]\tilde{\Psi}_U(R_x)= \sum_{i,j} R_{t_i^*\theta (E_{\L(\Gamma_0)}(s_ixs_j^*)){t_j}}.
$$ 
Viewed as map from $\L(\Gamma)$ onto $\L(\Gamma)$
$($via the identification of $L_x$ with $R_{x^*}$ through $J\cdot J)$ the formula becomes
$$
\Psi_U(x)=
[\Gamma:\Gamma_0]\tilde{\Psi}_U(x)= \sum_{i,j=1}^n   t_i^*\theta (E_{\L(\Gamma_0)}(s_ixs_j^*)){t_j}, 
\quad x\in \L(\Gamma).
$$
\end{thm}

\begin{proof}
Fix $\gamma$ in $\Gamma$. We will first determine a formula for $UR_\gamma U^*$.
We use the fact $s_i$ are a system of representatives for right cosets for $\Gamma_0$
in $\Gamma$, so that $\Gamma=\bigcup \Gamma_0 s_i$.

Hence for every $\gamma$ in $\Gamma$, and $i\in\{1,2,\ldots,n\}$
there exists a permutation $\pi_\gamma$ of $1,2,\ldots,n$ and
an element $\theta_i(\gamma)$ in $\Gamma_0$ such that 
$$
s_i \gamma = \theta_i (\gamma) s_{\pi_\gamma(i)},
$$
and hence, $\theta_i (\gamma)=s_i \gamma s_{\pi_\gamma(i)}^{-1}$. One
other  way to write this expression is
$$
\theta_i (\gamma)= \sum_j E_{\L(\Gamma_0)}(s_i\gamma s_j^{-1}). \leqno(3.2)
$$
Then for an arbitrary basis element $\gamma_1 t_i$ in $\overline{\ell^2(\Gamma_1)t_i}$,
$\gamma_1\in\Gamma_1$, we have 
$$
(U R_\gamma U^*) (\gamma_1 t_i) = U R_\gamma \theta^{-1}(\gamma_1)s_i 
= U \theta^{-1}(\gamma_1)s_i \gamma = 
U \theta^{-1}(\gamma_1) \theta_i(\gamma) s_{\pi_\gamma(i)}.
$$
Since $\theta^{-1}(\gamma_1)\theta_i(\gamma)$ belongs to $\Gamma_0$ this is
further equal to
$$ 
\theta(\theta^{-1}(\gamma_1)\theta_i(\gamma))t_{\pi_\gamma(i)}=
\gamma_1 \theta(\theta_i(\gamma)) t_{\pi_\gamma(i)}.
$$
Hence
$(U R_\gamma U^*) (\gamma_1 t_i) = \gamma_1 \theta(\theta_i(\gamma)) t_{\pi_\gamma(i)}$.
With the terminology from Lem\-ma~5, it follows that the restriction of $U R_\gamma U^*$
to $\overline{\ell^2(\Gamma_1)t_i}$ is exactly  $V_{t_i,t_{\pi_\gamma(i)}}^{\theta(\theta_i(\gamma))}$,
which is a partial isometry whose initial space is exactly 
$\overline{\ell^2(\Gamma_1)t_i}$. Since the space $\overline{\ell^2(\Gamma_1)t_i}$
are pairwise orthogonal it follows that
$$U R_\gamma U^* =\sum_{i=1}^n V_{t_i,t_{\pi_\gamma}(i)}^{\theta(\theta_i(\gamma))}.$$
Hence by Lemma 5 it follows that $ E ^{\L(\Gamma_1)'} _{\L(\Gamma)'} (U R_\gamma U^*)$
is equal to the right convolutor by
$$
\frac{1}{[\Gamma:\Gamma_0]} \sum_{i} t_i^*\theta (\theta_i(\gamma))t_{\pi_\gamma(i)}.
$$
By formula (3.2), this turns out to be
$$
\frac{1}{[\Gamma:\Gamma_0]}\sum_{i,j} t_i^*\theta (E_{\L(\Gamma_0)}(s_i \gamma s_j^{-1}))t_j.
$$
By linearity it then follows that
$$ 
[\Gamma:\Gamma_0]\tilde{\Psi}_U(R_x)= \sum_{i,j}R_{t_i^*\theta (E_{\L(\Gamma_0)}
(s_ixs_j^{-1}))t_j},\quad R_x\in \cR(\Gamma).
$$
Passing from $\cR(\Gamma)$ to $\L(\Gamma)$, ($R_x$ being mapped into $L_{x^*}$)
this is then (after switching the indices $i$ and $j$) the completely positive map  taking $L_{x^*}$  into 
$$
 \sum_{i,j} L_{t_i^*\theta (E_{\L(\Gamma_0)}(s_i x^*s_j^{-1}))t_j}
$$
and thus as, a map on $\L(\Gamma)$ this is the completely positive, unital map on $\L(\Gamma)$
$$
[\Gamma:\Gamma_0]\tilde{\Psi_U}(x)=  \sum_{i,j} t_i^*\theta (E_{\L(\Gamma_0)}
(s_ixs_j^{-1}))t_j,\quad \gamma\in\Gamma.
$$
If we use the conjugation map ${}^{\overline{\;\;\;}}$ on $\L(\Gamma)$, this map becomes
$$
[\Gamma:\Gamma_0]\overline{\tilde{\Psi}}_U(\overline{x})=  \overline{t}_i^*\theta (E_{\L(\Gamma_0)}
(s_i\overline{x}s_j^{-1}))\overline{t}_j
$$ 
or
$$
[\Gamma:\Gamma_0]\overline{\tilde{\Psi}}_U(x)=  \overline{t}_i^*\theta (E_{\L(\Gamma_0)}
(s_i x s_j^{-1}))\overline{t}_j
$$ 
for $x$ in $\L(\Gamma)$. 
\end{proof}

We note here that the result in this section are in fact true in a much more general context, (see the Appendix 1) which also explains why the statements remain true in the presence of a two-cocycle on $G$.

\vskip6pt

 {\bf Remark}. {\it Let $M$ be a type II$_1$ factor with unital trace $\tau$. Let $P_0, P_1$ two subfactors of equal, finite, integer index in $M$. Assume that $U$ is a unitary in
 $B(L^2(M, \tau)$ that maps, by conjugation, the $II_1$ factor $P_0$ onto $P_1$ (that is $\text{Ad\ } U(P_0)=U(P_0)U^\ast=P_1$). Let $\tilde {\theta}$ be the automorphism from $P_0$ onto $P_1$ induced by $\text{Ad\ } U$. Let $M$ act on $L^2(M, \tau)$ and denote the commutants of the corresponding algebras by
 $M', P_0', P_1'$. Let $\tPsi_U$ be the completely positive, unital map on
 $M'$ defined by 
 $$  \tPsi_U(m')=E^{P_1'}_{M'}(m'),\quad m' \in M'.$$
 Let $s_i, i=1,2,\ldots, [M: P_0]$ be a (left) Pimsner Popa basis for $P_0$ in
 $M$ (a left orthonormal $P_0$ module basis for $M$ over $P_0$).
 Let $t_i= U(s_i), i=1,2,\ldots, [M: P_0]=[M:P_1]$ Thus $t_i, i=1,2,\ldots, [M: P_1]$ is  a (left) Pimsner Popa basis for $P_1$ in
 $M$. Then the following formula holds true for $\tPsi_U$. Let $x$ be an element in $M$ and denote by $R_x\in M'$ be the right convolutor by $x$. Then:
 $$
 [M: P_1]\tPsi_U(R_x)=\sum_{i,j=1}^{[M: P_0]} R_{t_i^\ast\tilde {\theta} (E^M_{P_0}(s_i^\ast x s_j))t_j}.
 $$
 }
 \begin{proof} This is almost contained in the previous proof. The only more general fact that is needed is that in general, if $a,b\in M$ are two $P_1$ orthonormal elements, 
 (e.g. $E_{P_1}(aa^*)=E_{P_2}(bb^*)=1$) then if $V_{a,b}\in P_1'$ is the isometry from $L^2(P_1)a$ onto $L^2(P_1)b$ mapping $p_1a$  into $p_1b$, $p_1\in P_1$, then 
 $$E^{P_1'}_{M'}(V_{a,b})=([M: P_1])^{-1}R_{a^\ast b}.\qedhere$$
 \end{proof}

 {\bf Remark.} {\it
Let $M=\L (G, \varepsilon)$, $P_0=\L(\Gamma_{\sigma^{-1}},\varepsilon)$, $P_1=\L(\Gamma_{\sigma},\varepsilon)$
  and $U=\pi(\sigma), \sigma \in G$, where $\pi$ is a projective unitary representation of $G$ on $\ell^2(\Gamma)$, extending the left regular (projective) representation of $\Gamma$. Let  $\theta=\theta_{\sigma}$ be the group morphism from $\Gamma_{\sigma^{-1}}$ onto $\Gamma_{\sigma}$ defined by $\theta(\gamma_0)=\sigma\gamma_0\sigma^{-1}$, $\gamma_0\in \Gamma_{\sigma^{-1}}$.
Let 
$$\chi(\sigma,\gamma_0)=\frac{\varepsilon(\sigma\gamma_0\sigma^{-1},\sigma)}{\varepsilon(\sigma,\gamma_0)},
\quad \sigma\in G,\ \gamma_0\in\Gamma_{\sigma^{-1}}.$$
Then $\tilde{\theta}$ is related to $\theta$ by the formula}
$$\tilde{\theta}(\gamma_0)=\chi(\sigma,\gamma_0)\theta(\gamma_0),\quad \gamma_0\in \Gamma_{\sigma^{-1}}.$$

\section{The type II$_1$ representation for the Hecke algebra\\ of a pair
$\Gamma \subseteq G$, when the regular representation\\ of $\Gamma$ may
be unitarily extended to $G$}

In this section we consider the case of an almost normal subgroup $\Gamma$
of a countable discrete group $G$.  We assume that $G$ has the property that there exists a (projective) unitary representation
$\pi:G \to \U(\ell^2(\Gamma))$ that extends the left (projective) regular representation
of $\Gamma$. In this case, as noted before, for every $\sigma$ in $G$, the
groups $\Gamma_\sigma=\Gamma \cap \sigma \Gamma\sigma^{-1}$ and 
$\Gamma_{\sigma^{-1}} = \Gamma \cap \sigma^{-1} \Gamma \sigma$ have equal indices. Let
$\H_0=\H_0(G,\Gamma)$ which we will also denote as $\H(\Gamma\setminus G /\Gamma)$, be the Hecke algebra of the pair $\Gamma \subseteq G$. All  the proofs  in this section remain valid in the presence of a group 2-cocycle on $G$, which restricts to the group $\Gamma$. We present the proofs for the case when no cocycle is present,  and refer for the general case to  Appendix 1.

We recall from [Krieg], that $\H(\Gamma \setminus G/\Gamma)$ is simply the linearization
of the algebra of double cosets of $\Gamma$ in $G$. The product formula is as follows: let
$\sigma_1,\sigma_2$ be elements of $G$
$$
[\Gamma \sigma_1 \Gamma][\Gamma \sigma_2 \Gamma] = \sum c(\sigma_1,\sigma_2,z)[\Gamma z \Gamma],
$$
where $[\Gamma z \Gamma]$ runs over the space of double cosets of $\Gamma$ contained
in $\Gamma \sigma_1 \Gamma \sigma_2\Gamma$. The multiplicity $c(\sigma_1,\sigma_2,z)$
is computed by the formula
$$
c(\sigma_1,\sigma_2,z) = \#\{\Gamma \theta_2 \mid \Gamma \theta_2 \subseteq \Gamma
\sigma_2 \Gamma \hbox{ s.t. } (\exists) \theta_1 \hbox{ in } 
\Gamma \sigma_1 \Gamma \hbox{ with } z = \theta_1 \theta_2\} \leqno(3)
$$
(see [Krieg], formula on page 15).

Moreover, $\H(\Gamma\setminus G /\Gamma)$ acts on the vector space
of left cosets $\ell^2(\Gamma/G)$, which has as a basis the set $\{\Gamma s\}$
of left cosets representatives for $\Gamma$ in $G$.

The formula of the action is for $g,h \in G$, 
$$
[\Gamma g\Gamma] [\Gamma h] = \sum_{\Gamma g_i \subseteq \Gamma g \Gamma} \Gamma g_i h.
$$
This $*$-representation is called ([BC], [CM], [Tz]) the left regular representation of the Hecke algebra on 
$\ell^2(\Gamma\setminus G)$ and is denoted by $\lambda_{\Gamma\setminus G}$.

Consequently, the above formula reads as
$$
\lambda_{\Gamma\setminus G} ([\Gamma g\Gamma])( [\Gamma h]) = 
\sum_{\Gamma g_i \subseteq \Gamma g \Gamma} \Gamma \sigma g_i h,
$$
where $\Gamma g_i$ are a system of representatives for left cosets
of $\Gamma$ that contained in~$\Gamma g\Gamma$.

The Hecke algebra comes with a natural multiplicative homeomorphism
$\ind : \H(\Gamma\setminus G /\Gamma) \to \C$ which is defined by the requirement
that
$$
\ind[\Gamma g\Gamma] = \# \hbox{ right cosets of $\Gamma$ in $\Gamma g \Gamma$} = \card
[\Gamma :\Gamma_g].
$$

The space of cosets has a natural  Hilbert space structure defined by imposing the
condition that the representatives of cosets $[\Gamma g]$, $g\in G$ 
are an orthonormal basis in $\ell^2(\Gamma\setminus G)$.

The reduced Hecke von Neumann algebra $\H$ is  the von Neumann subalgebra of $B(\ell^2(\Gamma \setminus G))$ generated 
by the left multiplication with elements in $\H(\Gamma\setminus G /\Gamma)$ (the weak closure). By $\H_{\rm red}(\Gamma\setminus G /\Gamma)$ we will denote the reduced C*-Hecke algebra which is the normic closure of  $\H(\Gamma\setminus G /\Gamma)$. 
These algebras are the weak (respectively the norm)  closure of the algebra generated by the image of $\lambda_{\Gamma \setminus G}$.
Note that this algebras come with a natural state $\varphi=\omega_{\Gamma,\Gamma}$ which
is simply
$$
\varphi(x) = \langle x[\Gamma],[\Gamma]\rangle.
$$
In particular,
$$
\varphi([\Gamma g\Gamma]^{\ast}[\Gamma g\Gamma]) = \ind[\Gamma g\Gamma]
$$
If for all $g$ in $G$, the subgroups $\Gamma_g$ and $\Gamma_{g^{-1}}$ have equal
indices in $\Gamma$ then $\varphi$ is a trace, and the reduced C$^{\ast}$ algebra $\H_{\rm red}(\Gamma\setminus G /\Gamma)$ is  obtained through the GNS construction from the trace $\varphi $ on $\H(\Gamma\setminus G /\Gamma)$.
(Note that $\H(\Gamma\setminus G /\Gamma)$ has  involution
$[\Gamma\sigma \Gamma]^* = [\Gamma \sigma^{-1}\Gamma]$ and hence the Hecke algebra is a 
$*$-algebra.)

\begin{prop}\label{continuity}

Recall that the generators of the Hecke algebra of $G = \PGL_2(\Q_+)$ over $\Gamma=\PSL_2(\Z)$
are then the cosets of the form $\alpha_{p^k} =\Gamma \sigma_{p^k} \Gamma$, with
$\sigma_{p^k}  = \begin{pmatrix} 1 & 0 \\ 0 & p^k 
\end{pmatrix}$, where $p\geq 2$ runs over the prime numbers and $k$ is a natural number. 

Consequently, the spectrum of $\alpha_p$ in the
reduced $C^*$-Hecke algebra  of $G = \PGL_2(\Q_+)$ over $\Gamma=\PSL_2(\Z)$ is  exactly $[-2\sqrt p, 2 \sqrt p]$. 

In particular if $\zeta$ is an eigenvector for the classical Hecke operator $T_p$ with eigenvalue $c_p$, let $\aleph_p$ be the corresponding character induced by $\zeta$ on the algebra 
generated by the double cosets $\alpha_{p^k}$. It follows that $c_p$ belongs to the interval 
$[-2\sqrt p, 2 \sqrt p]$ if and only if $\aleph_p$ extends to a continuous character of the $C^*$-algebra generated by the $\alpha_{p^k}, k\geq 1$ in the reduced $C^* $-Hecke algebra $\H_{\rm red}$.
\end{prop}

\begin{proof}
 Fix $p\geq 3$ be a prime number.  Let $N=(p-1)/2$ and let $F_N$  be the free group with N generators. Let $\chi_k\in\L(F_N) $ be the sum of words of length $k$, $k\geq 0$. It is proved in [Py] that the algebra generated by the selfadjoint elements 
$\chi_k, k\geq 0$ is abelian, and that the spectrum of $\chi_1$ is exactly $[-2\sqrt p, 2 \sqrt p]$. Moreover the recurrence relations  for $\chi_k$ are the same as the one for $\alpha_{p^k} $, and hence we have 
an algebra morphism mapping $\alpha_{p^k} $ into $\chi_k$. Since this morphism is trace preserving, we actually obtain an isomorphism of $C^*$-algebras. It is easily seen that this is also valid for $p=2$.
\end{proof}

We can now state the main result of this section. In particular, 
this proves that if $G$ has a unitary representation on $\ell^2(\Gamma)$ that extends
the left regular representation, then $\H_{\rm red}(\Gamma\setminus G /\Gamma)$ and $\H$ embeds in
a natural way into~$\L(G)$.

\begin{thm}\label{heckerep}
Let $G$ be a discrete group with an almost normal subgroup $\Gamma$. Assume
that $G$ admits a unitary representation $\pi$ on $\ell^2(\Gamma)$ that extends
the left regular representation of $\Gamma$ on $\ell^2(\Gamma)$. Let $e$ be the neutral element of $\Gamma$,
viewed as on element of the Hilbert space $\ell^2(\Gamma)$.

For $\theta$ in $G$, we use the scalar product on $\ell^2(\Gamma)$ to define
$$
t(\theta) = \langle \pi(\theta) e,e\rangle.
$$
This is a specific matrix coefficient of the representation $\pi$.

For $\a = [\Gamma \sigma \Gamma]$ a double coset in $\H(\Gamma\setminus G /\Gamma)$ 
define
$$
t^\a = \sum_{\theta\in\a} t(\theta)\cdot \theta.
$$
Then $t^\alpha$ is an element of  $\ell^2(\Gamma g \Gamma)\subseteq  \ell^2(G)$
and the map $\rho$
$$
\a \to t^\a, \quad \a=[\Gamma g \Gamma] \in \H(\Gamma\setminus G /\Gamma)
$$
extends by linearity and continuity to a unital $*$ normal isomorphism $\rho$ from the von Neumann algebra
$\overline{\H_{\rm red}(\Gamma\setminus G /\Gamma)}^w=\H$
into $\L(G)$. The restriction of $\tau_{\L(G)}$ to the image of $\H$ correspond to the
state $\omega_{[\Gamma],[\Gamma]}$ on the Hecke algebra.

For $c=[\Gamma s]$ a coset in $\ell^2(\Gamma/G)$, define
$$
t^c = \sum_{\theta \in \Gamma s} t(\theta) \theta.
$$
Then $t^c \in \ell^2(\Gamma s)$, and the family $t^c$, where $c$ runs
over the space of left cosets of $\Gamma$ in $G$ is an orthonormal system generating a Hilbert space
$K$.
Then $K$ is a reducing space for the representation $\rho$. The restriction of the
representation $\rho$ of $\H$ to $K$ is unitarily equivalent to the left representation
$\lambda_{\Gamma\setminus G}$ of $\H_{\rm red}(\Gamma\setminus G /\Gamma)$ on $\ell^2(\Gamma\setminus G)$,
by the unitary that maps $t^c$ into the coset $c\in \ell^2(\Gamma \setminus G)$.

As explained in the Appendix 1, this construction obviously extends to the case of a projective representation of $G$.

$($Note that by replacing $t(\theta)$ by $\overline{t(\theta)}$, for all $\theta\in G$, the results remain valid, since this corresponds to taking  the conjugation map on the group algebra$)$.
\end{thm}

\begin{rem}
Note that, in particular, the theorem implies that the following properties
hold true.

For all $a_1=[\Gamma \sigma_1 \Gamma]$, $a_2=[\Gamma \sigma_2 \Gamma]$ 
double cosets of $\Gamma$ in $G$

a) $t^{a_1} t^{a_2} = \sum_{\Gamma z\Gamma \subseteq \Gamma \sigma_1 \Gamma \sigma_2 \Gamma}
c(a_1,a_2,z)t^{\Gamma z\Gamma}$.

b) For all double cosets $\Gamma\sigma\Gamma$ we have 
$$
(t^{\Gamma\sigma\Gamma})^* = t^{\Gamma\sigma^{-1}\Gamma}.$$

c) If $a=[\Gamma\sigma\Gamma]$, and $c = [\Gamma s]$ is a 
coset then
$$
t^a \cdot t^c = \sum_{\Gamma g_i \subseteq \Gamma\sigma\Gamma} t^{[\Gamma g_i s]},
$$
where $\Gamma g_i$ runs over a set of representatives for left cosets of $\Gamma$ that 
are contained in $\Gamma\sigma\Gamma$.

d) For every coset $c = \Gamma s$, $\|t^c\|_2^2 =1$ and $\{t^c\}$, where $c$ runs
over cosets of $\Gamma$, is an orthonormal basis.

Moreover, the following additional properties 1) through 9) hold true.

\ 

1) $\|t^{\Gamma\sigma\Gamma}\|_2^2 = \tau((t^{\Gamma\sigma\Gamma})^* t^{\Gamma\sigma\Gamma})
=
\ind [\Gamma\sigma\Gamma] .$

2) If $a_1 = [\Gamma\sigma_1\Gamma]$, $a_2=[\Gamma\sigma_2\Gamma]$ are two different
double cosets, then for all $\gamma$ in $\Gamma$
$$E ^{\L(G)}_{\L(\Gamma)} (t^{a_1} \gamma t^{a_2})=0.$$
(In particular, $t^{a_1}, t^{a_2}$ are orthogonal.)

3) If $a=[\Gamma\sigma\Gamma]$ then 
$$
E^{\L(G)}_{\L(\Gamma)} (t^a (t^a)^*) = \ind a.
$$

4) For all $\xi, \eta$ in $\ell^2(\Gamma)$ and $a=[\Gamma\sigma\Gamma]$
$$
\eta t^a \xi^* = \sum_{\theta\in [\Gamma\sigma\Gamma]} \langle \pi (\theta)
\overline{\xi}, \overline{\eta} \rangle \theta,
$$
where $\overline{\xi}, \overline{\eta}$ are the images of $\xi, \eta$ to
the conjugation map $\overline{\sum \xi_\gamma \gamma} = \sum \overline{\xi_\gamma} \gamma$.

5) If $s_i$ is a system of representatives for right cosets $\Gamma_{\sigma^{-1}}$ in $\Gamma$, 
so that $\Gamma\sigma\Gamma$ is as a set the disjoint union of $\Gamma \sigma s_i$ (since
$\Gamma = \bigcup \Gamma_{\sigma^{-1}} s_i$) then
$$
t^{\Gamma\sigma\Gamma} = \sum_{i=1}^{[\Gamma:\Gamma_{\sigma^{-1}}]} t^{\Gamma\sigma s_i}.
$$

6) If $\sigma$ in $G$ commutes with $\Gamma$, then $t^{\Gamma\sigma\Gamma}$ is
simply a multiple of $\sigma$ as an element of $\L(G)\subseteq \ell^2(G)$.

7) The representation $\pi$ can be recovered from the coefficients $t(\theta)$, $\theta$ in $G$.
Indeed, for all $\theta$ in $G$, $\gamma$ in $\Gamma$
$$
\pi(\theta)\gamma = \sum_{\gamma_1} t(\gamma_1^{-1} \theta \gamma) \gamma_1.
$$

In particular, $\pi(\sigma)e$ as an element of $\ell^2(\Gamma)$ is equal to
$\sigma \cdot \overline{t^{\sigma^{-1} \Gamma}}$ and hence
$$
(\pi (\sigma)e)^* = \overline{t^{\Gamma \sigma}} \cdot \sigma^{-1}.
$$

Recall that if $x = \sum x_\gamma \gamma$ is an element of $\L(\Gamma)$, then
$\overline{x} = \sum \overline{x}_\gamma \gamma$.

8) Let $\Gamma s, \Gamma t$ be two left cosets of $\Gamma$ in $G$. Let $A_{\Gamma s,\Gamma t}$ be the subset
of $\Gamma$ defined by $A_{\Gamma s,\Gamma t} = \Gamma \cap s^{-1} \Gamma t$.

Let $\a_{\Gamma s,\Gamma t}$ be the projection from $\ell^2(\Gamma)$ onto the Hilbert space
generated by the elements in $A_{\Gamma s,\Gamma t}$. 
In particular, $\gamma$ belongs to $A_{\Gamma s,\Gamma t}$ is equivalent
to $\a_{\Gamma s,\Gamma t}(\gamma) \neq 0$ (and hence $\a_{\Gamma s,\Gamma t}(\gamma)=\gamma$)
and this is further equivalent to the fact that there exist $\theta$ in $\Gamma$ such that
$$
s \gamma = \theta t \quad (\gamma = s^{-1} \theta t).
$$
Then, for $x$ in $\L(\Gamma)$, 
$$E ^{\L(G)} _{\L(\Gamma)} (t^{\Gamma s} x (t^{\Gamma t})^*) = t^{\Gamma s}
\a_{\Gamma s,\Gamma t}(x) (t^{\Gamma t})^*.
$$

9)Let $\Gamma\sigma\Gamma$ be a double coset in $G$, and let $(s_i)_{i=1}^{[\Gamma:\Gamma_{\sigma^{-1}}]}$ 
be a set of representatives for left $\Gamma_{\sigma^{-1}}$ cosets of $\Gamma_{\sigma^{-1}}$
in $\Gamma$ (that is $\Gamma = \bigcup \Gamma_{\sigma^{-1}} s_i$, so that $\Gamma\sigma\Gamma
= \bigcup \Gamma \sigma s_i$). For $\gamma$ in $\Gamma$, let $\pi_\gamma$ be the permutation of
$\{1,2,\ldots, [\Gamma:\Gamma_{\sigma^{-1}}]\}$ defined by 
the requirement that for $i$ in $\{1,2,\ldots, [\Gamma:\Gamma_{\sigma}]\}$, $\pi_\gamma(i)$
is the unique element of $\{1,2,\ldots, [\Gamma:\Gamma_{\sigma}]\}$, such that 
there exists $\theta$ in $\Gamma_{\sigma^{-1}}$ with $s_i \gamma = \theta s_{\pi_\gamma(i)}$
(in particular, $\theta = s_1 \gamma s^{-1}_{\pi_\gamma(i)} \in \Gamma_{\sigma^{-1}}$).

Then
$$
[\Gamma:\Gamma_{\sigma}] E ^{\L(G)}_{\L(\Gamma)} (t^{\Gamma \sigma s_i} \gamma(t^{\Gamma \sigma s_j})^*)=
t^{\Gamma \sigma s_i} \a_{\Gamma \sigma s_i,\Gamma \sigma s_j} (\gamma)
(t^{\Gamma \sigma s_j})^*
$$
is different from $0$, if and only if $j = \pi_\gamma(i)$, in which case it is equal to
$$t^{\Gamma \sigma s_i} \gamma(t^{\Gamma \sigma s_{\pi_\gamma(i)}})^*.$$
This is equivalent to fact that $\gamma$ belongs to 
$A_{\Gamma \sigma s_i,\Gamma \sigma s_j}$
which is equivalent to the fact that there exists $\theta$ in $\Gamma$ such that
$(\sigma s_i)\gamma = \theta(\sigma s_j)$.
\end{rem}

To prove the remark, we will first prove the following lemma, which is the main computational 
tool for all these equalities.

\begin{lemma}
For all $\theta_1, \theta_2$ in $G$ the following equality holds:

{\rm 1)} $\overline{t(\theta_1)} = t(\theta_1^{-1})$;

{\rm 2)} $\sum_{\gamma \in \Gamma} t(\theta_1 \gamma) t(\gamma^{-1} \theta_2) = t(\theta_1\theta_2)$.
\end{lemma}

\begin{proof}
Clearly 
$$\overline{t(\theta_1)} = \overline{\langle \pi(\theta_1)e,e \rangle_{\ell^2(\Gamma)}} =
\overline{\langle e, \pi(\theta_1^{-1})e\rangle} =
\langle \pi(\theta_1^{-1})e,e \rangle = t(\theta_1^{-1}).
$$
To prove the second property note that 
$$
t(\theta_2\theta_1) = \langle \pi(\theta_1)e, \pi(\theta_2^{-1})e \rangle
$$
which by property 7) (that we will prove below) is 
$$
\Big\langle \sum_{\gamma_1} t(\gamma_1^{-1} \theta_1)\gamma_1,
\sum_{\gamma_2} t(\gamma_2^{-1} \theta_2)\gamma_2 \Big\rangle
= 
$$
$$
=
\sum_{\gamma_1} t(\gamma_1^{-1} \theta_1)
\overline{t(\gamma_1^{-1} \theta_2^{-1})}=
\sum_{\gamma} t(\theta_2\gamma_1)t(\gamma_1^{-1} \theta_1).
$$

\end{proof}

The proof of property 7) is as follows: Fix $\theta\in G$, $\gamma\in\Gamma$. Then
$$
\pi(\theta) \gamma = \sum_{\gamma_1}\langle \pi(\theta)\gamma,\gamma_1 \rangle \gamma_1 =
\sum_{\gamma_1} \langle \pi(\gamma^{-1}_1 \theta \gamma)e,e \rangle \gamma_1
= \sum t(\gamma^{-1}_1 \theta \gamma)\gamma_1.
$$

We now start the proof of Theorem 7.

The most relevant properties are a), c) that we will prove first.

To prove property a) let $a_1 = [\Gamma \sigma_1 \Gamma]$, $a_2 = [\Gamma \sigma_2 \Gamma]$
be two double cosets in $\H(\Gamma\setminus G/\Gamma)$. Then
$$
t^{a_1}\cdot t^{a_2} = \sum_{\theta_1 \in \Gamma\sigma_1\Gamma \atop
\theta_2 \in \Gamma\sigma_2\Gamma}
t(\theta_1)t(\theta_2) \cdot \theta_1\theta_2
$$
and hence this is equal to
$$
\sum_{z\in\Gamma\sigma_1\Gamma\sigma_2\Gamma} z\bigg(
\sum_{\theta_1 \in \Gamma\sigma_1\Gamma, \,
\theta_2 \in \Gamma\sigma_2\Gamma \atop \theta_1\theta_2 = z}
t(\theta_1)t(\theta_2) \bigg).
$$
To identify the coefficient 
$$
\sum_{\theta_1 \in \Gamma\sigma_1\Gamma, \,
\theta_2 \in \Gamma\sigma_2\Gamma \atop \theta_1\theta_2 = z}
t(\theta_1)t(\theta_2) \leqno(4)
$$
for any $z\in G$, that belongs to $\Gamma\sigma_1\Gamma\sigma_2\Gamma$, we consider
$$
A_z = \{(\theta_1,\theta_2) \in \Gamma\sigma_1\Gamma \times \Gamma\sigma_2\Gamma \mid
\theta_1\theta_2 = z\}.
$$

Clearly, the group $\Gamma$ acts on $A_z$, the action of $\gamma$ on an element
$(\theta_1,\theta_2)$ being
$$
\gamma(\theta_1,\theta_2) = (\theta_1 \gamma^{-1}, \gamma \theta_2).
$$
It is obvious that this is a free action of $\Gamma$. Let $\O$ be the space
of orbits of $\Gamma$. Each orbit is of the form $\{(\theta_1 \gamma^{-1}, \gamma \theta_2)
\mid \gamma \in \Gamma\}$,
with the action of $\Gamma$ being bijective.
It follows by property 2) of Lemma 9 that for every orbit $o$ in $\O$
$$
\sum_{(\theta_1,\theta_2)\in o} t(\theta_1)t(\theta_2) = t(z).
$$
Hence the coefficient in formula (4) is $n(z)t(z)$, where $n(z)$ is the
number of orbits of $\Gamma$ for the given action on $A_z$.

We consider the following map $\Phi$ from $\O$ into the space of cosets of $\Gamma$
in $G$. 
If $o\in \O$, is defined as $o=\{(\theta_1\gamma,\gamma^{-1}\theta_2)\mid \gamma \in \Gamma\} \subseteq
A_z$
for some $\theta_1 \in \Gamma \sigma_1 \Gamma$, $\theta_2 \in \Gamma \sigma_2 \Gamma$,
(with the necessary property that $\theta_1\theta_2 = z$) then we define
$$
\Phi(o)=\Gamma \theta_2.
$$

Clearly, this map is well defined. 

Moreover, the image lies in the set $M=M(\sigma_1,\sigma_2, z)$ of cosets of $\Gamma y$ in $G$
that verify that there exists $x$ in $\Gamma \sigma_1 \Gamma$ with $xy =z$. (This is
the set defining the coefficient $c(\sigma_1, \sigma_2, z)$ in formula (3).

Now, clearly $\Phi$ is injective since if $o' = \{(\theta'_1,\gamma^{-1},\gamma\theta'_2) \mid 
\gamma \in \Gamma\}$ is another orbit in $A_z$, such that
$$
\Phi(o') = \Phi (o)
$$
then it follows that
$$
\Gamma \theta_2 = \Gamma \theta'_2.
$$
But this implies $\theta'_2=\gamma_0 \theta_2$ for some $\gamma_0$ in $\Gamma$.
Since $\theta_1\theta_2=\theta'_1\theta'_2=z$, this
implies that $\theta'_1 = \theta'_2 \gamma_0^{-1}$
and hence that $o$ and $o'$ are the same orbit.

Thus the number $n(z)$ in formula (4) is $c(\sigma_1, \sigma_2, z)$, and
since this only depends of the double coset of 
$\Gamma z \Gamma$ and not of the individual value of $z$,
this proves that in the product
$t^at^c$ the element $t^{\Gamma z \Gamma}$ appears with coefficient
$c(\sigma_1, \sigma_2, z)$.

This completes the proof of property a).

We now prove property c). Let $a=\Gamma\sigma\Gamma$ be a double coset and let 
$c=\Gamma s$ be a left coset of $\Gamma$ in $G$. We want to determine $t^at^c$.
Then
$$
t^at^c = \sum_{\theta\in \Gamma\sigma\Gamma,\, g\in\Gamma s} 
t(\theta)t(g)\theta g
=\sum_{z \in \Gamma \sigma \Gamma\Gamma g}
z \bigg( \sum_{\theta\in \Gamma \sigma \Gamma, \, g\in \Gamma s\atop
\theta g = z} t(\theta)t(g)\bigg). \leqno(5)
$$

Let $(r_a)_{a=1}^n$, with $n =[\Gamma:\Gamma_{\sigma^{-1}}]$,
be a set of representatives for right cosets of $\Gamma_{\sigma^{-1}}$ in $\Gamma$.
Then $\Gamma = \bigcup_{a=1}^n \Gamma_{\sigma^{-1}}r_a$ as
a disjoint union. Since $\sigma \Gamma_{\sigma^{-1}}\sigma^{-1} = 
\Gamma_\sigma \subseteq \Gamma$
it follows that 
$$
\Gamma\sigma \Gamma = \bigcup_a \Gamma\sigma \Gamma_{\sigma^{-1}} r_a
= \bigcup_{a=1}^n \Gamma_{\sigma} r_a.
$$

Clearly, this is also a disjoint union, since if $\gamma_1 \sigma r_a = \gamma_2 \sigma r_b$
with $\gamma_1, \gamma_2$ in $\Gamma$, then it follows that
$$
r_b r_a^{-1} = \sigma^{-1}(\gamma_2^{-1}\gamma_1)\sigma
$$
and hence since $r_a r_b^{-1}$ belongs to $\Gamma$.
It follows that $\sigma^{-1} (\gamma_2^{-1}\gamma_1)\sigma$
belongs to $\sigma^{-1} \Gamma \sigma \cap \Gamma$. Hence
$r_b r_a^{-1}$ belongs to $\Gamma_{\sigma^{-1}}$ or
$r_b$ belongs to $\Gamma_{\sigma^{-1}}r_a$. But this implies
$r_a=r_b$, since the $r'_as$ were a set of representatives. We decompose the set
$\Gamma \theta \Gamma \times \Gamma s$ as the reunion
$
\bigcup_{a=1,2,\ldots,n} \bigcup_{\gamma_1\in \Gamma} A_{\gamma_1, a},
$
where $A_{\gamma_1, a}$ is the set  $\{(\gamma_1 \sigma r_a \gamma, \gamma^{-1}s)\mid \gamma \in \Gamma\}$.
Note that the sets $A_{\gamma_1, a}$ are disjoint.

Indeed, if $A_{\gamma_1, a} \cap A_{\gamma_2, b} \neq \emptyset$,
then there exists $\gamma',\gamma'' \in \Gamma$ such that 
$$
(\gamma_1 \sigma r_a \gamma', (\gamma')^{-1} s) = 
(\gamma_2 \sigma r_b \gamma'', (\gamma'')^{-1} s) 
$$
but this implies that $\gamma'=\gamma''$ and hence this implies that
$$
\gamma_1 \sigma r_a =\gamma_2 \sigma r_b.
$$
Since as we have shown before the union
$\Gamma = \bigcup_{c=1}^n \Gamma \sigma r_c$ is disjoint
it follows that $r_a=r_b$ and hence that $\gamma_1=\gamma_2$.

By formula (5) we thus have

 $$
 t^at^c=\sum_{\gamma_1,a}\sum_{\gamma\in\Gamma}t(\gamma_1 \sigma r_a\gamma)t(\gamma^{-1}s)
 \gamma_1 \sigma r_a\gamma\gamma^{-1}s.
 $$
By Lemma 9, this is further equal to
$$
\sum_{\gamma_1,a}t(\gamma_1 \sigma r_as)\gamma_1\sigma r_as
=\sum_{a}\bigg(\sum_{\gamma_1}t(\gamma_1 \sigma r_as)\gamma_1\sigma r_as\bigg)
=\sum_{a}t^{\Gamma\sigma r_as}
$$
which is exactly
$$
\sum_{\Gamma_z \subseteq \Gamma \sigma \Gamma} t^{\Gamma z s},
$$
where the sum runs over right cosets of $\Gamma$ contained in 
$\Gamma \sigma \Gamma$.

We now prove property d) in Remark 8.
Let $c=\Gamma s$, $d=\Gamma t$ be two cosets of $\Gamma$
in $G$.

Then $\langle t^c,t^d \rangle_{\ell^2(G)}$ is equal to
$$
\Big\langle \sum_{\gamma_1\in \Gamma} t(\gamma_1 s) \gamma_1 s, 
\sum_{\gamma_2\in \Gamma} t(\gamma_2 t) \gamma_2 t\Big\rangle_{\ell^2(G)}
=  \sum_{\gamma_1,\gamma_2\in \Gamma} t(\gamma_1 s) 
\overline{t(\gamma_2 t)} \langle \gamma_1 s, \gamma_2 t \rangle.
$$

If the cosets $\Gamma s$ and $\Gamma t$ are disjoint then
this is clearly 0. Otherwise, if $s=t$ then this is further equal to
$$
\sum_{\gamma_1\in \Gamma} t(\gamma_1 s) \overline{t(\gamma_1 s)}
= \sum_{\gamma\in \Gamma} t(\gamma_1 s)t(s^{-1}\gamma_1)=
\sum_{\gamma\in \Gamma} t(s^{-1}\gamma_1)t(\gamma_1 s)=
t(s^{-1} s)=1
$$
again by Lemma 9.

This completes the proof of properties a), b), c), d) from Remark 8.

We now proceed to the proof of Theorem 7.

By properties a), b) it is then obvious that the map $\Phi$ from $\H(\Gamma\setminus G/\Gamma)$
into $\L(G)$ defined by $\Phi([\Gamma\sigma\Gamma])=t^{\Gamma\sigma\Gamma}$
and then extended by linearity is $*$ homeomorphism.

Because of properties c), d) the map $V$ which maps $t^{\Gamma s}$ into the
coset $\Gamma s$ in $\ell^2(\Gamma\setminus G)$ is a unitary operator.
Moreover, $\Phi(\H(\Gamma\setminus G/\Gamma))$ invariates
$K$, so the projection $P_K$ from $\ell^2(G)$ onto $K$ belongs to the commutant
of the algebra $\H_0 = \Phi(\H(\Gamma\setminus G/\Gamma))$.

Moreover, by property d) and because of the definition of the left action
$\lambda_{\Gamma/G}$ of $\H(\Gamma\setminus G/\Gamma)$ on 
$\ell^2(\Gamma\setminus G)$ it follows that
$$
U(P\Phi(a)P)U^* = \lambda_{\Gamma/G}(a)
$$
for all $a$ in $\H(\Gamma\setminus G/\Gamma)$.

Moreover, $$\tau_{\L(G)}(\Phi(a))=\omega_{\Gamma,\Gamma}(\lambda_{\Gamma/G}(a)) =
\langle \Phi(a)e,e\rangle =  
\langle P\Phi(a)Pe,e\rangle.$$
Here we use the fact that $e$ (the unit of $\Gamma\subseteq G$) belongs to $K$, as $t^\Gamma=e$.

To conclude the fact that $\Phi$ is an isomorphism from 
$\H_{\rm red}(\Gamma\setminus G/\Gamma)$ into $\L(G)$ we need the following lemmata
that summarizes the properties we obtained so far

\vskip6pt

{\bf Lemmata.} {\it Let $M$ be a finite von Neumann algebra with finite
faithful trace $\tau$. Let $N_0$ be a unital $*$-subalgebra of $M$ that contains
the unit. Assume that there exist a projection $P$ onto a subspace $K$
of $L^2(M,\tau)$ that contains $1$, and such that $P$ commutes with $N_0$.
Let $B_0 = P_0 N_0 P_0$, and let $B$ be the von Neumann algebra
generated by $B_0$ in $B(K)$. Assume that $\tau=\omega_{1,1}$ is a faithful  
state on $B$. 

Then the reduction map, which maps $n_0 \in N_0$ into $p_0n_0p_0$ extends
to a von Neumann algebra isomorphism from $N=\{N_0\}''$ onto $B$.
}

\begin{proof}
Indeed $\Phi$ becomes a unitary from $L^2(N,\tau)$ onto $L^2(B,\omega_{1, 1})$
which then implements the isomorphism from $N$ onto $B$.

This concludes the proof of the fact that $\Phi: \H(\Gamma\setminus G/\Gamma)$
extends to a von Neumann algebras isomorphism, from 
$\H_{\rm red}(\Gamma\setminus G/\Gamma)$ into $\H=\{\Phi(\H(\Gamma\setminus G/\Gamma))''\}$
because of the unitary $U$ that intertwines the left regular representation of $\H$
with the restriction of $\Phi$ to $\ell^2(\Gamma\setminus H)$ (which generates 
$\H_{\rm red}(\Gamma\setminus G/\Gamma)$ with the representation $a \to P_K \Phi(a)P_K$.
This concludes the proof of Theorem 7. $\hfill\Box$

We now proceed to the proof of the properties 1)--8) in Remark 8 (since 7) was already
proven).

We start with property 2). Assume that $a_1=\Gamma \sigma_1 \Gamma$, 
$a_2=\Gamma \sigma_2 \Gamma$ are two double cosets such that
$E ^{\L(G)}_{\L(\Gamma)} (t^{a_1} \gamma(t^{a_2})^*)$
is different from 0 for some $\gamma$ in $\Gamma$.

The terms in $t^a \gamma (t^{a_2})^*$ are sums of multiples
of elements of the form\break
$(\gamma_1 \sigma_1 \gamma_2)\gamma(\gamma_3 \sigma_2^{-1} \gamma_4)$, with
$\gamma_1, \gamma_2,\gamma_3,\gamma_4 \neq 0$ and hence if
$E ^{\L(G)}_{\L(\Gamma)} (t^{a_1} \gamma(t^{a_2})^*)$
is different from 0, it follows that there exists 
$\gamma_1, \gamma_2,\gamma_3,\gamma_4$ and $\theta$ in $\Gamma$
such that
$$
(\gamma_1 \sigma_1 \gamma_2)\gamma(\gamma_3 \sigma_2^{-1} \gamma_4) = \theta.
$$
Hence $\sigma_2 = (\gamma_4^{-1} \theta^{-1} \gamma_1)\sigma_1(\gamma_2 \gamma \gamma_3)$
and hence $\Gamma\sigma_2 \Gamma = \Gamma\sigma_1\Gamma$ or $a_1=a_2$.

This proves property 3) and also proves that
$$
\tau_{\L(G)}((t^{\Gamma\sigma_1\Gamma})^*(t^{\Gamma\sigma_2\Gamma}))=0
$$
if $\Gamma\sigma_1 \Gamma \neq \Gamma\sigma_2\Gamma$.

To prove the remaining part of property 1), note that by property 5)
(which is obvious since the sets $\Gamma\sigma s_i$ are disjoint) we have that 
$$
t^{\Gamma\sigma\Gamma} = \sum_{\Gamma z \subseteq \Gamma \sigma \Gamma} t^{\Gamma z}.
$$
Since we know that for different cosets $\Gamma z_1, \Gamma z_2, t^{\Gamma z_1}$ and
$t^{\Gamma z_2}$ are orthogonal, it follows that
$$
\|t^{\Gamma\sigma\Gamma}\|_2^2 = \tau_{\L(G)} (t^{\Gamma\sigma\Gamma}t^{\Gamma\sigma\Gamma})
= \sum_{\Gamma z \subseteq \Gamma\sigma\Gamma} \langle \Gamma z,\Gamma z \rangle
= \sum_{\Gamma z \subseteq \Gamma\sigma\Gamma}1,
$$
and this is exactly the number of left cosets in $\Gamma\sigma\Gamma$, 
which is $\ind[\Gamma\sigma\Gamma]$.

Property 3) is now a consequence of property 1). Indeed, as we have
proven in property 1), for every double coset $a=[\Gamma\sigma\Gamma]$,
we have that $t^a(t^a)^*$ is the sum
$$
\sum_{\Gamma z \Gamma\subseteq \Gamma \sigma \Gamma \sigma^{-1} \Gamma }
c(\sigma,\sigma^{-1}, z) t^{\Gamma z \Gamma}.
$$
Hence 
$$
E ^{\L(G)}_{\L(\Gamma)} (t^a (t^a)^*)=c(\sigma,\sigma^{-1},e) t^{\Gamma e}.
$$
But $t^{\Gamma e}$ is just the identity.

On the other hand, if we apply the trace $\tau$ into the previous relation,
and since $E$ preserves the trace it follows that
$$\ind a = \tau (t^a (t^a)^*) = \tau(
E ^{\L(G)}_{\L(\Gamma)} (t^a (t^a)^*)=c(\sigma,\sigma^{-1},e) .
$$
Hence $c(\sigma,\sigma^{-1},e)=\ind a$ where $a=[\Gamma\sigma\Gamma]$
and hence
$$
E ^{\L(G)}_{\L(\Gamma)} (t^a (t^a)^*)= (\ind a).
$$

We now proceed to the proof of property 4). By bilinearity it is sufficient to prove this
property for $\xi = h_1$, $\eta=h_2$, where $h_1,h_2$ are two elements in $\Gamma$.

Hence we have to prove that for $a=[\Gamma\sigma\Gamma]$
$$
h_1 t^a h_2 = \sum_{\theta \in \Gamma\sigma\Gamma} \langle \pi(\theta)h_2,h_1 \rangle \theta,
$$
i.e.,
$$
t^a = \sum_{\theta \in \Gamma\sigma\Gamma} \langle \pi(h_1^{-1}\theta h_2)e,e \rangle h_1^{-1}\theta h_2.
$$
Doing a change of variable $\theta' = h_1^{-1} \theta h_2$
this equality becomes   the definition of $t^a$.

Finally, property 6) follows from the fact that in this case $\pi(\sigma)$
commutes with $\Gamma$ on $\ell^2(\Gamma)$ so it must be a scalar $\lambda$. Hence
$$
t ( \gamma_1 \sigma \gamma_2 ) = \langle \pi(\gamma_1 \sigma \gamma_2)e,e\rangle
= \lambda  \langle \gamma_1 \gamma_2 e,e \rangle
$$
which is different from 0, if and only if $\gamma_1\gamma_2=e$.

But in this case $\Gamma \sigma\Gamma$ is simply $\Gamma \sigma$ and hence
$$
t^{\Gamma\sigma\Gamma} = t^{\Gamma \sigma} = \sum_{\gamma \in \Gamma} t(\gamma \sigma)
\gamma \sigma = t(\sigma) \sigma = \lambda \sigma.
$$
For property 7) note that
$$
\pi(\sigma)e = \sum_{\gamma \in \Gamma} \langle \pi(\sigma)e,\gamma \rangle \gamma
= \sum_\gamma t(\gamma^{-1} \sigma)\gamma.
$$
Hence
$$
\sigma^{-1}(\pi(\sigma)e) = \sum_{\gamma \in \Gamma} t(\gamma^{-1}\sigma)\sigma^{-1}\gamma
=\sum_{\gamma \in \Gamma} \overline{t(\sigma^{-1}\gamma)} \sigma^{-1}\gamma =
\overline{t^{\sigma^{-1} \Gamma}}.
$$
Taking the adjoint we obtain
$$
(\pi(\sigma)e)^* \sigma = \overline{t^{\Gamma\sigma}}.
$$

We now prove property 8).

Let $\Gamma s, \Gamma t$ be two left cosets as in the statement. Let $\gamma$
be any element in $\Gamma$. Then $E^{\L(G)}_{\L(\Gamma)} (t^{\Gamma s} \gamma(t^{\Gamma t})^*)$ is different
from 0, if and only if there exists $\gamma_1, \gamma_2$ and $\theta$ in $\Gamma$
such that
$$
\gamma_1 s \gamma t^{-1} \gamma_2 = \theta
$$
which is equivalent to
$$
\gamma= s^{-1}(\gamma_1\theta \gamma_2^{-1})t = s^{-1}(\gamma')t,
$$
where $\gamma'$ belongs to $\Gamma$.

Thus $E^{\L(G)}_{\L(\Gamma)}(t^{\Gamma s} \gamma(t^{\Gamma t})^*)$ is different from 0, if and only if
$\gamma$ belongs to $\Gamma \cap s^{-1} \Gamma t = A_{\Gamma s, \Gamma t}$.
But this gives exactly that 
$$
E^{\L(G)}_{\L(\Gamma)}(t^{\Gamma s} \gamma(t^{\Gamma t})^*) = \a_{ \Gamma s, \Gamma t}(\gamma)
$$
which by linearity proves the statement of property 8).

Note that $\a_{ \Gamma s, \Gamma t}$ is the zero projection if $\Gamma\cap s^{-1} \Gamma t$
is void.

To prove property 9) we use property 8).
$\a_{\Gamma \sigma s_i, \Gamma \sigma s_j}(\gamma)$ is different from 0, if and only if
$\gamma$ belongs to $\Gamma \cap (\sigma s_i)^{-1} \Gamma (\sigma s_j) =
\Gamma \cap s_i^{-1} \sigma^{-1} \Gamma \sigma s_j$
for $\gamma$ in~$\Gamma$. 

So $\a_{\Gamma \sigma s_i, \Gamma \sigma s_j}(\gamma)$ is different from 0, if and only if
there exists $\theta$ in $\Gamma$ such that
$$
\gamma = s_i^{-1} (\sigma^{-1}\theta \sigma)s_j \quad\hbox{ (or } \sigma s_i \gamma = \theta\sigma s_j)
$$
or
$$
s_i \gamma s_j^{-1} = \theta' = \sigma^{-1} \theta \sigma.
$$
Hence $\theta$ belongs to $\Gamma_{\sigma^{-1}}$ and
$s_i \gamma = \theta' s_j$ so $j$ must be equal to $\pi_\gamma(i)$.
\end{proof}

\section{The representation of the Hecke operators for $\Gamma \subseteq G$\\ on
the type II$_1$ von Neumann algebra $\L(G)$}

This section contains the main technical result of the paper. In Section 2 we obtained 
an explicit formula for abstract Hecke operator.

In this section we prove that the algebra consisting of completely positive
maps representing the Hecke operators has a lifting to $\L(G)$. This lifting
is similar to the dilation of a semigroup of completely positive maps, as explained in the introduction.
It relies on the representation for the Hecke algebra given in the
previous section. The result is a formula that does not involve in its expression any 
choice of a system of representatives.

The main theorem of this paper is the following.

 Note that
all the result in this section remain valid in the presence of a group 2-cocycle on $G$, which restricts to the group $\Gamma$. We will assume that all the partial automorphisms of $\Gamma$, Ad\ $ \sigma $, $\sigma \in G$ are $\varepsilon$ preserving, (see also Appendix 1).

\begin{thm}\label{homeo}
Let $G$ be a discrete group and $\Gamma\subseteq G$ an almost normal subgroup.
Assume that $G$ admits a unitary representation $\pi$ on $\ell^2(\Gamma)$ that
extends the left regular representation. For a coset $[\Gamma \sigma \Gamma]$ let
$\tilde{\Psi}_{\sigma}=\tilde{\Psi}_{[\Gamma\sigma\Gamma]}$ be the abstract Hecke operator, associated with the unitary
$\pi(\sigma)$, 
$$
\tilde{\Psi}_\sigma(x) = 
E ^{\L(\Gamma_\sigma)'} _{\L(\Gamma)'}
(\pi(\sigma)x \pi(\sigma)^*)
$$
for $x$ in $\cR(\Gamma)$. We identify $\cR(\Gamma)$ with $\L(\Gamma)$
via the canonical anti-isomorphism and hence consider $\tilde{\Psi}_\sigma$ as a map from $\L(\Gamma)$
into $\L(\Gamma)$.

Let $\rho : \H(\Gamma\setminus G/\Gamma) \to \L(G)$ be the
representation of the Hecke algebra constructed in
the previous section, so that 
$$
\rho([\Gamma\sigma\Gamma]) =
\sum_{\theta \in [\Gamma\sigma\Gamma]} \langle \pi(\theta) e,e,\rangle_{\ell^2(\Gamma)}
\theta = t^{\Gamma\sigma\Gamma}
$$
for $a=[\Gamma\sigma\Gamma]$ a double coset.

Then for $x$ in $\L(\Gamma)$, $$\Psi_\sigma (x)=[\Gamma:\Gamma_\sigma] \tilde{\Psi_\sigma} (x) =  
[\Gamma:\Gamma_\sigma]E ^{\L(G)}_{\L(\Gamma)} 
(\overline{\rho(a)} x \overline{\rho(a)}^*).$$

Note that in particular $\Psi_\sigma$ depends only on the coset $\Gamma\sigma\Gamma$.

This formula is  a dilation formula, for the ``pseudo-semigroup" of completely positive maps $\Psi_\sigma$, in the sense of the corresponding theory for semigroups of completely positive maps ([Ar]).

\end{thm}

{\bf Remark.} In Appendix 4 we are constructing a two variable version of the Hecke operators. One starts with a representation of  the groupoid $(G\times G^{\rm op}) \rtimes K$ on a Hilbert space $V$  ($K$ is the profinite completion of $\Gamma$). By restricting to 
$\Gamma \times\Gamma$ invariant vectors in $V$, one obtains a new representation of the Hecke algebra associated to $\Gamma \subseteq G$. In Example \ref{Heckeas2}, we prove that the above construction is a particular realization of this new model for the Hecke operators.

{\bf Remark.} By using the anti-linear isomorphism ${}^{\overline{\;\;\;}} : \L(G) \to \L(\Gamma)$
defined by $\sum x_\gamma \gamma \to \sum \overline{x}_\gamma \gamma$, where  $x_\gamma$ are  complex numbers, the formula for
$\Psi_\sigma(x)$ becomes
$$[\Gamma:\Gamma_\sigma]\tilde{\Psi}_\sigma(x)=[\Gamma:\Gamma_\sigma]E ^{ \L(G)} _{\L(\Gamma)}
(\rho(a) x \rho(a)^*),\ x \in \L(G),\ a=[\Gamma\sigma\Gamma],\ \sigma \in G.$$

If $\Gamma = \PSL_2(\Z)$, $G = \PGL_2(\Q_+)$ , by Proposition 4,  
$\Psi_\sigma=[\Gamma:\Gamma_\sigma]\tPsi_{\sigma}$ is unitary equivalent
to the Hecke operator associated with $\Gamma\sigma\Gamma$ on Maass forms. 

In the next proposition
we will prove that, as in the classical case
$$
[\Gamma:\Gamma_{ \sigma_1}][\Gamma:\Gamma_{\sigma_2}] \tilde{\Psi}_{\sigma_1}  \tilde{\Psi}_{\sigma_2} =  
\sum_{\Gamma z \Gamma \subseteq \Gamma \sigma_1 \Gamma \sigma_2 \Gamma}
c(\sigma_1, \sigma_2,z)[\Gamma:\Gamma_z]  \tilde{\Psi}_z.
$$

Recall that $\Gamma_{\sigma^{-1}}=\Gamma\cap \sigma^{-1}\Gamma\sigma$, $\Gamma_{\sigma}=\Gamma\cap \sigma 
\Gamma\sigma^{-1}$ and $s_i$ is a system of left representatives for left cosets for $\Gamma_{\sigma^{-1}}$
in $\Gamma$, that is $\Gamma=\bigcup_{i=1}^n \Gamma_{\sigma^{-1}}s_i$, where 
$n = [\Gamma:\Gamma_\sigma] = [\Gamma:\Gamma_{\sigma^{-1}}]$. Let $t_i=\pi(\sigma)s_i$,
which is a $\L(\Gamma_{\sigma^{-1}})$ orthonormal family of vectors in $\ell^2(\Gamma)$ (that is
$E_{\L(\Gamma_{\sigma^{-1}})}(t_it_j^*)=\delta_{ij}$). 

Moreover, $\theta_\sigma  : \Gamma_{\sigma^{-1}} \to \Gamma_{\sigma}$ is the isomorphism
implemented by $\sigma$, defined by $\theta_\sigma (\gamma)=\sigma\gamma \sigma^{-1}$ for $\gamma$
in $\Gamma_{\sigma^{-1}}$. In particular, for $y$ in $\L(\Gamma)$ we have~that
$$
E^{\L(\Gamma)}_{\L(\Gamma_\sigma)}(\sigma x \sigma^{-1})= \sigma E_{\L(\Gamma_{\sigma^{-1}})}(x) \sigma^{-1}.\leqno(6)
$$

In Proposition 6 we proved that $\tPsi_\sigma(x)$ is given by the following formula, for $x$ in $\L(\Gamma)$:
$$
[\Gamma\sigma\Gamma] \tPsi_\sigma(x) = 
\sum_{i,j=1}^{[\Gamma:\Gamma_\sigma]} t_i^* \theta_\sigma E^{\L(\Gamma)}_{\L(\Gamma_{\sigma^{-1}})}
(s_i x s_j^{-1})t_j.
$$
By linearity we may assume that $x$ is equal to $\gamma=L_\gamma$ a group element in $\L(\Gamma)$.
Let $\pi_\gamma$ be the permutation of the set $\{1,2,\ldots,[\Gamma:\Gamma_{\sigma^{-1}}]\}$,
determined by the requirement that 
$$
s_i \gamma = \theta_i(\gamma) s_{\pi_\gamma(i)},
$$
where $\theta_i(\gamma)=s_i \gamma s_{\pi_\gamma(i)}^{-1}$ belongs to $\Gamma_{\sigma^{-1}}$.
Then
$$
\tPsi_\sigma(\gamma) = \frac{1}{[\Gamma:\Gamma_\sigma]}
\sum_{i=1}^{[\Gamma:\Gamma_\sigma]} t_i^* \theta_\sigma(s_i \gamma s_{\pi_\gamma(i)}^{-1})
t_{\pi_\gamma(i)}.
$$
Because $\theta_\sigma$ is conjugation by $\sigma$ this is further equal to
$$
\frac{1}{[\Gamma:\Gamma_\sigma]}
\sum_{i=1}^{[\Gamma:\Gamma_\sigma]} t_i^* \sigma s_i(\gamma) s_{\pi_\gamma(i)}^{-1}
\sigma^{-1} t_{\pi_\gamma(i)}. \leqno(7)
$$

By property 7) in Remark 8
$$
t_i^* = (\pi(\sigma)s_i)^* = (\pi(\sigma s_i)e)^*
$$ 
is equal to 
$$
\overline{t^{\Gamma\sigma s_i}}(\sigma s_i)^{-1}
$$
and hence
$$
t_i^* \sigma s_i = \overline{t^{\Gamma\sigma s_i}} \quad
\hbox{for } i=1,2,\ldots,n.
$$

Consequently, combining this with formula 7) it follows that
$\tPsi_\sigma(\gamma)$ is further equal to
$$
\frac{1}{[\Gamma:\Gamma_\sigma]}
\sum_{i=1}^{[\Gamma:\Gamma_\sigma]}
\overline{t^{\Gamma\sigma s_i}}
\gamma ( \overline{t^{\Gamma\sigma s_{\pi_\gamma(i)}}})^*. \leqno(8)
$$

Note that
$
E_{\L(\Gamma)}^{\L(G)} (\overline{t^{\Gamma\sigma s_i}}
\gamma ( \overline{t^{\Gamma\sigma s_j}})^*)
$
is equal to
$
\delta_{i,\pi_{\gamma}(i)}\overline{t^{\Gamma\sigma s_i}}
\gamma ( \overline{t^{\Gamma\sigma s_{\pi_\gamma(i)}}})^*.
$

Indeed, a term of the form $\overline{t^{\Gamma\sigma s_i}}
\gamma \overline{t^{\Gamma\sigma s_j}}$ contains various terms of the form
$a(\gamma_1 \sigma s_i \gamma s_j^{-1} \sigma^{-1} \gamma_2)$.
Then $E ^{\L(G)}_{\L(\Gamma)}$ of such a term is different from 0 if and
only if there exists a $\theta$ in $\Gamma$ such that
$$
\gamma_1 \sigma s_i \gamma s_j^{-1} \sigma^{-1} \gamma_2 = \theta
$$
which implies that
$$
\sigma (s_i \gamma s_j^{-1}) \sigma^{-1} = \gamma_1^{-1}\theta\gamma_2
$$
and hence
$$
s_i \gamma s_j^{-1} = \sigma^{-1}(\gamma_1^{-1}\theta\gamma_2).
$$
Thus, $s_i \gamma s_j^{-1} = \theta_1$ for some $\theta_1$ in $\Gamma_{\sigma^{-1}}$
and hence $s_i \gamma =\theta_1 s_j$. But this by the definition of the 
permutation $\pi_\gamma$ implies that $j=\pi_\gamma(i)$.

Thus the equality (8) might be continued as
$$
\sum_{i,j=1}^{[\Gamma:\Gamma_\sigma]}
E^{\L(G)}_{\L(\Gamma)}
(\overline{t^{\Gamma\sigma s_i}}
\gamma  \overline{t^{\Gamma\sigma s_j}}).
$$
But $t^{\Gamma\sigma\Gamma} = \sum_{i=1}^n t^{\Gamma\sigma s_i}$ and
hence this is further equal to
$$
E_{\L(\Gamma)}^{\L(G)}
(\overline{t^{\Gamma\sigma \Gamma}}
\gamma  \overline{t^{\Gamma\sigma \Gamma}}).
$$
By linearly this gives the required formula for $\tPsi_\sigma$.

It is well known that the Hecke operators on Maass forms (or cusp forms)
give a representation for the Hecke algebra.

This is also true for the abstract Hecke operators, and we prove
this, directly from the formula in the preceding theorem.

\vskip6pt

{\bf Proposition.} {\it The map $\a=[\Gamma\sigma\Gamma] \to \Psi_\a= [\Gamma:\Gamma_\sigma] \tPsi_\a$
described in the previous theorem is a $*$ morphism from $\H(G\setminus \Gamma/G)$
into the algebra of bounded operators on $\ell^2(\Gamma)$.
If $a_1 = \Gamma \sigma_1 \Gamma$, $a_2 = \Gamma \sigma_2 \Gamma$
are two double cosets with multiplication rule
$$
a_1 a_2 = \sum_{\Gamma z \Gamma \subseteq \Gamma \sigma_1 \Gamma \sigma_2
\Gamma} c (\sigma_1,\sigma_2, z) \Gamma z \Gamma.
$$
Then for all $x$ in $\L(\Gamma)$}
$$
[\Gamma:\Gamma_{\sigma_1}][\Gamma:\Gamma_{\sigma_2}] \tilde{\Psi}_{\sigma_1} \tilde{\Psi}_{\sigma_2} =  
\sum_{\Gamma z \Gamma}
c (\sigma_1,\sigma_2, z) [\Gamma:\Gamma_z] \tilde{\Psi}_z
$$
and hence
$$
[\Gamma:\Gamma_{\sigma_1}][\Gamma:\Gamma_{\sigma_2}] E(\overline{t^{a_1}}
E(\overline{t^{a_2}} x (\overline{t^{a_2}})^* )(\overline{t^{a_1}})^*)
=
$$
$$
= \sum_{\Gamma z \Gamma \subseteq \Gamma \sigma_1 \Gamma \sigma_2
\Gamma} 
c (\sigma_1,\sigma_2, z) [\Gamma:\Gamma_z]
E^{\L(G)}_{\L(\Gamma)}(\overline{t^{\Gamma z \Gamma}} x \overline{t^{\Gamma z \Gamma}}).
$$

\begin{proof}
To do this we need first to formulate another variant for the
formula of $\Psi^{\Gamma \sigma \Gamma}$, for $\Gamma \sigma \Gamma$ a double
coset of $\Gamma$ in $G$.

Let $(s_i)_{i=1}^{[\Gamma:\Gamma_{\sigma^{-1}}]}$ be a system of representatives for
right cosets for $\Gamma_{\sigma^{-1}}$ in $\Gamma$, that is $\Gamma$
is the disjoint union of $\Gamma_{\sigma^{-1}} s_i$, $i=1,2,\ldots,[\Gamma:\Gamma_{\sigma^{-1}}]$.
Then for each $\gamma$ in $\Gamma$ there exists a permutation $\pi_\gamma$ of the
set $\{i=1,2,\ldots,[\Gamma:\Gamma_{\sigma^{-1}}]\}$ such that for each $i$, there exists
$\theta_i(\gamma)$ in $\Gamma_{\sigma^{-1}}$ with the property that
$$
s_i \gamma = \theta_i (\gamma) s_{\pi_\gamma}(i).
$$
We proved that
$$
[\Gamma:\Gamma_\sigma] \tPsi_\sigma(\gamma)=
\sum_{i=1}^{[\Gamma:\Gamma_{\sigma^{-1}}]}
\overline{t^{\Gamma\sigma s_i}}
\gamma ( \overline{t^{\Gamma\sigma s_{\pi_\gamma(i)}}})^*. \leqno(9)
$$

By property 9) in Remark 8, let $\a_{ij} = \a_{\Gamma\sigma s_i,\Gamma\sigma s_j}$
be the projection from $\ell^2(\Gamma)$ onto the Hilbert space generated by
$$\Gamma\cap (\sigma s_i)^{-1}\Gamma(\sigma s_j)=\Gamma \cap s_i (\sigma \Gamma \sigma^{-1})s_j.$$
Combining properties 8) and 9) it follows that
$$
[\Gamma:\Gamma_\sigma] \Psi_\sigma(x)=
\sum_{i,j=1}^n \overline{t^{\Gamma\sigma s_i}}
\a_{\Gamma\sigma s_i,\Gamma\sigma s_j}(x) 
\overline{t^{\Gamma\sigma s_j}}. 
$$

Let now $a_1 = [\Gamma \sigma_1 \Gamma]$, $a_2 = [\Gamma \sigma_2 \Gamma$], be two
double cosets in $G$, for which we want to compute the composition
$$
[\Gamma:\Gamma_{\sigma_1}][\Gamma:\Gamma_{\sigma_2}]  \tilde{\Psi}_{a_2} \circ \tilde{\Psi}_{a_1}.
$$

Assume that $s_i$ are representatives for left $\Gamma_{\sigma_1^{-1}}$ cosets 
in $\Gamma$ (that is $\Gamma =$\break $\bigcup_{i=1}^{[\Gamma:\Gamma_{\sigma_1}]} s_i \Gamma_{\sigma_1^{-1}}$
 and similarly assume that $r_\a$, $\a = 1,2,\ldots, [\Gamma:\Gamma_{\sigma_2^{-1}}]$
are representatives for $\Gamma_{\sigma_2^{-1}}$ left cosets, that is 
$\Gamma = \bigcup \Gamma_{\sigma_2^{-1}} r_\a$.

Recall that by property c) in Remark 8, we have that for all 
$i = 1,2,\ldots,$ $[\Gamma:\Gamma_{\sigma_1}]$
$$
t^{\Gamma\sigma_2 \Gamma}t^{\Gamma\sigma_1 s_i} = \sum_{a=1}^{[\Gamma:\Gamma_{\sigma_2}]}
t^{\Gamma \sigma_2 r_a \sigma_1 s_i}. \leqno(11)
$$

Let $\pi_\gamma$ be the permutation associated to the cosets $\Gamma_{\sigma_1^{-1}} s_1,
\Gamma_{\sigma_1^{-1}} s_2,\ldots$ as in Remark 8 and 9.

Then by using property 10) for $[\Gamma\sigma_1 \Gamma]$ we obtain that for every
$\gamma$ in $\Gamma$, we have~that
$$
[\Gamma:\Gamma_{\sigma_2^{-1}}][\Gamma:\Gamma_{\sigma_1^{-1}}] 
\tPsi_{\sigma_2} (\tPsi_{\sigma_1}(\gamma))=
$$
$$
=
\sum_{i,j=1}^{[\Gamma:\Gamma_{\sigma_1}]}
[\Gamma:\Gamma_{\sigma_2^{-1}}][\Gamma:\Gamma_{\sigma_1^{-1}}]
E^{\L(G)}_{\L(\Gamma)}
\overline{t^{\Gamma\sigma_2 \Gamma}}\,\overline{t^{\Gamma\sigma_1 s_i}}
\a_{\Gamma\sigma s_i,\Gamma\sigma s_j} (\gamma)
(\overline{t^{\Gamma\sigma_1 s_j}})^* (\overline{t^{\Gamma\sigma_2 \Gamma}})^* 
$$
which by using the equality (11) is further equal to 
$$
\sum_{a,b=1}^{[\Gamma:\Gamma_{\sigma_2^{-1}}]}
\sum_{i,j=1}^{[\Gamma:\Gamma_{\sigma_1^{-1}}]}
[\Gamma:\Gamma_{\sigma_2^{-1}}][\Gamma:\Gamma_{\sigma_1^{-1}}]
E^{\L(G)}_{\L(\Gamma)}
(\overline{t^{\Gamma \sigma_2 r_a \sigma_1 s_i})
\a_{\Gamma\sigma s_i,\Gamma\sigma s_j}} (\gamma)
(\overline{t^{\Gamma \sigma_2 r_b \sigma_1 s_j})^*}). \leqno(12)
$$

As noted in property 9) of Remark 8, $\a_{\Gamma\sigma s_i,\Gamma\sigma s_j}(\gamma)$
is different from 0 if and only if $\gamma \in \Gamma \cap s_1^{-1} \sigma^{-1}\Gamma\sigma s_j$
which is equivalent to the fact that there exists $\theta$ in $\Gamma$ such that
$
\sigma s_i\gamma = \theta \sigma s_j.
$

Moreover, still as a consequence from property 8) in Remark 8 it follows that
a term in the sum (12) is different from 0 if and only if there exists $\theta'$ in $\Gamma$
such that 
$$
(\sigma_2 r_a \sigma_1 s_i)\gamma = \theta'(\sigma_1 r_b \sigma_1 s_j).
\leqno(13)
$$
But $j$ was determined by the fact that 
$$\sigma_1 s_i = \theta \sigma_1 s_j \quad \hbox{for some }\theta \hbox{ in }\Gamma.
\leqno(14)
$$
From (14) we deduce that 
$$
\sigma_2 r_a \sigma_1 s_i\gamma = \sigma_2 r_a \theta \sigma_1 s_j
$$
and using (13) we deduce that
$$
\sigma_2 r_a \theta \sigma_1 s_j = \theta' \sigma_1 r_b \sigma_1 s_j
$$
and hence
$$
\sigma_2 r_a \theta = \theta' \sigma_1 r_b.
$$
Hence $b$ and $\theta'$ are uniquely determined by $\theta$ and $a$ and
hence by $a,i$ and $\gamma$.

Thus there exists a bijection $\a_\gamma=(\a_\gamma^1,\a_\gamma^2)$
of the set $\{1,2,\ldots, [\Gamma:\Gamma_{\sigma_2^{-1}}]\}\times
\{1,2,\ldots, [\Gamma:\Gamma_{\sigma_1^{-1}}]\}$
which to every pair $(a,i)$ associates 
the unique $b=\a_\gamma^1(a,i)$, $j=\a_\gamma^2(a,i)=\pi_\gamma(i)$
for which the term starting with 
$t^{\Gamma \sigma_1 r_a \sigma_2 s_i}$ in the sum (12) remains
non zero after applying $E ^{\L(G)}_{\L(\Gamma)}$. Moreover, this bijection has the 
property that for all $(a,i)$ in $\{1,2,\ldots, [\Gamma:\Gamma_{\sigma_2^{-1}}]\}\times
\{1,2,\ldots, [\Gamma:\Gamma_{\sigma_1^{-1}}]\}$
we have that there exists $\theta'$ in $\Gamma$ 
such that
$$\sigma_2 r_a \sigma_1 s_i\gamma = \theta'
\sigma_2 r_{\a^1_{\gamma(a,i)}} \sigma_1 s_{\pi_\gamma(i)}.
$$

Thus,
$
[\Gamma:\Gamma_{\sigma_2^{-1}}]
[\Gamma:\Gamma_{\sigma_1^{-1}}]
\tPsi_{\sigma_2}\tPsi_{\sigma_1}(\gamma)
$
is equal to
$$
\sum_{a=1}^{[\Gamma:\Gamma_{\sigma_2^{-1}}]}
\sum_{i=1}^{[\Gamma:\Gamma_{\sigma_1^{-1}}]}
\overline{t^{\Gamma \sigma_2 r_a \sigma_1 s_i}}
\gamma\Big(\overline{t^{\Gamma \sigma_2 r_{\a^1_{\gamma(a,i)}} \sigma_1 s_{\pi_\gamma(i)}}}\Big)^*, \leqno(15)
$$
where $\a_\gamma = (\a^1_\gamma,\a^2_\gamma)$ is a bijection.

On the other hand, we know that
$$
t^{\Gamma\sigma_2 \Gamma}t^{\Gamma\sigma_1 \Gamma}
= \sum_{\Gamma z \Gamma \subseteq \Gamma \sigma_1 \Gamma \sigma_2
\Gamma} c (\sigma_2,\sigma_1, z) t^{\Gamma z \Gamma},
$$
where the multiplicities $c (\sigma_2,\sigma_1, z)$
are strictly positive integer numbers that come from
the algebra structure of the Hecke algebra of double cosets.

Moreover, as we have seen above 
$$
t^{\Gamma\sigma_2 \Gamma}t^{\Gamma\sigma_1 \Gamma}
= 
\sum_{a=1}^{[\Gamma:\Gamma_{\sigma_2^{-1}}]}
\sum_{i=1}^{[\Gamma:\Gamma_{\sigma_1^{-1}}]}
t^{\Gamma \sigma_2 r_a \sigma_1 s_j}.
$$
Hence the enumeration of left cosets in 
$[\Gamma\sigma_2 \Gamma][\Gamma\sigma_1 \Gamma]$ is
$\Gamma \sigma_2 r_a \sigma_1 s_i$, 
$a=1,2,\ldots,$ $[\Gamma:\Gamma_{\sigma_2^{-1}}]$, 
$i=1,2,\ldots,[\Gamma:\Gamma_{\sigma_1^{-1}}]$.

This enumeration will contain for each coset $[\Gamma z \Gamma] \subseteq
[\Gamma\sigma_1 \Gamma][\Gamma\sigma_2 \Gamma]$ exactly $c(\gamma_1,\gamma_2,z)$
sets of $[\Gamma:\Gamma_{z^{-1}}]$ cosets, that together constitute of
$\Gamma z \Gamma$.

The contribution of any such group in the sum (15), 
will be one copy of\break 
$E^{\L(G)}_{\L(\Gamma)}(t^z \gamma(t^z)^*)$.

But this proves exactly that
$$
[\Gamma:\Gamma_{\sigma_1^{-1}}][\Gamma:\Gamma_{\sigma_2^{-1}}]E^{\L(G)}_{\L(\Gamma)}
(\overline{t^{\Gamma\sigma_2 \Gamma}} E(\overline{t^{\Gamma\sigma_1 \Gamma}} \gamma (\overline{t^{\Gamma\sigma_1 \Gamma}})^*)
(\overline{ t^{\Gamma\sigma_2 \Gamma}})^*)
$$
is
$$
\sum_{\Gamma z \Gamma \subseteq \Gamma \sigma_1 \Gamma \sigma_2
\Gamma} \a (\sigma_2,\sigma_1, z)
[\Gamma:\Gamma_z]
 E^{\L(G)}_{\L(\Gamma)}(\overline{t^z} \gamma (\overline{t^z})^*).
$$
By linearity this proves our result.
\end{proof}

In concrete examples, it might happen that we have the unitary representation
$\pi$ of $G$ on a Hilbert space $H$, and that we know that
$\pi|\Gamma$ is unitarily equivalent to the left regular representation,
but without knowing precisely the structure of the intertwiner realizing the unitary equivalence.

Hence, it would be useful to proceed with the construction of the elements 
$t^{\Gamma\sigma \Gamma}$, but starting just with a cyclic vector $\eta$
(which automatically is separating) which is not necessary a trace vector, in the Hilbert space of the representation
of $\pi$.

So, in this case we would start with
$$
\tilde{t}^{\Gamma\sigma\Gamma} = \sum_{\theta\in \Gamma\sigma\Gamma}
\langle \pi(\theta)\eta,\eta\rangle.
$$

For example, in the case of $\PSL_2(\Z)$ represented on the space $H_{13}$ ([GHJ])
by Perelmov ([Pe], see also [KL]) we know that the evaluation
vector at any given point in $\mathbb{H}$ is cyclic.
Then the $\tilde{t}^{\Gamma\sigma\Gamma}$ might have an 
easier expression.

To exemplify we replace $\PSL_2(\R)$ by $SU(1,1)$. Hence the upper half plane
gets replaced by the unit disk, and $\PSL_2(\Z)$ gets replaced by a
discrete subgroup  $\Gamma_0$ of $SU(1,1)$. Let $\eta$ be the evaluation
vector at 0, so $\eta$ becomes the constant function $1$ and $\langle \pi(\theta)\eta,\eta\rangle_{H_{13}}$
is clearly easy to compute (since $\pi(\gamma)1$ is a multiple of the evaluation
vector at $\gamma_0$).

In the next lemma we prove that the family of ``deformed'' $\tilde{t}^{\Gamma\sigma\Gamma}$
might be used to compute $\Psi_{\sigma}$.

\begin{prop}
Let $\eta$ be a cyclic separating vector in $\ell^2(\Gamma)$ and let, for $\sigma$ in~$G$,
$$
\tilde{t}^{\Gamma\sigma\Gamma} = 
\sum_{\theta\in \Gamma\sigma\Gamma}
\langle \pi(\theta)\eta,\eta\rangle \theta.
$$

Let $x = (\eta^* \eta)^{1/2}$ which is invertible at least in the affiliated algebra
of unbounded operators. Then $\xi = x^{-1/2} \eta$ is a cyclic trace vector,
and hence by Remark~$8$, property $4)$,
$$
t^{\Gamma\sigma\Gamma} = \overline{x}^{-1/2} \tilde{t}^{\Gamma\sigma\Gamma}
x^{-1/2}
$$
and hence for $y$ in $\L(\Gamma)$,
$$
[\Gamma:\Gamma_\sigma]\tPsi_\sigma(y)=[\Gamma:\Gamma_\sigma]E^{\L(G)}_{\L(\Gamma)}\big(\overline{x}^{-1/2}\overline{ \tilde{t}^{\Gamma\sigma\Gamma}}
(\overline{x}^{-1/2}y  x^{-1/2})
(\overline{ \tilde{t}^{\Gamma\sigma\Gamma}})^* x^{-1/2}\big).
$$
\end{prop}

\begin{proof} This is now obvious.
\end{proof}

There is a very simple way to compute the element $x$ in the 
preceding lemma from the matrix coefficients $\langle \pi(\theta)\eta,\eta\rangle$,
$\gamma \in \Gamma$. This is certainly well known to specialist, but
for completeness we include the exact result here.

\vskip6pt

{\bf Lemma.} {\it Let $\eta$ in $\ell^2(\Gamma)$ be given. Assume
we know the element $A = \sum \langle \gamma\eta,\eta\rangle \gamma^{-1}$.
Then $\xi = (A^*A)^{-1/2}\eta$ is a cyclic trace vector in $\ell^2(\Gamma)$.}

\begin{proof}
Indeed, 
$$A=\sum_\gamma \tau(\eta^*\gamma\eta) \gamma^{-1} =
\sum_\gamma \tau((\eta\eta^*) \gamma^{-1})\gamma =
\sum_\gamma \langle \eta\eta^*, \gamma\rangle \gamma = \eta\eta^*.
$$
Hence $(\eta\eta^*)^{-1/2} = A^{-1/2}$ which is invertible since $\eta$ is
cyclic and separating. Then
$(\eta\eta^*)^{-1/2} \eta$ is a unitary, that is (as a vector) a cyclic trace vector.
\end{proof}

\section{Complete positivity multipliers properties\\ for eigenvalues for a joint
eigenvector\\ of the Hecke operators}

In this section we derive further consequences, from the relations derived
in the previous chapter, regarding the relative position
in $\L(G)$ of the algebra $\L(\Gamma)$, and the von Neumann algebra $\H \subseteq \L(G)$ generated by the
$\overline{t}^{\a}$'s $\a$ running in the space of double cosets of $\Gamma$ in $G$.
This  will also work for the image of $\H$, through the canonical conjugation anti-isomorphism
in $\L(\Gamma)$.

To avoid cumbersome notations we will use $\rho(a)$ and $t^a$ for $\overline{\rho(a)}$
and $\overline{t^a}$ for $a=[\Gamma\sigma\Gamma]$ double cosets.

Let $\D$ be the von Neumann subalgebra in $B(\ell^2(G))$ generated by the operators
of left and right multiplication with elements in $\H$, that is by $L_\a = L_{t^\a}$, 
$R_\a = R_{t^\a}$, the left and right convoluters by the elements $t^\a \in \H$,
that are associated to double cosets $\a = [\Gamma\sigma\Gamma]$, $\sigma$ in $G$.

From an algebra point of view, the algebra $\D$ is isomorphic to $\H \otimes \H$, but, when taking
closures, this might be false (e.g., see the action of the algebra $\D$ on the
vector 1 (the unit of $G$), viewed as a vector in $\ell^2(G)$ (see [Po])).

Let $P$ be the projection from $\ell^2(G)$ onto $\ell^2(\Gamma)$. Then, by property
2) in Remark~8, it follows that 
$$
PL_\a R_\beta^* P=0
$$
unless $\a = \beta$ in which case
$$
PL_\a R_\a^* P=[\Gamma:\Gamma_\sigma]\tilde{\Psi}_\a, \quad \a=[\Gamma\sigma\Gamma].
$$

Then $\D$ has the following remarkable property:
$$
(P\D P)(P\D P) \subseteq (P\D P)
$$
and hence $P\D P$ is an algebra.

Moreover, the algebra $t_\a \to P (L_\a R_\a^*) P$, 
$\a = [\Gamma\sigma\Gamma]$, $\sigma \in G$
extends to a $*$-algebra homeomorphism $\Phi$ from the $*$-algebra generated by the
$t_\a$'s, into $P\D P$.

Although we do not know the structure of the action of the algebra $\D$ on a vector
$\xi$ in $\ell_2(\Gamma)$, that is different from 1,
we can still derive some conclusion in the case when the unit vector $\xi$ in $\ell_2(\Gamma)$
is a joint eigenvalue for all the $[\Gamma:\Gamma_{\sigma^{-1}}]\tilde{\Psi}_\a$'s of eigenvalue
$\lambda(\a)$, $\a = [\Gamma\sigma\Gamma]$ running over all double cosets of $\Gamma$.

Let $K$ be the Hilbert subspace of $\ell^2(G)$ generated by $\H \xi \H$.

The fact that $\xi$ is a norm 1 eigenvector for all the $\tPsi_\a$'s, and $\a$ a double coset,
implies that
$$
[\Gamma:\Gamma_\a]\tPsi_\a(\xi) =[\Gamma:\Gamma_\a] E_{\L(\Gamma)} (t_\a \xi t_\beta^*) = \delta_{\a\beta}
\lambda(\a) \xi
$$
for all double cosets $\a,\beta$ of $\Gamma$ in $G$, and hence
$$
\tau(t_\a \xi (t_\beta^*) \xi^*) = \delta_{\a\beta} \lambda(\a).
$$
(Here $\xi$ is a norm 1 eigenvector for $[\Gamma:\Gamma_\sigma]\tPsi_\a$,
of eigenvalue $\lambda(\a)$, with $\a = \Gamma\sigma\Gamma$, $\sigma$ in $G$.)

We	note the following consequence of these considerations.

\begin{lemma}
Let $\xi$ be a norm $1$ joint eigenvector for the maps
$[\Gamma:\Gamma_{\sigma^{-1}}]\tPsi_\a = [\Gamma:\Gamma_{\sigma^{-1}}]E_{\L(\Gamma)} (t^\a \cdot (t^\a)^*)$
on $\ell_2(\Gamma)$, of eigenvalue $\lambda (\a)$ where $\a$ is the double coset
$\Gamma\sigma\Gamma$, $\sigma$ in~$G$.

Recall that $\H$ is the von Neumann algebra generated by all the $t_\a$'s. Then
$$
E_\H (\xi^* t_\a \xi) = \frac{\lambda_\a}{[\Gamma:\Gamma_{\sigma^{-1}}]} t_\a
$$
for all double cosets $\a =\Gamma\sigma\Gamma$, $\sigma$ in $G$.
\end{lemma}

\begin{proof} Let $\eta$, in $L^2(\H,\tau_G)$, be the vector $E_\H (\xi^* t_\a \xi)$.

Then for all cosets $\beta = \Gamma\sigma_1\Gamma$, $\sigma_1$ in $G$
we have that
$$
\langle \eta, t_\beta \rangle_{\ell^2(G)} = \tau_{\L(G)} (\xi^* t_\a \xi t_\beta^* )
= \tau_{\L(G)} (t_\a \xi t_\beta^* \xi^*) =
\tau_{\L(G)} ( E_{\L(\Gamma)} (t_\a \xi t_\beta^*) \xi^*) 
$$
and this is 0  unless, $\a=\beta$, case in which the quantity above is further
equal to
$$
\tau_{\L(G)} (\lambda_\a\xi \xi^*)  = \lambda_\a.
$$
Thus $\eta$ is a vector in $L^2(\H, \tau_{\L(G)})$ which verifies that
$\langle \eta,t^\beta\rangle$ is 0 unless $\a=\beta$ case in which
$\langle \eta,t^\a\rangle = \lambda_\a$.

Since as proven in Remark 8, $\{t^\a\}$ is an orthogonal basis for
$L^2(\H, \tau_{\L(G)})$ implies that (again by Remark 8)
$$
\eta = \frac{\lambda_\a}{\|t_\a\|_2^2}t_\alpha =\frac{\lambda_\a}{ [\Gamma:\Gamma_{\sigma^{-1}}]}t_\alpha,
$$
if $\a = [\Gamma\sigma\Gamma]$.
\end{proof}

This observation has the following important corollary

\begin{cor} 
Let $G$ be a discrete group and $\Gamma\subset G$ an almost normal
subgroup. Assume that $G$ admits a unitary representation $\pi$ that
extends the left regular representation of $\Gamma$ on $\ell^2(\Gamma)$. For
$\a = [\Gamma\sigma\Gamma]$ a double coset of $\Gamma$ in $G$, 
let $\tilde{\Psi}_\a$ be the completely positive map on $\L(\Gamma)'$ defined by
the formula 
$$
\tPsi_\a(x) = 
E ^{\L(\Gamma_\sigma)'}_{\L(\Gamma)'}
(\pi(\sigma) x \pi(\sigma^{-1})).
$$

Let $\xi$ in $\ell^2(\Gamma)$ be a joint eigenvector of eigenvalue
$\lambda(\a)$, for all the completely positive
linear maps $[\Gamma:\Gamma_{\sigma^{-1}}]\tPsi_\a$, $\a = [\Gamma\sigma\Gamma]$,
$\sigma$ in $G$.

Consider the linear map $\Phi_0$ on $\H(\Gamma\setminus G/\Gamma)$
(the linear span of double cosets) defined~by
$$
\Phi_0 (\a) = \frac{\lambda(\a)}{\ind \a} \a.
$$
Here $\a = [\Gamma\sigma\Gamma]$, runs over all double cosets $\Gamma\sigma\Gamma$
of $\Gamma$ in $G$, and $\ind \a= [\Gamma:\Gamma_{\sigma^{-1}}]$.

Then $\Phi_0$ extends to completely positive linear $\Phi_\lambda$ map on $\H$.

In particular, the sequence $\left(\frac{\lambda(\a)}{\ind \a}\right)_{\a = [\Gamma\sigma\Gamma],\,
\sigma\in G}$ is a completely positive bounded multiplier of the Hecke's double cosets algebra.
\end{cor}

\begin{proof}
The extension of the map $\Phi_0$ is the map $\Phi$ on $\H$ defined by
$\Phi(x) = E_\H(\xi^* x \xi)$.

But this is clearly completely positive. 
\end{proof}

\begin{cor}
Let $\Gamma\subset G$ be an almost normal subgroup as above. Let $\Delta$ be the map from the Hecke algebra $\H_0$
into $\H_0\otimes\H_0$, defined by
$$
\Delta([\Gamma\sigma\Gamma])=\frac{1}{[\Gamma:\Gamma_{\sigma}]}[\Gamma\sigma\Gamma]\otimes[\Gamma\sigma^{-1}\Gamma]
$$
for $\sigma$ in $G$. (We may also extend $\Delta$ to the reduced $C^*$-Hecke algebra $\H$.)

Then $\Delta$ is positive. In the terminology of Vershik {\rm ([Ve])} where this is proved for finite $G$,
the algebra $\H_0$ (with basis $[\Gamma\sigma\Gamma]$) is 2-positive.
\end{cor}

\begin{proof}
We have to verify that if $p$ is positive in $\H_0$ then $\Delta(p)$ is positive.
Since $\H_0$ is a commutative algebra, it is sufficient to prove that if $\chi_\lambda$ is
a character of $\H_0$, then $({\rm Id} \otimes \chi_{\lambda})$ ($\Delta(p))\in\H_0$ is positive).
 
But obviously $({\rm Id} \otimes\chi_{\lambda})$ is the previous map $\Phi_{\lambda}$, which
is positive for all $\lambda$ corresponding  to values in the spectrum of $[\Gamma\sigma\Gamma]$
in the reduced $C^*$-algebra.
\end{proof}

\begin{rem} In the case of $G=\PGL_2(\Z[\frac1{p}])$, $\Gamma = \PSL_2(\Z)$, $p\geq 3$, as we observed before, the 
Hecke algebra $\H_0$ is isomorphic to the radial algebra in the free group with $N=\frac{p-1}{2}$ generators.
The results of [Py], [DeCaHa], also prove that $\Phi_\lambda$ is a completely positive map on $\H$,
for $\lambda$ in the interval $[-(p+1),(p+1)]$. So, we cannot exclude values of $\lambda$ by this method,
in the case of $\PSL_2(\Z)$.

However, we have the following additional property of the map $\Phi_\lambda$,
that is derived from the representation of the primitive structure of the Hecke algebra.
\end{rem}

\begin{prop} Let $\tilde{C}$ be 
the vector space of sets of the form $[\sigma_1 \Gamma \sigma_2]$, $\sigma_1,\sigma_2\in G$. 
We let $\cC(G,\Gamma)$ be the vector space obtained from $\tilde{C}$
by factorizing at the linear relations of the form
$$
\sum[\sigma_1^i\Gamma\sigma_2^i]=\sum[\theta_1^j\Gamma\theta_2^j]
$$
if $\sigma_s^i, \theta_r^j$ are elements of $G$, and the disjoint union $\sigma_1^i \Gamma \sigma_2^i$
is equal to the disjoint union~$\theta_1^j \Gamma \theta_2^j$.  Let $\xi$ be an eigenvector (see Appendix 2). Then there exists a bilinear map 
$\chi:\cC(G,\Gamma)\times \cC(G,\Gamma)\to\C$ such that $\chi$ recovers the value of the eigenvector, that is $\chi|_{\H_0\times\H_0}$
is defined by $\chi([\Gamma\alpha\Gamma],[\Gamma\beta\Gamma])=\delta_{\alpha\beta}
\frac{\lambda([\Gamma\alpha\Gamma])}{[\Gamma:\Gamma_{\alpha}]}$
and $\chi$ is positive in the following sense
$$
\sum\lambda_{i_2i_3}\overline{\lambda_{i_1i_4}}([\sigma_{i_1}\Gamma\sigma_{i_2}]),
([\sigma_{i_3}\Gamma\sigma_{i_4}])\geq 0.
$$
\end{prop}

\begin{proof} Indeed we define 
$$\chi ([\sigma_1\Gamma\sigma_2],[\sigma_3\Gamma\sigma_4])=
\tau(t^{\Gamma\sigma_{i_2}}\xi t^{\sigma_{i_3}\Gamma}t^{\Gamma_{i_4}}\xi^*t^{\sigma_{i_1}\Gamma}).
$$
Note that $\chi$ has a second positivity property coming from the inequality
$$
\tau\Big(\xi\Big(\sum \eta_{i_3}\overline{\eta_{i_4}}t^{\sigma_{i_3}\Gamma}t^{\Gamma_{\sigma_{i_4}}}\Big)\xi^*
\Big(\sum\theta_{i_1}\overline{\theta_{i_2}}t^{\sigma_{i_1}\Gamma}t^{\Gamma_{\sigma_{i_2}}}\Big)\Big)\geq 0
$$
for all complex numbers $\eta_i, \theta_j$.
\end{proof}

\begin{rem}
It is not clear if the completely positive map, for values of $\lambda$ outside $[-2\sqrt{p}, 2\sqrt{p}]$ would
have such an extension $\chi$.
\end{rem}

\

\section{The  structure of the crossed product algebra, modulo the compact operators, of left and right convoluters in $\PGL_2(\Z[\frac1p])$ acting
on $\ell^2(\PSL_2(\Z))$, $p$ prime number}

In this section $G$ will be the discrete group $\PGL_2(\Z[\frac1p])$ and
$\Gamma=\PSL_2(\Z)$. By $\varepsilon$ we denote the 2 group cocycle on $G$ 
with values in $\pm 1$ introduced in Chapter 2
(corresponding to the projective representation $\pi_{13}$ for $\PSL_2(\R)$ on~$H_{13}$).

We will prove an  extension of the usual Akemann-Ostrand property ([AO]),
that asserts the $C^*$-algebra generated by left and right convolution of $\Gamma$ on $\ell^2(\Gamma)$,
is isomorphic modulo the compact operators to the reduced $C^*$-algebra $C^*_{{\rm red}}(\Gamma\times \Gamma^{\rm op})$.
(Here $\Gamma^{\rm op}$ is the group $\Gamma$ considered with the opposite multiplication, so that we have
a natural action of $\Gamma\times \Gamma^{\rm op}$ on $\Gamma$.)

We will extend this result to the (partial) action of $G \times G^{\rm op}$ on $\Gamma$ and identify the structure of the crossed product algebra in the quotient, modulo the compact operators.

 As a consequence,
and since the representation we constructed  in Chapter 5 for the Hecke algebra, (giving unitarily equivalent operators to the classical Hecke operators),
takes values into the $C^*$-algebra generated by left and right convolutors from $G\times G^{\rm op}$,
and characteristic functions of cosets of modular subgroups acting $\ell^2(\Gamma)$, we can compute the essential spectrum 
of the classical Hecke operators.

Let $\mathbf Z_p$ be the $p$-adic integers and $K$ be the compact group $\PSL_2(\mathbf Z_p)$. Note that $K$ is
totally disconnected and that $\Gamma$ is dense in $K$. Let $\mu_p$ be the normalized 
Haar measure on $K$.

We will use the following embedding of the algebra of continuous functions on $K$ into 
$\B(\ell^2(\Gamma))$. To each function $f$ in $C(K)$ we associate the diagonal multiplication
operator on $\ell^2(\Gamma)$ by the restriction of $f$ to $\Gamma \subseteq K$.

In this way, $C(K)$ is identified with the commutative $C^*$-subalgebra $X_\Gamma$ of $\ell^\infty(\Gamma)$
generated by characteristic functions of left cosets (equivalently right) of modular subgroups (we have to add to the generators of $X_\Gamma$  a continuous sign function to separate the points).
Thus $C(K)=X_{\Gamma}\subseteq\ell^{\infty}(\Gamma)$.

The Haar measure $\mu_p$ on $K$ then correspond to the state (trace) on $X_\Gamma$ that associates
to the characteristic function of a coset $s\Gamma_ \sigma$ of a modular subgroup $\Gamma_\sigma$ of $\Gamma$ the value
$\frac1{[\Gamma:\Gamma_\sigma]}$. Note the group
$\tilde{G} = G\times G^{\rm op}$ (where $G^{\rm op}$ is the group $G$ with opposite multiplication)
acts as a groupoid on $\Gamma$ and hence it acts also on $K$, by the formula 
$$
(g_1 \times g_2^{\rm op})(\gamma) = g_1 \gamma g_2, \quad g_1 \times g_2^{\rm op} \in \tilde{G}, \
\gamma \in g_1^{-1}\Gamma g_2^{-1}\cap\Gamma\subseteq\Gamma.
$$

If we take  into account also the cocycle $\varepsilon$  (thus working with $\L(G,\varepsilon)$ instead of $\L(G)$, the formula of the action of $(g_1 \times g_2^{\rm op})$ on $\gamma$ is modified by the factor $\varepsilon(g_1,\gamma) \varepsilon (\gamma, g_2)$.

The domain of $g_1 \times g_2^{\rm op}$ is $\Gamma \cap g_1^{-1} \Gamma g_2^{-1}$. This shows that only elements of the form $g_1 \times g_2^{\rm op}$ with $g_1,g_2$ belonging
to the same double coset of $\Gamma$ in $G$ will have a nontrivial domain. Because $g_1,g_2$ are
in the same double coset the action is measure preserving. Hence we can construct the reduced and maximal 
crossed product algebra
$$
\A = C^*_{\rm red} ((G \times G^{\rm op})\rtimes C(K)),\quad  \A_{\text{max}}=C^*((G \times G^{\rm op})\rtimes C(K))
$$
To construct the reduced crossed product algebra we use the canonical trace $\tau_p$ on the algebraic crossed product
$(G \times G^{\rm op})\rtimes C(K)$ induced by the $G \times G^{\rm op}$ invariant measure
$\mu_p$ on $K$.

We have consequently a covariant representation of the crossed product
$$C^* ((G \times G^{\rm op})\rtimes C(K))$$ which comes from the embedding of $C(K)$ into
$\B(\ell^2(\Gamma))$ described above, and by representing elements in $(G \times G^{\rm op})$ as left and respectively right convolutors.

Indeed, let $\theta: G \times G^{\rm op} \to \B(\ell^2(\Gamma))$
be the representation (by partial isometries) of $G \times G^{\rm op}$ by left and right 
convolutions on $\Gamma$.
Then $\theta$ is compatible (equivalent) with respect to the action of $G \times G^{\rm op}$
on $K$ and hence we get in this way a covariant representation of the $C^*$-algebra
$\A_{\text{max}}=C^* ((G \times G^{\rm op})\rtimes C(K))$ into $\B(\ell^2(\Gamma))$. We denote the  $C^*$-algebra generated by the image of this representation 
by $\B$. Thus $\B$ is generated as a $C^*$-algebra by $\theta(G \times G^{\rm op})$
and $C(K)$.
By the results of Akemann-Ostrand ([AO]; see also [Co]) this algebra contains the compacts operators
$\K(\ell^2(\Gamma))$.

Note that the algebras $\B$ is in fact a corner (under the projection represented by the
characteristic function of $\Gamma$, in the larger crossed product algebra, of the group
$G \times G^{\rm op}$ acting on $\ell^2(G)$).

\vskip6pt
We formulate now our main result, which proves that the quotient algebra, modulo the ideal of compact operators, is isomorphic to the reduced, groupoid crossed product algebra.

\begin{thm}\label{AOlocal}[{Local Akemann-Ostrand property for $\PGL_2(\Z[\frac1p])$}] Let $p$ be a prime number and let $G$= $\PGL_2(\Z[\frac1p])$. 
Let $\B$ be the $C^*$-algebra generated by left and right convolutors in $C^*_{\rm red} (G)$ and by the image of the algebra $C(K)$  acting by multiplication operators on $\ell^2(\Gamma)$. 

Let $\A_0=\B/K(\ell^2(\Gamma))$ be the projection in the Calkin algebra of the algebra $\B\subseteq\B(\ell^2(\Gamma))$ 
Then $\A_0$ is canonically isomorphic to the $C^*$ algebra $\A$, the reduced, $C^\ast$- groupoid crossed product of $G \times G^{\rm op}$ acting on $K$, with respect to the invariant Haar measure on $K$.

 The statement remains valid if instead of the $C^*$-algebra $C^*_{\rm red} (G)$, we consider the skewed crossed product $C^*$-algebra, in which    the canonical $\Z_2$ valueded,  2-group cocycle on $PSL(2,\Z)$ also intervenes.
 \end{thm}

\begin{proof} First we give an outline of the proof.

The reduction of the general case with a $\Z_2$ valued 2-cocycle, to the case when no cocycle is present, will be done in the Remark \ref{z2}.

 In the appendices 4,5,6 we prove a reduction procedure that reduces the analysis of essential states (states induced by the representation in the Calkin algebra) to the analysis of states that are concentrated  in the identify fiber of $K$ (e.g. limits of averaging sets, contained in families of normal subgroups, with trivial intersection,as explained bellow). This reduction procedure  is  valid for more general inclusions $\Gamma\subseteq G$. 
 
 Then, to analyze the specific states, concentrated in the fiber at $e$ we use properties specific to the dynamics of the action, by conjugation with elements in the group  $G$, on  the subgroups of $\Gamma$. One essential property of the inclusion $\Gamma=\PSL(2,\Bbb Z)\subseteq G=\PSL(2, \Bbb Z[1/p])$  is the following.  Let $\mathcal S_0$ be the subset of stabilizer groups (except the stabilizer of the identity), of the action of $G$, by conjugation on $\Gamma$. Then the groups in $\mathcal S_0$ are amenable. Moreover, the cosets of the groups in $\mathcal S_0$ are  asymptotically disjoint (that is they have finite intersections). 
 
 This will be used to prove  that in the realization of the essential states in infinite measures, acted by $G\times G^{\rm op}$, by measure preserving  transformation (the essential states are measuring displacement by translations in the group), the actions may be assumed to be free, with trivial stabilizers.

The proof is organized in the following steps:
 We prove in Appendix 6 that the analysis of the behavior of essential  states on the crossed product $C^\ast$ algebra may be reduced to states that are a convex combination of limits of mutually singular averaging sets of points in $\Gamma$. 

Using the elements of  Loeb measure theory developed in   ([Lo]), it follows that is sufficient to analyze states that are realized as the  measure of the displacement, due to   the measure preserving  action of the groupoid $(G\times G^{\rm op}) \rtimes K$, of a fixed finite measure subset $F$ in an infinite measure space $(\Y,\nu)$ (more precisely the space $\Y$ is a $G\times G^{\rm op}$-equivariant,  measurable bundle   over $K$ and we are computing $\nu (g_1Fg_2\cap F)$, for $ g_1, g_2 \in G$).

 If the groups are exact we may further reduce to the case when $F$ is a $\Gamma$ wandering set, whose translates by the groupoid action cover $\Y$. By using the action of $\Gamma$, we prove in Appendix 5,  by using $\Gamma$-equivalent subsets ([Ng]), that we may substitute $F$ with a subset that "sits" in the fiber of the neutral element of $K$. 
 
 Consequently, it is sufficient to analyze essential states that are obtained by limits of averaging sets contained in a family of normal subgroups shrinking to the identity (with trivial intersection). The state now becomes equivalent (through $\Gamma$ translations) to  a state concentrated on the $C^\ast$ algebra $C^\ast( G)$ 
($G$ is identified with the subgroup $\{g\times g^{-1}| g\in G\}$ of $(G\times G^{\rm op})$. 

In the  Appendix 6, we prove that if certain conditions of  temperedness (in the sense of the tempered Koopman representations in [Ke]) on the state  on  $C^\ast( G)$, constructed  above,  are verified, then the original state, given by the measure of the displacements of the set $F$ under the action of the groupoid 
$(G\times G^{\rm op}) \rtimes K$ is continuous with respect to the 
$C^{\ast}_{\rm red}((G\times G^{\rm op}) \rtimes C(K))$ topology, and hence so are the essential states. This is outlined in Corollary \ref {reduction}. 

What remains to be verified is that the state on $C^\ast(G)$ given by limits of averaging sets, shrinking to the identity, is continuous with respect to the $C^\ast_{\rm red}(G)$ topology (and is a limit of states having support in a finite reunion of double cosets of $\Gamma$ in $G$). This last statement (specific to $\PSL(2,\Bbb Z)\subseteq \PGL(2, \Bbb Z[1/p]))$ is proved in the statements \ref{conj1} through \ref{Sinfinit} bellow, in this chapter, by using the dynamics of the action by conjugation with elements in the group $G$ on subgroups of $\Gamma$.

We now start the exposition of the proof of Theorem \ref{AOlocal}.
On $\B(l^2(\Gamma))$ we consider the essential  states (forms) which are states (forms) that factorize to the Calkin algebra
$$
\mathcal{Q}(l^2(\Gamma)) = \B(l^2(\Gamma))/ \K(l^2(\Gamma)).
$$
By the original work of Calkin ([Ca]), to describe the essential states (forms), we let $\omega$ be any free ultrafilter on $\N$, and let $\mathop{\lim}\limits_{n \to \omega}$ be the corresponding ultrafilter limit on bounded sequences.

Then the essential forms are obtained as follows. Let $\xi = (\xi_n)$, $\eta = (\eta_n)$ be sequences of unit vectors in $l^2(\Gamma)$, converging weakly to $0$. For $X$ in $\B(l^2(\Gamma))$ define
$$
\varphi_{\xi, \eta, \omega}(X) = \mathop{\lim}\limits_{n \to \omega}\langle X \xi_n, \eta_n \rangle.
$$

Then the forms (respectively states) of the form $\varphi_{\xi, \eta, \omega}$ (respectively $\varphi_{\xi, \xi, \omega}$) exhaust the space of essential forms (respectively states) on $\B(l^2(\Gamma))$ (that is the forms (respectively states) that vanish on $\K(l^2(\Gamma))$).

In the Appendix 6 (Theorem \ref{points}), we prove that the analysis of these essential states, from the point of view of the topology induced on $\B$, is further reduced to the case when the sequence of vectors $\xi = (\xi_n)$ is of the following form:
Let $(A_n)_n$ be a family of finite sets in $\Gamma$, that eventually avoid any fixed, finite subset of $\Gamma$. Then consider $\xi_A = ((\card A_n)^{-1/2} \chi_{A_n})_n$, where by $\chi_{A_n}$ we denote the characteristic function of the set $A_n$.

Then it will be sufficient, to determine the topology on the $C^*$-algebra $\B$ induced by essential states of the form $\xi_{A}$.

The most general case corresponds to a countable family $A^s = (A_n^s)$, $s \in \N$ of such sets, disjoint (after $s \in \N $) for any fixed $s$ and $\xi = \mathop{\sum}\limits_{s} \dfrac{1}{2^s} \chi_{A_n^s}$, where the states $\chi_{A^s}$ are singular to each other as explained bellow.

By using the Loeb measure construction, we construct a probability measure space ($\cC_{\omega}(A)$, $\mu_{\omega, A}$), where $\cC_{\omega}(A)$ is the infinite product of the $A^{ s}_{n}$ and the probability measure $\mu_{\omega, A}$ is the ultrafilter limit of the normalized counting measure. In particular if $B_n \subseteq A_n$, $n \in \N$ is a sequence of subsets, letting $\cC_{\omega}((B_n))$ be the subset of all sequences $(a_n)_n \in \cC_{\omega}(A)$ that eventually belong to $B_n$, then
$$
\mu_{\omega, A}(\cC_{\omega}((B_n))) = \mathop{\lim}\limits_{n \to \omega} \dfrac{\card B_n}{\card A_n}.
$$

Then $\cC_{\omega}(A)$ is obviously a $C(K)$ module, simply by defining, for a coset $s\Gamma_0$ of a  $\Gamma_0$ a modular subgroup of $\Gamma$,  the action to be $\chi_{\overline{s\Gamma_0}}\cC_{\omega}(A) = \cC_{\omega}((A_n \cap s\Gamma_0)_n)$.

Moreover we may construct an (infinite) measure space acted by measure preserving transformations by $G \times G^{\rm op}$ as follows:

Let $\Y = \Y_{\omega, A}$ (which we will denote simply $\Y_\omega$ when no confusion is possible) be the reunion (as subsets of $\Gamma^{\aleph_0}$) of the sets $g_1\cC_{\omega}(A)g_2^{-1}$. Thus
$$\Y_{\omega, A} = \mathop{\bigcup}\limits_{g_1, g_2 \in G} \cC_{\omega}(\Big(g_1(\Gamma_{g_1^{-1}, g_2} \cap A_n){g_2}\Big)_n).$$

Because we are taking the counting measures, the corresponding ultralimit measures do coincide on the overlaps and hence we obtain a measure $\nu_{\omega}$ on $\Y$ that is invariant to the (partial) action of $G \times G^{\rm op}$. Note that $\Y$ remains a $C(K)$ module, and in fact in this way we obtained a measure space $(\Y, \nu_{\omega})$ acted by (partial) measure preserving transformations  of $G \times G^{\rm op}$. Hence  we obtain  a unitary Koopmann representation of $C^{\ast}((G \times G^{\rm op})\rtimes C(K))$ on $L^2(\Y, \nu_{\omega})$. Note that the absence of Folner sets automatically implies that $\nu_{\omega}(\Y)$ is infinite.

The goal of this section is to prove that this Koopmann representation is continuous with respect to the $C_{\rm red}^{\ast}((G \times G^{\rm op}) \rtimes C(K))$ norm and that it verifies the additional  conditions  (FS1), (FS2) of Theorem \ref{quotienthecke}.
We will then apply Corollary \ref{reduction}.

We  have thus  to analyze the states $\varphi_{\omega, A}$ on $C^{\ast}((G \times G^{\rm op}) \rtimes C(K))$ which on an element of the form $(g_1, g_2)\chi_{s \Gamma_0} \in C^{\ast}((G \times G^{\rm op}) \rtimes C(K)) $, where $s\Gamma_0$ is a coset of a modular subgroup of $\Gamma$, $g_1, g_2 \in \Gamma$ take the value
$$
\varphi_{\omega, A}((g_1, g_2)\chi_{s\Gamma_0}) = \varphi_{\omega, A}((g_1, g_2)\chi_{s\Gamma_0 \cap g_1^{-1}\Gamma g_2}) =
$$
$$
= \nu_{\omega}(g_1(s\Gamma_0 \cap g_1^{-1}\Gamma g_2)F g_2^{-1} \cap F) =
\mu_{\omega, A}(g_1(s\Gamma_0 \cap g_1^{-1}\Gamma g_2)F g_2^{-1} \cap F)=$$
$$= \mathop{\lim}\limits_{n \to \omega} \dfrac{\card (g_1(A_n \cap s\Gamma_0 \cap g_1^{-1}\Gamma g_2)g_2^{-1} \cap A_n)}{\card A_n}.
$$

In the above, the meaning of the notation $CF$, where $C$ is a coset of a modular group in $\Gamma$, is precisely the
$\cup\{ cF \ c \in {\overline C}\}$, where $\overline C$ is the closure in $K$ of $C$.

If more generally we have a family $A^s=(A_n^s)_n$, $s\in \Bbb N$,  such that the measures $\mu_{\omega, A^s}$ are all singular, then the state corresponding to the vector $\xi = \mathop{\sum}\limits_{s} \dfrac{1}{2^s}\xi_{A^s}$ will be in fact a direct sum of infinite measure space $(\Y^s, \nu_{\omega}^s)$, so the proof may be reduced to the case of a single family.

Note that because of the arguments in [Ra6], we may assume that $G \times G^{\rm op}$ acts freely on $\Y_\omega$.

In the Appendix 6, Corollary \ref{reduction}, we prove that we may further reduce the analysis of the continuity of the states $\varphi_{\omega, A}$, to the following more particular situation.  Recall that the group $\Gamma(p^{n})$ is the kernel of the surjection $\PSL(2, \Z) \to \PSL(2, \Z_{p^n}), n\in \N$. Then let  $\Gamma_n = \Gamma(p^{s_n}),n \in \N$, where $s_n$ is strictly increasing sequence of integers. Then we may further assume,  that $A_n \subseteq \Gamma(p^{s_n}),n\in \N$.

In this case $\varphi_{\omega, A}$ is simply a state on $C^{\ast}(G)$ (since it vanishes on $(g_1, g_2)f$, if $g_2 \neq g_1$) and the statement to be proved is that $\varphi_{\omega, A}$ is continuous with respect to the $C^{\ast}_{\rm red}(G)$ norm (more precisely that it verifies the conditions of Theorem \ref {quotienthecke}).

Denote by $F = \cC_{\omega}((A_n)_n)$, which is  finite  measure subset of $\Y_{\omega, A}$. Then 
$$
\varphi_{\omega, A}(g) = \nu_{\omega}(gF g^{-1}\cap F) = \mu_{\omega,A}(gF g^{-1}\cap F) = \mathop{\lim}\limits_{n \to \omega}\dfrac{\card (g A_n g^{-1} \cap A_n)}{\card A_n}.
$$

Note that in this context, the particular choice of the sets $(A_n)_n$ implies that $(\Y, \nu_{\omega})$ is an infinite measure space, acted by bijective, measure preserving transformations of $G$. We  prove bellow that the Koopmann representation is tempered (in the sense of Kechris ([Ke])),  i.e. that the representation is continuous with respect to the $C^{\ast}_{\rm red}(G)$ norm.

We will do this by proving that $\Y$ inherits a finer module structure  over the Borel $\ast$ - algebra generated by characteristic functions of subgroups of ${\PSL_2(\mathbf Z_p)}$.

For simplicity we denote the positive definite function $\varphi_{\omega, A}$ by $\varphi_A$ and the measure $\mu_{\omega,A}$ by $\mu_A$.   The positive definite function $\varphi_A$ is then computed by the formula
$$
\varphi_A (g)=\mu(gF\cap F), g \in G
$$

We  prove  that $\varphi_A$ belongs to $C^*_{\rm red}(G)$. Although we are not using the following remark, we note that ultimately,
to prove the Ramanujan--Pettersson conjecture (for the essential spectrum)
we are interested in the positive definite functional associating to the double coset
$[\Gamma\sigma\Gamma]$ the sum
$$
\Psi_A ([\Gamma\sigma\Gamma])=  \sum_{g\in[\Gamma\sigma\Gamma]}\varphi_A(g)
|\langle\pi_{13}(g)I,I\rangle|^2.
$$
It is easy to see that $\Psi_A$ is a positive definite functional on the reduced Hecke algebra of double cosets if  the positive definite function 
$\varphi_A(g)
|\langle\pi_{13}(g)I,I\rangle|^2$
 is a positive definite function on  $C^*_{\rm red}(G)$.
This might be useful for other groups $\Gamma, G$.

For $x$ in $\Gamma$ denote by $O^{\Gamma}_x$
(respectively $O^{G}_x$) the orbit of $x$, under the conjugation action,
by $\Gamma$ (respectively by $G$).

It is obvious that for $g\in\Gamma\sigma\Gamma$, $g O_x g^{-1}\cap\Gamma_n$ is non void and only if
$O_x$ intersects 
$x\in \Gamma_n\cap g^{-1}\Gamma_n g$, i.e., $O_x$ intersects
$\Gamma_n\cap \sigma^{-1}\Gamma_n \sigma$.
Thus for $e\in\{1,2,\ldots\}$ we may consider 
$$
A^e_n=\{a\in A_n \mid O^{\Gamma}_a
\hbox{ does not intersect } \Gamma_n \cap ({\sigma}_{p^{e+1}})^{-1}\Gamma_n {\sigma}_{p^{e+1}}\}.
$$




We let $F^e$ be the subset of $\cC_{\omega}(A)$ defined by
$F^e=\cC_{\omega}((A_n^e)_n)$, $e=0,1,2,\ldots$.
Let $F_{\infty}$ be defined by the formula
$$
F_{\infty}=\cC_{\omega}(A)\setminus\bigg[\bigcup\limits_e F^e\bigg].
$$
Then $F_{\infty}$ consists of all sequences $(a_n)_n$ in $\cC_{\omega}(A)$
such that for every integer $k$, the set $\{n\mid \O_{a_n}^{\Gamma}$
intersects $\Gamma_n \cap (\sigma_{p^k})^{-1} \Gamma_n \sigma_{p^k}\}$
is cofinal in $\omega$.

The sets $\bigcup\limits_e F^e$ and $F_{\infty}$ have disjoint $G$ orbits in $\Y_{\omega}$.

Moreover, $(\Gamma\sigma_{pf}\Gamma)F^e\cap F^e$ is non-void only if $f\leq e$.

Thus the states
$$
g\to\frac{1}{\mu(F_e)}\langle g(F^e),F^e\rangle
$$
are $C^*_{\rm red}(G)$ continuous and  verify the conditions of Theorem \ref{quotienthecke}
(here we use the fact that  the Akemann--Ostrand property holds true for the group $\Gamma$ ([AO], [Oz])), 
for all finite $e$. Hence the same holds true for the 
 state corresponding to the reunion
$\bigcup\limits_{e\in \N}F^e$, which consequently is  is a $C^*_{\rm red}(G)$  continuous state.

It remains to analyze the state corresponding to $F_\infty$. 
To prove that the state on $G$ corresponding to displacement of $F_\infty$ is continuous with respect to the $C^{\ast}_{\rm red}(G)$ norm and verifies the hypothesis (FS1) of Theorem \ref{quotienthecke}, we introduce the following definition. 

Afterwards we will prove that the conditions in the next definition hold true for the action of $G$ on $\Y_{\omega}$. Note that condition (FS2)  follows from N. Ozawa's papers ([Oz1], [Oz]). The proof of Theorem \ref{AOlocal} will the be completed by applying Corollary \ref{reduction} in Appendix 6.

\vskip6pt

\begin{defn} \label{conj1} 
Let $\mathcal{M}\S$ be the G-equivariant, $\ast$ Borel algebra of Borel functions on $\PGL_2(\Q_p)$ generated by the characteristic functions of conjugates by elements in $G$, of the group $K =PSL_2(\mathbf Z_p)$ (Recall that $\mathbf Z_p$ are the p-adic integers).   Then $\mathcal{M}\S$ contains all intersections $K \cap g K g^{-1}$, $g \in G$ (and infinite intersections of type the above). We let $G$ act on $\mathcal{M}\S$ by conjugation. By $\mathcal{M}\S\cap K$, we denote the Borel algebra obtained by intersecting all the sets in $\mathcal{M}\S$ with $K$.

Let $(\Y, \nu)$ be an infinite, measure space  and assume that the group $G = \PGL_2(\Z[\frac{1}{p}])$ acts by measure preserving transformations on $\Y$. We also fix  
 $F$  a finite measure subset of $\Y$, that is $\Gamma$ - wandering (i.e. $\nu(\gamma F \cap F) = 0$ for $\gamma \neq e, \gamma \in \Gamma$). We also assume that $\mathop{\bigcup}\limits_{g \in G}g F = \Y$.
 
We will say that the $G$ system $(\Y, \nu)$ is quasi-expanding ,if $\Y$ has a $\mathcal{M}\S$ module structure, that is $G$ - equivariant and that verifies the normalizing property $\chi_K(\Gamma F) = \Gamma F$.

\end{defn}

\ 

We will prove bellow (Lemma \ref{ms}) that the infinite measure space $(\Y_{\omega}, \nu_{\omega, A})$, constructed as above, starting with   the Loeb space $\cC_{\omega}((A_n))$, has the property in the previous definition.

We first prove a "nesting" property for the subgroups, whose characteristic function generate $\mathcal{M}\S$.

\begin{lemma}
For $g$ in $G$ let $K_g$ be the subgroup of $K$ given by $K_g = gKg^{-1}\cap K$. Then $K_g$ is uniquely determined by the coset $s\sigma_{p^e}\Gamma$ to which $g$ belongs.

Moreover there exists an order preserving equivalence between the cosets of $\Gamma\sigma_{p^e}$, $e \geq 1$ in $\Gamma$ and such subgroups: if $g$ belongs to $s\sigma_{p^e}\Gamma$, and $s\Gamma_{\sigma_{p^e}}$ is contained in $s_1\Gamma_{\sigma_{p^{e-1}}}$ then for any $g_1$ in $s_1\sigma_{p^{e-1}}\Gamma$ we have $K_g \subseteq K_{g_1}$.
\end{lemma}

\begin{proof}
This is equivalent to the corresponding property of the subgroups $\Gamma_g = g \Gamma g^{-1}\cap \Gamma$ of $\Gamma$ and this property is almost tautological.
Indeed $s\Gamma_{\sigma_{p^{e}}}s^{-1} = \Gamma_{\sigma_{p^e}}$ for $s \in \Gamma, e \geq 1$. 
On the other hand for $g$ in $G$, $\gamma \in \Gamma$ we
 have $\Gamma_{ g\gamma} = \Gamma_g$.

If $s$ belongs to $\Gamma_{\sigma_{p^{e-1}}}$ then $s\Gamma_{\sigma_{p^{e-1}}}s^{-1} = \Gamma_{\sigma_{p^{e-1}}}$ and hence
$$
\Gamma_{\sigma_{p^{e}}} = s\Gamma_{\sigma_{p^{e}}}s^{-1} \subseteq s\Gamma_{\sigma_{p^{e-1}}} = \Gamma_{\sigma_{p^{e-1}}}.
$$

\end{proof}

Because of the "nesting" property, it follows that any infinite intersection of sets in $\mathcal{M}\S$, reintersected with $K$, will contain a reunion of infinite intersections of the form $K \cap K_{s_1\sigma_{p^1}} \cap \ldots \cap K_{s_e\sigma_{p^e}} \cap \ldots$ where $s_e\Gamma_{\sigma_{p^e}}\subseteq s_{e-1}\Gamma_{\sigma_{p^{e-1}}}$ for all $e \geq 1$.

But such a decreasing sequence of cosets corresponds to a coset in $K/K_{\infty}$, where $K_{\infty}$ consists of the lower triangular matrices in K, that is the subgroup of  matrices of the form
$\left(
\begin{array}{cc}
a & 0 \\[3mm]
c & d
\end{array}
\right)$
in $K$.

Hence the intersection is determined uniquely by an element in the projective space $P^1(\mathbf Z_p^2)$.

Indeed if $s_e = 
\left(
\begin{array}{cc}
x_e & y_e \\[3mm]
z_e & t_e
\end{array}
\right)$, modulo the scalars, then $s_e\Gamma_{\sigma_{p^e}}$ is determined by $(y_e, t_e) \in P^1(\Z_{p^e}^2)$, and the nesting condition
$$
s_e\Gamma_{\sigma_{p^e}} \subseteq s_{e-1}\Gamma_{\sigma_{p^{e-1}}}
$$

\noindent corresponds to the fact that for $e \geq 1$, $(y_e, z_e) \equiv (y_{e-1}, z_{e-1})$ in $P^1(\Z_{p^{e-1}} ^2)$.

We analyze now the structure of infinite intersections.

\begin{lemma}
Denote the infinite intersection, described above, corresponding to $(y, t) \in P^1(\mathbf Z_p^2)$ by $K_{(y, t)}$.

Given 2 distinct points $(y_1, t_1)$ and $(y_2, t_2)$ in $P^1(\mathbf Z_p^2)$, the intersection $K_{(y_1, t_1)}\cap K_{(y_2, t_2)}$ will reintersect a third $K_{(y_3, t_3)}$, with $(y_3, t_3)$ in $P^1(\mathbf Z_p^2)$, different from the previous two, in the trivial element.
\end{lemma}

\begin{proof}
By left translations by elements in $K$, we may assume that $(y_1, t_1) = (0, 1)$ in $P^1(\mathbf Z_p^2)$ and thus $K_{0,1} = K\cap  \mathop{\bigcap}\limits_{e \geq 1} K_{\sigma_{p^e}} = K_{\infty}$.

Assume that 
$\left(
\begin{array}{cc}
x_2 & y_2 \\[3mm]
z_2 & t_2
\end{array}
\right)$
is a representative in $K=\PSL_2(2, \mathbf Z_p)$ of the coset of $K/K_{\infty}$ represented by $(y_2, t_2) \in P^1(\mathbf Z_p^2)$.
 
Then $K_{y_2, t_2}$ is 
$\left(
\begin{array}{cc}
x_2 & y_2 \\[3mm]
z_2 & t_2
\end{array}
\right)
K_{\infty}
\left(
\begin{array}{cc}
x_2 & y_2 \\[3mm]
z_2 & t_2
\end{array}
\right)^{-1}$ 
and hence $K_{(0, 1)} \cap K_{(y_2, t_2)}$ consists of all elements 
$\left(
\begin{array}{cc}
a & 0 \\[3mm]
c & d
\end{array}
\right)$
in $K_{(0, 1)}$ such that
$$
\left(
\begin{array}{cc}
x_1 & y_1 \\[3mm]
z_1 & t_1
\end{array}
\right)^{-1}
\left(
\begin{array}{cc}
a & 0 \\[3mm]
c & d
\end{array}
\right)
\left(
\begin{array}{cc}
x_1 & y_1 \\[3mm]
z_1 & t_1
\end{array}
\right)
$$.

This condition becomes in $\mathbf Z_p$
$$
y_1t_1(a-d) = y_1^2c.
$$

Thus, if $(0, 1) \neq (y_1, t_1)$ in $P^1(\mathbf Z_p^2)$, the intersection $K_{(0, 1)} \cap K_{(y, t)}$ is, :
$$
\left\{ 
\left(
\begin{array}{cc}
a & 0 \\[3mm]
c & d
\end{array}
\right) 
\in \PSL_2(2, \mathbf Z_p) \mid t_1(a-d) = y_1c\right\}.
$$

Clearly this can reintersect $K_{(0, 1)} \cap K_{(y_2,t_2)}$ in a non-trivial element if and only if $(y_2, t_2) = (y_1, t_1)$ in $ P^1(\mathbf Z_p^2)$. 
\end{proof}

In the following we describe the $\mathcal{M}\S$ module structure on the measure space $(\Y_{\omega}, \nu_{\omega})$ introduced at the beginning of the proof the theorem. Recall that the group $\Gamma(p^{n})$ is the kernel
 of the surjection $\PSL(2, \Z) \to \PSL(2, \Z_{p^n})$.

\begin{lemma}\label{ms}

For a family of a subgroups $H_n$ of $\Gamma$ let $\cC_{\omega}((H_n))$ consist of all sequences $(\gamma_n)_n$, such that $\gamma_n$ belongs to $H_n$ eventually, with respect to the ultrafilter $\omega$.

Let $s_n$ be a strictly increasing sequence of natural numbers and let $\Gamma_n = \Gamma(p^{s_n})$. Then $(\Gamma_n)$ is a decreasing sequence of normal subgroups of $\Gamma$, with trivial intersection.

Let $\mathcal{M}\S_{\omega}((\Gamma_n))$, which, for simplicity, when no confusion is possible, we denote by $\mathcal{M}\S_{\omega}$, be the countable Borel algebra (of functions on $\Gamma^{\aleph_0}$) generated by the characteristic functions of $\cC_{\omega}((\Gamma_n)_n)$ and their conjugates
$$
\cC_{\omega}((g\Gamma_ng^{-1})_n) = g\cC_{\omega}((\Gamma_n)_n)g^{-1}, g \in G.
$$

Then there exists a $\ast$ homeomorphism from $\mathcal{M}\S_{\omega}$ onto
$\mathcal{M}\S$, uniquely determined by the following requirements. 

(a). The morphism  is $G$ - equivariant and  maps the characteristic function of $\cC_{\omega}((\Gamma_n))$ into the characteristic function of $K = \PSL(2, \mathbf Z_p)$. 

(b). For unicity purposes, we require that the kernel  of the above $\ast$ homeomorphism, restricted to subgroups of $\cC_{\omega}((\Gamma_n))$, which are mapped into to subgroups of $K$, is   contained in the space of the characteristic function of
$$
\mathop{\bigcap}\limits_{e \geq 1}\cC_{\omega}((\Gamma(p^{s_n+e}))_n).
$$

Then the space $(\Y_{\omega}, \nu_{\omega})$ has a canonical, $G$ - equivariant, $\mathcal{M}\S_{\omega}$ module structure, and $\chi_{\cC_{\omega}((\Gamma_n))}(\Gamma F) = \Gamma F$. 

\end{lemma}

\begin{proof}
Since every intersection of finite index subgroups is again a finite index subgroup, it follows that if $(H_n^s)_n$ , $s \in \N$ is an infinite collection of decreasing sequences of finite index subgroups of $\Gamma$, then $\mathop{\bigcap}\limits_{s}\cC_{\omega}((H_n^s)_n)$ is always non trivial, as it contains
$$
\cC_{\omega}((H_n^1 \cap H_n^2 \cap \ldots \cap H_n^n)_n).
$$

Hence the only problem in establishing the homeomorphism from $\mathcal{M}\S_{\omega}$ onto
$\mathcal{M}\S$ will consist in determining the kernel.

To do this observe that the nesting property proven for the subsets $K_g, g \in G$ also holds true for the groups
$$
A_g = \cC_{\omega}((\Gamma_n \cap g\Gamma_ng^{-1})_n) = \cC_{\omega}((\Gamma_n))\cap g\cC_{\omega}((\Gamma_n)_n)g^{-1}.
$$
Indeed it is obvious that if $g$ belongs to $s\sigma_{p^e}\Gamma$, then $A_g$ depends only on $s\sigma_{p^e}$. Indeed $A_{g\gamma} = A_g$ for all $g \in G, \gamma \in \Gamma$ since computing $A_{g\gamma}$ corresponds to conjugate $\Gamma_n$ by $\gamma$, but the conjugate is  again $\Gamma_n$, since the subgroups $\Gamma_n$ are normal.

We also have to prove that if $[s_e\Gamma_{\sigma_{p^e}}]$ is contained in $[s_{e-1}\Gamma_{\sigma_{p^{e-1}}}]$, where $s_e, s_{e-1} \in \Gamma$, $e \geq 1$ then $A_{s_e\sigma_{p^e}} \subseteq A_{s_{e-1}\sigma_{p^{e-1}}}$.

It is obvious that
$$
A_{sg} = sA_{g}s^{-1}.
$$
Hence to prove the inclusion it is sufficient to assume that $s$ belongs to $\Gamma_{\sigma_{p^{e-1}}}$ and to prove that $sA_{\sigma_{p^e}}s^{-1} \subseteq A_{\sigma_{p^{e-1}}}$.
But if $s \in \Gamma_{\sigma_{p^{e-1}}}$ then $s\sigma_{p^e} = \sigma_{p^{e}}\T'$ for some $\T'$ in $\Gamma$ and hence
$$
s\sigma_{p^{e-1}}\Gamma_n(\sigma_{p^{e-1}})^{-1}s^{-1} \cap \Gamma_n = \sigma_{p^{e-1}}\T'\Gamma_n(\T')^{-1}(\sigma_{p^{e-1}})^{-1} \cap \Gamma_n =
$$
$$
= \sigma_{p^{e-1}}\Gamma_n(\sigma_{p^{e-1}})^{-1} \cap \Gamma_n.
$$
Thus $sA_{\sigma_{p^{e-1}}}s^{-1} = A_{\sigma_{p^{e-1}}}$ and hence, since $A_{\sigma_{p^e}} \subseteq A_{\sigma_{p^{e-1}}}$ (by the choice we made for the groups $\Gamma_n$) it follows that $sA_{\sigma_{p^e}}s^{-1} \subseteq A_{\sigma_{p^{e-1}}}$.

Thus, as in the case of subgroups in $\mathcal{M}\S$, any infinite intersection of subgroups in $\mathcal{M}\S_{\omega}$, when intersected with $\cC_{\omega}((\Gamma_n)_n)$, will contain a reunion of infinite intersections of the form
$$
\cC_{\omega}((\Gamma_n)_n) \cap A_{s_1\sigma_p} \cap \ldots \cap A_{s_e\sigma_{p^e}} \cap \ldots \leqno(\ast)
$$

\noindent where $[s_1\Gamma_{\sigma}] \supseteq [s_2\Gamma_{p^2}] \supseteq \ldots \supseteq [s_e\Gamma_{\sigma_{p^e}}]$, and $s_e \in \Gamma, e \geq 1$.

Again this will depend only on a coset of a point in $P^1(\mathbf Z_p^2)$ that in turn determines a coset of $K/K_{\infty}$.
We denote the infinite intersection in formula ($\ast$) corresponding an element $(y, t) \in P^1(\mathbf Z_p^2)$ (which in turns corresponds to $[s_1\Gamma_{\sigma}] \supseteq [s_2\Gamma_{\sigma_{p^2}}] \supseteq \ldots \supseteq [s_e\Gamma_{\sigma_{p^e}}]$) by $K_{(y,t)}^{\omega}$.

We will verify the same property of intersection for this class of subgroups as the one holding for the  for subgroups in $\mathcal{M}\S$. We check that $K_{(y_1,t_1)}^{\omega} \cap K_{(y_2, t_2)}^{\omega} \cap K_{(y_3, t_3)}^{\omega}$ is contained in the kernel of the morphism from $\mathcal{M}\S_{\omega} \cap \cC_{\omega}((\Gamma_n))$ onto $K = \PSL(2, \mathbf Z_p)$.

Indeed to check this we may assume that $(y_1, t_1) = (0, 1)$ in $P^1(\mathbf Z_p^2)$.
Thus assume representatives for $(y_2, t_2), (y_3, t_3)$ are 
$\left(
\begin{array}{cc}
x_2 & y_2 \\[3mm]
z_2 & t_2
\end{array}
\right)$
and
$\left(
\begin{array}{cc}
x_3 & y_3 \\[3mm]
z_3 & t_3
\end{array}
\right)$
and 
$$
K_{(y_1, t_1)}^{\omega} = K_{(0, 1)}^{\omega} = \cC_{\omega}((\Gamma_n)) \cap \mathop{\bigcap}\limits_{e\geq 1}\cC_{\omega}(\Gamma_n \cap \sigma_{p^e}\Gamma_n\sigma_{p^e}^{-1}).
$$
Assume that $[s_e^j\Gamma_{\sigma_{p^e}}]$ are the decreasing sequence of cosets that determine $K_{(y_j, t_j)}^{\omega}$, and thus we may assume $s_e^j =
\left(
\begin{array}{cc}
x_e^j & y_e^j \\[3mm]
z_e^j & t_e^j
\end{array}
\right)$,
$e \geq 1$, $j = 1,2$, where the sequence $(y_e^j, t_e^j)$ in $P^1(\Z_{p^e}^2)$ represents $(y_j, t_j)$ in $P^1(\mathbf Z_p^2)$.

Then $K_{0,1}^{\omega} \cap K_{(y_j, t_j)}^{\omega}$, for $j = 1,2$, by the same computations that we have  performed for the subgroups of $\PSL(2, \mathbf Z_p)$, consists of  the group of sequences:
$$
\left\{ 
\left(
\begin{array}{cc}
a_n & b_n \\[3mm]
c_n & d_n
\end{array}
\right) 
\in \cC_{\omega}((\Gamma_n)_n) \mid b_n \equiv 0, y_j^et_j^e(a_n - d_n) \equiv (t_j^e)^2, {\rm mod} \; p^{s_n + e}. 
\right\}
$$

Because $(y_e^j, t_e^j)_e$, $j = 1,2$ in the $p$ - adic completion correspond to different elements in $(y_j, t_j)$ in $P^1(\mathbf Z_p^2)$, the triple intersection will be contained in
$
\left\{ 
\left(
\begin{array}{cc}
a_n & b_n \\[3mm]
c_n & d_n
\end{array}
\right) 
\in \cC_{\omega}((\Gamma_n)_n) \mid a_n \equiv d_n \; c_n \equiv 0 \; ({\rm mod} \; p^{s_n + e - f}), b_n \equiv 0 \; ({\rm mod} \; p^{s_n + e})
\right\}
$, where $f$ depends on which power of $p$ divides $(y_j, t_j)$.

Replacing $e$ by $e+f$, when necessary, this is contained in the required kernel.

To complete the proof we note that because of this argument, the only  non-trivial intersections of subgroups in $\mathcal{M}\S_{\omega}\cap K$ are the intersections $K_{(y_1, t_1)}^{\omega} \cap K_{(y_2, t_2)}^{\omega}$ which may also be intersected by finite intersection of the form $$\mathop{\bigcap}\limits_{i = 1}^{r}\cC_{\omega}((\Gamma_n \cap g_i\Gamma_ng_i^{-1})_n),$$ where $g_1, g_2, \ldots, g_r$ belongs to $G$.

The $\mathcal{M}\S_{\omega}$ structure on $\Y_{\omega}$ is now simply the appurtenance relation, defined by the fact that points in $\Y_{\omega}$ are sequences $(a_n)_n$ in $\Gamma$. Thus  the characteristic function of $\cC_{\omega}((H_n)_n)$ will multiply $(a_n)_n$ by 1 (or 0) if $\{ a_n \in H_n \}$ is cofinal in $\omega$ (respectively is not cofinal).

\end{proof}

In the above terminology the set $\Gamma F_{\infty}$ is contained in $$\mathop{\bigcap}\limits_{e \geq 1}(\mathop{\bigcup}\limits_{i}\cC_{\omega}\Big( s_i^e\sigma_{p^e}\Gamma_n(\sigma_{p^e})^{-1}(s_i^e)^{-1} \cap \Gamma_n\Big)),$$ where $s_i^e$ are the coset representatives for $\Gamma_{\sigma_{p^e}}$ in $\Gamma$, for $e \geq 1$.

The characteristic function in $\mathcal{M}\S$, acting identically on $\Gamma F_{\infty}$, via the module structure, is the characteristic function of the set:
$$
S_{\infty} = \mathop{\bigcap}\limits_{e}(\mathop{\bigcup}\limits_{i} K_{s_i^e\sigma_{p^e}}).
$$

We have thus proved that modulo the trivial element of $K$, The set $S_{\infty}$ is a reunion of sets of the form  $K_{(y_1, t_1)} \cap K_{(y_2, t_2)} \cap \mathop{\bigcap}\limits_{i =1}^{r} K_{g_i}$ where $(y_1, t_1), (y_2, t_2)$ are distinct points in $P^1(\mathbf Z_p^2)$ and $g_1, g_2, \ldots , g_r$ belong to $G$.

The similar statement holds true in $\mathcal{M}\S_{\omega}$ (modulo the kernel).

Note that, (in the terminology  introduced at the end of the proof above),
 $\Gamma F$ is indeed contained in $\cC_{\omega}((\Gamma_n))$, because $F$ is contained in $\cC_{\omega}((\Gamma_n)_n)$ and all subgroups in $(\Gamma_n)_n$ are normal.
 
 \begin{cor}\label{Sinfinit}
There exists a $\mathcal{M}\S \cap K$, $G$ - equivariant module structure on $\Gamma F$.

Moreover, let $\mathbf{m}\mathbf{s}$ be the Gelfand spectrum of $\mathcal{M}\S \cap K$, and let $p$ be the corresponding, $G$ - equivariant, projection from $\mathcal{M}\S \cap K$ onto $\mathbf{m}\mathbf{s}$. 
Then $p(F_{\infty}) \subseteq S_{\infty}$. Hence the dynamics of the action of the group $G$, on $F_{\infty}$ (e. g. the precise movement of the subsets of $F_{\infty}$ that are brought back to into $F_{\infty}$ by the action of $G$) is determined by the dynamics of the (conjugation) action of $G$ on $\mathcal{M}\S \cap K$.
Moreover $$S_{\infty} = \mathop{\bigcup}\limits_{(y, t) \in P^1(\mathbf Z_p^2)}K_{(y, t)} = \mathop{\bigcup}\limits_{s \in K/K_{\infty}} sK_{\infty}s^{-1}.$$

The only possible intersections of subgroups in $\mathcal{M}\S \cap K$ yielding a nontrivial intersections (that is not equal to identity)  are $$K_{(y, t)} \cap K_{(y_1, t_1)} \cap K_{g_1} \cap \ldots \cap K_{g_r},$$ where $(y, t), (y_1, t_1)$ are distinct elements in $P^1(\mathbf Z_p^2)$ and $g_1, g_2, \ldots, g_r$ are elements in $G$.

Moreover $g(K_{(y, t)})g^{-1} \cap K$ is non-trivial if and only if $(y, t)$ corresponds to $s \in K/K_{\infty}$ with the property that $K_{y, t} = sK_{\infty}s^{-1}$. In this case necessary $g$ is of the form $s\sigma_{p^e}$ for some $\sigma_{p^e}$, $e \geq 1$.

In addition $\sigma_p(K_{(0,1)}\cap K_{(y, t)})$ is $K_{0, 1} \cap K_{y, pt} \cap \sigma_p(K_{0,1})$.

\end{cor}

\begin{proof}
The fact that there exists such a bimodule structure follows from the previous Lemma.

The only part of the statement that was not yet proved is the statement about $g(K_{y, t})g^{-1}\cap K$.

To prove this we may assume that $(y, t) = (0, 1)$ in $P^1(\Z_p^2)$ and hence we are analyzing the set
$$
L_{\infty} = K \cap g(K \cap \sigma_p K\sigma_{p}^{-1} \cap \ldots \cap \sigma_{p^e}K(\sigma_{p^e})^{-1}\cap \ldots)g^{-1}.
$$

But, unless $g$ is of the form $s\sigma_{p^e}$ for some $e \geq 1$, the intersection is trivial. In the non-trivial case  the intersection is
$$
L_{\infty} = s\sigma_{p^e}(K_{\infty})(\sigma_{p^e})^{-1}s^{-1}.
$$

The last computation is trivial.

\end{proof}

 Hence we observe that the subset of $S_{\infty}$, that is brought back into $S_{\infty}$ by the action by elements in the group $G$, is $$\mathop{\bigcup}\limits_{\gamma \in \Gamma/K_{\infty} \cap \Gamma} \gamma K_{\infty}\gamma^{-1}.$$

To conclude the proof of Theorem \ref{AOlocal}, we observe that, by using the same arguments as  Proposition \ref{shrink} (Appendix 6) we may replace $F_{\infty}$ by a measurable subset $F_{\infty, 1}$, which $\Gamma$ equivalent to $F_{\infty}$ ([Ng]), and such that $F_{\infty, 1} \subseteq K_{\infty}$.

In this case the only movements by $G$ that bring back pieces of $F_{\infty, 1}$ are those implemented by $K_{\infty}\cap G$, which is an amenable group. Moreover because of the last statement of Corollary \ref{Sinfinit}, this action verifies the conditions of Theorem  \ref{quotienthecke} (Appendix 4).

The  remaining case is the analysis of the case in which the set $\Gamma F$ has a part sitting in the kernel of the morphism from $\mathcal{M}\S_{\omega}$ onto $\mathcal{M}\S$.

But the orbits of $G$ through points in the kernel are returning to the kernel, and hence the dynamics under the action of $G$  of the subset $p^{-1}(e)\cap F$ (where $p$ is the projection from Corollary \ref{Sinfinit}) might be analyzed separately. But this subset of $F$ corresponds to a finer selection of the groups $\Gamma_n$. If we require that the original set $F$ has effective mass in $(\Gamma_n)$ (i.e. that the sequence of subgroups $(\Gamma_n)_n$ is minimal for $F$), then we may proceed by transfinite induction on smaller sets of normal subgroups shrinking to $e$.

This completes the proof of Theorem \ref{AOlocal}.

\end{proof}

\  


Although this is not needed for the proof, we note that if the sets $(A_n)_{n\in N}$ are equidistributed in the coset representatives, so that the measure $\mu_A$ is the Haar measure $\mu_p$ on $K$ then one can obtain an explicit formula for the essential states.

\begin{prop}We use the notations from the previous theorem.
Let $B_t=\{g \mid g \in \PSL_2(\R), \|g\|_2 \leq t\}$ be the hyperbolic ball in $\PSL_2(\R)$ of radius $t$. Because of 
the work of  Gorodnik and Nevo ([GoNe]), see also [EM], [DRS]), it follows that the sets $\Gamma_t=\Gamma\cap B_t$ are equidistributed in the cosets of modular subgroups in $\Gamma$.  We let the sets $A_n$ be defined by the formula $A_n=\Gamma\cap B_n$. Then the measure $\mu_A$ (constructed in the previous theorem) induces the Haar measure on $K$.  Moreover the state $\phi_p$ on $\A_{\text{max}}$ corresponding to this choice of the
sets $A_n$ is given by the formula
$$\phi_p=\sum\limits_{(g_1\times g_2)\in G\times G^{\rm op}} F(g_1,g_2)
\chi_{\overline{\Gamma \cap g_1^{-1}\Gamma( g_2)^{-1}}} d\mu_p\times d\mu_p \ (g_1\times g_2).$$
Here $\chi_{\overline{\Gamma \cap g_1^{-1}\Gamma( g_2)^{-1}}}$ is the characteristic function of the closure of $\Gamma \cap g_1^{-1}\Gamma( g_2)^{-1}$ in $K$ and 
 $F$ is a numerical, positive definite function on $G\times G^{\rm op}$, depending only on $\|g_1\|_2, \|g_2\|_2$, of the order of
$$\frac{\ln \|g_1\|_2+\ln \|g_2\|_2}{ \|g_1\|_2\cdot\|g_2\|_2}. $$

More precisely 
$F(g_1,g_2)$ is the asymptotic displacement of the  family of well
rounded balls $B_t$ (as in [GoNe]), that is 
$$
F(g_1,g_2)=\lim_{t\to\infty} \frac{{\rm vol}(B_t \cap g_1B_tg_2)}{{\rm vol}(B_t)}
$$
(volumes computed with respect Haar measure on $G$).

Then $\phi_p$ is state on the reduced C*-crossed product. Indeed if $\Psi$ is the completely  positive map on $C^*(G)$  mapping an element $g\in \Gamma\sigma\Gamma$ into $\frac {1}{[\Gamma : \Gamma_\sigma]} g$, then viewing $\Psi \otimes \Psi$ as a map on 
$\A_{\max}$, then $\phi_p\circ  (\Psi \otimes \Psi)^\epsilon$ is  square summable for any
$\epsilon >0$ (see the $L^{2+\epsilon}$ summability criteria in 
([DeCaHa]).

\end{prop}

\begin{proof}Because the points in $A_n$ are equidistributed in cosets it follows that the measure $\mu_{\omega, A}$ from the Theorem \ref{AOlocal} is absolutely continuous with respect to the Haar measure $\mu_p$ on $K$.
It follows that for every $g=(g_1,g_2)$ in $G\times G^{\rm op}$, there exists   for $(g_1\times g_2)\in G\times G^{\rm op}$, there exists a density $\theta_{A,g^{-1}A}$ a measurable function $K$, computing the displacement:
 $$\phi_A(g_1, g_2)=\int_K \theta_{A,g^{-1}A}\d \mu_p.$$ Moreover 
 $\theta_{A,g^{-1}A}$  is equal to the limit, for  $\Gamma_{\sigma_{p^e}}s$ a modular subgroup coset with closure $K(p^e,s)$, of the following expression. 
 $$
 \frac{1}{\mu_p(K(p^e,s))} \int_{K(p^e,s)}\theta_{A,g^{-1}A}(\omega)\d\mu_p(\omega)=
 \lim _{t\to \infty} \frac  {\text {card}(g_1\Gamma_tg_2\cap \Gamma_t\cap \Gamma_{\sigma_{p^e}}s)}
 {\text{card}(\Gamma_t \cap  \Gamma_{\sigma_{p^e}}s)}.$$
 
 Since the sets $g_1B_tg_2\cap B_t$, $g_1,g_2\in \PSL_2(\R)$, $t>0$, ([GoNe]) are well rounded, it follows that this is equal to 
 $$ \lim _{t\to \infty} \frac  {\text {card}(g_1\Gamma g_2\cap \Gamma \cap (\Gamma_{\sigma_{p^e}}s)\cap (B_t\cap g_1B_tg_2 ))}
 {\text{card}(\Gamma \cap B_t \cap  \Gamma_{\sigma_{p^e}}s)}.$$
 
 For a large exponent $e$, the above quantity is non-zero, if and only if the coset $\Gamma_{\sigma_{p^e}}s$ is contained in
  $g_1^{-1}\Gamma g_2^{-1}\cap \Gamma$, and hence it follows, by ([GoNe]), that 
 $\theta_{A,g^{-1}A}$ is given by a constant density with respect to $\mu_p$, supported on the closure in $K$ of
 $g_1^{-1}\Gamma g_2^{-1}\cap \Gamma$, of weight 
 $$F(g_1,g_2)=\lim _{t\to \infty} \frac  {\text {vol\ }  (B_t\cap g_1B_tg_2 )} { \text {vol\ } B_t}.$$
 Here $vol$ stands for the volume computed with respect to Haar measure on $\PSL_2(\R)$.
 Note that $F$ is in itself a positive definite function on $\PSL_2(\R)\times \PSL_2(\R)^{\rm op}$, but we are only interested in values of $F$ at
 $(g_1\times g_2)\in G\times G^{\rm op}$, whenever $g_1,g_2$ determine the same double coset of $\Gamma$ in $G$ (so that   $g_1^{-1}\Gamma g_2^{-1}\cap \Gamma$ is non-void).

 To finish the proof we have to find the order of growth of $F$. To do this we switch to $SU(1,1)$ instead of 
 $\PSL_2(\R)$. Assume $g_1,g_2\in SU(1,1)$ are given by:
 $$g_1= \begin{pmatrix} x & y \\ \overline {y}  &\overline {x} \end{pmatrix}, \quad g_2 =
 \begin{pmatrix} s & r \\ \overline {r}  &\overline {s} \end{pmatrix}.$$
 Since $B_t$ is invariant to left and right multiplication by unitaries, it follows that $F(g_1,g_2)$ only depends on
 $\|g_1\|_2, \|g_2\|_2$ and hence we may assume that the numbers $x,y,s,t$ are all positive. 
 
 We have to compute the relative  volume (with respect to the volume of $B_t$), as $t$ tends to infinity, of the intersection $g_1B_t\cap B_t g_2^{-1}$. Using the $KAK$ decomposition of $SU(1,1)$, and the corresponding Haar measure, we have to compute the volume of the set
 $$\Big\{\begin{pmatrix} a & b \\ \overline {b}  &\overline {a} \end{pmatrix} \in \text{SU\ }(1,1) \mid |xa+y\overline b| \leq t,
 |as+br| \leq t\Big\}.$$
 
 We denote $|a|=\rho $, $a =\rho \exp{i\theta_1}$, $b=|b| \exp {i\theta_2}$. Since we are interested only in the asymptotic ratio of volumes as $t$ tends to infinity, we may substitute $|b|=\sqrt {|a|^2-1}$ with $|a|$ and we may 
 replace the Haar measure on $\text{SU\ }(1,1)=KAK$, $\text{d}\mu_{\text{SU\ }(1,1)}= \text{d}k_1\text{d} |a| \text{d} k_2=\text{d}k_1 (\cosh^2\alpha )\text{d} \alpha \text{d} k_2$ (where $|a|=\cosh \alpha$) with 
 $ \text{d}\theta_1 \rho \text{d} \rho  \text{d} \theta_2$ (since we are interested only is asymptotic relative size of volumes).
 
 Hence the formula for $F(g_1,g_2)$ is 
 
 $$\int _{-\pi}^{\pi}\int _{-\pi}^{\pi}\int _0^{\text {min}(\frac{1}{|x\exp {i\theta_1}+y|}, \frac{1}{|s\exp {i\theta_2}+r|})}\rho \text {d} \rho\text {d} \theta_1 \text {d} \theta_2,$$ which up to a constant is
 $$\frac {1}{x^2 s^2} \int _{-\pi}^{\pi}\int _{-\pi}^{\pi}  \text {min}
\left (\frac{1}{|\exp {i\theta_1}+\frac{y}{x}|^2}, \frac{1}{|\exp {i\theta_2}+\frac {r}{s}|^2}\right)\text {d} \theta_1 \text {d} \theta_2.$$
 Denote $\alpha= y/x$  and $\beta=r/s$ and note that these two quantities are of the order of respectively $1/x^2$ and $1/s^2$ . Using arclenght approximation it follows that the integral is of the order of 
 $$\frac {1}{x^2 s^2}\int _{-1}^{1} \int _{-1}^{1}\text {min}
 \left(\frac{1}{\alpha^2+\theta_1^2}, \frac{1}{\beta^2+\theta_2^2}\right)\text {d} \theta_1 \text {d} \theta_2.$$
The result follows then by a straightforward computation.
 \end{proof}
 

In the following remark we describe a method to avoid the use of the cocycle $\varepsilon$ from the projective representation, by passing to a $\Z_2$ cover.

 \vskip6pt

\begin{rem}\label{z2}
{\it Assume that $\Gamma \subseteq G$ is an almost normal subgroup of $G$. We assume that 
$G$ is presented in the following way:
(Here we assume $\Z_2$ is mapped into the center of $\tilde{G}$.)
$$
\begin{array}{ccccccccc}
0 & \longrightarrow & \Z_2 & \longrightarrow & \tilde{G} & \longrightarrow & G & \longrightarrow & 0\\
&      &  ||  & &  \cup & & \cup \\
0 & \longrightarrow & \Z_2 & \longrightarrow & \tilde{\Gamma} & \longrightarrow & \Gamma & \longrightarrow & 0
\end{array}
$$

Let $u$ be the  image of the non-identity element of $\Z_2$ in $\tilde{G}$
(or which is the same, in  $\tilde{\Gamma}$). Then we assume that $u$ is central element in $\tilde{G}$.

In the group algebra of $\tilde{G}$ let $P=1-u$, which is a projection (corresponding to the negative part of $u$).
We consider the reduced algebra $\A_P=P\L (\tilde{G})P\supseteq P\C(\tilde{G})P$.
$(P\C(\tilde{G})P$ is like the group algebra of $\C(\tilde{G})$ modulo the identity $u=-P$,
$P$ being the identity of the reduced algebra.)

A similar construction the one in the preceding chapters, can be done in this setting, as follows.

Let $H_P=L^2(\A_P,\tau_P)$, where $\tau_P$ is the reduced trace. The group $\tilde{\Gamma}$
acts by left and right convolutors
$L_{\tilde{\gamma}}$, $R_{\tilde{\gamma}}$, $\tilde{\gamma}\in\tilde{\Gamma}$ on $H_P$ and we
obviously have
$L_{u\tilde{\gamma}}=-L_{\tilde{\gamma}}$, $R_{\tilde{\gamma} u}=-R_{\tilde{\gamma}}$ $(\tilde{\gamma}\in\tilde{\Gamma})$.
Assume $\tilde{\sigma}\in\tilde{G}$ and let
$\tilde{\Gamma}_{\tilde{\sigma}}=\tilde{\sigma}\tilde{\Gamma}\tilde{\sigma}^{-1}\cap\tilde{\Gamma}$.
Then every $\tilde{X}$ in $B(L^2(\A_P,\tau_P))$, such that
$\tilde{X}L_{\tilde{\gamma}_0}=L_{\tilde{\sigma}\tilde{\gamma}_0\sigma^{-1}}\tilde{X}$
for $\tilde{\gamma}_0\in\tilde{\Gamma}_{\sigma^{-1}}$
will give raise to a completely positive map $\Psi_{\tilde{X}}$ obtained from the following diagram
$$
\begin{array}{ccc}
(P\L(\tilde{\Gamma}_{\tilde{\sigma}^{-1}})P)' & \stackrel{X^*\cdot X}{\longrightarrow} &  (PL(\tilde{\Gamma}_{\tilde{\sigma}})P)' \\[4pt]
\phantom{aa}\supseteq & \circlearrowright \phantom{a}&\hspace{-0.3cm} \swarrow {}_E \\[4pt]
&  (P\L(\tilde{\Gamma})P)'&
\end{array}.
$$
Here the commutants are computed in the Hilbert space $H_P$ and  $E$ is the canonical conditional expectation.

If we start with a representation $\tilde{\pi}$ of $\tilde{G}$ on $L^2(\A_P,\tau_P)$ extending the left
regular representation of $\tilde{\Gamma}$ on $L^2(\A_P,\tau_P)$ (thus $\pi(u)=-1$), then we can
construct as before 
$$
\tilde{T}^{[\tilde{\Gamma}\tilde{\sigma}\tilde{\Gamma}]}=
\frac12\bigg(\sum_{\theta\in[\tilde{\Gamma}\tilde{\sigma}\tilde{\Gamma}]}
\overline{\langle\pi(\theta)P,P\rangle}\,\theta
\bigg).
$$
The factor $\frac12$ is needed because when reducing by $P$ the terms
$\langle\pi(\theta)P,P\rangle\, \theta$ and $\langle\pi(\theta u)P,P\rangle\, \theta u$ 
correspond to the same term.

 Assuming that $\tilde{X}=\tilde{\pi}(\tilde{\sigma})$, $\tilde{\sigma}\in \tilde{G}$,  by using the identification $P\L(\tilde{\Gamma})P=P\cR(\tilde{\Gamma})P$ and the same construction as in the Appendix 1,
we obtain that
$\Psi_{\tilde{X}}$ is unitarily equivalent to 
$$
\tilde{\Psi}_{[\tilde{\Gamma}\tilde{\sigma}\tilde{\Gamma}]}(x)=
E_{P\L(\tilde{\Gamma})P}^{P\L(\tilde{G})P}
\Big(P\tilde{T}^{[\tilde{\Gamma\sigma\Gamma}]}P  x  P\tilde{T}^{[\tilde{\Gamma\sigma\Gamma}]}P \Big),
\quad x\in P\L(\tilde{\Gamma})P,
$$
where $[\tilde{\Gamma}\tilde{\sigma}\tilde{\Gamma}]=[\tilde{\Gamma\sigma\Gamma}]$
is a double coset.

Choosing a system of representatives for the elements $\tilde{\gamma}P$, $\tilde{\gamma} \in \tilde{\Gamma}$
amounts  to give a cocycle $\varepsilon$, and working with $\L(\Gamma,\varepsilon)$ instead
of $\A_P$, and hence also the operators $\tilde{T}^{[\tilde{\Gamma}\tilde{\sigma}\tilde{\Gamma}]}$
are unitarily equivalent to the classical Hecke operators where
$\tilde{G}=\PGL_2(\Z[\frac1p])$,  $\tilde{\Gamma}=\PSL_2(\Z)$.
}
\end{rem}


We now return to the context of Theorem \ref{heckerep}. The previous remark shows that we may always switch from the skewed algebra with cocycle to a reduced algebra of the cover group $\tilde{G}$. 
In chapter 5, we proved that the Hecke algebra $\H_0$ of double cosets $[\Gamma \sigma\Gamma]$ of $\Gamma$ in $G$
 admits a $\ast$-representation into $\L(G,\varepsilon)$, by mapping a double coset $[\Gamma \sigma\Gamma]$ into 
 $$t^{[\Gamma \sigma\Gamma]}=\sum_{\theta \in \Gamma \sigma\Gamma}\overline{ \langle \pi_{13}(\theta)e,e}\rangle_{13} \theta \in l^2(\Gamma \sigma\Gamma)\cap \L(G),$$
 where $e$ is the identity element of $G$.
 
 (For simplicity, from now on we do the notational substitution consisting in changing the coefficients $t(\theta)$ with $\overline{t(\theta)}$, $\theta \in G$). Here $e$ is a vector in the Hilbert space $H_{13}$ that is  a cyclic trace vector for the von Neumann algebra generated by $\pi_{13}(\Gamma)$, which is isomorphic to $\L(G, \varepsilon).$

We will apply this theorem to the representation of the completely positive maps
$
\Psi_\sigma (x)=[\Gamma:\Gamma_\sigma]E_{\L(\Gamma,\varepsilon)}^{\L(G, \varepsilon)}
\big(t^{\Gamma\sigma\Gamma} x t^{\Gamma\sigma\Gamma}\big)$, $x\in \L(\Gamma, \varepsilon).
$
We  analyze the spectrum of the maps $\Psi_\sigma$ modulo the compact operators.
We require then that  the convolutors $t^{\Gamma\sigma\Gamma}$
are in the reduced $C^*$-algebra $C^*_{\text{red}}(G, \varepsilon)$. To obtain this requirement, we  prove that 
there exists a choice for the cyclic trace vector $\xi$ in $H_{13}$ such that $t^{\Gamma\sigma\Gamma}$
belong to the reduced $C^*$-algebra $C^*_{\text{red}}(G, \varepsilon)$.

Note that changing $\xi$ into $u\xi$, where $u$ is a unitary in $\L(\Gamma, \varepsilon)$,  changes $t^{\Gamma\sigma\Gamma}$ into 
$u^*t^{\Gamma\sigma\Gamma}u$. We are proving that the orbit
$\big\{u^*t^{\Gamma\sigma\Gamma}u\mid u\in\U(\L(\Gamma, \varepsilon))\big\}$ intersects the reduced $C^*$-algebra.

\begin{lemma}
With the notations from Proposition 4, there exists a choice of the cyclic trace vector $\xi$ in $H_{13}$ used in the construction of the elements  $t^{\Gamma\sigma\Gamma}$, such that for all
double cosets $[\Gamma\sigma\Gamma]$, the elements $t^{\Gamma\sigma\Gamma}$ belong to the reduced 
$C^*$-algebra $C_{\rm red}^*(G,\varepsilon)$. 
\end{lemma}

{\it Proof.} Consider the space $\H_{13}$ of positive functions on $\PSL_2(\R)$ that are obtained
as matrix coefficients from elements $\eta$ in $H_{13}$ (that is $\varphi : G \to \C$ belongs
to $\H_{13}$ if there exists $\eta$ in $H_{13}$ such that $\varphi(g) = \langle \pi_{13}(g) \eta,\eta \rangle$,
$g$ in $\PSL_2(\R)$.

Obviously, $\H_{13}$ is a cone closed to infinite convergent sums. Indeed if $(\eta_i)$
is a family of vectors in $H_{13}$, $\sum\|\eta_i\|^2<\infty$, each determining the positive
functional $\varphi_i$. Consider the Hilbert subspace $L$ of $H_{13}\otimes \ell^2(I)$
generated by $\bigoplus\limits_{i \in I} \pi(g)\eta_i$. This space is
obviously invariant to the action of $G$. Since $\pi_{13}$ is irreducible $\pi(g)|_L$ is a multiple of 
the representation $\pi_{13}$ and because we have a cyclic vector, it is unitary equivalent to 
$\pi_{13}$. The vector $\eta = \bigoplus\limits_{i \in I} \eta_i$ will then determine the positive
definite function on $G$ defined by the formula
 $$\sum\limits_i \varphi_i(g) =\sum\limits_i \langle \pi_{13}(g) \eta_i,\eta_i \rangle, g \in G.$$

In the sequel we denote $\pi_{13}$ simply by $\pi$.

As it was noted in the list of properties of $t^{\Gamma\sigma\Gamma}$, this is equal to 
$$\sum\limits_{g\in [\Gamma\sigma\Gamma]}\overline{ \langle \pi(g) \xi,\xi \rangle }g.$$
 
If $\varphi_\eta(g) = \overline{\langle \pi(g)\eta,\eta \rangle}$, $g\in\PSL_2(\R)$ is determined by the vector 
$\eta$, then for $a$ in $L^1(\L(\Gamma),\tau)$ the vector $\pi(a)\eta$ (note that 
$\pi|_\Gamma$ extends from $G$ to a representation of $\Gamma$ on $H_{13}$ to a representation
of $\L(\Gamma, \varepsilon)$) will determine a functional $\varphi_a$, that has the property that
$$
\varphi_a | _{\PGL_2(\Q)} = a^* \varphi a.
$$

We are looking to find a positive functional in $\H_{13}$ that has the property that 
$\varphi | _{\PGL_2(\Q)}$ belongs to the reduced $C^*$-algebra, and such that moreover
$\varphi$ is implemented by a trace vector (as we have seen in Chapter 3, this is equivalent
to the pseudo-multiplicative property
$$
\varphi(g_1g_2) = \sum_{\gamma \in \Gamma} \varphi(g_1\gamma)\varphi(\gamma^{-1}g_2), \quad
g_1,g_2\in \PSL_2(\R).
$$

To find such a $\varphi=\varphi_\xi$ is therefore sufficient to find a vector $\xi$ such that the corresponding
positive functional has the following properties:

1) the restriction of $\varphi_\xi$ to $\Gamma\sigma\Gamma$ determines an element in
$C_{\rm red}^*(\PGL_2(\Q)_+, \varepsilon)$;

2) $\varphi_\xi|_\Gamma$ is invertible in $C_{\rm red}^*(\Gamma)$.

Indeed if we found such a vector $\xi$ then we are done because the vector $\xi_0=
\pi((\varphi_\xi|_\Gamma)^{-1})\xi$ is a trace vector.

Moreover, let $$\Psi(g)= \langle \pi(g) \xi_0,\xi_0 \rangle$$
and let $t_0^{\Gamma\sigma\Gamma}= \sum\limits_{g\in \Gamma\sigma\Gamma}\Psi(g)g$.
Then $t_0^{\Gamma\sigma\Gamma}= (\varphi_\xi |_\Gamma)^{-1/2} t^{\Gamma\sigma\Gamma}
(\varphi_\xi |_\Gamma)^{-1/2}$, where $t^{\Gamma\sigma\Gamma}$ correspond to $\varphi_\xi |_{\Gamma\sigma\Gamma}$
and hence are in $C_{\rm red}^*(G, \varepsilon)$ and thus belongs to  $C_{\rm red}^*(G, \varepsilon)$.

Hence to conclude the proof it is sufficient to construct a vector with properties 1), 2). By
Jolissaint estimates, it is sufficient to take a fast decreasing vector for the group $G$, such that
$\varphi_\xi |_\Gamma$ is invertible.

We now use a result by  in  [BH] (proof of Theorem A1) which says that given $x\geq 0$, $x\neq 0$ in
$C_{\rm red}^*(\Gamma, \varepsilon)$ there exists unitaries $\gamma_1,\ldots,\gamma_n$ in $\Gamma$ such
that $\sum \gamma_i x \gamma_i^{-1}$ is invertible.

Let $\xi$ be a vector in $H_{13}$, generating a positive definite function on G,  which  has rapidly decreasing coefficients ([Jo]). For example we may take the vector of evaluation
at $0$ in the model of the unit disk).

Then we construct the functional $\varphi_\xi$ and use the above mentioned result in [BH], to replace 
$\varphi_\xi$ by $\sum \gamma_i^{-1} \varphi_\xi \gamma_i = \Psi _0$.

Then $\Psi_0$ corresponds to the vector $\frac1{\sqrt n} ( \oplus\, \pi(\gamma_i)\xi)$
which is a vector generating a positive definite function with  rapidly decreasing coefficients. Consequently, by construction,  $\Psi_0$ is invertible and the inverse belongs to the $C^\ast$-algebra.

(See also [Ra5], where it is proved that the elements $t^{[\Gamma\sigma]}$, $\sigma \in G$ are a Pimsner Popa basis, and thus bounded).
\qed
\

\

The algebraic mechanism, that is implicit  in  the fact that the linear application, mapping a double coset $[\Gamma\sigma\Gamma]$ into the
completely positive map $\Psi_\sigma$ on $\L(\Gamma)$,  is a   *-algebra morphism (constructed in Chapter 5),  is summarized as follows:
 \vskip6pt

\begin {lemma}\label{algebra} \it Let $\tilde \A$ be the  $\ast-\Bbb C$-algebra generated by $\L(G)\otimes \L(G)^{\rm op}$  and the characteristic functions $\chi_C$, where  $C$ runs through the cosets, in $G$, of the subgroups $\Gamma_{\sigma}, \sigma \in G$, subject to the relation

 $$(g_1\otimes g_2)\chi_C=\chi_{g_1C(g_2)^-1}(g_1\otimes g_2),$$
  for  $g_1\otimes g_2$ in $G\otimes G^{\rm op}$. Let $\tilde \A_0$ be the the subalgebra $\chi_{\Gamma} \tilde \A \chi_{\Gamma}$, with unit  $\chi_{\Gamma}$. Note that $\tilde \A_0$ is a weakly dense sub algebra in the reduced, von Neumann  crossed product algebra of the measure  preserving, grupoid action of $G\otimes G^{\rm op}$ on the space $K$ (as is the $C^{\ast} $ algebra $\A$ from the statement of Theorem \ref {AOlocal}).  

  Then  map $\Phi$ from the Hecke algebra
 $\H_0=\H_0(G,\Gamma)$ into $\tilde \A_0$, defined by
$$\Phi([\Gamma\sigma\Gamma])= \chi_{\Gamma}( t^{[\Gamma\sigma\Gamma]}\otimes t^{[\Gamma\sigma\Gamma]})\chi_{\Gamma},$$
is a $*$-algebra morphism (see also [Ra5]).

\end{lemma}
 
\vskip6pt

\begin{proof} This obtained by passing to the quotient, modulo the compacts the fact that the map
taking the double coset $[\Gamma\sigma\Gamma]$ into the completely positive application $\Psi_\sigma$, defined by
$$\Psi_\sigma (x)=[\Gamma:\Gamma_\sigma]E_{\L(\Gamma)}^{\L(G)}
\big(t^{\Gamma\sigma\Gamma} x t^{\Gamma\sigma\Gamma}\big), \quad x\in \L(G, \varepsilon),$$
is a $*$-algebra morphism.

\end {proof}

In the next theorem, by using the identification proved  in Theorem \ref {AOlocal}, 
of the algebra $\A_0$ with  the reduced   $C^*$-algebra groupoid crossed product of $G \times G^{\rm op}$ acting on $K$, we prove that the linear  application  mapping a coset $\Gamma\sigma\Gamma$, $\sigma \in G$ into the class,  in the Calkin algebra, of the completely positive map $\Psi_\sigma$, extends to an isometric embedding from the reduced $C^\ast$ algebra $\H_{\rm red}$ into the reduced $C^\ast$-algebra $C^\ast_{\rm red}((G\times G^{\rm op})\ltimes C(K))$. We are now proving that  this latest map, is in fact the linear application    $\Phi$ constructed in Lemma \ref{algebra}.

\begin{thm}\label{main1}  Let
 $G =  \PGL_2(\Z[\frac 1p])$, $\Gamma = \PSL_2(\Z)$. Let   $\Pi_Q$ be  the canonical projection  from $\B(\ell^2(\Gamma))$ onto the Calkin algebra $Q(\ell^2(\Gamma))$.
 
 Then, the  $*$-algebra morphism, constructed in Chapter 5, which maps a double coset $[\Gamma\sigma\Gamma]$ into the
completely positive map $\Psi_\sigma$ on $\L(\Gamma,\varepsilon)$ given by  the Stinespring
 dilation formula 
$$
\Psi_\sigma(x) =  [\Gamma:\Gamma_{\sigma}] E_{\L(\Gamma,\varepsilon)}^{\L(G, \varepsilon)}
\big(t^{\Gamma\sigma\Gamma}x t^{\Gamma\sigma\Gamma}\big), \quad x\in \L(\Gamma,\varepsilon)
$$
extends, when composing with the canonical projection $\Pi_Q$  to an isomorphism  
  from the reduced $C^*$-Hecke algebra $\H-{\rm red}$ into the the Calkin algebra  $Q(\ell^2(\Gamma))$, mapping
  the double coset
  $\Gamma\sigma\Gamma$ into $\Pi_Q (\Psi_\sigma)$ for $\sigma $ in $G$.
  
  Here we use implicitly the fact that for all $\sigma $ in $G$, the linear continuous map
  $\Psi_\sigma$, which is defined a priori on $\L(\Gamma, \varepsilon)$, extends to a bounded operator on $\ell^2(\Gamma)$.
  
  \end{thm}

{\it Note.} We proved in Chapter 4 (Lemma \ref{continuity}, see also [Ra3]) that the validity of the estimates of the Ramanujan-Petersson conjecture
is equivalent to the continuity, with respect to the reduced $C^*$-algebra associated to the Hecke
algebra, of the map taking $\Gamma\sigma\Gamma$ into $\Psi_{\sigma}$. Hence we prove that the Ramanujan-Petersson conjecture holds true, modulo the compact operators (that is for the essential spectrum of the Hecke operators) in the
case $G = \PGL_2(\Z\big[\frac{1}{p}\big])$, $\Gamma = \PSL_2(\Z)$, for every prime number $p$.

\vskip6pt

In fact, a trivial application of classical Fredholm theory gives the following corollary.

\begin{cor}\label{main2}
For every prime number $p$, the essential spectrum of the classical Hecke operator  $T_p$, acting on Maass forms, is contained in the interval $[-2\sqrt p, 2\sqrt p]$, predicted by the Ramanujan Petersson conjectures. In particular, given an open interval containing $[-2\sqrt p, 2\sqrt p]$,  there are at most a finite number of possible exceptional eigenvalues lying outside this interval.

Note that as a corollary of the proof we reprove that the continuous part of the spectrum (corresponding
to Eisenstein series) also verifies the  Ramanujan-Petersson estimates. 
(See also the paper of P. Sarnak ([Sa])
where a  distribution formula for  the exceptional values is computed).  

Let $\Gamma_0=\Gamma_0(p^n)$, $n\geq 1$ be a modular subgroup  of $\PSL(2, \Bbb Z)$. We replace  in the above computations, the projective unitary representation $\pi_{13}$ by the projective unitary representation $\pi_t$, where $t$ is determined by the condition $[(t-1)/12]= \frac {1}{[\Gamma:\Gamma_0]}$. 

Then using the matrix coefficients of the representation $\pi_t$ restricted to $\PGL(2,\Bbb Z[1/p])$, the above methods prove that the essential norm for Hecke operator on the  $\Gamma_0$-invariant Maass forms,  associated to the double coset $\Gamma_0 \sigma\Gamma_0$, $\sigma \in \PGL(2,\Bbb Z[1/p])$ is equal to the norm of the convolution operator by the coset $\Gamma_0 \sigma\Gamma_0$ on the Hilbert space 
 $\ell^2(\Gamma_0\backslash \PGL(2,\Bbb Z[1/p]))$ (this norm is  by definition, is the norm of the double coset $\Gamma_0 \sigma\Gamma_0$, viewed   as en element of the reduced $C^\ast$-Hecke algebra
 $\H_{\rm red}(\Gamma_0\backslash \PSL(2,\Bbb Z[\frac{1}{p}])/ \Gamma_0))$.

\end{cor}

We  note that  the existence of a spectral gap bellow the eigenvalue $[\Gamma:\Gamma_{\sigma}]$  (corresponding to the eigenvector 1) of $\Psi(\sigma)$ is equivalent to the existence to a spectral gap in the sense of [Po2] (that is to the fact that a sequence in $\L(\Gamma, \varepsilon)$ that asymptotically commutes  
with $ t^{\Gamma\sigma_p\Gamma}$, $p$ a prime number,  should be an asymptotically scalar sequence).

 \vskip6pt

\begin{proof}(of Corollary \ref{main2})  The corollary follows from the Theorem \ref{main1}. Indeed for every prime $p\geq 2$ let  
$\sigma_{p}  = \begin{pmatrix} 1 & 0 \\ 0 & p
\end{pmatrix}$ and let $\alpha_p=[\Gamma\sigma_p\Gamma]$ be the corresponding double coset. By
Theorem 33 it follows that the essential spectrum of $\Psi_{\sigma_p}$ is equal to the spectrum of the
double coset $\alpha_p$ as a selfadjoint convolutor in the reduced $C^*$-Hecke algebra. By the Lemma \ref{continuity}, the spectrum of $\alpha_p$ is the interval $[-2\sqrt p, 2\sqrt p]$. Consequently the essential spectrum of $\Psi_{\sigma_p}$ is the interval $[-2\sqrt p, 2\sqrt p]$. 
By the  Proposition \ref{identification}, the classical Hecke operators $T_p$
are unitarily equivalent  (modulo a rescaling of the Hilbert space) to the completely positive map
$\Psi_{\sigma_p}$ acting on $\ell^2(\Gamma)$. Hence the essential spectrum of $T_p$ is
 $[-2\sqrt p, 2\sqrt p]$ and hence by Fredholm theory the discrete spectrum can only accumulate at the endpoints of the interval. 
 
 The last  part of the statement follows  from the fact that, by the dimension formula in ([GHJ]), we have that the Murray von Neumann dimension dim$_{\{\pi(\Gamma_0)\}'''}H_t=1$. Hence the construction (in Chapter 4) of the representation of  the reduced $C^\ast$-Hecke algebra $\H_{\rm red}(\Gamma_0\backslash \PSL(2,\Bbb Z[\frac{1}{p}])/ \Gamma_0)$ could be done also in this case. The local Akemann Ostrand property remains valid for finite index subgroups of $\PSL(2,\Bbb Z)\subseteq \PSL(2,\Bbb Z[\frac{1}{p}])$ , and this proves the last part of the statement.

 \end{proof}

\vskip6pt

\begin{proof}(of Theorem \ref{main1}).  By definition, for $\sigma \in G$, the operator $\Psi_\sigma$ belongs to 
the algebra $\B$,  which we recall that it is  the  $C^*$-algebra generated by $\chi_\Gamma L_{g_1} R_{g_2}\chi_\Gamma
\in \B(\ell^2(\Gamma))$, $g_1,g_2 \in G$, and $C(K) \subseteq \B(\ell^2(\Gamma))$
(by $\chi_{\Gamma}$ we denote the characteristic function of $\Gamma$ viewed as a multiplication operator on $\ell^2(G)$). 

 Taking its image into the Calkin algebra, and using the identification of the quotient algebra from Theorem \ref{AOlocal}, the only thing that remains to  be proved is that   the map 
$$
[\Gamma\sigma\Gamma]\to \chi_\Gamma \big(t^{\Gamma\sigma\Gamma}\otimes t^{\Gamma\sigma\Gamma}\big)\chi_\Gamma
\in \A= C_{\text{red}}^*((G\times G^{\rm op}) \rtimes C(K)).
$$
extend to a continuous isomorphism from the reduced, $C^\ast$-Hecke algebra $\H_{\rm red}$ into $\A$. But this is a trace preserving map when endowing $\A_0$ with the crossed product trace coming from the Haar measure on $K$. Hence the above map preserves moments and thus is an isomorphism. 
\end{proof}

\section*{Appendix 1\\
A construction of abstract Hecke operators on II$_1$ factors}

In this appendix we start with a pair of isomorphic subfactors of a given type II$_1$ factor.
We define the analogue of the first step of the Jones's basic construction for such a data, 
which is a correspondence between spaces of intertwiners and von Neumann bimodules
over the initial II$_1$ factor (see also [FV] for a related approach).

We then analyze the Connes' fusion for  these bimodules and prove a multiplicativity property 
for the associated completely positive maps, which generalizes the construction in Chapters 2,3.

\begin{defn}
 Let $M$ be a type II$_1$ factor and let $N_0, N_1\subseteq M$ be two
subfactors of equal index and $\theta:N_0\rightarrow N_1$ an isomorphism.
We denote by $I_\sigma\subseteq B(L^2(M,\tau))$ the linear space of all 
$X:L^2(M, \tau)\rightarrow L^2(M, \tau)$ such that
\begin{equation*}
Xa=\theta(a)X\hbox{ for all } a\in N_0.
\end{equation*}

Note that if $M=\L(\Gamma)$, $N_0=\L(\Gamma_{\sigma^{-1}})$ and
$N_1=\L(\Gamma_\sigma)$ then $\theta$ is $\sigma\cdot\sigma^{-1}$, viewed
as a map on $\Gamma_{\sigma^{-1}}$ into $\Gamma_\sigma$, and extended to the group algebra.

Also if $X$ belongs $I_{\sigma}$ then obviously $Y^*$ belongs to $I_{\sigma^{-1}}$.
$I_\sigma$ plays the role of the commutant algebra of a subfactor, in the case $N_0=N_1$ and if $\theta$ is the identity.

The following construction is a measure for the obstruction for $\sigma$ being implemented by an 
automorphism of $M$.
\end{defn}

\begin{defn}
 Let $X,Y$ in $I_\sigma$. Then $X\cdot Y^*$ maps $N'_0$ into $N'_1$
(e.g.  $X a Y^*$ belongs to $N'_1$ for all $a$ in $N'_0$),
and hence we have the following diagram
$$
\begin{array}{ccc}
N'_0 & \mathop{\longrightarrow}\limits_{X\cdot Y^*} &  N'_1\\[3pt]
\phantom{aaaa}_{\rm inc }\nwarrow  & & \hskip-4pt\swarrow   {}_{E_{M'}^{N'_1}}\\[3pt]
 & M'  &
\end{array}
$$
where $E_{M'}^{N'_1}$ is the canonical expectation from $N'_1$ onto $M'$
(both $N'_1$, $M'$ are II$_1$ factors, and the commutants are computed in the algebra 
$\B(L^2(M,\tau))$.
Denote $\Psi_{X,Y^*}$ the composition, which is thus a linear map from $M'$ into $M'$.

Thus the formula for $\Psi_{X,Y^*}$ is
$$
\Psi_{X,Y^{*}}(m')=E_{M'}^{N'_1}(Xm'Y^*),\quad m'\in M'.
$$

Note that if $X=Y$, then $\Psi_{X,X^\ast}$ is a completely positive map.
As explained in Chapter 2, $\Psi_{X,Y^\ast}$ is a generalization of the Hecke operators.

The analysis of  the maps $\Psi_{X,Y^\ast}$ is  a method to measure how far is $\theta$ 
from being implemented by an internal automorphism.
Indeed if $\sigma$ was the restriction of  an automorphism of $M$, then $\sigma$ would
implement an unitary $U$ on $L^2(M,\tau)$ which in turn
would have the property that $UM'U^*=M'$ and hence the completely positive map $\Psi_{U,U^*}$  would be simply an automorphism of~$M'$.
\end{defn}

We develop the analogy with the Jones' basic construction. In Jones's basic construction, for an inclusion of algebras with traces, $B\subseteq A$,
the first algebra in the basic construction is $Ae_BA$ (see [GHJ])
(as an $A\otimes A^{\rm op}$ bimodule) and it isomorphic to the algebra $B'\subseteq \B(L^2(A,\tau))$.

In our situation, we want to get an abstract definition of the $\Gamma\times\Gamma^{\rm op}$
bimodule $\ell^2(\Gamma\sigma\Gamma)$, starting from $\theta_\sigma : \Gamma_{\sigma^{-1}} \to
\Gamma_\sigma$
defined by $\theta_\sigma(\gamma)=\sigma\gamma\sigma^{-1}$. 
\begin{defn}\label{cocycle}
If a cocycle $\varepsilon$ is present on $G$ (coming from a projective unitary representation of $G$) then $\theta$ is replaced with the automorphism $\tilde {\theta}$, constructed at the 
end of Chapter 3. In this case the bimodule $\ell^2(\Gamma\sigma\Gamma)$ is identified with a subspace of
$\L(G,\varepsilon)$. In particular if $u_g, g\in G$ is the canonical basis of $\L(G,\varepsilon)$, the the bimodule
structure of $\ell^2(\Gamma\sigma\Gamma)$ over $\Gamma\times\Gamma^{\rm op}$  is so~that
$$\gamma_1u_\sigma\gamma_2= \chi_{\varepsilon}(\gamma_1, \sigma,\gamma_2)u_{\gamma_1\sigma\gamma_2}, 
\quad \gamma_1, \gamma_2 \in \Gamma,\ g \in G.$$
The coefficient $\chi_{\varepsilon}(\gamma_1, \sigma,\gamma_2)$ is determined by the cocycle $\varepsilon$.

\end{defn}

\begin{defn}
Let $N_0,N_1\subseteq M$ and let $\theta:N_0\rightarrow N_1$  an isomorphism.
(which should correspond respectively to $\Gamma_{\sigma^{-1}}, \Gamma_\sigma\subseteq\Gamma$
and $\theta_\sigma{(\gamma)}= \sigma\gamma\sigma^{-1}$ in the group case,  with the above amendment if a cocycle $\varepsilon$ is present). The bimodule generalizing for the pair of isomorphic subfactors, the commutant in the Jones's construction, is 
 the Hilbert space closure of $M\sigma M=M\sigma M^{\rm op}$
where $\sigma$ is a virtual element with the property $\sigma n_0\sigma^{-1}=\theta(n_0)$
or $\sigma^{-1} n_1 \sigma=\theta^{-1}(n_1)$ for $n_0$ in $N_0$, $n_1$ in $N_1$.

Here the element $m \sigma m'$ is the tensor product $m \otimes m'$, where
$m \otimes m'$ belongs to $M\otimes M^{\rm op}$, and the scalar product is
$$
\langle m \otimes m',a \otimes a'\rangle =\tau(a^*m\theta (E_{N_0}((a')^* \odot m'))
$$
for all $m,a$ in $M$, $m',a'$ in $M^{\rm op}$. Here $\odot$ stands for the product in $M^{\rm op}$
that is $x\odot y=yx$.

Thus the formula for the scalar product is
$$
\langle m \otimes m',a \otimes a'\rangle =\tau(a^*m\theta (E_{N_0}(m'(a')^*).
$$

Note that the above formula could also  be used to  define an $M$-left Hilbert module structure on $M\sigma M^{\rm op}$. 
\end{defn}

{\it Proof of the consistency of the definition}.
We have to prove that the definition is consistent with the formal  definition of $M\sigma M$,
which is equal as a vector space to $M\otimes M^{\rm op}$.

Thus we have to verify that 
$$mn_1\sigma m'=m\sigma(\sigma^{-1}n_1\sigma)m'=m\sigma\theta^{-1}(n_1)m'=m\sigma(m'\odot\theta^{-1}(n_1))$$
for all $m,m'$ in $M$, $n_1$ in $N_1$.

Thus we have to verify that
$
\langle m n_1 \otimes m'-m \otimes\theta^{-1}(n_1) m', a \otimes a'\rangle 
$
is zero for all $m,m',a,a'$ in $M$, $n_1$ in $N_1$.
But
\begin{gather*}
\langle m n_1 \otimes m', a\otimes a'\rangle =
\tau(a^*mn_1\theta(E_{N_0}(m'(a')^*)=\\
=\tau(a^*m\theta(\theta^{-1}(n_1))\theta(E_{N_0}(m'(a')^*)))=\\
=\tau(a^*m\theta(\theta^{-1}(n_1))(E_{N_0}(m'(a')^*)))=
\tau(a^*m\theta(E_{N_0}(\theta^{-1}(n_1)m'(a')^*))).
\end{gather*}

Here we use the fact that $E_{N_0}$ is a conditional expectation
and that $\theta^{-1}(n_1)$ belongs to~$N_0$.

Note that the scalar product corresponds exactly to the Stinespring dilation of the
completely positive map $m\rightarrow \theta(E_{N_0}(m))$ viewed
as a map from $M$ with values into $N_1\subseteq M^{\rm op}$.

\begin{rem}
 Without going into the complication of using  the definition of   $M^{\rm op}$, which is only
needed to have positivity of the scalar product, we could simply say that $M\sigma M$ 
is the Hilbert space completion, of the bimodule defined by the relation
$$
mn_1\sigma m'=m\sigma\theta^{-1}(n_1)m'
$$
for all $m$, $m'$ in $M$, $n_1$ in $N$ and $\theta$ is implemented formally by $\sigma$.

Then the scalar product
$\langle m\sigma m', a\sigma a'\rangle$ is formally trace of
$(a')^*\sigma^{-1}a^*m\sigma m'$
which, by the trace property, is equal to the trace of 
$a^*m\sigma(m'(a')^*)\sigma^{-1}$ and is formally equal to 
$\tau(a^*m\theta(E_{N_0}(m'(a')^*))).$

We define an anti-linear isomorphism between the intertwiner space and
the bimodule as follows.
\end{rem}

\begin{defn} For $X$ in $I_{\sigma}$ (that is $Xn_0=\theta(n_0)X$ for all $n_0 \in N_0$) we
associate to $X$ a canonical element in $M\sigma M$, where as above, the element $\sigma$
virtually implements $\theta$ (that is $mn_0\sigma m'=m\sigma\theta(n_0)m'$,
for all $m$, $m'$ in $M$, $n_0$ in $N$).

Then the anti-linear map $X\rightarrow\theta(X)\in L^2(M\sigma M)$ 
is defined by the relation
$$
\langle m\sigma m',\theta(X)\rangle=\tau(X(m')m)
$$
for all $m$, $m'$ in $M$.
\end{defn}

{\it Proof} (of the consistency of the definition).
We have to check that with this definition $X(n_0m)=\theta(n_0)X(m)$
or by taking a trace against on element $m'$ that 
$$
\tau(X(n_0m)m')=\tau(X(m)m'\theta(n_0)).
$$
By using the above definition of $\theta(X)$ this comes to
$$
\langle m'\sigma n_0m,\theta(X)\rangle=
\langle m'\theta (n_0)\sigma m,\theta(X)\rangle
$$
which is obviously true from the definition of the bimodule property of $M\sigma M$.

\begin{cor} With the notations introduced above,
assume that $s_i$ a left Pimsner Popa orthonormal basis for $N_0$ in $M$.
Consequently,  $M$ is as left $N_0$ bimodule the ($N_0$-orthogonal) sum of $N_0 s_i$.

Then $X(n_0s_i)=\theta(n_0)X(s_i)$, for all $n_0 \in N_0$.
Denote by $t_i=X(s_i)$. Then the $t_i$ are a $N_1$ Pimsner-Popa orthonormal basis
for $N_1$ in $M$.

Moreover, the formula for $\theta(X)$ is in this case
$$
\theta(X)=\sum_i t_i^*\sigma s_i.
$$
\end{cor}

{\it Proof.} Note that the decomposition
$M\sigma M^{\text{op}}=\bigcup[M\sigma s_i]$ is orthogonal.
Hence we may assume that assume $\theta(X)=\sum_i x_i\sigma s_i$.

The relation between $\theta(X)$ and $X$ is
$$
\langle m_0\sigma m_1,\theta(X)\rangle= \tau(X(m_1),m_0)
$$
and hence
$$
\langle X(m_1), m_0\rangle=
\langle m_0^*\sigma m_1,\theta(X)\rangle.
$$
Hence taking $m_1=s_i$ we obtain
$$
\langle t_i, m_0\rangle=
\langle X(s_i),m_0\rangle=\langle m_0^*\sigma s_i,\theta(X)\rangle=
\langle m_0^*\sigma s_i,x_i\sigma s_i\rangle.
$$

Hence we get that for all
$m_0$ in $M=\L(\Gamma)$ we have that
$$
\langle t_i, m_0\rangle=
\langle m_0^*,x_i\rangle
$$
or that $\tau(t_i m_0^*)=\tau(x_1^*m_0^*)$
and hence that $t_i=x_i^*$.

Hence
$$
\theta(X)=\sum_i(X(s_i))^*\sigma s_i.
$$
Another corollary is the explicit formula for $\theta(X)$
in the case we have that $\gamma_1X\gamma_2=X(\gamma_1\sigma\gamma_2)$,
$\gamma_1,\gamma_2 \in \Gamma$.

\begin{cor} We assume that we are in the case of a group $G$ with two-cocycle $\varepsilon$, as described in Definition 37. Let $\sigma$ be an element in $G$. Let $\tilde{\theta}$ be the corresponding isomorphism from
$N_0=\L(\Gamma_{\sigma^{-1}}, \varepsilon)$ onto
$N_1=\L(\Gamma_{\sigma}, \varepsilon)$.
Let $X$ be in $I_{\tilde{\theta}}$. Denote $\gamma_1X\gamma_2$ by $X(\gamma_1\sigma\gamma_2)$.
Then
\begin{gather*}
\theta(X)=\sum_{\alpha\in\Gamma\sigma\Gamma}\overline{\chi_{\varepsilon}(\gamma_1,\sigma,\gamma_2)(\langle X(\gamma_1\sigma\gamma_2)I,I\rangle)}, 
\end{gather*}
where $I$ is the unit element (or more generally a trace vector) in $\L(\Gamma, \varepsilon)$.

In particular, if $\pi$ is a projective unitary  representation
of $G$, with 2-cocyle $\varepsilon$, extending the left regular representation, and $\sigma \in G$, $X=\pi(\sigma)$, then
\begin{gather*}
\theta(\pi(\sigma))=\sum_{\alpha\in\Gamma\sigma\Gamma}\overline{(\langle X(\alpha)I,I\rangle)}\alpha, 
\end{gather*}

\end{cor}

\begin{proof} Again we use the formula $\langle m_0\sigma m_1,\theta(X)\rangle=
\tau(X(m_1)m_0)$ and hence
$$
\langle X(a),b\rangle=\langle b^*\sigma a,\theta(X)\rangle
$$
or
$$
\langle \theta(X),b^*\sigma a \rangle=\langle b,X(a)\rangle.
$$
Thus, using the notations from Definition \ref{cocycle}, we obtain:

$$
\langle 
\theta(X),u_{\gamma_0\sigma\gamma_1}\rangle=\chi_{\varepsilon}(\gamma_0,\sigma,\gamma_1)\langle \theta(X),\gamma_0 u_\sigma \gamma_1\rangle=$$
$$\chi_{\varepsilon}(\gamma_0,\sigma,\gamma_1)\langle\gamma_0^{-1},X(\gamma_1)\rangle=
\chi_{\varepsilon}(\gamma_0,\sigma,\gamma_1)\langle I,\gamma_0 X\gamma_1 I\rangle=$$
$$=\chi_{\varepsilon}(\gamma_0,\sigma,\gamma_1)\overline{\langle (\gamma_0X\gamma_1 I,I\rangle}.
$$
The second part of the statement follows from the fact that since $\pi$ is a projective unitary representation extending the left regular representation of $\Gamma$, we have
$$\gamma_1\pi(\sigma)\gamma_2=\chi_{\varepsilon}(\gamma_0,\sigma,\gamma_1)\pi(\gamma_1\sigma\gamma_2), 
\quad\gamma_1, \gamma_2 \in \Gamma.$$
\end{proof}

The isometrical property of the map $\theta$ from intertwiners into bimodules is described
in the next proposition.

First, we define an $M$-valued pairing $\mathcal{P}$
from $M\sigma M \times M\sigma^{-1} M$ into $M$ as follows

\begin{defn}
There is a well defined projection $\mathcal{P}:M\sigma M\times M\sigma^{-1} M$
into $M$, defined by the formula
$$\mathcal{P}((m_0\sigma m_1)(m_2\sigma^{-1}m_3))=$$
$$=\mathcal{P}((m_0\sigma m_1),(m_2\sigma^{-1}m_3))= m_0\theta(E_{N_0}
(m_1m_2))m_3.
$$
Indeed, this is the $M$-component of the Connes' fusion product
$$(M\sigma M)\mathop{\otimes}\limits_M (M\sigma^{-1} M).$$
\end{defn}

\begin{prop}
Let $N_0, N_1\subseteq M$ and $\theta$
an isomorphism from $N_0$ into $N_1$, virtually  implemented
by $\sigma$.
Fix $m'=R_m\in M'$, for $m$ in $M$ be the right convolutor by $m$.

Then for all $X,Y$ in $I_{\sigma}$ we have that
$$
E_{M'}^{N'}(Xm'Y^*)=\mathcal{P}(\theta(Y)m\theta(X^*))
$$
for all $m'=R_m$ in $M'$.
\end{prop}

\begin{proof}
The proof is essentially that fromTheorem 22.
and won't be repeated here. Note that by linearity here we can assume
 simply that
$\theta(X)=s\sigma r$, $\theta(Y)=s_1\sigma_1 r_1$, for unitaries
$s$, $s_1$, $r$, $r_1$ in $M$.
\end{proof}

From here on we work (for simplicity) only in the case
$G$, $\Gamma$ and $\sigma_1,\sigma_2$ partial automorphism
of $\Gamma$ reduced by elements in $G$, but we maintain the generality
of the choice $X$, $Y$.
We assume that we are given  $\varepsilon$ a 2-cocycle on $G$, preserved by all $\sigma$'s and all algebras are group
algebras with cocycle.

\begin{defn}
Fix an element $\sigma'$ in $[\Gamma\sigma\Gamma]$.
Let two orthogonal projections $\mathcal{P}_{\sigma'\Gamma}$ 
and $\mathcal{P}_{\Gamma{\sigma'}}$  be the projections on $\ell^2[\sigma'\Gamma]$
and $\ell^2[\Gamma{\sigma'}]$ respectively.
For $\alpha$ in $\ell^2(\Gamma\sigma\Gamma)$ we denote $\alpha|_{\sigma'\Gamma}$
or ${}_{\Gamma{\sigma'}}|\alpha$ the projection $\mathcal{P}_{\sigma'\Gamma}(\alpha)$
and $\mathcal{P}_{\Gamma{\sigma'}}(\alpha)$.

We now prove various formulas of the multiplication of $\theta(X),\theta(Y)$
where $X$, $Y$ are in the intertwiners set $I_{\sigma_1}$, $I_{\sigma_2}$
for various $\sigma_1$, $\sigma_2$.

\end{defn}

The multiplicativity property for $\theta$ is then as follows:

\begin{prop}
We assume that $G$ is a discrete group containing $\Gamma$ almost normal.
For $\sigma$ in $G$ denote $\Gamma_\sigma = \Gamma \cap \sigma \Gamma {\sigma^{-1}}$.

Assume $\varepsilon$ is a cocycle on $G$ coming from a projective representation $\pi$ of $G$.
Let $\sigma_1,\sigma_2$ in $\Gamma$ and $X,Y$ in $I_{\sigma_1}$, $I_{\sigma_2}$ respectively. Assume $X=\pi(\sigma_1)$ and $Y=\pi(\sigma_2)$.

We consider the algebra $M = \L(G,\varepsilon)$, $N = \L(\Gamma,\varepsilon)$, and by
$N_\sigma=\L(\Gamma_\sigma,\varepsilon|_{\Gamma_{\sigma}})$ we denote  the corresponding subalgebras for $\sigma \in G$.
Denote the basis of $\L(G,\varepsilon)$ by $u_g$, and note that 
$u_{g_1}u_{g_2} =\varepsilon(g_1,g_2)u_{g_1g_2}$. 

We have:

$(1)$ The coefficient of $\alpha \in [\Gamma{\sigma_1} \Gamma{\sigma_2}\Gamma]$ in $\theta(X)\theta(Y)$
is given by the formula
$$
\sum\overline{\langle(r_1Xr_2)(r'_2Yr_3)I,I\rangle}
\varepsilon(r_1\sigma_1r_2,r'_2\sigma_2r'_3),
$$
where the sum runs over all $r_1,r_2,r_2',r_3$ in $\Gamma$ such that 
$(r_1\sigma_1r_2)(r_2'\sigma_2r_3)=\alpha$, with no repetitions of the type
$[(r_1\sigma_1(r_2\gamma)][(\gamma^{-1}r_2')\sigma_2r_3]$
allowed.

Note that if $\pi$ is a representation of $G$ extending the left regular representation
and $X = \pi(\sigma_1)$, $Y=\pi(\sigma_2)$ then the summand becomes
$\overline{\langle\pi(\alpha)I,I\rangle}$.

$(2)$ For all $s$ in $\Gamma$, $\sigma_1,\sigma_2$ in $G$, $X$ in $I_{\sigma_1}$, $Y$ in $I_{\sigma_2}$
$$
\theta(X){}_{[\Gamma{\sigma_2]}}\Big|{\theta(Y)}=
\sum{}_{[\Gamma{\sigma_1}s_j\sigma_2]}\Big|{\theta(Xs_jY)},
$$
where $s_j$ are a system of coset representatives for $\Gamma_{\sigma_1^{-1}}$ in $\Gamma$.
\end{prop}

\begin{proof}
It is clear that (2)  is a  consequence of formula (1); consequently, we will only prove (1).

First note that the following identity
$$
\varepsilon(\sigma_1,\sigma_2)=\varepsilon(\sigma_1\gamma,\sigma_2\gamma^{-1})
\varepsilon(\sigma_1,\gamma)\varepsilon(\gamma^{-1}\sigma_2)
$$
is a consequence of the projectivity property of a representation $\pi$ having $\varepsilon$ as a 
cocycle.

Indeed, just expand in two ways
$$
\pi(\sigma_1\sigma_2)=\pi((\sigma_1\gamma)(\gamma^{-1}\sigma_2)).
$$
Recall that
$$
\theta(X)=\sum_{\theta=\gamma_1\sigma\gamma_2\in\Gamma \sigma_1\Gamma}
\overline{\langle\pi(\gamma_1\sigma_1\gamma_2)I,I\rangle}u_{\theta},
$$
$$
\theta(Y)=\sum_{\theta=\gamma'_2\sigma_2\gamma_3\in\Gamma\sigma_2\Gamma}
\overline{\langle\pi(\gamma'_2\sigma_2\gamma_3)I,I\rangle}u_{\theta}.
$$

We want to compute the coefficient of $u_\alpha$ in $\theta(X) \theta(Y)$, where $$\alpha = (\gamma_1 \sigma_1 \gamma_2)
(\gamma_2' \sigma_2 \gamma_3).$$

We will compute the sum of  all the terms corresponding to non-allowable repetitions. Denote
$\sigma = \gamma_1 \sigma_1 \gamma_2$, $\sigma' = \gamma_2' \sigma_2 \gamma_3$.
Since $$u_{\sigma r} u_{r^{-1}\sigma'} = \varepsilon(\sigma r, r^{-1}\sigma')u_{\sigma\sigma'},$$
the  sum of this coefficients will be
\begin{gather*}
\sum_r\varepsilon(\sigma r,r^{-1}\sigma')\overline{\langle\pi(\sigma r)I,I\rangle}\,
\overline{\langle \pi(r^{-1}\sigma')I,I\rangle}=\\
=\sum_r\varepsilon(\sigma r,r^{-1}\sigma')\varepsilon(\sigma ,r)\varepsilon( r^{-1},\sigma')
\overline{\langle \pi(r)1,\pi(\sigma)^*I\rangle}\,
\overline{\langle(\pi(\sigma')I,\pi(r)I\rangle}=\\
=\sum_r\varepsilon(\sigma,\sigma')\langle(\pi(\sigma)^*I,\pi(r)I\rangle
\overline{\langle(\pi(\sigma')I,\pi(r)I\rangle}
\end{gather*}
which, since $\pi(r)I, r\in \Gamma$ is an orthonormal basis, is equal to
\begin{gather*}
\varepsilon(\sigma,\sigma')\langle(\pi(\sigma)^*I,\pi(\sigma')I\rangle
=\varepsilon(\sigma,\sigma')\langle I,(\pi(\sigma)\pi(\sigma')I\rangle=\\
=\varepsilon(\sigma,\sigma')\overline{\langle\pi(\sigma)\pi(\sigma')I,I\rangle}=
\overline{\langle\pi(\sigma \sigma')I,I\rangle}.
\end{gather*}

This completes the proof of formula (1), and the other two are simple consequences.
\end{proof}

Using  formula (3) we obtain  a generalization of the composition formula for the completely positive  maps from Chapter 5.

\begin{prop}
Let $\sigma_1,\sigma_2$ be elements in $G$, and $A,B$ in $I_{\sigma_1}$, $C,D$ in $I_{\sigma_2}$.
Let $N_{\sigma_j}= \L(\Gamma_{\sigma_j})$, $j=1,2$.

Let $I_{\sigma_1,\sigma_2}=\{\sigma_3\mid[\Gamma\sigma_3\Gamma]\subseteq[\Gamma\sigma_1\Gamma\sigma_2\Gamma]\}$
and let $X_{\sigma_3}, Y_{\sigma_3}$ be the $I_{\sigma_3}$ intertwiners that are obtained by taking products
of the form $Ds_iB$, $Cs_iA$, where $s_i$ is a system of representatives for $\Gamma_{\sigma_2}$.

Let $\Psi_{AB} =[\Gamma:\Gamma_{\sigma_1}] E^{\L(G,\varepsilon)}_{\L(\Gamma,\varepsilon)} (\theta(A)\,\cdot\, \theta(B))$,
$\Psi_{CD} = [\Gamma:\Gamma_{\sigma_2}]E^{\L(G,\varepsilon)}_{\L(\Gamma,\varepsilon)} (\theta(C)\cdot \theta(D))$.
Then
$$
\Psi_{CD}\circ\Psi_{AB} = \sum N_{\sigma_1\sigma_2}^{\sigma_3}
\Psi_{X_{\sigma_3},Y_{\sigma_3}},
$$
where $\sigma_3$ runs over $I_{\sigma_1,\sigma_2}$, and $N_{\sigma_1,\sigma_2}^{\sigma_3}$
are the multiplicities.
\end{prop}

\begin{proof}
Fix $u$ a unitary in $M'$. Note that $u\theta(X) u^*$ denoted by $\theta_u(X)$ has the same properties as
$\theta(X)$, as it is obtained by using the cyclic vector $u$ instead of  the unit vector  $1$ in the matrix coefficient computations for the map $\theta$.

Then $E(\theta(A)u \theta(B)) = E(\theta(A) \theta_u(B))u$ where $E = E^{\L(G,\varepsilon)}_{\L(\Gamma,\varepsilon)}$. Let  $s_i$ be a system of representatives 
for $\Gamma_{\sigma_1^{-1}}$ in $\Gamma$.

Applying the condition expectation we obtain that
$$
E(\theta(A) \theta_u(B))u=
\sum_i 
\big [\theta(A) |_{\Gamma{\sigma_1} s_i}\big]\big[ {}_{s_i^{-1}\sigma_1^{-1}\Gamma}\Big|{( \theta_u(B))\big]u}\;.
$$
Apply $\theta(C),\theta(D)$ to the right and left, taking $r_j$ a system of representations
for $\Gamma_{\sigma_2^{-1}}$ in $\Gamma$, we get by  using formula (3) in the preceding statement that formula the following expression for 
$C(\Psi_{AB}(u))D^*$:
$$
C(\Psi_{AB}(u))D^*=
\sum_{j,k,i}
\big[\theta(CA) \big|_{\Gamma_{\sigma_2}r_j\sigma_1 s_i}\big]\big[ {}_{s_1^{-1}\sigma_1^{-1}r_j^{-1}\sigma_2^{-1}\Gamma}\Big|{ \theta_u(BD)\big] u}.
$$
When applying $E_{\L(\Gamma,\varepsilon)}^{\L(G,\varepsilon)}$, only the terms with  $j=k$ will remain  in the above formula. The conclusion follows 
from the fact that the cosets $[\Gamma\sigma_2 r_j\sigma_1s_i]$
when grouped into double cosets will make a list of   the double cosets in the product $[\Gamma\sigma_1\Gamma][\Gamma\sigma_2\Gamma]$,
with multiplicities taken into account.
\end{proof}

\section*{Appendix 2\\
A more primitive structure of the Hecke algebra}

Behind the structure of the Hecke algebras of double cosets of an almost
normal subgroup $\Gamma$ of $G$ (discrete and countable) there exists in fact  a more natural pairing operation
between left and right cosets, which in fact gives all the information about the multiplication structure and the embedding of the Hecke algebra. We refers to  this structure as to  a "primary structure" of the Hecke algebra.

We prove here that our construction in chapters 2,3  is in fact a  representation
of the primary structure of the Hecke algebra.

First, we describe this primitive structure of the Hecke algebra.

Let $\H_0=\C(\Gamma\setminus G/\Gamma)$ be the algebra of double cosets, which
is represented either on $\ell^2(\Gamma\setminus G)$ or $\ell^2(G/ \Gamma)$ (by left or right convolution).

\begin{defn} The "primary structure" of   the Hecke algebra. (This is an operator system in the sense of Pisier ([Pi])). Let $\tilde{C}$ be 
the vector space of sets of the form $[\sigma_1 \Gamma \sigma_2]$, $\sigma_1,\sigma_2\in G$. 
We let $\cC(G,\Gamma)$ be the vector space obtained from $\tilde{C}$,
by factorizing by the subspace generated by  the linear relations of the form
$$
\sum[\sigma_1^i\Gamma\sigma_2^i]=\sum[\theta_1^j\Gamma\theta_2^j]
$$
if $\sigma_s^i, \theta_r^j$ are elements of $G$, and the disjoint union $\sigma_1^i \Gamma \sigma_2^i$
is equal to the disjoint union~$\theta_1^j \Gamma \theta_2^j$.

Then there exist a natural bilinear pairing  
$\C(\Gamma\setminus G) \times \C(G/\Gamma)\to \cC(G,\Gamma)$
extending  the usual product of the Hecke algebra. (Note that the Hecke algebra of double cosets is contained in $\cC(G,\Gamma)$. We obtain a natural isomorphism, by considering the tensor product $\C(\Gamma\setminus G) \otimes_{\H_0} \C(G/\Gamma)$ and extending the bilinear map to the tensor product). 
\end{defn}

We prove that our construction in Chapter 4, beyond proving a representation of 
the Hecke algebra (and of its subjacent left and right Hilbert space module) in $\L(G)$
it also gives a representation of the more primitive structure described above.
The proof of the  following theorem is contained in what we proved in  Chapter 4.

\begin{thm}\label{primary}
Let $G$ be a countable discrete group and let $\Gamma \subseteq G$ be an almost normal subgroup.
Assume that there exists a projective representation $\pi$ with cocycle $\varepsilon$ of $G$,
which, when restricted to $\Gamma$ is unitarily equivalent with the left regular
representation $\lambda_{\Gamma,\varepsilon}$ of $\Gamma$ on $\ell^2(\Gamma)$.
For $\sigma$ in $G$, let $t^{\Gamma\sigma}=(t^{\sigma^{-1}\Gamma})^*$ be the $\L(\Gamma_\sigma)$-unitary
element (that is $E_{\L(\Gamma_{\sigma,\varepsilon})}((t^{\Gamma\sigma})^* t^{\Gamma\sigma})=1$) constructed in 
Chapter 4. 

Moreover, we proved in Chapter 4 that the elements $(t^{\Gamma\sigma})$, where $\Gamma\sigma$ runs over a system of representatives of right cosets of $\Gamma$ in $G$,
form a Pimsner-Popa basis for $L(\Gamma) \subseteq L(G)$.

Then the map $\Phi : \ell^2(\Gamma/G)\to \L(G,\varepsilon)$ mapping $\Gamma\sigma$ into $t^{\Gamma\sigma}$
along with its dual $\tilde{\Phi} : \ell^2(G\setminus \Gamma)\to \L(G,\varepsilon)$ mapping (mapping
$\sigma\Gamma$ into $t^{\sigma\Gamma}$) extends to a representation of the "primary" structure, by defining:
$\Phi^2 (\sigma_1\Gamma\sigma_2) = t^{\sigma_1\Gamma\sigma_2} = \sum_{\theta\in [\sigma_1\Gamma\sigma_2]}
\overline{\langle\pi(\theta)I,I\rangle}\theta$.

In particular, $t^{\Gamma_\sigma}$ is determined by the following identity:
$$
\sum (t^{\Gamma{\sigma_i^1}})^*(t^{\Gamma{\sigma_i^2}})=\sum(t^{\Gamma\theta_i^1})^*(t^{\Gamma\theta_i^2})
$$
if the disjoint union $\bigcup \sigma_1^i \Gamma\sigma_2^i$ is equal to the disjoint union
$\bigcup \theta_j^1 \Gamma\theta_j^2$. Moreover, $t^{\Gamma\sigma}$ is $L^2(\Gamma_\sigma) \cap \L(G)$.
\end{thm}

\begin{rem}
There exists a remarkable pairing involving $\Phi^2$, which is defined by the following formula:
$$\chi([\sigma_1\Gamma\sigma_2],[\sigma_3\Gamma\sigma_4])=
\tau(\tilde\Phi([\sigma_1\Gamma\sigma_2])\tilde\Phi([\sigma_3\Gamma\sigma_4])).$$

It is easy to compute that 
$$\chi([\sigma_1\Gamma\sigma_2],[\sigma_3\Gamma\sigma_4])=\sum_{\theta\in \sigma_1\Gamma\sigma_2 \cap \sigma_3\Gamma\sigma_4} |t(\theta)|^2.
$$

Moreover, $\chi$ has special  positivity properties that  derive from  the fact that $\tau$ is a trace ($\chi$ is a cyclic Hilbert space product in the sense of [Ra4]).
$$
\chi([\sigma_1\Gamma\sigma_2],[\sigma_3\Gamma\sigma_4]) = 
\tau(t^{\sigma_1\Gamma}(t^{\sigma_2\Gamma})^* t^{\sigma_3\Gamma}(t^{\sigma_4\Gamma})^*),
$$
\end{rem}

\begin{proof}
The proof of the representation $\tilde\Phi$ is a straight consequence of the 
identity (proved in Chapter 4, see also the preceding Appendix for the cocycle 
$\varepsilon$).

$$
t(\theta_1\theta_2)=\varepsilon(\theta_1,\theta_2)\sum_{\gamma\in\Gamma} t(\theta_1\gamma)
t(\gamma^{-1} \theta_2)
$$
which implies
$$
(t^{\Gamma{\sigma_1}})^* t^{\Gamma{\sigma_2}}= t^{\sigma_1\Gamma\sigma_2}.
$$
\end{proof}

The existence of the representation $\Phi^2$ is equivalent the
existence of the unitary representation, extending to $G$ the left regular representation of $\Gamma$, as explained bellow.

\begin{prop}
Assume that there exists  a representation (as in Theorem \ref{primary}) $\Phi^2, \Phi, \tilde\Phi$ of  $\cC(G,\Gamma)$,
$\C(\Gamma\setminus G)$, $ \C(G/\Gamma)$,
 into the algebra
$\L(G,\varepsilon)$. Assume that $a^{\Gamma\sigma} = \Phi(\Gamma\sigma)$ is a Pimsner-Popa basis
for $L(\Gamma)\subseteq L(G)$ such that in addition $\Phi([\Gamma\sigma])$ belongs to $\ell^2([\Gamma\sigma])$. 

We also assume the following property, which is implicit in the statement of Theorem \ref{primary}. Let $\beta, \alpha$ be cosets, $\alpha$  of the form $\sigma_1\Gamma\sigma_2$, $\sigma_1,\sigma_2 \in G$ and $\beta$=$\sigma_3\Gamma$ or $\Gamma\sigma_3$, $\sigma_3 \in G$. Let $I$ be their intersection and $P_I$ the projection from $\ell^2(G)$ onto the subspace $\ell^2(I)$. 
 We assume that  $P_I(\Phi^2(\alpha))=P_I(\Phi (\beta))$, for all $\alpha, \beta$, as above.

Then there exists a projective unitary representation $\pi$ of $G$ onto $\ell^2(\Gamma)$, extending the left regular representation
of $\Gamma$ with cocycle $\varepsilon$, on $\ell^2(\Gamma)$. Moreover, $\pi$ is projective with cocycle~$\varepsilon$. Through the construction in Theorem \ref{primary}, the representation $\pi$ corresponds to the representation $\Phi^2$ in the hypothesis of this statement.

We also assume that  
$a^{\Gamma}$ is the identity element in the group $G$.
 Then the condition $a^{\sigma_1\Gamma} a^{\Gamma\sigma_2}= \Phi^2 (\sigma_1\Gamma\sigma_2)$ implies the above Pimsner-Popa basis condition.  \end{prop}

\begin{proof}  Denote the basis of $\L(G,\varepsilon)$ by $u_g$, and note that 
$u_{g_1}u_{g_2} =\varepsilon(g_1,g_2)u_{g_1g_2}$, $g_1,g_2\in G$. 

Let $\sigma$ an element of $G$, $s_i$ a set of representatives for $\Gamma_{\sigma^{-1}}$ in $\Gamma$. Then
we define
$$
\pi(\sigma) s_i = \varepsilon (\sigma,s_i)[t^{\Gamma\sigma s_i}(\sigma s_i)^{-1}]^*\in\L(\Gamma, \varepsilon) .
$$

Then the fact that $\pi(\sigma)$ is a representation follows form the identity $a(\theta_1\theta_2)=
\sum a(\theta_1\gamma)a(\gamma^{-1}\theta_2) \varepsilon(\sigma_1,\sigma_2)$.
Here $a(\theta)$ is the $u_\theta$ coefficient of $a^{[\Gamma\sigma]}$. The identity is a consequence of the
fact that $a^{\sigma_1\Gamma}a^{\Gamma{\sigma_2}}$ depends only on the set $\sigma_1\Gamma\sigma_2$ and
of the fact that $a^{\Gamma\sigma_1\Gamma}a^{\Gamma\sigma_2\Gamma}=\sum_{\sigma_3} N_{\sigma_1\sigma_2}^{\sigma_3}
a^{\Gamma\sigma_3\Gamma}$, where $N_{\sigma_1\sigma_2}^{\sigma_3}$ are the multiplicities from the Hecke
algebra structure. The fact that $\pi(\sigma)$ is a unitary follows from the last condition in the statement.
\end{proof}

We also note that the free $\C$-algebra generated by left or right cosets, subject admits a canonical $C^*$-representation (in fact a representation into $\L(G, \varepsilon)$, in the above terms.

\begin{thm}Let $\A (G,\Gamma)$ the free ${}^* -\C$-algebra generated by all the cosets   $[\Gamma\sigma]$, $\sigma \in  G$, and their adjoints
($[\Gamma\sigma]^\ast=[\sigma^{-1} \Gamma]$, subject to the relation $$
\sum[\sigma_1^i\Gamma][\Gamma\sigma_2^i]=\sum[\theta_1^j\Gamma][\Gamma \theta_2^j]
$$
if $\sigma_s^i, \theta_r^j$ are elements of $G$, and the disjoint union $\sigma_1^i \Gamma \sigma_2^i$
is equal to the disjoint union~$\theta_1^j \Gamma \theta_2^j$. Note that the above relation corresponds exactly to the fact that the Hecke algebra
of double cosets is a  canonical subalgebra of $\A (G,\Gamma)$ , with the trivial embedding mapping a double coset into the formal sum of its left or right cosets (using representatives).

Then we have that 
the ${}^*- \C$-algebra $\A (G,\Gamma)$  admits at least one unital $C^*$ algebra representation into $\L(G, \varepsilon)$.
\end{thm}
\begin{proof}
This is a trivial consequence of the relation described above , by mapping the coset $[\Gamma\sigma]$, $\sigma \in  G$ into $t^{[\Gamma\sigma]}.$
\end{proof}


\section*{Appendix 3\\
Properties of the "square root state" of the state measuring the displacement of a fundamental domain}

Let $\mathbb H$ be the upper half plane and let $F$ be a fundamental domain for the action of the
group $\Gamma=\PSL_2(\Z)$. Let $\mu$ be the canonical  $\PSL_2(\Z)$ invariant measure on $\mathbb H$.

Let $\varphi$ be the positive state on $G=\PGL_2(\Z[\frac1{p}])$ defined by
$$
\varphi(g)=\mu(gF\cap F),\quad g\in G.
$$

Consider $\mathcal F$ be the set of states on $G$ that are obtained as
follows. Let $(F_i)_{i=1}^n$ be a partition of $F$ with measurable sets. Let 
$\varphi_{ij}(g)=\mu(gF_i\cap F_j)$.

For every family  $(\xi_i)_{i=1}^n$ in $\C$ of scalars, consider the state
$$
\sum\xi_i\xi_j\varphi_{ij}
$$

The set $\mathcal F$ is the collection of all such states.

Given any state $\varphi_1$ on $G$ such that the restriction of 
$\varphi_1$ to any coset $\Gamma s$ belongs to $\ell^2(\Gamma s)$,
we define a state $\theta(\varphi_1)$ on the Hecke
algebra of $\Gamma$ in $G$ by
$$
\theta(\varphi_1)(\Gamma\sigma\Gamma)=\sum_{\theta\in\Gamma\sigma\Gamma}\varphi_1(\theta).
$$
To prove the Ramanujan-Petersson conjecture one should prove that
$\theta(\varphi_1)$ is  continuous on the reduced Hecke algebra for any $\varphi_1$ in  $\mathcal F$ as above
(with $\sum \xi_i \mu(F_i)=0$).

Our approach is based on the existence of a ``square root'' of the state of the type $\mu(gF\cap F)$
as above.

Assume that $G$ is an abstract discrete group, $\Gamma$ is a discrete (infinite subgroup), $X$ 
is an infinite measure space with measure $\mu$ and assume that $G$ acts on $X$ by preserving the measure.
Also we assume that $F\subseteq X$ is subset of measure $1$, that is a fundamental 
domain for $\Gamma$ (in particular, we assume that $X=\bigcup\limits_{\gamma\in\Gamma}\gamma F$).

Let as before $\varphi_X$ be the positive definite function $G$, defined by
$$
\varphi_X(g)=\mu(gF\cap F), \quad g\in G.
$$

Then $\varphi_X|_{\Gamma}$ is zero, unless we evaluate at the neutral element.

We assume that $\varphi_X$ has a square root, that is there exists a positive definite
state on $G$ such that
$$
|\varphi_0(g)|^2=\varphi_X(g),\quad g\in G.
$$
Here we may also assume, with no loss in the conclusion that $\varphi_0$ is positive definite on the group algebra of $G$ twisted by a cocycle.

Then $ \varphi_0$ has built in very strong algebraic identifies, that may be
derived as follows.

The fact that $\sum\limits_{\gamma\in \Gamma}\varphi_{X}(\gamma g)=1$
for all $g\in G$, implies that in the $GNS$ representation 
$(H_{\varphi_0},\pi_{\varphi_0},\xi_{\varphi_0})$ of the state
$\varphi_0$ (see, e.g., [Dix]) we have that
$\pi\xi_{\varphi_0}$ belongs to closed linear square of 
$\pi(\gamma)\xi_{\varphi_0}$, $\gamma\in\Gamma$.

In particular, $H_0$, the Hilbert closure of
$\pi(\Gamma)\xi_{\varphi_0}\subseteq H_{\varphi_0}$ is invariant under $G$.
Moreover, the vectors $\{\pi(\gamma)\xi_{\varphi_0}\mid\gamma\in\Gamma\}$ are an orthonormal basis for
$H_0$.

If we apply the Parseval Identity, with respect to this basis of $H_0$ we obtain the following identity
(with $\pi=\pi_{\varphi_0}\mid_{H_0}$)
\begin{gather*}
\varphi_0(g_1g_2)=\langle\pi(g_1g_2)\xi_{\varphi_0},\xi_{\varphi_0}\rangle
=\langle\pi(g_2)\xi_{\varphi_0},\pi(g_1^{-1})\xi_{\varphi_0}\rangle=\\
=\sum_{\gamma\in\Gamma}\langle\pi(g_2)\xi_{\varphi_0},\pi(\gamma)\xi_{\varphi_0}\rangle
\overline{\langle\pi(g_1^{-1})\xi_{\varphi_0},\pi(\gamma)\xi_{\varphi_0}\rangle}=\\
=\sum_{\gamma\in\Gamma}\varphi_0(\gamma^{-1}g_2)\overline{\langle\xi_{\varphi_0},\pi(g\gamma)\xi_{\varphi_0}\rangle}
=\sum_{\gamma\in\Gamma}\varphi_0(g\gamma)\varphi_0(\gamma^{-1}g_2).
\end{gather*}
Thus we obtain that for all $g_1,g_2\in G$ we have that
$$
\varphi_0(g_1g_2)=\sum_{\gamma\in\Gamma}\varphi_0(g\gamma)\varphi_0(\gamma^{-1}g_2)\leqno (*)
$$
Thus if write for a coset $g\Gamma$ of $\Gamma$ in $G$
$$
\varphi_0\mid_{g\Gamma}=\sum_{\theta\in g\Gamma}\overline{\varphi_0(\theta)}\theta
$$
then the identity $(*)$ has as consequence that
$$
\varphi_0\mid_{g_1\Gamma}*\varphi_0\mid_{\Gamma g_2}=
\varphi_0\mid_{g_1\Gamma g_2}\leqno(**)
$$
where $\varphi_0\mid_A=\sum\limits_{\theta\in A}\overline{\varphi_0(\theta)}\theta.$

The identity $(**)$,  in the case of an almost normal subgroup $\Gamma$ of $G$,  it is  proven in this paper to
be equivalent to the fact that the map
$$
[\Gamma\sigma\Gamma]\to\varphi_0\mid_{\Gamma\sigma\Gamma}\in\ell^2(\Gamma\sigma\Gamma),\quad
\sigma\in G
$$
is a representation of Hecke algebra.

If the state $\varphi_0$ is represented as the matrix coefficient of the representation $\pi$ of $G$ on $H_0$,
that is $\varphi_0(g)=\langle\pi(g)\xi_{\varphi_0},\xi_{\varphi_0}\rangle$,
then we may replace the states in the set $\mathcal F$ (that are supposed to
be obtained from vectors orthogonal to the constant function), by the states
$$
g\to\varphi_0(g)\overline{\varphi_0(\gamma^{-1}g\gamma)},\quad \gamma\in \Gamma.\leqno (***)
$$

For $\gamma \in \Gamma$ we denote  by $\pi^\gamma$ the representation $g\to\pi(\gamma^{-1}g\gamma)$ on $H_0$.
Then the analysis of the states in $\mathcal F$ may be replaced by the matrix coefficients
of $\pi \otimes \overline{\pi^\gamma}$, evaluated at the vector
$\xi_{\varphi_{_X}} \otimes \xi_{\varphi_{_X}}$.


In our case the square root state is provided by Jones' theorem, as the discrete series
representation $\pi_{13}$ of $\PSL_2(\Z)$ verifies all of these conditions 
Moreover, the identities $(**)$ give a ``double''
representation of the Hecke algebra, which allows to analyze the matrix
coefficients in $(***)$  as $\gamma \to \infty$.


\

\section*{Appendix 4. A two variable version of the Hecke operators}

In this appendix we are constructing  a new type of  representation of the Hecke algebra.

Let $\Gamma\subseteq G$ be a countable discrete group with an infinite,
almost normal subgroup $\Gamma$. We assume that we are given a directed net $\S$
(closed to intersections) of finite index subgroups of $\Gamma$,
that contains a family of normal subgroups
shrinking to the neutral element, and such that for any $\sigma$ in $G$,
the subgroup $\Gamma_{\sigma}=\sigma\Gamma\sigma^{-1}\cap\Gamma$
belongs to $\S$. Let $K$ be the profinite completion of $\Gamma$
with respect to $\S$. For $g_1,g_2$ in $G$, let 
$\Gamma_{g_1^{-1},g_2}=g_1^{-1}\Gamma g_2\cap\Gamma$
and let $K_{g_1^{-1},g_2}$ be the closure of this coset in $K$. 
Let $\chi_{g_1^{-1},g_2}=\chi_{\overline{\Gamma}_{g_1^{-1},g_2}}$ be
the characteristic function of this coset, viewed as a multiplication operator on $C(K)$.
Then there exists a partial action of $G\times G^{\rm op}$ on $K$, 
defined by $(g_1, g_2)k=g_1\,k\, g_2^{-1}$, for $g_1,g_2$ in $G$, 
$k$ in $K_{{g_1^{-1},g_2}}$.

By $C^*((G\times G^{\rm op})\rtimes C(K))$ we denote the canonical, full groupoid crossed product
$C^*$-algebra associated to this action. Since the Haar measure on $K$,
is invariant under the (partial) action of $G\times G^{\rm op}$,
we also have a reduced groupoid crossed product $C^*$-algebra,
denoted as $C^*_{\rm red}((G\times G^{\rm op})\rtimes C(K))$.

  One starts with a representation of  the groupoid $(G\times G^{\rm op}) \rtimes K$ on a Hilbert space $V$. By restricting to 
$\Gamma \times\Gamma$ invariant vectors in $V$, one obtains a new type of representation of the Hecke algebra. In Example \ref{Heckeas2}, we prove that the  construction in Chapter 5,  is a particular realization of this new representation for the Hecke operators.

We introduce the following definition.

\begin{defn}\label{bivariablehecke}
Let $H$ be a countable discrete group, and consider a unitary representation of $H$ on
a Hilbert space $V$. We denote the action (representation) by $h\cdot v\in V$,
for $h\in H$, $v\in V$. We assume that there exists a Hilbert subspace
$W\subseteq H$, such that $hw_1\perp w_2$ for $h\in H$, $h \neq e$
and $w_1,w_2\in W$ (such a property for $W$ will be called $H$-wandering).
We also assume that $W$ is $H$-generating, that is $V$ is the closure of the span 
$\bigcup\limits_{h\in H}hW$.

We define $V^H(W)$, as the subspace of $H$-invariant vectors on $V$ (with
respect to $W$), consisting of the subspace of the densely defined, $H$-invariant 
functionals on $V$. We  identify the space $V^H(W)$ with the the space of formal sums 
$\sum\limits_{h\in H}hw$, $w\in W$. 
It is an obvious Hilbert space (isomorphic to $W$), with scalar product,
for $w_1,w_2\in W$, defined by
$$
\left\langle\sum _{h_1\in H}h_1w_1,\sum_{h_2\in H} h_2w_2\right\rangle_{V^H(W)}
=\sum _{h_1\in H}\langle h_1w_1,w_2\rangle_{V}=\langle w_1,w_2\rangle_W.
$$

This formalism will be useful for our description of the Hecke operators.
It is obvious that if $(\Y,\nu)$ is a measure space, and $H$ acts by measure
preserving transformations on $\Y$, with a fundamental domain $F$, then
$L^2(\Y,\nu)^H(L^2(F,\nu))$ is obviously isomorphic to $L^2(F,\nu)$
and also to $L^2(\Y^H,\nu)$; the Hilbert space of $H$-invariant functions on $\Y$, with
(Pettersson style) scalar product, defined by the following formula: for $f,g\in L^2(\Y^H,\nu)$
$$
\langle f,g\rangle_{L^2(\Y^H,\nu)}=\int_{F}f\overline{g}d\nu.
$$

Obviously, this scalar product is independent on the choice
of the fundamental domain $F$ for $H$ in $\Y$.

Note that with this definition, the Hilbert space
$H_{12}$, which is acted by the unitary representation of
$\Gamma=\PSL_2(\Z)$ via $\pi_{12}$, does not have $\Gamma$-invariant
vectors, since it doesn't have a wandering subspace, (because the Murray von Neumann
dimension $\dim_{\Gamma}H_{12}<1$, by Jones's formula [GHJ]).
However the modular form $\Delta$ gives a $\Gamma$-invariant,
densely defined functional on $H_{12}$.

With this definition we can describe a new approach to the 
Hecke operators, on bivariant functions, in the presence 
of a unitary representation of $C^*((G\times G^{\rm op})\rtimes C(K))$.
We assume that this representation has a $\Gamma\times\{e\}$
wandering, generating subspace.
\end{defn}

\begin{thm}\label{2heck}
Let $G$ be a countable discrete group, and let $\Gamma$ be an infinite, almost normal subgroup.
We assume that $\S$ is a family of finite index subgroups of $\Gamma$, directed downward, and containing all  the subgroups of the form
$\Gamma_g=g\Gamma g^{-1} \cap \Gamma$, for $g$ in $G$.

Let $K$ be the profinite completion of $\Gamma$ with respect to this
family of subgroups.
Let $V$ be a Hilbert space endowed with a unitary action (representation) 
of the full $C^*$-algebra
$C^*((G\times G^{\rm op})\rtimes C(K))$.

The action of $C(K)$, will
be denoted simply $f\cdot v$, for $f$ in $C(K)$, $v\in V$.
Let $\chi_{g_1^{-1},g_2}$ be the characteristic function of the closure in $K$ of the intersection $g_1^{-1}\Gamma g_2\cap \Gamma$, for $g_1, g_2\in G$.
Then the range in $V$ of the projection $\chi_{g_1^{-1},g_2}$ will be the domain for the partial isometry on $V$ defined by the action of $(g_1,g_2)\in G\times G^{\rm op}$.

We denote the action of $(g_1,g_2)\in G\times G^{\rm op}$ on a vector
$v$ in $\chi_{g_1^{-1},g_2}V$ by 
$$
(g_1,g_2)v=(g_1,g_2)(\chi_{g_1^{-1},g_2}v)=
g_1(\chi_{g_1^{-1},g_2}v)g_2^{-1},\quad g_1,g_2\in G,\ v\in V.
$$

Assume that the unitary representation of $\Gamma\times \{e\}$ on $V$,
obtained by restriction of the action of $G\times G^{\rm op}$,
admits a $\Gamma\times\{e\}$ wandering  subspace $W_0$
such  that the translations by elements in $\Gamma\times\{e\}$ of 
 $W_0$ cover $V$.

Let $\sigma$ be an arbitrary element in $G$, and assume that
$[\Gamma\sigma\Gamma]$ is the disjoint union of the cosets $s_i\sigma\Gamma$,
$s_i\in \Gamma$, $i=1,2,\ldots,[\Gamma:\Gamma_{\sigma}]$.
We use the obvious extension of the action of $C(K)$ on
$V^{\Gamma\times \{e\}}(W_0)$, which maps a $\Gamma$ invariant vector into a vector that is invariant with respect to a smaller subgroup.

For $v$ in $V^{\Gamma\times  \{e\}}(W_0)$, we define
$$
\Pi(\sigma)v=\sum\limits_is_i\sigma(\chi_{\sigma^{-1},\sigma}v)\sigma^{-1}.
$$
We have the following: if the vector  $v$ has the expression
$
v=\sum\limits_{\gamma\in\Gamma}\gamma w,
$
for some $w\in V^{\Gamma\times \{e\}}(W_0)$, then
$$
\Pi(\sigma)\bigg(\sum_{\gamma}\gamma w\bigg)=\sum_is_i\sigma\bigg(\sum_{\gamma\in\Gamma}
\chi_{\sigma^{-1},\sigma}\gamma w\bigg)\sigma^{-1}=
\sum_{\theta\in\Gamma\sigma\Gamma}\theta(\chi_{\theta^{-1},\sigma}w)\sigma^{-1}.
$$
Moreover, in this case, $\Pi$ is a unitary representation of $G$ on $V^{\Gamma\times  \{e\}}(W_0)$.\end{thm}

\begin{obs}\label{doublefund}
Assume in addition that the unitary action of $\Gamma\times\Gamma$ on $V$,
obtained by restriction of the action of $G\times G^{\rm op}$,
admits a $(\Gamma\times\Gamma)$ wandering, generating subspace $W$
(and thus $W_0=\overline{{\rm Sp}(W\gamma\mid\gamma\in \Gamma)}$ is a $(\Gamma\times \{e\})$-wandering,
generating, subspace of $V$).
Since the representation $\Pi$ acts on $(\Gamma\times \{e\})$ equivariant
vectors on $V$ and $W_1 =\overline{{\rm Sp}(\gamma W\mid\gamma\in \Gamma)}$ is wandering Hilbert subspace for $\{e\}\times\Gamma$, we
can define a unitary representation of the Hecke operators on 
$(V^{\Gamma\times \{e\}}(W_0))^{\{e\}\otimes\Gamma}(W_1)=V^{\Gamma\times\Gamma}(W)$
by defining for a double coset $[\Gamma\sigma\Gamma]$ of $G$, and for $v$ in $V^{\Gamma\times\Gamma}(W_0)$,
$$
T([\Gamma\sigma\Gamma])(v)=\sum_{i,j}s_i\sigma\big(\chi_{\sigma^{-1},\sigma}v\big)\sigma^{-1}s_j.
$$
If $v$ is given as 
$\sum\limits_{\gamma_1,\gamma_2\in\Gamma}\gamma_1w\gamma_2^{-1}$, for
$w\in W$,  the formula has the expression:
$$
\sum_{\theta_1,\theta_2\in\Gamma\sigma\Gamma}\theta_1\big(\chi_{\theta_1^{-1},\theta_2}w\big)\theta_2^{-1}.
$$
We will prove in  Appendix 7 that in fact the Hecke operators, that we introduced in Chapter 3, are of this form.
\end{obs}

\begin{obs}\label{formula2heck}
In the context of Theorem \ref{2heck},  let  $w$ be a vector in $W$, and let $v=\sum\limits_{\gamma\in\Gamma}\gamma w$.
Let $\sigma$ be an element in $G$ and assume that $\Gamma\sigma\Gamma=\bigcup s_i\sigma\Gamma$
where $s_i\in\Gamma$ are coset representatives.
Let $K_i$ be a the closure of $\overline{s_i\Gamma_{\sigma^{-1}}}$ in $K$.
Note that $\Gamma=\bigcup\limits_i s_i\Gamma_{\sigma^{-1}}$.
Consider $w_i=\chi_{K_i}w$. Let 
$$
w_0=\sum s_i^{-1}(\chi_{K_i}w).
$$
Then $\chi_{\sigma^{-1},\sigma}(w_0)=w_0$ and
$$
\sum_{\gamma\in\Gamma}\gamma w=\sum_{\gamma\in\Gamma}\gamma w_0.
$$
The formula for the representation $\Pi(\sigma)$ becomes
\begin{gather*}
\Pi(\sigma)\bigg(\sum_{\gamma\in\Gamma}\gamma w\bigg)=\sum_is_i\sigma\bigg(\chi_{\sigma^{-1},\sigma}
\Big(\sum_{\gamma\in\Gamma}\gamma w_0\Big)\bigg)\sigma^{-1} \\ =
\sum_is_i\sigma\bigg(\sum_{\gamma\in\Gamma_{\sigma^{-1}}}\gamma w_0\bigg)\sigma^{-1}
=
\sum_{\gamma\in\Gamma}\gamma\sigma w_0\sigma^{-1}=\sum_{\theta\in[\Gamma\sigma]}
\theta w_0 \sigma^{-1}.
\end{gather*}
\end{obs}

\begin{proof}
(of Theorem \ref{2heck}).
To prove that $\Pi$ is a representation of $G$ on $V^{\Gamma\times \{e\}}(W_0)$,
take $\sigma_1,\sigma_2\in G$, and
a vector $v=\sum\limits_{\gamma\in\Gamma}\gamma w$ in $V^{\Gamma\times \{e\}}(W_0)$.

We want to prove that
$$
\Pi(\sigma_1)\Pi(\sigma_2)\bigg(\sum_{\gamma\in\gamma}\gamma w \bigg)=\Pi(\sigma_1\sigma_2)\bigg(\sum_{\gamma\in\Gamma}\gamma w\bigg).
$$

By the Observation \ref{doublefund}, letting the subgroup
$$L=L(\sigma_1,\sigma_2)=\sigma_2^{-1}
(\sigma_1^{-1}\Gamma\sigma_1)\sigma_2\cap\sigma_2^{-1}\Gamma_{\sigma_1}\sigma_2\cap\Gamma,
$$
one can replace $w$ in the previous formula, by another vector in $w_0$ such that $\chi_Lw_0=w_0$.
By the formula in Observation \ref{formula2heck}, we have that
$$
\Pi(\sigma_1)v=\Pi(\sigma_1)\bigg(\sum_{\gamma\in\Gamma}\gamma w_0\bigg)=
\sum_{\gamma\in\Gamma}\gamma(\sigma_1w_0\sigma_1^{-1})
$$
and hence 
$\Pi(\sigma_2)(\Pi(\sigma_1)v)=\Pi(\sigma_2)\Big(\gamma\sum\limits_{\gamma\in\Gamma}
\sigma_1w_0\sigma_1^{-1}\Big)$.

But $\chi_{\sigma_2^{-1}\sigma_2}(\sigma_1w_0\sigma^{-1})=\sigma_1w_0\sigma_1^{-1}$
by our assumption, and hence this is equal for 
$\sum\limits_{\gamma\in\Gamma}\gamma\sigma_2\sigma_1w_0\sigma_1^{-1}\sigma_2^{-1}$
which is the formula for
$\Pi(\sigma_2\sigma_1)\Big(\sum\limits_{\gamma\in\Gamma}\gamma w_0\Big)$.

To verify that $\Pi$ is a unitary, it is this sufficient to check that
$\Pi(\sigma)^*=\Pi(\sigma^{-1}),$
i.e., to check for all $\sigma$ in $G$, $w_1$,$w_2$ in $W_0$ we have
$$
\Big\langle\Pi(\sigma)\Big(\sum_{\gamma\in\Gamma}\gamma w_1\Big), 
\sum_{\gamma\in\Gamma}\gamma w_2\Big\rangle_{V^{\Gamma\times \{e\}}(W)}=
\Big\langle\sum_{\gamma\in\Gamma}\gamma w_1, 
\Pi(\sigma^{-1})\sum_{\gamma\in\Gamma}\gamma w_2\Big\rangle_{V^{\Gamma\times \{e\}}}
$$

Using the formula from Observation \ref{formula2heck}, and the definition of the scalar product
on $\Gamma\times \{e\}$ invariant vectors, 
and replacing $w_1$,$w_2$ with vectors in the image of $\chi_{\sigma^{-1},\sigma}$ and
$\chi_{\sigma,\sigma^{-1}}$ respectively, this is then equivalent to
$$
\Big\langle\sum_{\gamma}\gamma \sigma w_i \sigma^{-1},w_2\Big\rangle_{V}=
\Big\langle w_1, 
\sum_{\gamma}\gamma \sigma^{-1}w_2\sigma\Big\rangle_{V}
$$
which holds true because of the unitarity of the action of $G$.
Note that the matrix coefficients of representation $\Pi$ are of the form
$$
\Pi_{w_1,w_2}(\sigma)=\sum_{\theta\in\Gamma\sigma\Gamma}
\Big\langle\theta(\chi_{\theta^{-1},\sigma} w_1)\sigma^{-1},w_2\Big\rangle,
\quad w_1,w_2\in W,\ \sigma\in G.
$$

The formula from Observation \ref{doublefund}, is an obvious consequence of the formula
for the Hecke operators on $\Gamma$-invariant vectors.
\end{proof}

We analyze now the case when the Hilbert $V$ is a Hilbert space of $L^2$ - functions on an infinite measure space $\Y$ and the representation on $V$ comes from the Koopmann representation a groupoid action of $(G \times G) \rtimes K$.

\begin{thm}\label{quotienthecke} Let $(\Y,\nu)$ be an infinite measure space and let
 $V=L^2(\Y,\nu)$. We assume that we have  a module action of $C(K)$  on $L^2(\Y,\nu)$, (that is we assume that we are given
a projection $\pi:\Y\to K$), and we  assume that $G\times G^{\rm op}$ has a groupoid action on
$\Y$. Here we are given a partial action of $G\times G^{\rm op}$ on $\Y$, denoted by 
$g_1 yg_2^{-1}$, defined if $y\in\Y$,  $g_1,g_2 \in G$ and $\pi(y)$ belongs to
$$
K_{g_1^{-1},g_2}=\overline{\Gamma_{g_1^{-1},g_2}}=
\overline{g_1^{-1}\Gamma g_2\cap\Gamma},
$$
the closure of $g_1^{-1}\Gamma g_2\cap\Gamma$ being computed in $K$.

Since $\pi(g_1y g_2^{-1})=g_1\pi(y)g_2^{-1}$, for all $g_1,g_2 \in G, y\in \Y$, it follows that this action gives a unitary representation
$C^*((G\times G^{\rm op})\rtimes C(K))$. The representation is unitary if the action of
$(g_1 ,g_2) \in  G\times G^{\rm op}$ from $\pi^{-1}(K_{g_1^{-1},g_2})$ onto $\pi^{-1}(K_{g_1,g_2^{-1}})$
is a measure preserving transformation (on $\Y$).
Note that the representation of $C^*((G\times G^{\rm op})\rtimes C(K))$ on $L^2(\Y,\nu)$
is  the Koopman representation ([Ke]) associated to the groupoid action of
$(G\times G^{\rm op})\rtimes K$.

Assume that there exists a $\Gamma \times\{e\}$ - fundamental domain $F_1$ in $\Y$. In particular the quotient space $\Gamma \setminus \Y$, with the induced quotient measure $\nu^{\Gamma\setminus \Y}$ is isomorphic to the measure space $(F_1, \nu|_{F_1})$. 

Through the construction from Theorem \ref{2heck}, the Koopmann unitary representation of $C^\ast((G \times G)\rtimes C(K))$ on $L^2(\Y, \nu)$ gives rise to a representation $\Pi$ of $G$ on $L^2(\Gamma \setminus \Y)$. There exists  a canonical measure preserving action $\alpha$ of $G$ on $\Gamma \setminus \Y$, whose associated Koopmann representation is exactly the representation $\Pi$ of $G$ on $L^2(\Gamma \setminus \Y)$.

We assume in additional the following set of conditions:

(FS1) There exists a finite measure subset $F_2$ of $\Gamma \setminus \Y$ whose translates, by $\a(g)$, $g \in G$, cover $\Gamma \setminus \Y$. Moreover assume that there exists an increasing sequence $E_n$ of measurable subsets of $F_2$, whose union is $F_2$ and assume there exist  natural numbers $e_n$, such that the states on $G$ defined  by the formulae 
$$
\psi_{E_n}(g) = \langle \; \Pi(g) \chi_{E_n}, \chi_{E_n} \rangle_{L^2(\Gamma\setminus \Y)} = \nu^{\Gamma\setminus \Y}(\a(g)(F_n)\cap F_n), \; g \in G
$$

\noindent have support in $\mathop{\bigcup}\limits_{e \leq e_n} [\Gamma\sigma_{p^e}\Gamma]$ (in particular this implies that the Koopmann representation of $G$ on $L^2(\Gamma \setminus \Y)$ is tempered (see e.g  [Ke] for the definition of the Koopman representation)).

 (FS2) We assume that restriction of the  groupoid action of $(G \times G) \rtimes K$ to $(\Gamma \times \Gamma) \rtimes K$ gives, through the Koopmann representation, a tempered representation (that is continuous with respect to the $C_{\rm red}^{\ast}((\Gamma \times \Gamma) \rtimes C(K))$ norm).

Then, if the above conditions FS1, FS2, hold true, the Koopmann representation of $C^{\ast}((G \times G) \rtimes C(K))$ on $L^2(\Y, \nu)$ is tempered (continuous with respect to the $C^{\ast}_{\rm red}((G \times G) \rtimes C(K))$ norm).

\end{thm}

\begin{proof}

Let $\Gamma\setminus K$ be the space of left orbits of $\Gamma$ in $G$.
There exists a canonical action of $G$ on $\Gamma\setminus K$, described as follows.

Fix $\sigma$ in $G$, and take an orbit $\Gamma k$ for some $k \in K$.
Define for $\sigma\in G$, the action of $G$ on the orbit $\Gamma k$ by the formula
$$
\alpha (\sigma)(\Gamma k)=
\bigcup\limits_is_i\sigma(\overline{\Gamma_{\sigma^{-1},\sigma}}\cap\Gamma k)\sigma^{-1},
$$
where $\Gamma\sigma\Gamma=\bigcup s_i\Gamma\sigma$, is the coset description.
Equivalently, if $k'=\theta k$ is such that
$\theta k\in\overline{\Gamma_{\sigma^{-1},\sigma}}$, then
$\alpha(\sigma)(\Gamma k)=\alpha(\sigma)(\Gamma\theta k)=\Gamma\sigma(\theta k)\sigma^{-1}$.

This can also be described as 
$$
\alpha(\sigma)[\Gamma k]=[\Gamma\sigma\Gamma] k\sigma^{-1},
$$
with the convention that, if $\Gamma\sigma\Gamma$ in the formula $\Gamma\sigma\Gamma k\sigma^{-1}$ is decomposed
into the cosets, then $\Gamma\sigma s_i k\sigma^{-1}$ is taken to be zero, if
$k$ does not belong to $\overline{\Gamma_{(\sigma s_i)^{-1},\sigma}}$
(i.e., if $k$ not in the corresponding domain).

Then the action of $\Pi(\sigma)$ on $\Gamma$ invariant
functions on $\Y$, is described in the same way.
If we identity the points of $\Gamma \setminus\Y$ with $\Gamma$-orbits,
$\Gamma y$, $y\in\Y$, then the projection $\pi$ induces a projection 
$\tilde{\pi}:\Gamma\setminus\Y\to\Gamma/K$.
The action of $\Pi(\sigma)$ on $\Gamma y$, $y\in\Y$, is described as follows:
choose $\theta\in\Gamma$ so that $y'=\theta y$ belongs to $\pi^{-1}(\overline{\Gamma_{\sigma^{-1},\sigma}})$
and let
$$
\Pi(\sigma)([\Gamma y])=\Pi(\sigma)[\Gamma y']=
\Gamma(\sigma y'\sigma^{-1})
$$
which again as above way be described as
$$
\Pi(\sigma)[\Gamma y]=
[\Gamma\sigma \Gamma]y\sigma^{-1}.
$$

Let $F_0, E_n^0$ be finite measure, $\Gamma \times\{e\}$ wandering subsets of $\Y$, that project in the quotient $\Gamma \setminus \Y$ into the sets $F_2$ and $E_n$ respectively. We may assume that $E_n^0$ is an increasing sequence of subsets of $\Y$, whose union is $F_0$.

The state $\varphi_{E_n^0}$ on $C^{\ast}((G \times G) \rtimes C(K))$, associated through the Koopmann representation to the characteristic function $\chi_{E_n^0}$, is
$$
\varphi_{E_n^0}((g_1, g_2)f) = \int_{\Y}g_1(f\chi_{E_n^0})g_2^{-1}\chi_{E_n^0} d\nu.
$$

\noindent These states converge weakly to the state $\varphi_{F_0}$ on $C^{\ast}((G \times G) \rtimes C(K))$ associated, by the same type of formula, to the characteristic function $\chi_{F_0}$.

Since $\psi_{E_n}$ has support in $\mathop{\bigcup}\limits_{e \leq e_n} \Gamma \sigma_{p^e}\Gamma$ it follows, by the construction of the representation of $G$, in Theorem \ref{2heck}, that $\varphi_{E_n}$ has support in $\mathop{\bigcup}\limits_{e \leq e_n}[\Gamma p^e \Gamma] \times [\Gamma p^e \Gamma] \subseteq G \times G$. Hence by the hypothesis, on the continuity with respect to the reduced $C^*$-algebra norm, of the representation $\Pi$ restricted to $C^{\ast}((\Gamma \times \Gamma) \rtimes C(K))$, it follows that $\varphi_{E_n}$ is a state on $C_{\rm red}^{\ast}((G \times G) \rtimes C(K))$. Hence $\varphi_{F_0}$ is a state on $C_{\rm red}^{\ast}((G \times G) \rtimes C(K))$. Since the translates of $F_0$ through $G \times G$ cover $\Y$, it follows that the associated Koopmann representation is continuous with respect to the $C_{\rm red}^{\ast}((G \times G) \rtimes C(K))$ norm.

\end{proof}

\vskip6pt

Appendix 5. Analysis of the correspondence between states on $C^{\ast}((G \times G) \rtimes C(K))$ and $C^{\ast}(G)$ through the Koopmann representation.

\

\vskip6pt

In this appendix we work in the hypothesis of Theorem \ref{quotienthecke}. Let $(\Y, \nu)$, together  with the action $(G \times G) \rtimes K$ on $\Y$ be as in the above mentioned theorem. Let $\Pi$ be the associated Koopmann representation of $G$ on $L^2(\Gamma \setminus \Y)$. We want to analyze the relation between the states on $C^{\ast}((G \times G) \rtimes C(K))$ induced by $\Gamma \times \{e\}$ wandering subsets $F_0$ of $\Y$ and the state on $C^{\ast}(G)$ induced (through Koopmann representation) by the image $\tilde{F}$ of the set $F_0$ in $\Gamma \setminus \Y$.

Obviously the continuity properties of the state on $C^{\ast}(G)$ do not change if we replace $F_0$ by a $\Gamma$ - equivalent ([Ng]) subset $F_n$ of $\Y$, and hence we  will "shrink" $F_n$ to the fiber $\pi^{-1}(e)$ of $\Y$ over $e\in K$. We obtain a sequence of states on $C^{\ast}((G \times G) \rtimes C(K))$ that weakly converge to a state on $C^{\ast}((G \times G) \rtimes C(K))$ which is "supported" at the neutral element $e$ of $K$, and which by restricting to $G$ (viewed as the subgroup $\{ (g, g^{-1}) \mid g \in G \}$ of $G \times G$) gives the state on $C^{\ast}(G)$ associated to $\tilde{F}$. The procedure is explicitly described in terms of cosets for a family of normal subgroups in $\S$, of $\Gamma$, shrinking to $e$.

Recall that $F_0\subseteq\Y$ is a $(\Gamma\times \{e\})$ wandering
subset. Let $F=\Gamma F_0$. Let $\tilde{F}=\Gamma\setminus F$ be the image of F in the quotient $\Gamma\setminus \Y$. Then $\tilde{F}$   is  a finite measure subset
of $\Gamma\setminus \Y$, with respect to the induced measure on the quotient.

Then we have the following two states on $C^*((G\times G^{\rm op})\rtimes C(K))$ and $C^* (G)$ respectively, defined  (extending by linearity) by the formula
$$
\varphi_{F_0}((g_1,g_2)f)=\int_{\Y}g_1f\chi_{F_0}g_2^{-1}\chi_{F_0} d\nu
$$
for $(g_1,g_2)\in G\times G^{\rm op}$, $f\in C(K)$  (here $(g_1,g_2)f$ is a generic element in $C^*((G\times G^{\rm op})\rtimes C(K))$). For $g$ in $G$ we define
$$
\varPsi_{\tilde{F}}(g)=\langle \Pi(g)\chi_{\tilde{F}},\chi_{\tilde{F}}\rangle_{L^2(\Gamma\setminus\Y)},
$$
where $\Pi$ is the representation introduced in Theorem \ref{quotienthecke}.
The precise relation between the two states is described in following definition.

\begin{defn}\label{defdiag}
Let $\theta$ be the linear map that to every state $\varphi$ on
$C^*((G\times G^{\rm op})\rtimes C(K))$, with positive coefficients (that is
$\varphi((g_1,g_2)f)\geq 0$ for all $(g_1,g_2)\in G\times G^{\rm op}$, and $f$
 a positive function on $K$, associates the functional an $\C(G)$, (which is then extended to a state on $C^* (G)$)  defined by the formula
 $$
 \theta(\varphi)(g)=\sum_{\theta\in\Gamma\sigma\Gamma}\varphi((\theta,g)\chi_{\theta^{-1},g}).
 $$
Then clearly, with the above notations, we have $\theta(\varphi_{F_0})=\varPsi_{\tilde{F}}$.  Moreover if $F_0'$ is $\Gamma$-equivalent to $F_0$ in the sense of ([Ng]), then $\theta(\varphi_{F_0})=\theta(\varphi_{F_0'})$.

\

\

In the following we describe an explicit  process that gives a formula    for the state $\varPsi_{\tilde{F}}$
constructed above, on 
$C^*(G)$, as the restriction to the diagonal C*-subalgebra generated by $\{(g,g^{-1}), g\in G\} 
\subseteq C^*((G\times G^{\rm op})\rtimes C(K))$, of a limit of states 
 of the  type $\varphi_{F_0}$, as above.
\end{defn}

\begin{prop} \label{diag}
Let $\Gamma\subseteq G$, and $K$ as above and let $(\Y,\nu)$
be an infinite measure space, so that $L^2(\Y,\nu)$ is the Hilbert space of  a Koopman unitary 
representation of $C^*((G\times G^{\rm op})\rtimes C(K))$.
Thus, we assume that we are given a surjective projection
$\pi:\Y\to K$ (corresponding to the action of $C(K)$ on $L^2(\Y,\nu)$),
which is $G\times G^{\rm op}$ equivariant, (that is
$
\pi(g_1 y g_2)$=$g_1\pi(y)g_2,$ for $(g_1,g_2)\in G\times G^{\rm op}$,
where $y$ belongs to $\pi^{-1}(K_{g_1^{-1},g_2})\subseteq\Y$). Moreover, we assume that the partial transformations of $G\times G^{\rm op}$ are measure preserving (and hence we assume that the Koopman representation is unitary).

Let $F_0$ be a subset of $\Y$ that is $\Gamma\times \{e\}$ wandering.
Let $F=\Gamma F_0$, $\tilde{F}=\Gamma\backslash F$
and let $\varphi_{F_0}$, $\Psi_{\tilde{F}}$
be the states on $C^*((G\times G^{\rm op})\rtimes C(K))$, and respectively $C^*(G)$,
introduced above, in the Definition \ref{defdiag}.
Thus $\Psi_{\tilde{F}}=\theta(\varphi_{F_0})$.

Then, there exists a sequence of states $\varphi_{F_n}$, for a suitable
choice of $\Gamma\times \{e\}$-wandering subsets $F_n$ of $\Y$, that are $\Gamma\times \{e\}$
equivalent to $F_0$ in the sense of [Ng], such that $\theta(\varphi_{F_n})=
\Psi_{\tilde{F}}$ for all $n$. Moreover the states $\varphi_{F_n}$ are converging
weakly to a state $\varphi_0$ on $C^*((G\times G^{\rm op})\rtimes C(K))$, $\theta(\varphi_0)=\Psi_{\tilde{F}}$. In addition $\varphi_0$ has 
 the property that $\varphi_0((g_1,g_2)f)$ is equal to
$$
\Psi_{\tilde{F}}(g_1)\delta_{g_1,g_2}f(e),\quad g_1,g_2\in G
$$
where $\delta_{g_1,g_2}$ is the Kronecker symbol.
\end{prop}

Before giving the proof of the proposition we make the following
observations describing the structure of the  $G$-measure space $\Gamma\setminus\Y$, (acted
by $G$, through  the transformations constructed in Theorem \ref{diag}).

\begin{obs}\label{adelichilbert}
Along with the subspace $L^2(\Gamma\setminus\Y, \nu_{\Gamma\setminus \Y})$, we consider the
Hilbert space $\H=L^2(\Y^{\rm ad},\nu^\Gamma)$ defined
as the profinite limit of $L^2(\Gamma_i\setminus\Y)$, after $\Gamma_i$ in $\S$. Here $(\Gamma_i)_i$ is a decreasing family of normal subgroups, with trivial intersection.
We let $\nu^\Gamma = \nu_{\Gamma\setminus \Y}$ be the induced measure on the quotient $\Gamma\setminus\Y$.

Then $\H$ is naturally acted by $C^*((G\times G^{\rm op})\rtimes C(K))$.
This is simply because multiplying a $\Gamma$-invariant
function with the characteristic function of the closure of
the coset $\overline{\Gamma_0 s}$, $s\in \Gamma$, gives a function that
is $\Gamma_0$ invariant.

Clearly, the fiber at $e$ of $L^2(\Y^{\rm ad},\nu^\Gamma)$ is $L^2(\Gamma\setminus\Y,\nu^{\Gamma\setminus \Y})$ and the restriction
of the action of $(G\times G^{\rm op})\rtimes K$ to the fiber at $e$ is exactly the
 unitary representation  $\Pi$ of $G$, that we have constructed in Theorem \ref{quotienthecke}, on $L^2(\Gamma\setminus\Y,\nu^{\Gamma\setminus \Y})$.
\end{obs}

\

\

The above Hilbert space is in fact the Hilbert space of
germs of $\Gamma$-invariant functions on $\Y$. Recall that the space $\Y$  admits a fibbering over the compact set $K$.

A more convenient description of such a space of germs
is obtained by considering an adelic completion of $\Y$. By using this representation we obtain an alternative description of the measure on $\Gamma\setminus\Y$ (represented as the
fiber at $e$ in the adelic description) in terms of $\Gamma \times \{e\}$ wandering subsets of $\Y$.

We describe this construction, assuming first only of the action of the group $\Gamma$.

\begin{prop}
Let $(\cX,\mu)$ be an infinite measure space, and assume that we are given an action 
of $\Gamma$ on $\cX$, by measure preserving automorphisms of $\cX$.
We denote the action of $\gamma \in \Gamma$ by $\gamma x$, for  $\gamma\in \Gamma$, $x\in \cX$.

Also, we are given an action of $C(K)$ on $\cX$, equivalently a projection 
$\pi:\cX \to K$, which is also $\Gamma$ equivariant (that is $\pi(\gamma x)=\gamma\pi(x)$,
$\gamma\in\Gamma$, $x\in \cX$).
Thus, via the Koopmann representation we have a representation of
$C^*(G\rtimes C(K))$ on $L^2(\cX,\mu)$.

Let $\cX^{\rm ad}$ be the measure space $K\times_{\Gamma}\cX$, where $\Gamma$ acts on the left
on both $K$ and $\cX$, and the equivalence relation is
$$
(k,x)\sim(\gamma k,\gamma x),\quad k\in K,\ x\in\cX,\ \gamma\in \Gamma.
$$

If $\Gamma$ admits a fundamental domain in $\cX$, then $\cX^{\rm ad}$ has a canonical
measure. Moreover, $\cX^{\rm ad}$ is again fibered over $K$, via
$\tilde{\pi}(k,x)=k^{-1} x$. Let $\cX^{\rm ad}_e$ be the fiber at $e$, 
with the induced measure $\mu_{\cX^{\rm ad}_e}$.
Thus
$$
\cX^{\rm ad}_e=\{(k,x)\in K\times_{\Gamma}\cX\mid k=\pi(x)\}.
$$

Note that every fundamental domain $F$ for $\Gamma$ in
$\cX$ is canonically isomorphic to $\cX^{\rm ad}_e$, simply
by mapping $F$ into $\tilde{F}=\{(\pi(f)^{-1},f)\mid f\in F\}$.
Clearly, this map is surjective.

Then for every two $\Gamma$ wandering subsets $G_1,G_2$ of $\cX$, we denote by
$\tilde{G}_1,\tilde{G}_2$ their image into $\cX^{\rm ad}_e$.

Let $\Gamma_n$ be a family of subgroups in $\S$. (Recall that $\S$ is the family
of finite index subgroups of $\Gamma$ used in order to construct  the profinite completion $K$.)
Let $(s^n_i)_i$ be a system of coset representatives for $\Gamma_n$ in $\Gamma$.

Then, we have that $\mu_{\cX_e}(\tilde{G}_1\cap\tilde{G}_2)$ is equal to the limit of the following increasing sequence
$$
\lim_{n\to\infty}\sum_{i,j}\mu\Big((s_i^n)^{-1}\big[\pi^{-1}(\overline{(s_i^n)\Gamma_n})\cap G_1\big]
\cap
(s_j^n)^{-1}\big[\pi^{-1}\overline{((s_j^n)\Gamma_n})\cap G_2\big]\Big).\leqno(*)
$$
\end{prop}

\begin{proof}
Indeed, the formula for the intersection of $\tilde{G}_1,\tilde{G}_2$ can be written obviously as
$$
\mu_{\cX_e}(\tilde{G}_1\cap\tilde{G}_2)=\mu\big(G_1\cap (\cup_{\gamma}\gamma G_2)\big)
=\sum_{\gamma\in\Gamma}\mu(G_1\cap\gamma G_2).\leqno(**)
$$

Formally, we disintegrate $G_1,G_2$ as measures over $K$ 
$$
G_i=\int_K^{\oplus}(\mu_{G_i})_k dk.
$$
We have that
$$
\mu_{\cX_e}(\tilde{G}_1\cap\tilde{G}_2)=
\iint\limits_{K^2}\big\langle(\mu_{G_1})_k, k^{-1}l(\mu_{G_2})_l\big\rangle,
$$
and translating this at the origin, this gives
$$
\iint\limits_{K^2}\big\langle k^{-1}(\mu_{G_1})_k, l^{-1}(\mu_{G_2})_l\big\rangle
dkdl.
$$
Here by the scalar product of two positive measure
we understand
$$
\langle\mu,\nu\rangle=\int\frac{d\mu}{d\nu}d\nu,
$$
(which could also be $\infty$).

Rigorously the proof is as follows: the sequence on the right hand side of formula $(*)$
is increasing, as the reunion is taking into account more and more intersections,
when $n$ is increasing.

Since the subgroups ($\Gamma_n$) are separating the points of $\Gamma$,
the family $\{(s_i^n)^{-1}$ $(s_j^n)\mid i,j,n\}$
is exhausting the points of $\Gamma$ and hence by formula $(**)$, the two
quantities in the statement are equal.
\end{proof}

The following observation is used only to clarify the relation between the $G$-system obtained as the fiber over $e\in K$ of the adelic system, and the $G$ system in the quotient $\Gamma\backslash\Y$ described in the Theorem \ref{quotienthecke} .

\begin{obs}
Assume that  $\Gamma, G,\Y,\nu$ are as at the beginning of this section with the left and right
action of $G\times G^{\rm op}$ on $\Y$. This action is equivariant with respect to the projection
$\pi:\Y \to K$.

Then  $\Y^{\rm ad}$ also admits a $(G\times G^{\rm op})\rtimes K$ action, defined as follows:

Recall that $\Y^{\rm ad}=K\times\Y$,
is defined by the equivalence
relation defined by requiring that $(\gamma k,\gamma y)$ is equivalent to $(k,y)$ $\gamma\in\Gamma$, $k\in K$, $y\in\Y$. 

Then $G\times G^{\rm op}$ acts as
$$
(g_1,g_2)(k,y)=(kg_1^{-1},yg_2^{-1}).
$$

The projection $\tilde{\pi}:K\times_{\Gamma}\Y$ is $\pi((k y))=k^{-1}\pi(y)$
and hence 
\begin{align*}
\tilde{\pi}((g_1,g_2))(k,y))& =\tilde{\pi}((kg_1^{-1},y g_2^{-1}))=
(kg_1^{-1})^{-1}\pi(y g_2^{-1})\\ & =
g_1k^{-1}\pi(y)g_2^{-1}=g_1(\pi(k,y))g_2^{-1}.
\end{align*}
Thus we have a new representation of $C^*((G\times G^{\rm op})\rtimes C(K))$
on $L^2(\Y^{\rm ad},\nu)$. The fiber at $e$ of this action corresponds to 
$L^2(\Y_e^{\rm ad},\mu_{\Y_e^{\rm ad}})$. Clearly the fiber at $e$, as a $G$ system,
is  the same as the fiber  at $e$ in the above construction,  by taking the profinite limit
after subgroups $\S$.

Thus adjoint  action of $G$ in the fiber at $e$ of $L^2(\Y^{\rm ad},\nu^{\Gamma})$   is then equivalent to the action  $\Pi$ of $G$ on
on $\Gamma\setminus\Y$, that we have constructed in the  Theorem \ref{quotienthecke}.
\

\

The advantage of the adelic formulation is the fact  that we obtain the explicit formula
$(*)$ which is used to compute  the measure  displacement function, by translations representation of $G$ in $(\Gamma\backslash \Y, \nu^{\Gamma\backslash \Y})$.
\end{obs}

We obtain consequently:

\begin{cor}\label{pieces}
With the previous notations, let $F_0 \subseteq \Y$ be a $\Gamma \times \{e\}$
wandering, measurable subset of $\Y$, of finite measure. Let $\tilde {F}$ be the projection of this set in the quotient space $\Gamma\setminus \Y$. Recall that the space $(\Gamma\backslash \Y, \nu_{\Gamma\setminus\Y})$ is acted by $G$ through the transformation $\Pi$ described in the last theorem of the previous appendix.
Then the measure $\nu_{\Gamma/\Y}\Big(\tilde {F}\cap \Pi(g)(\tilde {F})\Big)$ of the displacement of $F_0$ in the quotient  $\Gamma\setminus \Y$,
by elements in $G$, is given by the formula
\begin{gather*}
\lim_{n\to\infty}\sum_{i,j}\nu\bigg(\Big(\Big[(s_i^n)^{-1}\Big(\pi^{-1}(\overline{s_i^n\Gamma_n})\cap F_0\Big)\Big]
\cap\Big[g\Big(s_j^n)^{-1}\pi^{-1}(\overline{s_j^n\Gamma_n})\cap F_0)\Big)g^{-1}\Big]\bigg).
\end{gather*}
\end{cor}
\begin {proof}
This is essentially formula $(*)$.
The fact that for any fixed $g\in G$, we  obtain in the formula simply  conjugation by $g$, instead of the more  complicated expression for the action of $\Pi(g)$, is due to
the fact that the groups $\Gamma_n$ are normal in $\Gamma$ and due to  the fact then when the cosets $s_j^n\Gamma_n$ are very small (for large $n$) , that is, if the groups $\Gamma_n$  we started with (eventually for large $n$) are in the domain of the adjoint action by $g$ on $\Gamma$,
then the expression   for $\Pi (g)$ becomes, by the preceding observation, simply conjugation by $g$. 
\end {proof}

We now return to the proof of Proposition \ref{diag}.

\begin{proof} (Proposition \ref{diag})
 Let $(\Gamma_n)$ be finite index normal subgroups
shrinking to $e$ in $\S$ (the family of subgroups that defines the profinite completion
of $K$). Let $(\Gamma_n)$ be a family of coset representatives for $\Gamma_n$, $n\in \N$. Start with 
$F_0$ a $\Gamma \times 1$ wandering subset of $\Y$. Then take 
$$
F_n=\bigcup_{i}(s_i^n)^{-1}[\pi^{-1}(\overline{s_i^n\Gamma_n})\cap F_0].
$$
Then for $g_1,g_2\in \Gamma_n$ the value of $\varphi((g_1,g_2)\chi_{F_n})$ is zero unless  
$g_1^{-1}g_2$ belongs to $\Gamma_n$. Hence the state $\varphi_{F_N}$ on
$C^*((G\times G^{\rm op})\rtimes C(K))$ defined by the following formula, for $(g_1,g_2)\in G\times G^{\rm op},\ f\in C(K)$,
$$
\varphi_{F_N}((g_1,g_2)f)=\int_{\Y}g_1(f\chi_{F_n})g_2^{-1}\chi_{F_n}d\nu
$$
has support in the $\Gamma_n$-tubular neighborhood of the diagonal.

Then $\varphi_{F_n}$ converges weakly to a state concentrated $\varphi_{\infty}$ on the diagonal.
On the other hand since $\Gamma F_n=\Gamma F_m$ for all $n,m$, by the Observation 61, $\theta(\varphi_{F_N})=\theta(\varphi_M)$
for all $N,M$. Thus the states $\varphi_{F_N}$, induce,
through the map $\theta$ from Definition \ref{defdiag},
the same state on $C^*(G)$.
This state on $C^*(G)$ will thus be the restriction to the diagonal
$G\times G^{\rm op}$ of the state $\Psi_{\tilde{F}}=\theta(\varphi_{F_0})$ constructed in Proposition \ref{diag}.
\end{proof}

\section*{Appendix 6. Analysis of the essential states on $C^*((G\times G^{\rm op})\rtimes C(K))$ coming from the 
embedding into the Calkin algebra $Q(\ell^2(\Gamma))$}

We consider as in the previous section, $\Gamma\subseteq G$ a pair consisting of a 
discrete group and an  almost normal subgroup of 
the countable discrete group $G$. As before, we assume that we have a directed 
family $\S$ of finite index subgroups of $\Gamma$, that also contains
a family, shrinking to the identity, of normal subgroups $\Gamma_n$ of $\Gamma$.
Let $K$ be the profinite completion of $\Gamma$ with respect to $\S$.  By definition,
$C(K)$ is generated by characteristic functions of cosets of elements in $\S$, and hence
acts on $\ell^2(\Gamma)$. 

The left and right action of $G\times G^{\rm op}$ on $\ell^2(\Gamma)$
give the action of $G\times G^{\rm op}$ (the domain of $g_1,g_2$ is $\chi_{\overline{\Gamma_{g_1^{-1},g_2}}}\ell^2(\Gamma))$). Together, the left and right representations determine a representation of 
$C^*((G\times G^{\rm op})\rtimes C(K))$ on $\ell^2(\Gamma)$.
We want to analyze states 
on $C^*((G\times G^{\rm op})\rtimes C(K))$, 
which are obtained by composing the above representation, with the
projection onto
$\QC(\ell^2(\Gamma))= \B(\ell^2(\Gamma))/\K(\ell^2(\Gamma))$.

By Calkin [Ca], it sufficient to consider the essential states on
$B(\ell^2(\Gamma))$ of the form
$$
\omega_{\xi\zeta}=\omega_{(\xi_n),(\zeta_n)}(A)=
\lim_{n\to\infty}\langle A\xi_n,\zeta_n\rangle,\quad A\in B(\ell^2(\Gamma)),
$$
where
$\xi=(\xi_n)$, $\zeta=(\zeta_n)_n$ are sequences in $\ell^2(\Gamma)$,
weakly convergent to zero. Here the limit is after a free ultrafilter. It is sufficient (for continuity purposes),
by linearity, to consider states such that 
$\xi_n$, $n\in \N$ are vectors in $\ell^2(\Gamma)$
with finite support and positive coefficients.

We will prove, by using the Loeb measure construction [Lo], that all such states
are reconductible to states of the form 
$$
G\times G^{\rm op}\ni(g_1,g_2)\to\nu(g_1 Fg_2^{-1}\cap F),
$$
where $\nu$ is an infinite measure on an infinite measure space $\Y$, with an
action of $C(K)$ and an equivariant groupoid action of 
$G\times G^{\rm op}$, invariating the measure.

We may exclude suitable measurable sets from $\Y$, (corresponding to averaging sets of points concentrated in cosets of amenable subgroups) so that this
action becomes free (see [Ra6]).

Assuming $\Gamma$ is exact, it will also follow that we may assume that the
action of $\Gamma \times \{e\}$ (which is by construction continuous on 
$C_{\rm red}^*(\Gamma)$) has either a fundamental domain for $\Gamma$, or either has a fundamental domain for a coamenable quotient of $\Gamma$.
First we prove the representation result for the essential states.

\begin{thm}\label{points}
With the above notations, any state $\omega_{\xi,\xi}$ is a weak limit of states 
of the following, form described bellow.

There exist $(\Y,\nu)$ an infinite probability measure
space, with a surjective projection onto $K$ (thus $C(K)$ acts by
multiplication on $L^2(\Y,\nu)$), a measure preserving,  groupoid
action of $G\times G^{\rm op}$ on $\Y$, that is $G\times G^{\rm op}$ equivariant with respect to
$\pi$ (that is 
$\pi(g_1y g_2^{-1})=g_1\pi(y)g_2^{-1}$ if
$g_1,g_2\in G$, $\pi (y)\in \overline{g_1^{-1}\Gamma g_2\cap\Gamma}$, $y\in\Y$)
and a finite measure subset $F$ of $\Y$. Associated to this data,
we define a state $\varphi_0$ on $C^*((G\times G^{\rm op}) \rtimes C(K))$ as 
follows:

For $(g_1,g_2)\in G\times G^{\rm op}$, $\theta\in C(K)$, let
$$
\varphi_0((g_1,g_2)\theta)=\int_{\Y}g_1\theta\chi_F g_2^{-1}
\chi_F d\nu. \leqno(***)
$$

Then  the state $\omega_{\xi,\xi}|_{C^*((G\times G^{\rm op}) \rtimes C(K))}$
is a weak limit  of convex combinations of states of the form $f^* \varphi_0 f|_{C^*((G\times G^{\rm op}) \rtimes C(K))}$, with $\varphi_0$ as above,
where $f$ is a positive, measurable, square integrable function on $\Y$.

Therefore, the continuity problem for essential state on
$C^*((G\times G^{\rm op})\rtimes C(K))$ is reduced to states of the form $(***)$.

Moreover we may restrict to states $\varphi_0$ as above, so that, in addition, the Koopman representation of $C^*((G\times G^{\rm op})\rtimes C(K))$ into $B(L^2(\Y,\nu))$ is continuous with respect to norm inherited from the norm on the crossed product representation into 
$\mathcal Q(\ell^2(\Gamma))$,  of the
$C^\ast$-algebra $C^*((G\times G^{\rm op})\rtimes 
\ell^{\infty}(\Gamma))$. Here we view $C(K)$ as a subalgebra of $\ell^{\infty}(\Gamma)$. This corresponds to the fact that the states in the convex combinations are ultrafilter limit of states coming from averaging sets.
\end {thm}

\begin{proof}
Let $(\xi_n)_n\subseteq \ell^2(\Gamma)$ be a sequence weakly convergent to zero,
$\omega$ a free ultrafilter and $\omega_{\xi,\xi}$ the corresponding essential
states. 
We may assume that
$$
\xi^n=\sum_{a\in A_n}\lambda_n(a)a,
$$
where $A_n$ are finite subsets of $\Gamma$, and
$\lambda_n(a)$, $a\in A_n\geq 0$, are positive weights. Then
$\omega_{\xi,\xi}$ gives a Loeb measure $\mu_\lambda$ on
$\cC_\omega((A_n)_n)$. Here $\cC_\omega((A_n)_n)$ is the ultra-product of
the sets $A_n$. Note that on $\cC_\omega((A_n)_n)$ we also have the canonical
Loeb counting measure, that we will denote by $\mu_\omega=\mu_{\omega,(A_n)_n}$ (see also [Ra6]).

Because of $\aleph_{1}$-saturation ([Lo], [Cut]) and since we are
interested only in weak approximation, we may assume that the
support of $\mu_\lambda$ is $\cC_\omega((A_n)_n)$
(eventually by replacing the set $\cC_\omega((A_n)_n)$ with a subset of
the same type (Lemma 1.19, [Cut]).

For $M>0$, let 
$A_n^M=\left\{a\in A_n\mid\lambda_n(a)\leq\frac{M}{\card A_n}\right\}.$
For every $(\alpha_n)_n$,  positive sequence of numbers increasing to $\infty$
we let $$A_n^{\alpha}=\left\{a\in A_n\mid\lambda_n(a)>\frac{\alpha_n}{\card A_n}\right\}.$$

Then $\cC_\omega((A_n)_n)$ is the reunion of 
$\bigcup\limits_{M > 0}\cC_\omega((A_n^M)_n)$ and $\bigcup\limits_{\alpha}\cC_\omega((A_n^\alpha)_n)$,
where the second, directed, reunion runs over all positive increasing sequences 
$(\alpha_n)$. By 
$\aleph_{1}$ saturation,
it will be sufficient then to assume the case when the support of $\mu_\lambda$ 
is of the form
 $\cC_\omega((A_n^M)_n)$ for a sufficiently large $M$,
union with $\cC_\omega((A_n^\alpha)_n)$ for a sufficiently slow decreasing sequence $\alpha$.
(In fact, $\mu_\lambda$ is here decomposed  into a measure absolutely
continuous with the Loeb measure $\mu_\omega^M$ on $\cC_\omega(A^M)$, (where $A^M = (A_n^M)_n)$,
and another measure $\mu_\omega^\alpha$ supported on $\cC_\omega(A^\alpha)$
(where $A^\alpha = (A_n^\alpha)_n$). Note that $\mu_\omega^M$ and 
$\mu_\omega^\alpha$ are singular. Also, the total mass of $\mu_\omega^M$ is non zero, by our initial
assumption, that the support of $\mu_\lambda$ is $\cC_w((A_n)_n)$.)

We repeat this procedure by transfinite induction for $\cC_w(A^\alpha)$.

Because the mass of the measures is always non zero, this procedure will stop
after a countable number of iterations.

In this way we end up by writing the initial measure $\mu_\lambda$ in the form
$$
\mu_{\lambda}=\sum_{k=1}^{\infty}(f_k^n)d\mu_{\omega, (A_n^k)_n}
$$
where $f_k=(f_k^n),$ are measurable functions positive on $\cC_\omega(A)$ (and we also may assume
$f_k$ are step functions with finite values), 
and $\mu_{\omega, (A_n^k)_n}$ is singular with
respect $\sum\limits_{s>k}\mu_{\omega,(A_n^s)_n}$.

We take the measure $\mu_{\lambda}^0=\sum\frac{1}{2^k}\mu_{\omega, (A_n^k)_n}$
and by renormalizing the functions $f_n$ into $\tilde{f}_n=2^nf_n$, we get
$$
\mu_\lambda=\tilde{F}d\mu_{\lambda}^0.
$$
The measure $\mu_{\lambda}^0$ is extended to a $\Gamma$-invariant measure
on the countable union $\bigcup \gamma \cC_\omega(A)$. This is because the pieces  of $\mu_\lambda^0$, which are weighted
copies of $\mu_{\omega, (A_n^k)_n}$ are reciprocally singular (with the translates of 
$\mu_{\omega, (A_n^l)_n}$, $l > k$), being multiples of counting
measures (so that the computations for $\mu_\lambda^0 (\cC_\omega(A_n) \cap g \cC_\omega(A_n)$)
involve only the diagonal pieces 
$\mu_{\omega,(A_n^k)}\big(\cC_\omega(A_n^k)\cap g\cC_\omega(A_n^k)\big), g \in G).$
The required functions $f$ from the statement are then the square root of $F$. Since we will prove temperedness (continuity with respect to the $C^\ast$ reduced crossed product norm) for all the representations involving these states, it will be sufficient for proving continuity to consider only states of the type (***).
\end{proof}

The analysis of the essential states on $C^*((G\times G^{\rm op})\rtimes C(K))$
could be further reduced, by noting that we may only consider states with the property that the measure $\mu_\lambda^0$, from the proof of the preceding proposition,
is concentrated at the fiber at $e$ (the unit element of $K$).

\begin{prop}\label{shrink}
Let $\varphi_0$ be a state of the form $\sum \frac{1}{2^n}\mu_{\omega,(A_n^k)_n}$ 
as in the preceding theorem and let $(\Y,\nu)$ be the corresponding, associated
measure space, with the $(G\times G^{\rm op})\rtimes C(K)$ action and $F$ the finite
measure subset of $\Y$ whose displacements by $G$ 
compute $\varphi_0$.  With the above notations let $\Psi_0$ be the state on
$C^*(G)$, associated to $\varphi_0$, constructed in Proposition \ref{diag} (which can be
extended to $C^*((G\times G^{\rm op})\rtimes C(K))$.

 Then $\Psi_0$ is a state of the type considered
in the previous proposition (constructed as an ultrafilter limit of states associated to averaging sets of points), with the additional property that there exists a decreasing family of normal, finite subgroups, $\Gamma_n$ of $\Gamma$, with trivial intersection, such that  the finite sets $(A^k_n)_n)$  from  the construction of the measure space $\Y$ in the previous proposition, have the property that
$A^k_n \subseteq \Gamma_n $,  for all $k,n$.
\end{prop}

\begin{proof}
Fix a family ($\Gamma_n)_n$ of finite index, normal subgroups in $\S$, shrinking to $e$.
We fix an exhausting family $(G_n)_n$ in $G$,
with finite sets, and let $C(K_n)$ be the finite subalgebra of $C(K)$ generated by
characteristic functions of the closure of cosets of the group $\Gamma_n$. We know
from Proposition \ref{diag} that the state $\Psi$ is realized as a weak limit of $\varphi_{F_n}$, measuring
the displacement of $F_n$, (with the notations from the Proposition \ref{diag} ) under the action of $(G\times G^{\rm op})\rtimes K$, where
$$
F_n=\bigcup\limits_i(s_i^n)^{-1}(\pi^{-1}(\overline{s_i^n\Gamma_n})\cap F_0),
$$
where $s_i^n$ are coset representatives for $\Gamma_n$.

Then we replace the sets $A_n^k$ by sets $\widetilde{A}_n^k$,
that are obtained as follows
$$
\tilde{A}_n^k=\bigcup\limits_i(s_i^n)^{-1}\big(\pi^{-1}(s_i^n\Gamma_n)\cap A_{p_n}^k\big),
$$
where $p_n$ are chosen so large that the characteristic functions
of $(A_{p_n}^k)$, $k=1,2,\ldots,n$ and of $g A_{p_n}^k g^{-1}$, $g\in G_n$,
verify up to $\frac{\varepsilon}{2^n}$ the same measure of intersections relations
as the corresponding measure of intersections relations of the characteristic functions
of $\chi_F$, $\chi_{gFg^{-1}}$, $g\in G_n$ in relation to the characteristic functions
in $C(K_n)$. Then $(\tilde{A}_n^k)_n$ have support in $\Gamma_n$ eventually, and
letting $\varepsilon \searrow 0$, and using Corollary \ref {pieces}, we get that the state on $C^*((G\times G^{\rm op})\rtimes C(K))$  corresponding to the new family of sets $(\tilde{A}_n^k)_n$, $k\in N$, is the state $\Psi_0$.
\end{proof}

Note that one could give an alternative proof by arguing that $\Psi_0$ is still an
essential state.

\begin{rem}
As in [Ra6], we may assume that  the action is free on $\Y$. To do this we subtract the sets 
corresponding to fixed points of the action of $G \times G$
on $\Y$. The fact that the fixed point sets 
are permuted by the action of 
$G \times G$, implies that the state $\Psi_0$ obtained in this case is represented
as the state associated to the ultra limit of averaging sets $\cC_\omega(\tilde{A}_n^k)_n)$  shrinking to $e$, minus the reunion
of the ultra limit of averaging sets of the same type.
In either case, it follows that we can represent
the state $\Psi_0$ by the displacement ultra limit measure   of averaging  finite sets $(A_n^k)_n$ whose support is shrinking
to $e$, and such that the action of $G\times G^{\rm op}$ on $\Y$ is free.
\end{rem}

\

\

To apply the machinery that we developed in the preceding appendix   for the quotient $\Gamma\backslash\Y$, we prove the following
result which establishes the existence of a fundamental domain for the action of $\Gamma$ (or for  a coamenable quotient). This will be applied to the measure preserving action of the group
$\Gamma\times \{e\}$, (the left action) on the  measure  spaces, that we are using to    represent essential
states on $C^*((G\times G^{\rm op})\rtimes C(K))$.

\begin{lemma}\label{fd}

Let $\Gamma$ be a countable discrete group that is non-amenable, with infinite conjugacy classes and exact.
Let $\omega $ be a free ultrafilter on $\N$.
Let  $A=(A_n)_n$ be a family of finite subsets of $\Gamma$, that avoids eventually (with respect to the ultrafilter
$\omega$) any finite, initial subset of $\Gamma$. 
Let, as above,  $(\cC_{\omega}\Big((A_n)_n\Big), \mu_{\omega , (A_n)_n} )$ be the Loeb probability measure space associated to this data, where  $\mu_{\omega , (A_n)_n} $ the ultrafilter limit of the counting measure.

Let ($\Y_\omega,  \nu_{\omega , (A_n)_n})$  be the infinite measure space constructed as follows.
Let
$$\Y_\omega=\bigcup_{\gamma \in \Gamma} \gamma \cC_{\omega}((A_n)_n)= \bigcup_{\gamma \in \Gamma}   \cC_{\omega}((\gamma A_n)_n).$$
Since the restriction of the corresponding
Loeb ultrafilter limits of counting measures $\mu_{\omega , (A_n)_n}$ coincide on overlaps,  the measures $\mu_{\omega, (\gamma A_n)_n}$, $\gamma \in \Gamma$, define a  $\Gamma$ - invariant measure $\nu_\omega$ on $\Y_{\omega}$.
Note that the absence of Folner sets for the group $\Gamma$, implies that $\nu_{\omega}(\Y_{\omega}) = \infty$. The same arguments will apply for a countable reunion of such spaces, if the corresponding Loeb measures are mutually singular.

Consequently $\nu_{\omega}$ defines an infinite measure on $\partial^{\beta}(\Gamma)=\beta\Gamma\setminus c_0(\Gamma)$, where $\beta(\Gamma)$ is the Stone-Cech compactification of $\Gamma$. Moreover $\nu_{\omega}$ defines a semifinite trace of the algebra $C^{\ast}(\Gamma\rtimes L^{\infty}(\Y_\omega, \nu_\omega))$. This $C^\ast$ algebra is then   a crossed product $C^{\ast}$ representation of the Roe-algebra $C^{\ast}(\Gamma\rtimes l^{\infty}(\Gamma)) \subseteq \B(l^2(\Gamma))$. Because of the exactness assumption,  we have that the maximal crossed product C$^\ast$-algebra  $C^{\ast}(\Gamma\rtimes L^{\infty}(\Y_\omega, \nu_\omega))$ coincides with the C$^\ast$-algebra  $C_{\rm red}^{\ast}(\Gamma\rtimes L^{\infty}(\Y_\omega, \nu_\omega))$.

We assume in addition that $\Gamma$ admits only a countable subset $\S\A$ of infinite amenable subgroups, and that the distinct cosets for all the subgroups in this family have finite intersections.

Then there exists a disjoint splitting of $\Y_{\omega}$ into $\Gamma$ - invariant, measurable subsets of infinite measure (or zero measure) $\Y_I$ and $\Y_{II}$, and furthermore we have the disjoint splitting into $\Gamma$ - invariant, measurable subsets, $\Y_{II} = \mathop{\bigcup}\limits_{\Gamma_0 \in \S\A} \Y_{\Gamma_0}$, such that the following happens

1) The action of $\Gamma$ on $\Y_I$ has a finite measure fundamental domain in $\Y$,

2) For each $\Gamma_0$ there exists a subset $F_{\Gamma_0}$ of finite measure in $\Y_{\Gamma_0}$, such that $F_{\Gamma_0}$ is invariated by $\Gamma_0$ and the $\Gamma$ - system $\Y_{\Gamma_0}$ is isomorphic to $F_{\Gamma_0} \times \Gamma/\Gamma_0$ (where $\Gamma/\Gamma_0$ has the counting measure).

The second situation corresponds, after doing a rearrangement of the sets $(A_n)_n$, by $\Gamma$-transformations, to the case $$\Gamma_0\cC_{\omega}((A_n))=
\cC_{\omega}((A_n))=\cC_{\omega}((B_nx_n)),$$ where $B_n$ are Folner sets in $\Gamma_0$, and $x_n$ are elements in $\Gamma$, $n$ in a cofinal subset of $\omega$. Note that by  doing a rearrangement by $\Gamma $ transformations, doesn't change the topology on the crossed product $C^\ast$ algebra (see [Ng]). 
\end{lemma}

\begin{proof}
The weight $\nu=\nu_{\omega}$ is semifinite, and   $\Gamma$ acts by measure preserving transformations on $\Y=\Y_\omega$, which is a subspace of the spectrum  $\beta(\Gamma)\setminus c_0(\Gamma)$ of $l^{\infty}(\Gamma)$. It follows  the algebra $C^{\ast}_{\rm red}(\Gamma \rtimes L^{\infty}(\Y, \nu))$ is a representation of the Roe algebra $C^{\ast}(\Gamma\rtimes l^{\infty}(\Gamma)) \subseteq \B(l^2(\Gamma))$, which by exactness is nuclear.

We have a canonical semifinite trace on this algebra, which is the composition of the canonical, normal conditional expectation $E$ onto $L^{\infty}(\Y, \nu)$ with the measure (weight) on $L^{\infty}(\Y, \nu)$ given by $\nu$. We consider the Koopman unitary representation of the reduced $C^{\ast}$ - algebra $C^{\ast}_{\rm red}(\Gamma \rtimes L^{\infty}(\Y, \nu))$ on the Hilbert space $\H_{\nu}=L^2(\Y, \nu)$ associated to the semifinite trace $\nu$ (the representation is isometric because of nuclearity).

Let $M$ be the corresponding von Neumann algebra, which is necessary of semifinite type. Let $D = L^{\infty}(\Y, \nu)$ be the corresponding MASA in $M$, and let $E$ be the normal conditional expectation from $M$ onto $E$. Because of the infinite conjugacy classes condition on the group $\Gamma$, the center $\mathcal{Z}(M)$ is is contained in $D = L^{\infty}(\Y, \nu)$.

 We identify the algebra $\mathcal{Z}(M)$  with  the algebra $L^\infty(\mathcal Z, \nu_0)$, for some measure space $\mathcal Z$, for a canonical measure $\nu_0$.  In fact $L^\infty(\mathcal Z, \nu_0)$ is the $\Gamma$-invariant part of $L^{\infty}(\Y, \nu)$. The measure $\nu_0$ is defined simply by letting $\nu_0(\tilde F)=\nu (F)$, if $F$ is measurable subset of $\Y$, of finite measure and the characteristic function  $\chi_{\tilde F}$ is the central support  in $M$ of the projection $\chi_F$.

We denote by $\nu$, the semifinite, faithful weight on $M$ induced by $\nu_{\omega}$. Note that $M$ can only have type $I_{\infty}$ or hyperfinite type $II_{\infty}$ components. (the infiniteness is a consequence of the absence of Folner sets). Indeed, by the nuclearity of the algebra $C^{\ast}_{\rm red}(\Gamma \rtimes L^{\infty}(\Y, \nu))$, the type $II$ components are hyperfinite ([Co]).

We disintegrate $M$ over the center $\mathcal{Z}(M)$ and obtain fibers $M_z \supseteq D_z$, $z\in \mathcal Z$, with normal faithful conditional expectation $E_z : M_z \to D_z$ and $\nu_{z}$ a semifinite trace on $D_z$, giving a semifinite faithful trace on $M_z$, for $z \in \mathcal{Z}$, almost everywhere.

In the case of type $I$, which corresponds to $\Y_I$, because of the existence of a normal conditional expectation onto the algebra $D_z$, it follows that the algebras $D_z$ are maximal abelian, diagonal algebras. Hence any field of minimal projections is the characteristic function a fundamental domain for the action of $\Gamma$ (e.g. by Vitali's criteria [Za]).

In the case of type $II$, which corresponds to the $\Y_{II}$ part in the statement, the fact that there exists a conditional expectation from $M_z$ onto $D_z$,  and since $M_z$ is of type $II_{\infty}$ it follows that $M_z$ admits a splitting $N_z \otimes \B(l^2(I_z))$, where $N_z$ is a type $II$, (hyperfinite) factor, and $l^2(I_z)$ is the Hilbert space associated to a countable set $I_z$.

Moreover, since $D_z$ is maximal abelian and generated by finite projections, it follows that $D_z$ splits as $D_z^1 \otimes D_z^2$, in such a way that $D_z^1$ is a MASA in $N_z^1$ and $D_z^2$ is the maximal abelian diagonal algebra of $\B(l^2(I_z))$ associated to the basis indexed by $I_z$.

Let $\pi_z$ be the disintegration of the left regular representation of the group $\Gamma$ in $\H_{\nu}$.
Thus $\pi_z(\Gamma)'' = M_z$ and $\pi_z(\gamma)$ normalizes the algebra $D_z$ for every $\gamma$. Then necessary there exists a permutation $P_z(\gamma)$ of $I_z$, $P_z(\gamma) : I_z \to I_z$ such that if $(e_{i, j}^z)$ is the matrix unit of $\B(l^2(I_z))$ associated to the basis indexed by $I_z$, then there exists unitaries $u_i^z(\gamma)$, $i \in I_z$ in the normalizer $\mathcal N_{N_z}(D_z^1)$ such that $\pi_z(\gamma) = \mathop{\sum}\limits_{i \in I_z} u_i^z(\gamma)\otimes e_{i, P_z(\gamma)(i)}$ for all $\gamma \in \Gamma$.
 But then necessary the map $\gamma \to P_z(\gamma)$ into the permutation group of $I_z$ is a homeomorphism and hence there exists a subgroup $(\Gamma_0^z)$ of $\Gamma$ such that the index set $I_z$ is identified with the set of cosets $[s\Gamma_0^z]$ in $\Gamma / \Gamma_0$, $s\in \Gamma$. The identification is $\Gamma$ - invariant. Note that $\Gamma_0^z$ is necessary infinite, since otherwise we are back in the  case of type $I_\infty$. Moreover $P_z(\gamma)$, in this identification, is translation by $\Gamma$ on $\Gamma/ \Gamma_0^z$. Let $e_0^z$ in $\B(l^2(\Gamma/ \Gamma_0^z))$ be the projection corresponding to $e_{[\Gamma_0^z], [\Gamma_0^z]}$.

Then $e_0^z$ is fixed by $\pi^z(\gamma)$, $\gamma \in \Gamma_0^z$, and hence after identifying $N_z$ with $N_z \otimes e_0^z$, we have a representation $\pi_0^z(\gamma)$, $\gamma \in \Gamma_0^z$ of $\Gamma_0^z$ in $N_z$, such that the original representation is now the induced representation ${\rm Ind}_{\Gamma_0^z}^{\Gamma}(\pi_0^z)$ on $L^2(N_z, \nu_0^z) \otimes l^2(\Gamma/ \Gamma_0^z)$. ($\nu_0^z$ is the canonical trace on $M_z$).

Because in the original representation $E_z(\pi^z(\gamma)) = 0$, it follows that, if we denote by $\nu_0^z = \nu^z(e_0^z\cdot)$ the trace induced by $\nu$ on $N^z$, then $\nu_0^z(\pi_0^z(\gamma)) =0$ for all $\gamma \in \Gamma_0$. Moreover $\pi_0^z(\Gamma_0)'' = N_z$ and hence $N_z$ is isomorphic the type $II$, factor associated to the group $\Gamma_0^z$.

Since $N_z$ is hyperfinite, it follows that $\Gamma_0^z$ is amenable and infinite. Since $e_0^z$ is the projection in $D_z$ corresponding to $1\otimes e_{[\Gamma_0], [\Gamma_0]}$ it also follows  that the $G$ system $\Y_z$ is isomorphic to a $G$ - system $F^z \times \Gamma / \Gamma_0^z$, where $F^z$ is a probability measure space, that is $\Gamma_0^z$ invariant. Since we have a countable set of infinite amenable subgroups, the property (2) holds true.

Moreover considering, any of the $\Gamma$ - invariant components $\Y_{\Gamma_0}$ of $\Y$, for some $\Gamma_0$ in $\S\A$, it follows that the original set $F = \cC_{\omega}((A_n)_n)$, out of which the space $\Y_{\omega}$ was constructed, by $\Gamma$ - translations, is decomposed into pieces corresponding to cosets of $\Gamma/ \Gamma_0$.

We may divide the sets $(A_n)_n$ by working with large $n$, so that recomposing the corresponding pieces, and bringing back by $\Gamma$ translations, to subsets of  the form $\cC_{\omega}((A'_n)_n)$ is contained in  $F_{\Gamma_0}$,by  which we denote the set corresponding to the projections $e^z_{[\Gamma_0], [\Gamma_0]}$, $z \in \Y_{\Gamma_0}$. Then $\cC_{\omega}((A'_n)_n)$ is $\Gamma$ - equivalent to $\cC_{\omega}((A_n)_n)$ (in the sense of [Ng]) Note that alternatively, we may argue that by continuity and linearity, we may reduce the proof to states such that $\cC_{\omega}((A_n)_n)$ has already  this property.

Since $\cC_{\omega}((A'_n)_n)$ is contained $F_{\Gamma_0}$, which is $\Gamma_0$ invariant, it follows that there exists Folner sets $B_n$ in $\Gamma_0$ and $x_n$ in $\Gamma$, such that $A'_n\subseteq B_nx_n$, for $n$ in a cofinal set of the ultrafilter $\omega$.

\end{proof}

\begin{lemma}\label{coamenable}
With the assumption from the previous lemma, assume in addition that we have a larger discrete group $G$, such that $\Gamma$ is almost normal in $G$. Assume that $G$ is exact. In the setting of the previous lemma, consider the larger measure space $(\Y_{\omega}, \nu_{\omega, A})$, defined by:
$$
\Y_{\omega} = \mathop{\bigcup}\limits_{(g_1, g_2) \in G\times G^{\rm op}} \cC_{\omega}(g_1(A_n\cap \Gamma_{g_1^{-1}g_2})g_2^{-1}).
$$

Then the measure $\nu_{\omega, A}$ is invariant to the partial $G \times G^{\rm op}$ transformations on $\Y_{\omega, A}$.

We assume the following additional property on the group $G$: for every amenable subgroup $\Gamma_0$ of $\Gamma$, the normalizer $\mathcal N_{G}$ of $\Gamma_0$ in $G$ is amenable and $\mathcal N_{G}(\Gamma_0)x \cap (g\mathcal N_G(\Gamma_0)g^{-1})y$ is finite, for all $\Gamma_0$ in $\S\A$, $g$ in $G$ and $x, y$ in $\Gamma$. Fix $\Gamma_0$ in $\S\A$.

Then the state on $C^{\ast}_{\rm red}((G \rtimes G^{\rm op}) \rtimes C(K))$ corresponding to a family $(A_n)$ of the form $(B_nx_n)_n$, where $(B_n)$ is a family of Folner sets in $\Gamma_0$, is continuous with respect to the $C^{\ast}_{\rm red}((G \rtimes G^{\rm op}) \rtimes C(K))$ topology.
\end{lemma}

\begin{proof}

Indeed in this case the state on $C^{\ast}((G \times G) \rtimes C(K))$ corresponding to $\cC_{\omega}((B_n)_n)$, will have support on $(\mathcal N_G(\Gamma_0))\times G$. Since $\mathcal N_G(\Gamma_0)$ is amenable and $G$ is exact, the result follows.

\end{proof}

The following remark explains the mechanics of the previous argument in Lemma \ref{fd}.

\begin{rem}
Let $H$ be a discrete, exact group acting ergodical, and measure preserving
on the (infinite) measure space $(X,\mu)$.
Assume that the action of $H$ has a fundamental domain $F$. 
Let $\A=W^*_{\rm koop}(H,L^\infty(X))$ be the crossed product algebra
(representing $W^*(H\rtimes L^\infty(X))$) in the space
$L^2(X,\mu)$. This is the Koopman representation ([Ke]) of $H$ on $L^2(X,\mu)$.
Then the center of $\A$, $Z(\A)$ is canonically identified to $L^{\infty}(X)^{H}$ (the $H$-invariant functions in $L^\infty(X)$) and $\A$ is isomorphic to
$W^*(H\rtimes\ell^\infty(H))\otimes L^\infty(X)^H$
acting on $\ell^2(\Gamma)\otimes L^2(F,\nu)$.
Here $X^H$, the spectrum of $L^\infty(X)^H$ is measurably identified to $F$, and the algebra $W^*(H\rtimes\ell^\infty(H))$
is the Roe crossed product ([Br Oz]).
\end{rem}

\begin{proof}
This simply result by the identification of $L^2(X,\mu)$ with $\ell^2(\Gamma)\otimes L^2(F)$.
\end{proof}

We can conclude the study of the essential states on $C^*((G\times G^{\rm op})\rtimes C(K))$,
induced by  the  representation into the Calkin algebra.
More precisely, we have following corollary, which is used in the proof of Theorem \ref{AOlocal}, to reduce to the case of essential states on  $C^*((G\times G^{\rm op})\rtimes C(K))$ to the case of essential states that vanish outside the diagonal
$\{(g,g^{-1})| g \in G\}$ of the group $G\times G^{\rm op}$.

\begin{cor}\label{reduction} Let $\Gamma\subseteq G$, $\S$, $K$, as above.
Then  the continuity property, with respect the reduced C* norm on $C^*((G\times G^{\rm op})\rtimes C(K))$  of the states coming from the Calkin algebra representation on the C*-algebra
$C^*((G\times G^{\rm op})\rtimes C(K))$, is determined by the analysis of states on
$C^*(G)$, which are of the form
$\varphi(g) = \nu(g F\cap F)$, where $(\Y,\nu)$ is an infinite
measure space of the type described described bellow, and $G$ is acting by measure preserving transformations and freely on the space
$\Y$. Here $F$ is a set of finite measure in $\Y$.

The measured space $(\Y, \nu)$ is constructed as follows: Let $\omega$ be a free ultrafilter on $\N$.
 The initial data is  a family of normal subgroups $\Gamma_n$ in $\S$, with trivial intersection, and $(A_n^k)_{n\in \N} $ is a family (indexed by $k\in \N$) of disjoint (for every fixed $n\in \N $) and finite, subsets of $\Gamma_n$,  for $n,k$, that eventually avoid (in the ultrafilter $\omega$, after $n\in \N$) any given, finite subset of $\Gamma$. For $k \in \N$,  let ($\cC_{\omega,(A_n^k)_n}, $ $\mu_{\omega, (A_n^k)_n})$ be the associated Loeb probability measure space.  We may assume  (by Proposition \ref{shrink})
 that for every $k\in\N$, the Loeb counting measure
$\mu_{\omega,(A_n^k)_n}$ is singular to $\sum\limits_{s>k}\mu_{\omega, (A_n^s)}$. We let $\Y^k$  be the reunion of the  by the adjoint action of $(G\times G^{\rm op})$ on the probability measure spaces $ \cC_{\omega,(A_n^k)_n}$. We obtain a well defined family of measured spaces $(\Y^k, \nu^k)$, the measure $\nu_k$ being obtained by patching together the Loeb measures of the terms in the above reunion. Because of coincidence on the overlaps, the measures $\nu_k$ are $G$-invariant. Then $\Y$ is the direct sum of the spaces $\Y^k$  with G-invariant measure
$\nu=\sum\limits_{k\geq 1}\frac{1}{2^k}\nu_k$.

To obtain a free action, we subtract (as in [Ra6])
the Loeb spaces associated 
  to infinite sets of fixed
points in $\Gamma$, corresponding to amenable subgroups.

Then, if all the states on $C^*(G)$, obtained through this method
are continuous with respect 
$C_{\rm red}^*(G)$, and verify the additional assumptions (FS1), (FS2) of Theorem \ref{quotienthecke}, then the essential states on 
$C^*((G\times G^{\rm op})\rtimes C(K))\subseteq \QC(\ell^2(\Gamma))$
are continuous with respect to
$C_{\rm red}^*((G\times G^{\rm op})\rtimes C(K))$.
\end{cor}

\begin{proof} Indeed by Theorem \ref{points}, for  the analysis of the continuity properties of essential states  on the C* algebra  $C^*((G\times G^{\rm op})\rtimes C(K))\subseteq \QC(\ell^2(\Gamma))$ it is sufficient to consider the states $\phi_{F_0}$ (as in Definition \ref{defdiag})
 measuring the displacement by $(G\times G^{\rm op})\rtimes K$   of a finite measure subset $F_0$ in an infinite invariant, measure space $(\Y,\nu)$ constructed as in Theorem \ref{points} and 
acted by $(G\times G^{\rm op})\rtimes K$. Because of Lemma \ref{fd} we know that the restriction of the action of $G\times G^{\rm op}$ to  $\Gamma \times \{e\}$ admits a fundamental domain in $\Y$ (the case of type $II_\infty$,  in the Lemma \ref{fd} was analyzed in the Lemma \ref{coamenable} directly). Consider the associated  state  $\varPsi=\theta (\phi_{F_0})$ introduced in Definition \ref{defdiag}
 on $C^*(G)$ and computed in Proposition \ref{diag}. By  Proposition \ref{shrink} and Corollary \ref{pieces}, the state $\varPsi$ on  $C^*(G)$ is of the same form as in the statement of the corollary  that we are now proving. 
 
 The conclusion of the corollary now follows,
 because  of Theorem \ref{quotienthecke}, which asserts that the initial state $\phi_{F_0}$  on 
$C^*((G\times G^{\rm op})\rtimes C(K))$ is continuous with respect to the norm on $C_{\rm red}^*((G\times G^{\rm op})\rtimes C(K))$,  if and only if the state $\varPsi=\theta (\phi_{F_0})$ on $C^*(G)$ is continuous on $C^*(G)$    and verifies the additional conditions (FS1), (FS2) in Theorem \ref{quotienthecke}. These additional conditions will be proved   hold true in the proof of Theorem \ref{AOlocal}.

\end{proof}

\section*{Appendix 7. Examples}

In the following we present a few examples of the  construction in Appendix 5, 
of a $C^*((G\times G^{\rm op})\rtimes C(K))$ action on a Hilbert
space $V$, and we determine the corresponding action $\Pi$ of $G$.
First, we consider the reduced $C^*$-algebra case.

\begin{ex}
Let $V=L^2((G\times G^{\rm op})\rtimes K)$ be the Hilbert
space of the reduced groupoid crossed product $(G\times G^{\rm op})\rtimes K$. 
Then $V$ may be identified with $L^2(\cX,\mu)$,
where $\cX = \bigcup\{ \underline{g_1} k \underline{g_2} \mid g_1,g_2 \in G,\ k \in K_{g_1,g_2}\}$.

Here $\cX$, as topological space, is a direct sum of copies of pieces of $K$,
which are labeled by $g_1,g_2$, and denoted in the sequel by underlined $\underline{g_1},\underline{g_2} \in G$.
The measure is the one induced from the Haar measure of $K$.

The action of $C(K)$ on $V$ which  identified with $L^2(\cX,\mu)$ is described, by giving
an explicit formula for the projection $\pi : \cX \to K$, which is simply
$\pi( \underline{g_1} \, k \, \underline{g_2^{-1}}) =  g_1 \, k \, g_2^{-1}$, $g_1,g_2 \in G$, $k \in K_{g_1,g_2^{-1}}$.

This action is compatible with the partial action of $G \times G$ on $\cX$, because 
$(g_3, g_4)(\underline{g_1} \,k \,\underline{g_2^{-1}}) = \underline{g_3g_1}\, k\, \underline{g_2^{-1}g_4^{-1}}$.
 
Here one requires that $g_1 \, k \, g_2^{-1}$ belongs to
$K_{g_3 ,g_4^{-1}}$ and that $k$ should belong to $K_{g_1,g_2^{-1}}$.

To describe the representation $\Pi$ we have to describe $L^2(\Gamma\backslash\cX)$. Clearly,
the points of this space are of the form
$\Gamma \underline{g_1} \, k \, \underline{g_2^{-1}}$ for all $k$ in $K_{g_1^{-1},g_2}$ and
the measure is induced from Haar measure on $K$ (by ignoring the symbol $\Gamma$).

To describe the formula for $\sigma \in G$, for the action
$\Pi(\sigma)$ on $y\in \Gamma\backslash \cX$ one has to consider a well chosen representative for $y$.
We choose $\theta \in \Gamma$ such that 
$y' = \theta y$ has the property that $\pi(y')$ belongs to $\Gamma_{\sigma^{-1}}$.

Then by using the definition of the  the action of $G$ in the Theorem \ref{quotienthecke}, we have that
$\Pi(\sigma)(\Gamma y) = \Pi(\sigma)\Gamma y'=\Gamma \sigma y' \sigma^{-1}$.
It is obvious from the above formula that $\pi$ is a representation.

The fact that $\Pi$ is equivalent to the $C^*_{\rm red}$ representation
of $G$ can be seen as follows:
$L^2(\Y,\nu)$ is 
$
\bigoplus\limits_{\Gamma\theta_1\in\Gamma\backslash G,\, \theta_2\in G}\underline{\Gamma\theta_1}L^2(K_{\theta_1^{-1},\theta_2})\underline{\theta_2}
$.

When applying $\Pi(\sigma)$, for every $y = \underline{g_1} k \underline{g_2}$, there exist a selection
of $\theta$ in $\Gamma$ such that the result is
$$
\Pi(\sigma)(\Gamma \underline{g_1}k \underline{g_2^{-1}})=\Gamma\underline{\theta g_1} k\underline{g_2^{-1}}
$$
we may describe 
$$
\Pi(\sigma)(\Gamma \underline{g_1} k \underline{g_2})=[\Gamma\sigma\Gamma] \underline{g_1} k\underline{g_2^{-1}\sigma^{-1}}
$$
where $\Gamma\sigma\Gamma$ is a sum of cosets $\sum\limits_i\Gamma\underline{\sigma s_ig_1}k\underline{g_2^{-1}\sigma^{-1}}$
and automatically only one index $i$ in this sum, gives a non zero term.

Because of the label on the right hand side, this action has $K$ as fundamental domain,
and the action is $C^*_{\rm red}(G)$.

Clearly, the Hecke operators are described pointwise as mapping
$\Gamma \underline{g_1}k g_2^{-1}\Gamma$ $g_1,g_2\in G$ into 
$[\Gamma\sigma\Gamma][\Gamma g_1]k[g_2^{-1}\Gamma][\Gamma\sigma\Gamma]$
and taking into account $k$ makes that this sum is performed on a suitable selection
of a permutation $\pi_\sigma$ of the indices
$$
\sum\underline{\sigma\sigma s_i g_1}k\underline{g_2s_{\pi_{\sigma(i)}}^{-1}}\sigma^{-1}\Gamma.
$$
\end{ex}
\vskip6pt

\

\

A second example will be obtained by tensoring a given action of $C^*((G\times G^{\rm op})\rtimes C(K))$
with a representation in which the action of $C(K)$ is trivial. We will prove below that
the Hecke operators we constructed in Section 5 in Theorem \ref{homeo} and Theorem \ref{AOlocal}, are of this form. In this
way, this example gives another direct proof of the algebraic relations implying that the
operators $\Psi([\Gamma\sigma\Gamma])$ that  we have constructed in Section 3 are a representation of the
Hecke algebra.

Also, in this way the continuity of the action of the Hecke algebra (relative to the $C^*_{\rm red}$
topology) is reduced to the analysis of the continuity properties of the associated unitary
representation of $G$.

\begin{ex}\label{Heckeas2}
Let $\pi$ be a representation of $C^*((G\times G^{\rm op})\rtimes C(K))$ on $V$.
Assume that $\pi_0$ is a unitary representation of the discrete group $G$ on the
Hilbert space $H_0$. We assume that $\pi_0|_\Gamma$ is unitarily equivalent to the left
regular representation of $\Gamma$.

We consider the unitary representation of $C^*((G\times G^{\rm op})\rtimes C(K))$
on $\H = H_0 \otimes \overline{H_0} \otimes V$, in which the representation of
$G\times G^{\rm op}$ is mapping
$(g_1\times g_2)$ into $\pi_0(g_1)\otimes\overline{\pi_0}(g_2)\otimes\pi((g_1,g_2))$, $g_1,g_2\in G$,
and by letting $C(K)$ act as $1 \otimes 1 \otimes \pi$.

We fix $1$, a cyclic trace vector of $\Gamma$ in $H_0$. Then 
$1\otimes \overline{H_0} \otimes V$ is a generating, wandering subspace for the action of $\Gamma \times 1$, while
$1\otimes 1 \otimes V$ is a generating, wandering subspace for $\Gamma \times \Gamma^{\rm op}$
acting on $H_0\otimes H_0 \otimes V$.

Thus $\H^{\Gamma \times 1}(1\otimes H_0 \otimes V)$ is identified with 
$1\otimes H_0 \otimes V$ and $\H^{\Gamma \times \Gamma}$ is identified with
$1\otimes 1 \otimes V$.

Then, the unitary representation $\Pi$ of $G$, associated in Theorem \ref {2heck} to this unitary representation of 
$C^*((G\times G^{\rm op})\rtimes C(K))$, is acting on $\overline{H_0} \otimes V$ and is
described by the linear application $\Pi(\sigma)$, $\sigma\in G$, mapping the vector $\xi \times v$, $\xi\in\overline{H_0}$, $v\in V$, into
$$
\sum_{\theta\in\Gamma\sigma\Gamma}\big\langle\pi(\theta)1,1\big\rangle\overline{\pi(\theta)}
\xi\otimes\big(\theta\chi_{\Gamma\sigma\theta^{-1}}(v)\sigma^{-1}\big).
$$
If $\tilde{\lambda}$ is the representation of $G$ on $\overline{H_0} \otimes V$ defined
as $\sigma\to \pi_0(\sigma)\otimes\pi(\sigma\otimes 1)$, $\sigma\in G$
then let $T^{[\Gamma\sigma\Gamma]}$ be the image of $t^{[\Gamma\sigma\Gamma]}$,
considered in Section 3, via this representation.

Note that mapping $[\Gamma\sigma\Gamma]$ into $T^{[\Gamma\sigma\Gamma]}$ is a representation
of the Hecke algebra on the Hilbert space $\overline{H_0} \otimes V$.

Then $\Pi(\sigma) = T^{[\Gamma\sigma\Gamma]}(1\otimes\pi(1\otimes\sigma^{-1}))\in \B(H_0\otimes V)$.

The representation of the Hecke algebra associated to the representation $\Pi$, will
act on the Hilbert space $V$. For $\sigma$ in $G$, the Hecke operator associated to the coset
$[\Gamma\sigma\Gamma]$ will map a vector $v \in H_0$ into
$$
\sum_{\theta_1,\theta_2\in\Gamma\sigma\Gamma}
\langle\pi(\theta_1)1,1\rangle
\langle\overline{\pi(\theta_2)}1,1\rangle
\theta_1\chi_{{\theta_1}^{-1},\theta_2}(v)\theta_2^{-1}.
$$
\end{ex}

\begin{proof}
We have proved that the matrix coefficients for the representation $\Pi$ associated to the action
of the $C^*$-algebra $C^*((G\times G^{\rm op})\rtimes C(K))$  on Hilbert space $V$,
with $\Gamma$-wandering, generating space $W_1$ are, for 
$v = \sum\limits_{\gamma\in\Gamma}\gamma w$, $w\in W_1$,
$$
\Pi(\sigma)v=\sum_i s_i\sigma\chi_{\sigma^{-1},\sigma}(v)\sigma^{-1}=
\sum_i s_i\sigma\chi_{\sigma^{-1},\sigma}(v)w\sigma^{-1}.
$$
Hence, for $w_1,w_2 \in W$, we have
\begin{gather*}
\Big\langle\sum_i s_i\sigma\chi_{\Gamma_{\sigma^{-1}}} \Big(\sum\gamma_1w_1\Big){\sigma^{-1}},
\sum\gamma_2w_2\Big\rangle_{ V^{\Gamma\times 1}(W)} =\\
=\sum_{i,\gamma}\big\langle s_i\sigma\chi_{\sigma_1\sigma^{-1}}(\gamma w_1)\sigma^{-1},w_2\big\rangle=
\sum_{\theta\in\Gamma\sigma\Gamma}\langle\theta\chi_{\theta,\sigma^{-1}}(w_1)
\sigma^{-1},w_2\rangle.
\end{gather*}

We take two vectors $\xi_i \otimes v_i$, $i=1,2$, in $1\otimes \overline{H_0} \otimes V$
and identity these vectors with the corresponding $\Gamma \times 1$ invariant vectors 
$$
\sum_{\gamma\in\Gamma}\pi_0(\gamma)1\otimes\overline{\pi_0(\gamma)\xi_i}\otimes
\gamma v_i,\quad i=1,2.
$$
It follows that the matrix coefficients corresponding to this vectors are
\begin{gather*}
\sum_{\theta\in\Gamma\sigma\Gamma}
\langle\pi_0(\theta)\, 1\otimes\overline{\pi_0(\theta)\xi_i}\otimes
\theta\chi_{\theta,\sigma^{-1}}(v_1)\sigma^{-1}, 1\otimes\xi_2\otimes v_2\rangle=\\
=
\sum_{\theta\in\Gamma\sigma\Gamma}
\langle\pi_0(\theta)1,1\overline{\langle\pi_0(\theta)\xi_1,\xi_2\rangle}
\langle\theta\chi_{\theta,\sigma^{-1}}(v_1)\sigma^{-1}, v_2\rangle.
\end{gather*}
But these  are  the matrix coefficients for the unitary representation $\Pi(\sigma)$, $\sigma \in G$,
that was announced in the statement.

The same argument will then work for the formula of the Hecke operators associated to $\Pi$.

Assume now that the representation $V$ is $\ell^2(\Gamma)$ with the canonical action of
the left and right representation.

Then the representation $\Pi$ will have the matrix coefficients on vectors $\xi \otimes \gamma$,
$\xi \otimes \gamma$ equal to
$$
\sum_{\theta\in\Gamma\sigma\Gamma}
\langle\pi_0(\theta)1,1\rangle\overline{\langle\pi_0(\theta)\xi,\xi\rangle}
\langle\theta\chi_{\theta,\sigma^{-1}}(\gamma)\sigma^{-1}, \gamma\rangle.
$$

In this sum the only non zero terms are obtained if
$\theta\chi_{\theta,\sigma^{-1}}(\gamma)\sigma^{-1}=\gamma$,
i.e., $\theta=\gamma\sigma\chi_{\theta,\sigma^{-1}}(\gamma)$,
i.e., $\theta=\gamma\sigma\gamma^{-1}$. This last equality holds true only if $\gamma$ belongs to $\sigma\Gamma\sigma^{-1}\cap\Gamma$. Such a term
  would give a matrix coefficient of the type 
$\langle\pi(\sigma)1,1\rangle\overline{\langle\pi(\sigma)\gamma,\gamma\rangle}$,
i.e., 
$$
\langle(\pi\otimes\overline{\pi})(\sigma)1\otimes\gamma,1\otimes\gamma\rangle_{H_0\otimes H_0}.
$$

This corresponds to the fact that the Hecke operators associated to this data are
obtained from the (diagonal) representation of $G$
$$
\sigma\to(\pi\otimes\pi)(\sigma) \hbox{ on } H_0\otimes\overline{H_0}.
$$
Hence these are the Hecke operators on $H_0 \otimes \overline{H_0}$,excluding the
part that gives eigenvalue 1, which is the subspace generated by $\{\pi_0(\gamma)1\otimes\overline{\pi_0(\gamma)}1,\gamma\in\Gamma\}$.

The matrix coefficients for the Hecke operators will be obtained summing over
$\Gamma\sigma\Gamma$.

Thus, the diagonal matrix coefficients of the Hecke operators, evaluated at elements of the group $\Gamma$  have the formula, for $\gamma\in\Gamma\setminus\{e\}$,
$$
\sum_{\theta\in\Gamma\sigma\Gamma}
\langle\pi_0(\theta)1,1\rangle\overline{\langle\pi_0(\theta)\gamma,\pi_0(\theta)\gamma\rangle},\quad
\sigma\in G.
$$
\end{proof}


\begin{thebibliography}{DeCaHa}

\bibitem[Ar]{Ar}W. Arveson, William Noncommutative dynamics and $E$-semigroups. Springer Monographs in Mathematics. Springer-Verlag, New York, 2003.

\bibitem[AO]{AO}  C.A. Akemann and P.A. Ostrand, On a tensor product C*-algebra associated with the free group on two generators, J. Math. Soc. Japan 27 (1975), 589--599.


\bibitem[AD]{AD} C. Anantharaman-Delaroche,  Amenability and exactness for dynamical systems and their C*-algebras,
Trans. Amer. Math. Soc. 354 (2002), 4153--4178.

\bibitem[Ar]{Ar}Arveson, William Noncommutative dynamics and $E$-semigroups. Springer Monographs in Mathematics. Springer-Verlag, New York, 2003.


\bibitem[AR]{AR}  C. Anantharaman-Delaroche et J. Renault, Amenable groupoids (avec un appendice par E. Germain), Monographie de l'Enseignement Mathematique (Geneve), 36.

\bibitem[Bel]{Bel} Belavkin, V. P. Optimal quantum filtration of Markovian signals. (Russian) Problems Control Inform. Theory/Problemy Upravlen. Teor. Inform. 7 (1978), no. 5, 345--360.

\bibitem [Be]{Be} Berezin, F. A. Quantization. (Russian) Izv. Akad. Nauk SSSR Ser. Mat. 38 (1974), 1116--1175. 

\bibitem[BhPa]{BhPa} Bhat, B. V. Rajarama; Parthasarathy, K. R. Kolmogorov's existence theorem for Markov processes in $C^\ast$ algebras. K. G. Ramanathan memorial issue. Proc. Indian Acad. Sci. Math. Sci. 104 (1994), no. 1, 253--262.


\bibitem[BC]{BC}
Bost, J.-B.; Connes, A., Hecke algebras, type III 
factors and phase transitions with spontaneous 
symmetry breaking in number theory.  Selecta 
Math. (N.S.)  {\bf 1} (1995),  no. 3, 411--457.


\bibitem[BH]{BH} Bridson, Martin R.; de la Harpe, Pierre Mapping class groups and outer automorphism groups of free groups are $C\sp *$-simple. J. Funct. Anal. 212 (2004), no. 1, 195--205. 

\bibitem[Ca]{Ca}  Calkin, J. W. Two-sided ideals and congruences in the ring of bounded operators in Hilbert space. Ann. of Math. (2) 42, (1941). 839-873. 

\bibitem[Co]{Co} Connes, A. Classification of injective factors. Cases $II\sb{1},$ $II\sb{\infty },$ $III\sb{\lambda },$ $\lambda \not=1$. Ann. of Math. (2) 104 (1976), no. 1, 73--115.

\bibitem[CM]{CM} A. Connes, H. Moscovici Modular Hecke algebras and their Hopf symmetry. Mosc. Math. J. 4 (2004), no. 1, 67--109.

\bibitem[Cu]{Cu} N. Cutland, Loeb measure in practice, Recent Advances, Lecture Notes in Mathematics,
vol. 1751, Springer, 2001.

\bibitem[DeCaHa]{DeCaHa} De Canni\`ere, Jean; Haagerup, Uffe Multipliers of the Fourier algebras of some simple Lie groups and their discrete subgroups. Amer. J. Math. {\bf 107} (1985), no. 2, 455--500.

\bibitem[Do]{Do} R. Douglas, Banach algebra, Techniques in Operator Theory, Academic Press. 

\bibitem[DRS]{DRS} W. Duke, Z. Rudnick and P. Sarnak, Density of integer points on affine 
homogeneous varieties. Duke Math. J. {\bf 71} (1993), 143--179. 


\bibitem[EM]{EM}A. Eskin and C. McMullen, Mixing, counting and equidistribution in Lie 
groups. Duke Math. J. {\bf 71} (1993), 181--209. 

\bibitem[FV]{FV}  Pierre Fima, Stefaan Vaes, HNN extensions and unique group measure space decomposition of II$_1$ factors, preprint
arXiv:1005.5002


\bibitem[GeGr]{GeGr} Gelfand, I. M.; Graev, M. I.; Pyatetskii-Shapiro, I. Representation theory and automorphic functions. Saunders Co., Philadelphia, Pa.--London--Toronto, Ont. 1969

\bibitem[GHJ]{GHJ}
Goodman, Frederick M.; de la Harpe, Pierre; 
Jones, Vaughan F.R., Coxeter graphs and towers of 
algebras. Mathematical Sciences Research 
Institute Publications, {\bf 14}. Springer-Verlag, New York, 1989, x+288 pp.


\bibitem [GoNe]{GoNe} Gorodnik, Alexander; Nevo, Amos, The ergodic theory of lattice subgroups. Annals of Mathematics Studies, 172. Princeton University Press, Princeton, NJ, 2010.

\bibitem[Ha]{Ha} Rachel W. Hall, Hecke C*-algebras, Thesis,
Pennsylvania State University, 1996.

\bibitem[Hej]{Hej}
Hejhal, D.A.; Arno, S., On Fourier coefficients 
of Maass waveforms for ${\rm PSL}(2,\bf Z)$.  
Math. Comp.  {\bf 61} (1993),  no. 203, 245--267, 
S11--S16.

\bibitem [He1]{He1} Hejhal, Dennis A.
The Selberg trace formula for ${\rm PSL}(2,\,R)$. Vol. 2. 
Lecture Notes in Mathematics, 1001. Springer-Verlag, Berlin, 1983.

\bibitem [Jo]{Jo} P. Jolissaint. Rapidly decreasing functions in reduced C*-algebras of groups. Trans. Amer. Math. Soc. {\bf 317} (1990), 167--196.

\bibitem [Ke]{Ke} A.S. Kechis, Global aspects of ergodic group actions,
Amer. Math. Soc., Mathematical Survey and Monographs, vol. 160, 2010.

\bibitem [KL]{KL}Kelly-Lyth, D.
Uniform lattice point estimates for co-finite Fuchsian groups,
Proc. London Math. Soc. (3) {\bf 78} (1999), no. 1, 29--51. 


\bibitem[Krieg]{Krieg}
Krieg, Aloys Hecke algebras.  Mem. Amer. Math. Soc.  {\bf 87}  (1990),  no. 435,

\bibitem[Lo]{Lo} P.A. Loeb, Conversion from nonstandard to standard measure spaces
and applications in probability theory, Trans. Amer. Math. Soc. {\bf 211} (1975),
113--122.

\bibitem[Ng]{Ng} Dang-Ngoc Nghiem, D\'ecomposition et classification des syst\`emes dynamiques,
Bull. Soc. Math. France {\bf 103} (1975), 149--175.

\bibitem[Oz]{Oz} N.Ozawa: A Kurosh type theorem for type II. 1 factors, Int. Math. Res. Not. (2006), 
Volume 2006, 1--21.



\bibitem[Oz1]{Oz2} Ozawa, Narutaka, Solid von Neumann algebras. Acta Math. 192 (2004), no. 1, 111Ð117


\bibitem[Perel]{Perel}
Perelomov, A.M., Remark on the completeness of 
the coherent state system. (Russian)  Teoret. 
Mat. Fiz.  {\bf 6} (1971), no. 2, 213--224.

\bibitem[PP]{PP}
Pimsner, Mihai; Popa, Sorin, Entropy and index for 
subfactors.  Ann. Sci. \'Ecole Norm. Sup. (4)  {\bf 19} (1986),  no. 1, 57--106.

\bibitem[Pi]{Pi} Pisier, Gilles Introduction to operator space theory. London Mathematical Society Lecture Note Series, 294. Cambridge University Press, Cambridge, 2003.

\bibitem[Po]{Po}
Popa, Sorin, Notes on Cartan subalgebras in type 
${\rm II}\sb 1$ factors.  Math. Scand.  {\bf 57}  
(1985),  no. 1, 171--188.


\bibitem[Po1]{Po1}
S. Popa, Some properties of the symmetric enveloping algebras with applications to amenability and property T, Documenta Mathematica {\bf 4} (1999), 665--744.


\bibitem[Po2]{Po2} Popa, Sorin,
On the superrigidity of malleable actions with spectral gap. 
J. Amer. Math. Soc. {\bf 21} (2008), no. 4, 981--1000. 

\bibitem[PS]{PS}
Popa, Sorin; Shlyakhtenko, Dimitri, Cartan 
subalgebras and bimodule decompositions of ${\rm 
II}\sb 1$ factors.  Math. Scand.  {\bf 92} (2003), no. 1, 93--102.

\bibitem[Py]{Py} Pytlik, T., Radial functions on free groups and a decomposition of the regular representation into irreducible components. J. Reine Angew. Math. 326 (1981), 124--135.


\bibitem[Ra1]{Ra1}
R\u adulescu, Florin, The $\Gamma$-equivariant 
form of the Berezin quantization of the upper 
half plane.  Mem. Amer. Math. Soc.  {\bf 133} (1998),  no. 630.

\bibitem[Ra2]{Ra2}
R\u adulescu, Florin, Arithmetic Hecke operators 
as completely positive maps.  C.R. Acad. Sci. 
Paris S\'er. I Math.  {\bf 322}  (1996),  no. 6, 541--546.

\bibitem[Ra3]{Ra3} R\u adulescu, Florin,  Free group factors and Hecke operators, Notes taken by N. Ozawa, to appear.


\bibitem[Ra4]{Ra4} R\u adulescu, Florin,  Cyclic Hilbert Spaces,
Studies in Informatics and Control, Vol. 18, No. 1/2009, pp. 83--86.

\bibitem[Ra5]{Ra5} R\u adulescu, Florin,  A universal, non-commutative C*-algebra associated to the Hecke algebra of double cosets, preprint, arxiv


\bibitem[Ra6]{Ra6} F. R\u adulescu, On the representation of a discrete group $\Gamma$ 
with subgroup $\Gamma_0$ in the Calkin algebra of $\ell^2(\Gamma/\Gamma_0)$,
preprint, arxiv.

\bibitem[Ra7]{Ra7} F. R\u adulescu,  Conditional expectations, traces, angles between spaces and Representations of the Hecke algebras, preprint, arxiv


\bibitem [Sa] {Sa} Sarnak, Peter, Statistical properties of eigenvalues of the Hecke operators. Analytic number theory and Diophantine problems (Stillwater, OK, 1984), 321--331, Progr. Math., 70, Birkh\"auser Boston, Boston, MA, 1987.

\bibitem [Sa1] {Sa1} Sarnak, Peter, Some Applications of Modular Forms, Cambridge University Press, 1990.

\bibitem [Su] {Su} Sunder, V. S. An invitation to von Neumann algebras. Universitext. Springer-Verlag, New York, 1987. xiv+171 pp.

\bibitem[Tz]{Tz}
Tzanev, Kroum,
Hecke $C\sp *$-algebras and amenability.
J. Operator Theory {\bf 50} (2003), no. 1.

\bibitem [Ve]{Ve} Vershik, A.M. Krein's duality, positive 2-algebras, and the dilation of comultiplications. Funct. Anal. Appl. {\bf 41} (2007), no. 2, 99--114.

\bibitem[Za]{Za} P. Zakrzewski, On nonmeasurable selectors of countable group actions,
Fund. Math. {\bf 202} (2009), 281--294.

\end{thebibliography}
\end{document}